\documentclass[12pt]{article}
\usepackage{amssymb,amsmath,amsfonts,mathrsfs,upgreek} 

\usepackage{bm}

\usepackage{graphicx}

\usepackage[colorlinks=true]{hyperref}
\usepackage{comment}

\newtheorem{theorem}{Theorem}[section]
\newtheorem{corollary}{Corollary}[section]

\newtheorem{lemma}{Lemma}[section]
\newtheorem{proposition}{Proposition}[section]

\newtheorem{remark}{Remark}[section]
\newtheorem{notation}{Notation}
\newcommand{\diam}{{\rm diam}}
\newcommand{\io}{{\infty}}

\newcommand{\real}{ {\mathbb R}   }
\newcommand{\torus}{ {\mathbb T}   }
\newcommand{\integer}{ {\mathbb Z}   }
\newcommand{\complex}{ {\mathbb C}   }

\renewcommand{\Im}{\, {\rm Im}\,}
\renewcommand{\Re}{\, {\rm Re}\,}

\newcommand{\eproof}{\qed}

\newcommand\beq[1]{ \begin{equation}\label{#1} }
\newcommand{\eeq}{ \end{equation} }

\newcommand\beqa[1]{ \begin{eqnarray} \label{#1}}
\newcommand{\eeqa}{ \end{eqnarray} }
\newcommand{\beqano}{ \begin{eqnarray*} }
\newcommand{\eeqano}{ \end{eqnarray*} }
\newcommand\arr[1]{\left\{\begin{array}{l}#1\end{array}\right.}
\newtheorem{definition}{Definition}[section]
\newcommand\dfn[1]{ \begin{definition}\label{#1} \rm}
\newcommand\edfn{ \end{definition} }
\newcommand{\proof}{\par\medskip\noindent{\bf Proof\ }}
\newcommand\equ[1]{{\rm (\ref{#1})}}
\newcommand{\nl}{{\smallskip\noindent}}

\newcommand{\Giu}{{\bigskip\noindent}}
\newcommand{\noi}{{\noindent}}
\newcommand{\qed}{\hskip.5truecm
\vrule width 1.7truemm height 3.5truemm depth 0.truemm
\par\Giu}
\newcommand{\qedeq}{\hskip.5truecm
\vrule width 1.7truemm height 3.5truemm depth 0.truemm}

\renewcommand{\b }{\beta }
\newcommand{\s }{\sigma }
\newcommand{\ii }{{\rm i} }
\renewcommand{\d }{\delta }
\newcommand{\D }{\Delta}

\newcommand{\g }{\gamma}
\newcommand{\f }{\varphi}

\renewcommand{\l }{\lambda }

\renewcommand{\O }{\Omega }
\newcommand{\C}{\mathbb{C}}
\newcommand{\Z}{\mathbb{Z}}
\newcommand{\z }{\zeta }

\newcommand{\checco}{\eta }
\newcommand{\ch}{{\bm{\eta}} }

\newcommand{\Hpend}{\Hpend }

\newcommand{\cgot}{\mathfrak c}

\def\R{\mathbb R}
\def\T{\mathbb T}

\def\const{{\rm\, const\, }}
\def\dst{\displaystyle}

\def\meas{{\rm\, meas\, }}

\newcommand\notat[1]{ \begin{notation} \label{#1} 
 }
\newcommand\enotat{\end{notation}}

\newcommand\rem{\begin{remark} 
}
\newcommand\erem{\end{remark} 
}

\newcommand{\normadue}{|}

\newcommand\casi[2]{ \left\{  
\begin{array}{l}
 {#1}  \\
 {#2} 
 \end{array} \right.}

\renewcommand{\Hpend}{H_{\mathtt{mech}}}

\newcommand\modulo{|}

\newcommand\sa{\theta} 
\newcommand\Sa{\Theta} 
\newcommand\loge{\lambda}
\newcommand\GG{{\mathtt G}}
\newcommand\GGS{{\mathtt G^*}}
\newcommand\FO{{\bar{\mathtt G}}}  
\newcommand\Gm{{\mathcal G}}

\renewcommand{\checco}{\eta }

\newcommand\hol{{\bf B}}

\newcommand{\morse}{{\rm M}}  
\newcommand\act{I}  
\newcommand\ang{\f}  

\newcommand\Fit{\Phi_{\mathtt{mech}}}
\newcommand\Fiq{\Phi_{\mathtt{aa}}^i}

\newcommand\cFiq{\check\Phi_{\mathtt{aa}}^i}
\newcommand\Bu{{\mathcal B}}

\newcommand{\acci}{{\mathtt a}}

\newcommand{\HS}{{H^*}}

\newcommand{\xx}{{\theta}}

\title{\bf
Action-angle Variables for Generic 1D Mechanical Systems
\ \\ \ \\ 
\small{(Draft)}
}

\begin{document}

\author{ 
\footnotesize L. Biasco  \& L. Chierchia
\\ \footnotesize Dipartimento di Matematica e Fisica
\\ \footnotesize Universit\`a degli Studi Roma Tre
\\ \footnotesize Largo San L. Murialdo 1 - 00146 Roma, Italy
\\ {\footnotesize biasco@mat.uniroma3.it, luigi@mat.uniroma3.it}
\\ 
}

\date{February 4, 2017}

\maketitle

\vskip1.truecm

\begin{abstract}
\noindent
We consider a 1D mechanical system 
$$\bar {\mathtt  H}(\mathtt P,\mathtt Q)=\mathtt P^2+\FO(\mathtt Q)$$
 in  action-angle variable $(\mathtt P,\mathtt Q)$ 
 where $\FO$ is a $2\pi$-periodic analytic function 
 with non degenerate critical points. 
 Then, we consider a small analytic perturbation of $\bar {\mathtt  H}$  of the form 
 $$
  {\mathtt  H}^*(\mathtt P,\mathtt Q;\hat{\mathtt  P})=
 \mathtt P^2+\FO(\mathtt Q)+ \eta {\mathtt F} (\mathtt P,\mathtt Q;\hat{\mathtt  P})=:\mathtt P^2+
 \GGS(\mathtt P,\mathtt Q;\hat{\mathtt  P})\,, \qquad \eta\ll 1\ ,$$
 where the perturbed potential $\GGS$ may depend on the action
 $\mathtt P$ and also on parameters $\hat{\mathtt  P}$
 (``the adiabatic  actions''); indeed, this is the form of a finite dimensional
  mechanical system
 close to an exact simple resonance after averaging over fast angles and disregarding the exponentially small remainder, see \cite{BCnonlin}.
 \\
 Up to a finite number of separatrices and elliptic/hyperbolic points
 the phase space of  ${\mathtt  H}^*$ is divided into a finite number of open
 connected components foliated by invariant circles.
 On every connected component we perform a (Arnold--Liouville) symplectic 
 action-angle
 transformation which integrates the system.
 We give a complete and quantitative description of the analyticity
 properties of such integrating transformations, estimating,  in particular,
 how  such transformations differ from the integrating transformation for $\bar {\mathtt  H}$; compare
 Theorem \ref{glicemiak} below.

\end{abstract}

\tableofcontents

\section{Set up and notations}

\begin{itemize}

\item[\rm (i)]
{\bf Norms on finite dimensional vector spaces} \\
$|\cdot|$ denotes the standard Euclidean norm on $\complex^n$ and its subspaces. 
\\
For linear maps and matrices (which are identified),
$\normadue\cdot\normadue$ denotes the  ``operator norm'' 
$$\dst |A|=\sup_{u\neq 0} |Au|/|u|\ .$$ 
\\
 $|k|_{{}_1}$  denotes the 1-norm $\sum |k_j|$.
\\
$|M|_{{}_\infty}$, with $M$ matrix (or vector), denotes the maximum norm $\max_{ij}|M_{ij}|$ (or $\max_i |M_i|$).

\item[\rm (ii)] {\bf Open covers}\\
Given a set $D\subseteq \R^m,$ $r>0$ 
we denote by $D_r\subseteq\C^m$  the complex  open neighborhood of $D$ formed by points  $z\in \complex^m$ such that $|z-y|<r$, for some $y\in D$.

\nl
Given $s>0,$ we denote by $\torus^n_s$ the open complex neighborhood of\footnote{$\torus^n$ denotes the standard flat $n$--dimensional torus $\real^n/(2\pi \integer^n)$.} $\torus^n$ given by 
$$
\T^n_s:=\{ x\in\C^n\ \ :\ \  \max_{1\leq j\leq n}|\Im x_j|<s \}
/2\pi\Z^n
$$
\item[\rm (iii)] {\bf Norms of analytic functions}\\
Given a  real--analytic function\footnote{$f_k$ denotes Fourier coefficients.}  $f:  \T^n_s\to\C$,
$f( x)=\sum_{k\in\Z^n}f_k e^{\ii k\cdot x}$,  
we define the ``sup--Fourier norm''
\begin{equation}\label{enorme}
\|f\|_s:=\sup_{k\in \integer^n} |f_k| e^{|k|_{{}_1}s}<\infty\ .
\end{equation}
Analogously, if
$f:D_r\times \T^n_s\to\C$,
$f(y,x)=\sum_{k\in\Z^n}f_k(y) e^{\ii k\cdot x},$
we let
\begin{equation}\label{laga}
\|f\|_{r,s}
 := \sup_{k\in\Z^n}
 \big(
 \sup_{y\in D_r}|f_k(y)|e^{|k|_{{}_1}s}
 \big)
 \, .
\end{equation} 
If the (real) domain needs to be specified, we let:
\beq{romaostia}
\|f\|_{D,r,s}:=\|f\|_{r,s}\ .
\eeq
Given a bounded holomorphic function 
$f: \T^n_s\to\C$, 
we set
\begin{equation}\label{ionicoionico}
{\modulo}f{\modulo}_s:=\sup_{\T^n_s}|f|\,.
\end{equation}
Given a bounded holomorphic function 
$f:D_{r}\times \T^n_s\to\C^m$, or $f:D_{r} \to\C^m$ with $D\subseteq \R^n
$
we set
\begin{equation}\label{ionico}
{\modulo}f{\modulo}_{D,r,s}={\modulo}f{\modulo}_{r,s}:=\sup_{D_r\times \T^n_s}|f|\,,\ 
{\rm or, \ respectively,\ }\quad 
{\modulo}f{\modulo}_{D,r}={\modulo}f{\modulo}_{r}:=\sup_{D_r}|f|\ .
\end{equation}
Notice that the following relations between the two norms
$\|\cdot\|$ and ${\modulo}\cdot{\modulo}$ hold:
for $\s>0$, we have\footnote{
Indeed: $\dst \sum_{k\in\Z^n}e^{-|k|_{{}_1}\s}=\Big(
\sum_{k\in\Z}e^{-|k|_{{}_1}\s}\Big)^n
=
\Big(1+ 2\sum_{j\geq 1}e^{-j\s}\Big)^n
=
\Big(\frac{e^\s+1}{e^\s-1}\Big)^n
= \coth^n(\s/2).$}
\begin{equation}\label{battiato2}
\|f\|_{r,s}\leq {\modulo}f{\modulo}_{r,s}\leq \coth^n(\s/2)\|f\|_{r,s+\s}
\leq
(1+2/\s)^n
\|f\|_{r,s+\s}\,.
\end{equation}

\item[\rm (iv)] {\bf Scale of Banach spaces of real--analytic periodic functions}\\
For $s\ge0$, we denote by  $\hol_s^n$ the following Banach space of real--analytic functions on $\T^n$ with vanishing average\footnote{The bar denotes complex conjugate.}:  
$$
\hol_s^n:=\big\{f:\torus^n_s\to\complex\ {\rm s.t.}\  \|f\|_s<\io \ {\rm and}\ f_0=0\ ,\ \overline{f_k}=f_{-k}\big\}\ .
$$

\end{itemize}

\section{Standard form of parametrized 1D mechanical systems}

Let $r_0,R_0,s_0>0$, $\hat D\subset \real^{n-1}$ and consider  the Hamiltonian
\begin{equation}\label{spigola}
\HS:= P_n^2
+
\GGS(P,Q_n)\,,
\end{equation}
real--analytic for
\begin{equation} \label{sontuosa}
(P,Q_n)\in D_{r_0}\times\T_{s_0}\,,\qquad \text{where}\qquad
D:=
\hat D \times (-R_0 ,R_0 )\ .
\end{equation}
This Hamiltonian represents a  
1D system in the symplectic variables $(P_n,Q_n)\in\real\times\torus$ depending on the parameter $\hat P=(P_1,...,P_{n-1})\in\real^{n-1}$.

Now, we want to reduce \equ{spigola} to a {\sl parameterized 1D mechanical system} (i.e., with a potential independent of $P_n$). We do this via the following  ``normalization lemma'', where $\HS$ is considered as a function of $2n$ variables $(P,Q)$.

\begin{lemma}\label{francesco}{\bf (Standard form of parametrized 1D mechanical systems)}\\
Let $D=\hat D\times(-R_0,R_0), \hat D\subset \real^{n-1}$, $s_0$, $r_0, R_0>0$, $0<\hat s\le s_0/2$, $\FO\in\hol_{s_0}^1$ and let
$$
\HS(P,Q):=P_n^2 +    \GGS(P,Q_n)\,,\qquad {\rm with}\quad {\modulo}\GGS-\FO{\modulo}_{D,r_0,s_0}\leq \checco_0
$$
and 
\begin{equation}\label{urea}
\checco_0\leq \frac{r_0^2}{64}\cdot\min\Big\{ \frac{s_0 }{\pi}\ ,\  1\Big\}
\ .
\end{equation}
Then,  there exists a symplectic transformation 
\begin{equation}\label{flaviano}
\Fit\ :(p,q)\in\ \hat D_{r_0}\times (-R_0,R_0)_{r_0/2}\times
 \T^{n-1}_{\hat s}\times\T_{s_0}
 \ \mapsto\ (P,Q)\in 
 D_{r_0}\times
 \T^{n-1}_{\hat s+16\pi\checco_0/r_0^2}\times\T_{s_0}\,,
\end{equation}
of the form
\begin{equation}\label{rondo}
(P,Q)=\Fit(p,q) :\ 
\casi{\hat P =\hat p}
{P_n=
p_n+a_*(\hat p,q_n)
=p_n-
{\mathtt P}^*_n(\hat p)
+ \mathtt P(\hat p,q_n)
}
\ 
\casi{\hat Q=\hat q +b_*(\hat p,q_n)}
{Q_n=q_n}\qquad\qquad
\end{equation}
such that $\HS\circ\Fit(p,q)=:\Hpend(p,q_n)$ has the ``standard form''
\beq{pasqua}
\casi{\boxed{\Hpend(p,q_n)=\big(1+ b (p,q_n)\big) \big(p_n-{\mathtt P}^*_n(\hat p)\big)^2  
  + \Gm (\hat p,q_n )} \phantom{\dst \sum_{i}^1}}
{{\rm where} \quad
\Gm:= \GGS(\hat p, \mathtt P(\hat p,q_n),q_n )
+ (\mathtt P(\hat p,q_n))^2\ .}
\eeq
Furthermore, the following estimates hold:
\begin{equation}\label{urina3}
{\modulo}{\mathtt P}^*_n{\modulo}_{\hat D,r_0}\leq 2\checco_0/r_0
\leq r_0/8\,,\quad 
{\modulo}a_*{\modulo}_{\hat D,r_0,s_0}\leq 4\checco_0/r_0\,,
\quad 
{\modulo}b_*{\modulo}_{\hat D,r_0/2,s_0}
\leq 16\pi\checco_0/r_0^2\,,
\end{equation}
and
\begin{eqnarray}\label{straussbis}
&&
{\modulo}\Gm-\FO {\modulo}_{\hat D,r_0,s_0}
\leq 2 \checco_0
\,,
\qquad
{\modulo}b{\modulo}_*\leq
\frac{32}{r_0^2}\checco_0
\,,
\qquad
{\modulo}\partial_{p_n}
b{\modulo}_*\leq
\frac{64}{r_0^3}\checco_0
\,,
\nonumber
\\
&&
 {\modulo} p_n \cdot  b(p,q_n){\modulo}_*\leq
\frac{10}{r_0}\checco_0\,,
\qquad
 {\modulo} p_n \cdot  \partial_{p_n}b(p,q_n){\modulo}_*\leq
\frac{100}{r_0^2}\checco_0
\,,
\end{eqnarray}
where by 
\begin{equation}\label{wagner}
|\cdot|_*:=\sup_{\hat D_{r_0}\times(-R_0,R_0)_{r_0/2}\times
\T_{s_0}}|\cdot|\,. 
\end{equation}
\end{lemma} 
\proof
Recalling \eqref{A1}, \eqref{urea}
and noting that $\FO$ does not depend on 
$p_n,$  by Cauchy estimates
we have
\begin{eqnarray}\label{mosella}
&&\sup_{p_n\in (-R_0,R_0)_{3r_0/4}}
{\modulo}\partial_{p_n}\GGS(\cdot,p_n,\cdot){\modulo}_{\hat D,r_0,s_0}\leq
\frac{4}{r_0}\checco_0\,,
\nonumber
\\
&&\sup_{p_n\in (-R_0,R_0)_{3r_0/4}}
{\modulo}\partial^2_{p_n}\GGS(\cdot,p_n,\cdot){\modulo}_{\hat D,r_0,s_0}\leq
\frac{32}{r_0^2}\checco_0\,,
\nonumber
\\
&&\sup_{p_n\in (-R_0,R_0)_{3r_0/4}}
{\modulo}\partial^3_{p_n}\GGS(\cdot,p_n,\cdot){\modulo}_{\hat D,r_0,s_0}\leq
\frac{384}{r_0^3}\checco_0\,.
\end{eqnarray}
As one easily checks by \eqref{urea},
 the fixed point equation
\begin{equation}\label{olimpia}
\mathtt P(\hat p,q_n)
=
-\frac12 \partial_{p_n} \GGS\big( \hat p,\mathtt P(\hat p,q_n),q_n   \big)
\end{equation}
has 
a unique solution
 $\mathtt P=\mathtt P(\hat p,q_n)$ with
\begin{equation}\label{urina}
{\modulo}\mathtt P{\modulo}_{\hat D,r_0,s_0}\leq 2\checco_0/r_0\leq r_0/8\,.
\end{equation}
 Set
\begin{equation}\label{mosa}
{\mathtt P}^*_n(\hat p)
:= \langle \mathtt P(\hat p,q_n) \rangle\,,\qquad 
 a_*(\hat p,q_n ):=\mathtt P(\hat p,q_n)-\langle \mathtt P(\hat p,q_n) \rangle\,,
\end{equation}
where $ \langle \cdot \rangle$ denotes the average with respect to $q_n.$
This proves the first two estimates in
\eqref{urina3}.
 Let $\phi=\phi(\hat p,q_n)$ the unique function
satisfying $a_*=\partial_{Q_n} \phi$
with $\langle \phi\rangle=0$. 
By the second estimate in \eqref{urina3}
we get $|\phi|_{\hat D,r_0,s_0}
\leq 8\pi\checco_0/r_0.$
Set $b_*:=-\partial_{\hat p} \phi;$ then
the third estimate in \eqref{urina3} follows
by Cauchy estimates.
Now, let $\Fit$ be the symplectic transformation in \equ{rondo}
obtained by the
 generating function
$\hat p \cdot \hat Q+p_n Q_n+ \phi(\hat p,Q_n)$.
By the estimates on $a_*$ and $b_*$
in \eqref{urina3} and \eqref{urea}
it turns out that the canonical transformation
$\Fit$ is well defined with respect to the domains
in \eqref{flaviano}.
Then
$\Fit$ casts $\HS$ into
\begin{eqnarray*}
&&\Big(
p_n-
{\mathtt P}^*_n(\hat p)
+ \mathtt P(\hat p,q_n)
\Big)^2+
  \GGS(\hat p,p_n-
{\mathtt P}^*_n(\hat p)
+ \mathtt P(\hat p,q_n) ,q_n)
\\
&&=
\big(1+ b (p,q_n)\big) \big(p_n-{\mathtt P}^*_n(\hat p)\big)^2  
  + \Gm (\hat p,q_n ) \,,
\end{eqnarray*}
where $\Gm$ was defined in \eqref{pasqua}
and\footnote{Using \eqref{olimpia}.}
omitting, for brevity, the dependence on
$\hat p, q_n,$
\begin{equation}\label{uforobot}
b=\frac{
\GGS(\mathtt P+p_n- {\mathtt P}_n^*)
-
\GGS(\mathtt P)
-
\partial_{p_n}\GGS(\mathtt P) ( p_n- {\mathtt P}_n^*)
}
{( p_n- {\mathtt P}_n^*)^2}
=
\int_0^1 (1-t) \partial^2_{p_n}\GGS
\big(\mathtt P+t( p_n- {\mathtt P}_n^*) \big)dt
\,.
\end{equation}
Then the first estimate in \eqref{straussbis} follow by 
\eqref{pasqua},
\eqref{A1}, \eqref{urea}, \eqref{urina}.
Note that by \eqref{urina} and 
\eqref{mosa} if 
$p_n\in(-R_0,R_0)_{r_0/2}$
then 
$\mathtt P+p_n- {\mathtt P}_n^*
\in(-R_0,R_0)_{3r_0/4}.$
By \eqref{uforobot} and
\eqref{mosella} 
we get the second estimate in 
\eqref{straussbis}.
Using the first equality in \eqref{uforobot} and the first estimate
in
\eqref{mosella} we get
$$
 {\modulo} (p_n- {\mathtt P}_n^*) \cdot
 b(p,q_n){\modulo}_*\leq \frac8{r_0}\checco_0
$$
and, therefore, by \eqref{urina3}, the second inequality in\eqref{straussbis}
and \eqref{urea} we get the fourth estimate in 
\eqref{straussbis}.
By \eqref{uforobot}, and the last estimate in \eqref{mosella}
we get\footnote{Note that $\int_0^1(1-t)t\,dt=1/6.$}

Finally, by the first equality in \eqref{uforobot}, we have
\begin{equation*}
\partial_{p_n} b=\frac{
\partial_{p_n}\GGS(\mathtt P+p_n- {\mathtt P}_n^*)
-
\partial_{p_n}\GGS(\mathtt P)
}
{( p_n- {\mathtt P}_n^*)^2}
- \frac{2b}{ p_n- {\mathtt P}_n^*}
\,,
\end{equation*}
and, by \eqref{mosella} and the second estimate in
\eqref{straussbis},
$$
|(p_n- {\mathtt P}_n^*)\cdot 
\partial_{p_n} b(p,q_n)|_* 
\leq
\left|
\frac{
\partial_{p_n}\GGS(\mathtt P+p_n- {\mathtt P}_n^*)
-
\partial_{p_n}\GGS(\mathtt P)
}
{p_n- {\mathtt P}_n^*}
\right|_*
 + 2|b|_*
\leq \frac{96}{r_0^2}\checco_0\,.
$$
Then by \eqref{urina3}, the third inequality in
\eqref{straussbis} and \eqref{urea}
we get the last estimate in \eqref{straussbis} 
\eproof

\begin{remark}\label{allergia}
 In the oscillatory regime, where $q_n$ is not an angle,
 we could perform a canonical transformation 
 (of the form in \eqref{grillo} below)
 which translates $p_n'=p_n-{\mathtt P}^*_n(\hat p)$, so that $\Hpend $
simplifies a little bit, but such translation is impossible in the
rotational regime, where $q_n=q_n'$ is  an angle,
since the transformation of $\hat q = \hat q'-
{\mathtt P}^*_n(\hat p)q_n'$  is not $2\pi$-periodic in 
$q_n'.$
\end{remark}

\rem\label{colomba}
The  symplectic transformation in \equ{rondo} belongs to the group ${\mathfrak G}$ of symplectic transformation of the special form:
\begin{equation}\label{grillo}
\Phi: \ (p,q)\ \mapsto \  (P,Q):\quad 
\casi{\hat P=\hat p}{P_n=P_n^*(p,q_n)}
\qquad
\casi{\hat Q=\hat q+\hat q^*(p,q_n)}
{Q_n=Q_n^*(p,q_n)} 
\end{equation} 
where, in general,  $Q,q$ may belong either to $\T^n$ or to $\R^n$. 
Notice that restricting the relation
$d P\wedge dQ=dp\wedge dq$ onto the planes $\hat p=\hat P=\const$, one has that $dP_n\wedge dQ_n=dp_n\wedge dq_n$, i.e.: 

\nl
{\sl For every fixed $\hat p,$ the map
 $
 (p_n,q_n)\mapsto
 \big(P_n^*(p,q_n), Q_n^*(p,q_n)\big)$
 is also symplectic}.
\erem

\notat{colomba2}
For  a transformation $\Phi\in {\mathfrak G}$  
we let
$\check \Phi$ 
denote the $(n+1)$--dimensional map
\begin{equation}\label{islands}
\check\Phi: \ (p,q_n)\mapsto \ (P,Q_n):=\big( 
\hat p, P_n(p,q_n), q_n(p,q_n)
\big)\,.
\end{equation}
\enotat

\rem\label{colomba3} If $\Phi_i\in {\mathfrak G}$, then 
\begin{equation}\label{peggylee}
{(\Phi_1\circ\Phi_2)}\check{\phantom .}=\check\Phi_1\circ\check\Phi_2\,,
\end{equation}
and, furthermore,
\beq{Ventura}
\check\Phi(E)\times \torus^{n-1}= \Phi(E\times \torus^n)\ ,\qquad \forall\  \Phi\in{\mathfrak G}\ ,\qquad \forall\  E\subseteq \real^n\times \torus\ .
\eeq
Relation \equ{Ventura} implies, in particular, that,
{\sl for avery map $\Phi\in {\mathfrak G}$, the map $\check \Phi$ is volume-preserving}.
\erem

\section{Morse non-degenerate potentials}

\begin{definition}\label{morso}
Let $s_0, \morse, \b>0$ and let $\FO$ be a 
$2\pi$-periodic holomorphic function.
We say that $\FO$  is
$(\morse,\b,s_0)$--Morse--non--degenerate if
\beqa{goffredo}
&& {\modulo}\FO{\modulo}_{s_0}\leq \morse\,,\\
&&\ \nonumber\\
&&\label{ladispoli}
\min_{\sa\in\real} \ \big( |\FO'(\sa)|+|\FO''(\sa)|\big) 
\geq\b\,,
\quad
\min_{1\leq i<j\leq 2N } 
|\bar E_i-\bar E_j|
\geq \b\,, 
\eeqa
where
\beq{maieli}
\bar E_i:=\FO(\bar\sa_i)\ ,\qquad    1\le i\le 2N\ ,
\eeq
are the $2N$ {\sl distinct} critical values (``critical energies'') and  $\bar\sa_i$ the corresponding non--degenerate critical points of $\FO$. 
\end{definition}
Note that by \eqref{goffredo} and \eqref{ladispoli}
we have\footnote{The first estimate is obvious
since $|\bar E_i|\leq \morse;$ then the second estimate 
directly follows if $s_0\leq 1,$ otherwise it follows by the first inequality in \eqref{ladispoli} and Cauchy estimates, which
imply $\b\leq \frac{\morse}{s_0 }+\frac{2\morse}{s_0^2 }.$}
\begin{equation}\label{ocarina}
\frac{\morse}{\b}\geq \frac12\,,\qquad
\frac{\morse}{\b s_0}\geq \frac13\,.
\end{equation}

\nl
{\sl 
We may  assume that, up to translation, the unique absolute maximum of $\FO$ is attained at} 
$$\bar \sa_0:=\bar \sa_{2N}-2\pi=-\pi$$ 
Then, {\sl the relative strict non--degenerate minimum and maximum points
follow in alternating order}: 
\begin{equation}\label{recremisik}
\bar\sa_0:= -\pi < \bar\sa_1<\bar\sa_2<\ldots <\bar\sa_{2N-1}<\bar\sa_{2N}:=\pi\,,\quad 
\left\{\begin{array}{l}
\bar\sa_{2j} \ \ \ \ {\rm maximum\ points}\\ 
\bar\sa_{2j-1}\ {\rm minimum\ points}
\end{array}
\right.
\end{equation}
and:
\beqano
&& the\ odd\ energies\ \bar E_1,...,\bar E_{2N-1}\ { are\ the} \ N\ { (local)\ minimal\  energies}\\
&&  the\ even\ energies\  \bar E_2,...,\bar E_{2N} \ { are\ the} \ N\ {  (local)\ maximal\  energies} 
\eeqano
and 
\begin{equation}\label{tracollo}
\bar E_0:=\bar E_{2N}
\end{equation}
 is the unique  global maximum.

\nl
By Cauchy estimates we get 
$$
\max_\R |\partial_\sa^k \FO|\leq k! \morse/s_0^k\,,\qquad
{\modulo}\partial_\sa^k \FO{\modulo}_{\s}\leq k! \morse/(s_0-\s)^k\,.
$$
Note that by  \eqref{ladispoli}, \eqref{goffredo}
and Cauchy estimates  it follows that\footnote{The first estimates
directly follows by \eqref{goffredo}
and the second inequality in
\eqref{ladispoli}. The other two estimates follow by contradiction: otherwise: i)
in a point $\sa_0$ with $\FO''(\sa_0)=0$, one has,
choosing $\xx=\pm s_0$ according to the sign of
 $\FO'(\sa_0)/\FO(\sa_0),$
 that
 $|\FO(\sa_0+\xx)|\geq|\FO(\sa_0)|+\b|\xx|-\morse\xx^2/s_0^2
 \geq \b s_0-\morse>\morse,$ contradicting \eqref{goffredo} and, so, proving the second estimate in \eqref{harlock}; 
 ii)
 in a point $\sa_1$ with $\FO'(\sa_1)=0$, one has
 $|\FO(\sa_1+\xx)|\geq
 |\FO(\sa_1)+\FO''(\sa_1)\xx^2/2|-
 \morse|\xx^3|/s_0^3,$
 then, by \eqref{ladispoli},
 $\sup_{|\xx|<s_0}|\FO(\sa_1+\xx)|\geq
 \b s_0^2/2-\morse>\morse,$ contradicting \eqref{goffredo} and, so, proving the third estimate in \eqref{harlock}.
}
\begin{equation}\label{harlock}
\b \leq 2 \morse\,,\qquad
\b s_0\leq 2 \morse\,,\qquad
\b s_0^2\leq 4 \morse\,.
\end{equation}
Obviuosly the second assumption in 
\eqref{ladispoli} directly implies \eqref{harlock}.
In the particular important case in which $\FO$
is minus cosine  we can explicitly evaluate
\begin{equation}\label{lontra}
\FO(\sa)=-\cos \sa\quad\Longrightarrow\quad
\left\{
\begin{array}{l}
\morse=\cosh s_0\,,\ \
N=1\,,\ \ 
\bar\sa_1=0\,,\ \bar\sa_2=\pi\\
\bar E_1=-1\,,\ \bar E_2=1\,,\ \
\b=1
\end{array}\right.
\end{equation}

\begin{lemma}\label{grancereali}
We have
\begin{equation}\label{alcafone}
2\xx_*\leq \bar\sa_i-\bar\sa_{i-1}\leq 2\pi\,,
\qquad
N\leq \frac{\pi}{2\xx_*}\,,
\end{equation}
where
 \begin{equation}\label{draghetto}
\xx_*:=\sqrt{\frac{\b s_0^3}{3 \morse}}\,.
\end{equation}

\end{lemma}
\proof
Since $\FO$ is convex in $\bar\sa_{2j-1}$
and concave in $\bar\sa_{2j}$, there exists
$\bar\sa_{2j-1}< \bar\sa<\bar\sa_{2j}$
such that $\FO''(\bar\sa)=0.$
By \eqref{ladispoli} and Cauchy estimates
we get
$$
\FO'(\bar\sa+\xx)\geq \b-3 \morse s_0^{-3}\xx^2\,.
$$
Then $
\FO'(\bar\sa+\xx)>0,
$ when $|\xx|<\xx_*$,
which implies
$\bar\sa_{2j-1}\leq \bar\sa-\xx_*
<\bar\sa+\xx_* \leq \bar\sa_{2j},$
 proving the first estimate
in \eqref{alcafone} (from which the estimate
on $N$ directly follows).
\eproof

\begin{lemma}\label{parcoavventura}
For every $1\leq i\leq 2N$
\begin{eqnarray}
&&\label{harlock1}
|\FO'(\bar\sa_i+\xx)|\geq \frac{\b}{2}|\xx|\,,
\quad \forall\, \xx\in\R\,,\  |\xx|\leq \xx_\sharp
:=\frac{\b s_0^3}{6\morse}
\,,
\ \ \ \ \text{and}\ \ \ \ 
\\
&&\label{harlock2}
\frac{\b s_0^3}{6\pi \morse}\leq 1\,,
\\
&&\label{harlock3}
\min_{[\bar\sa_{i-1}+\xx_\sharp/2,\bar\sa_i-\xx_\sharp/2]}
|\FO'|\geq \frac{\b^2 s_0^3}{32\morse}\,.
\end{eqnarray}
\end{lemma}
\proof
We will consider only the case\footnote{The case
$i=2j-1$ is analogous.} $i=2j$.
For\footnote{The case $-\xx_\sharp\leq \xx\leq 0$ is analogous.} 
$0\leq \xx\leq \xx_\sharp$ we get by Cauchy estimates
$$
\FO'(\bar\sa_i+\xx)
\geq 
\FO''(\bar\sa_i)\xx-\frac{3\morse}{s_0^3}\xx^2
\stackrel{\eqref{ladispoli}}\geq
\b \xx -\frac{3\morse}{s_0^3}\xx^2\geq 
\frac{\b}{2}\xx\,.
$$
Noting  that, as in \eqref{alcafone},
we have 
$2\xx_\sharp\leq \bar\sa_i-\bar\sa_{i-1}\leq 2\pi,$
 we get \eqref{harlock2}.
 \\
Regarding the minimum of $\FO'$
in the interval\footnote{Note that
on such interval $\FO'>0$.}
$[\bar\sa_{i-1}+\xx_\sharp/2,\bar\sa_i-\xx_\sharp/2]$
if it is achieved at the endpoints then the inequality
holds, otherwise if it is achieved in an inner point $\sa_*$ then $\FO''(\sa_*)=0$ and, by \eqref{ladispoli},
$\FO'(\sa_*)\geq \b\geq \frac{\b^2 s_0^3}{32\morse}$
by \eqref{harlock2}.
\eproof
\\
Fix $1\leq j\leq N$ and consider a minimum point $\bar\sa_{2j-1},$
thanks to \eqref{ladispoli} the function $\FO$ is strictly increasing, resp. strictly decreasing, in the interval $ [\bar\sa_{2j-1},\bar\sa_{2j}],$
resp. $ [\bar\sa_{2j-2},\bar\sa_{2j-1}],$
then we can invert $\FO$ on the above intervals 
obtaining two functions
\begin{equation}\label{andria}
\bar\Sa_{2j} : [\bar E_{2j-1}, \bar E_{2j}]\to  [\bar\sa_{2j-1},\bar\sa_{2j}]
\qquad {\rm and}\qquad
\bar\Sa_{2j-1} : [\bar E_{2j-1}, \bar E_{2j-2}]\to [\bar\sa_{2j-2},\bar\sa_{2j-1}]
\end{equation}
such that
$$
\FO (\bar\Sa_i(E))=E\,,\qquad
\bar\Sa_i(\FO(\sa))=\sa
\,,\qquad \forall\, 1\leq i\leq 2N\,.
$$
Note that $\bar\Sa_i$ is increasing, resp. decreasing, if 
$i$ is even, resp. odd.
The functions $\bar\Sa_i$ have a holomorphic extension
as it is shown below.

\section{The perturbed potential and the analytic properties of its inverse}

\subsection{Perturbed potential}

\nl
Recalling \eqref{pasqua} and \eqref{straussbis},
we now consider a perturbation
$\Gm(\sa,\hat\act )$
 of the function $\FO(\sa)$
depending also on a parameter $\hat\act \in \hat D$
satisfying
\begin{equation}\label{ciccio}
  {\modulo}
  \Gm(\sa,\hat\act )-\FO(\sa)
  {\modulo}_{\hat D,r_0,s_0}
  \leq 
\checco\,,
\end{equation}
where $\checco$ is a small parameter that we 
assume to satisfy the condition
\begin{equation}\label{genesis}
\checco
\leq
\checco_\diamond
=\checco_\diamond(\morse,\b,s_0,r_0)
:=
\frac{\b^9  s_0^{15}}{ 2^{120}\morse^9}
\min
\left\{
r_0^2\,, \ \frac{r_0^3}{\sqrt\morse}\,,
\ 
\frac{\b^{45}s_0^{75}}{2^{321}\morse^{44}}
\right\}\,.
\end{equation}
 Note that by \eqref{ocarina} and \eqref{harlock}
\begin{equation}\label{badedas}
 s_*:=\min\{s_0,1\}\geq \max\left\{\frac{s_0 \b}{4\morse}\,, \frac{s_0^2 \b}{4\morse}
 \right\}\,.
\end{equation}
In particular, by \eqref{goffredo}
\begin{equation}\label{ciccio2}
  {\modulo}\Gm{\modulo}_{\hat D,r_0,s_0}
  \leq 
\morse+\checco\,.
\end{equation}
Moreover, by \eqref{ciccio} and Cauchy estimates
\begin{equation}\label{ciccio*}
  {\modulo}
  \partial_{\hat\act}\Gm(\sa,\hat\act )
  {\modulo}_{\hat D,r_0/2,s_0}
  \leq 
2\checco/r_0\,.
\end{equation}
By \eqref{ladispoli}, for $\checco\leq\checco_\diamond $ small enough, 
we can continue  the critical points $\bar\sa_i$ (defined in \eqref{recremisik}),
resp. critical energies $\bar E_i,$ of $\FO$ obtaining
critical points
$\sa_i(\hat\act ),$
resp. critical energies $E_i(\hat\act ),$
 of $\Gm(\cdot,\hat\act ),$
solving the implicit function equation
\begin{equation}\label{faciolata}
\partial_\sa \Gm(\sa_i(\hat\act ),\hat\act )=0
\end{equation}
 and then evaluating
 \begin{equation}\label{faciolata2}
\Gm(\sa_i(\hat\act ), \hat\act )=:E_i(\hat\act )\,,
\end{equation}
respectively.
Note that, by definition (recall \eqref{tracollo})
\begin{equation}\label{tracollo1}
E_0=E_{2N}\,.
\end{equation}
 Note also that
 $\sa_i(\hat\act ),$
and $E_i(\hat\act )$ are analytic functions of $\hat\act \in \hat D_{r_0}.$ More precisely we have the following

\begin{lemma}
 Assume that 
 \begin{equation}\label{peterparker}
\checco\leq\checco_\diamond \leq 
\frac{\b s_0^2}{16 }
\left(\frac{12 \morse}{\b s_0^2} +1\right)^{-1}
\end{equation}
(which is implied by \eqref{genesis}).
Then, for every $1\leq i\leq 2N,$ there exists
a holomorphic function
$\sa_i(\hat\act )$ with
\begin{equation}\label{octoberx}
{\modulo} \sa_i-\bar\sa_i{\modulo}_{\hat D, r_0}\leq 
2\checco/\b s_0\leq s_0/8\,,
\end{equation}
solving equation \eqref{faciolata}.
Moreover
\begin{equation}\label{octoberx*}
{\modulo} \partial_{\hat\act}\sa_i{\modulo}_{\hat D, r_0/2}\leq 
4\checco/\b s_0 r_0\,,
\end{equation}
\end{lemma}
\proof
We know that $\partial_\sa \FO(\bar\sa_i)=0 .$
Then, considering $\hat\act \in \hat D_{r_0}$
as a parameter, we want to find
$\chi_i=\chi_i(\hat\act ),$
with 
$$
|\chi_i|\leq \rho
:= 2\checco/\b s_0\leq s_0/2\,,
$$
solving the equation
\begin{equation}\label{sterco}
\partial_\sa \Gm(\hat\act ,\bar\sa_i+\chi_i)=0\,.
\end{equation}
Introducing the parameter
$\epsilon$ we want to solve
the equation
$$
\mathcal F(y,\epsilon)=
\mathcal F(y,\epsilon;\hat\act )=0\,,
\qquad
\text{where}
\quad
\mathcal F(y,\epsilon;\hat\act ):=
\partial_\sa\FO(\bar\sa_i+y)
+\epsilon
\partial_\sa G(\hat\act ,\bar\sa_i+y)\,,
$$
with $G:=\Gm-\FO,$
finding $y=y(\epsilon;\hat\act )$
for every $|\epsilon|\leq 1.$
Then $\chi_i=\chi_i(\hat\act ):=
y(1;\hat\act )
$
solves \eqref{sterco}
(recalling \eqref{ciccio}).
We are going to solve $\mathcal F=0$
by the Implicit  Function Theorem.
First note that
$$
\mathcal F(0,0)=\partial_\sa\FO(\bar\sa_i)
=0\,.
$$
Denote by $\g:=
1/\partial_y \mathcal F(0,0)
=
1/\partial_{\sa\sa}\FO(\bar\sa_i)$
and note that
by \eqref{ladispoli}
$$
|\g|\leq 1/\b\,.
$$ 
In order to apply 
a quantitative version of the 
Implicit  Function Theorem, we have to verify the following conditions:
$$
\sup_{|\epsilon|\leq 1}
|\mathcal F(0,\epsilon)|
\leq \frac{\rho}{2|\g|}
$$
and
$$
\sup_{|\epsilon|\leq 1\,,\,
|y|\leq \rho}
\left|
1-\g \partial_y \mathcal F(y,\epsilon)
\right|\leq\frac12\,.
$$
Since, by Cauchy estimates,
$$
\sup_{|\epsilon|\leq 1}
|\mathcal F(0,\epsilon)|
=
|\partial_\sa G(\hat\act ,\bar\sa_i)|
\leq\frac{\checco}{s_0}
\leq \frac{\rho}{2|\g|}\,,
$$
the first condition follows 
by  \eqref{ladispoli}.
Regarding 
the second condition 
we note that
$$
\partial_y \mathcal F(y,\epsilon)
=
\partial_{\sa\sa}\FO(\bar\sa_i+y)
+\epsilon
\partial_{\sa\sa}G(\hat\act ,\bar\sa_i+y)\,.
$$
Then, by Cauchy estimates
and noting that $|y|\leq \rho\leq s_0/2$,
we get
$$
\sup_{|\epsilon|\leq 1\,,\,
|y|\leq \rho}
\left|
1-\g \partial_y \mathcal F(y,\epsilon)
\right|\leq
|\g| \frac{6 \morse\rho}{(s_0/2)^3}
+|\g|\frac{2\checco}{(s_0/2)^2}
=\frac{8\checco |\g|}{s_0^2}
\left(\frac{12 \morse}{\b s_0^2} +1\right)
\leq
\frac12
$$
by \eqref{peterparker}.
\\
\eqref{octoberx*} follows \eqref{octoberx} and Cauchy estimates.
\eproof
Note that by \eqref{genesis} and  \eqref{harlock} we get
\begin{equation}\label{genesis2}
\checco_\diamond\leq \frac{\morse}{16} \,. 
\end{equation}
Then by \eqref{genesis2}
and \eqref{ciccio2}
\begin{equation}\label{ciccio3}
  {\modulo}\Gm{\modulo}_{\hat D,r_0,s_0}
  \leq 
2\morse\,.
\end{equation}
Recalling \eqref{faciolata2}, by \eqref{ciccio}, \eqref{octoberx} and Cauchy estimates
we get, for $\checco\leq \checco_\diamond,$
\begin{equation}\label{october}
 {\modulo}E_i-\bar E_i {\modulo}_{\hat D, r_0}
 \leq \left(\frac{4\morse}{\b s_0^2}+1\right) \checco
\stackrel{\eqref{harlock}}\leq 2 \checco \,.
\end{equation}
Then by Cauchy estimates
\begin{equation}\label{october*}
 {\modulo}\partial_{\hat\act}E_i{\modulo}_{\hat D, r_0/2}
 \leq 4 \checco/r_0 \,.
\end{equation}
By \eqref{alcafone},\eqref{draghetto},
\eqref{octoberx},
\eqref{genesis},
\eqref{october}, we note that 
$\sa_i(\hat\act )$ and 
$E_i(\hat\act )$ maintain the same order 
(w.r.t. $i$)
of $\bar\sa_i$ and $\bar E_i;$
moreover (recalling also \eqref{harlock})
\beqa{ladispoli3}
&& \inf_{\hat\act \in \hat D_{r_0}}
\min_{x\in\R}
\Big( |\partial_\sa \Gm(\sa,\hat\act )|+ 
|\partial_{\sa\sa}\Gm(\sa,\hat\act )|\Big)\geq \frac\b 2\,,\nonumber\\
&&\inf_{\hat\act \in \hat D_{r_0}}
\min_{i\neq j}|E_j(\hat\act )-E_i(\hat\act )|\geq \frac{\b }{2}\,.
\eeqa
Finally by \eqref{octoberx}, \eqref{alcafone}, \eqref{genesis}, \eqref{harlock}
and \eqref{draghetto}
we get
\begin{equation}\label{octobery}
\inf_{\hat\act \in \hat D_{r_0}}|\sa_i(\hat\act )-\sa_{i-1}(\hat\act )|\geq \xx_*
=\sqrt{\frac{\b s_0^3}{3 \morse}}\,.
\end{equation}
 
\begin{lemma}\label{anna}
 There exists a function $\check\Gm(\sa,\hat\act )$ holomorhic
 in $\{|\sa|<\sa_\diamond\}\times \hat D_{r_0}$ 
 where
 \begin{equation}\label{fusek}
\sa_\diamond:=\frac{\b s_0^3}{2^9 \morse}
\stackrel{\eqref{harlock}}\leq
\frac{s_0}{2^7}\,,
\end{equation}
 such that
 $$
 \partial_\sa \Gm\big(\sa_i(\hat\act )+\sa,\hat\act \big)
 =\sa \check \Gm(\sa,\hat\act )\,.
 $$
 Moreover $\check \Gm(0,\hat\act )=\partial_{\sa\sa}\Gm\big(\sa_i(\hat\act ),\hat\act \big)$ and
 $$
 \sup_{\{|\sa|<\sa_\diamond\}\times \hat D_{r_0}}
 \frac{1}{|\check \Gm|}
 \leq
 \frac{4}{\b}
 \,.
 $$
\end{lemma}
 \proof
For brevity we skip to write the immaterial dependence on
$\hat\act .$ By Taylor expansion we get
 $$
 \check \Gm(\sa)=
 \partial_{\sa\sa}\Gm(\sa_i)+\sa \int_0^1
 (1-t)  \partial_{\sa\sa\sa}\Gm(\sa_i+t \sa)dt\,.
 $$
 Note that the above expression is well posed since, by \eqref{octoberx},
 $|\Im \sa_i|\leq s_0/8$ and $|\sa|<\sa_\diamond\leq s_0/2^7.$
 Then, by Cauchy estimates and \eqref{ciccio3}, we also get
 $$
  \sup_{\{|\sa|<\sa_\diamond\}\times \hat D_{r_0}}
  |\partial_{\sa\sa\sa}\Gm(\sa_i+t \sa)|
  \leq
  \frac{2^7 \morse}{s_0^3}\,.
 $$
 Then by \eqref{ladispoli3} we have that, uniformly on
  $\{|\sa|<\sa_\diamond\}\times \hat D_{r_0}$,
 $$
 |\check \Gm |\geq \frac{\b}{2}-\sa_\diamond  \frac{2^7 \morse}{s_0^3}
 =\frac{\b}{4}\,.
 $$
 \eproof

\begin{lemma}\label{paradisecity}
For every real $\hat\act \in\hat D,$
there are no more critical points of
$\sa\to\Gm(\sa,\hat\act )$ than 
$\sa_1(\hat\act ),\ldots,\sa_{2N}(\hat\act );$ 
namely
if $\sa^\sharp\in(\sa_{2N}(\hat\act )-2\pi,\sa_{2N}(\hat\act )]$ satisfies
$\partial_\sa\Gm(\sa^\sharp,\hat\act )=0,$
then $\sa^\sharp=\sa_i(\hat\act )$ for some
$i=1,\ldots,2N.$
\end{lemma}
 \proof
We assume by contradiction
that there exists 
 $\sa^\sharp$ satisfying
$\partial_\sa\Gm(\sa^\sharp,\hat\act )=0,$
with $\sa_{i-1}(\hat\act )<\sa^\sharp<\sa_i(\hat\act )$ for some
$i=1,\ldots,2N.$
 By Lemma \ref{anna} we have that, for every
 $j=1,\ldots,2N,$
\begin{equation}\label{welcometothejungle}
|\sa^\sharp-\sa_j(\hat\act )|\geq \sa_\diamond\,.
\end{equation}
Since, by \eqref{ciccio} and Cauchy estimates
we have
$$
|\FO'(\sa^\sharp)|\leq \checco/s_0
\stackrel{\eqref{genesis}}\leq
\b/4\,,
$$
we get $|\FO''(\sa^\sharp)|\geq \b/2.$
Then there exist $\bar\sa^\sharp$ with
$|\bar\sa^\sharp-\sa^\sharp|\leq 2\checco/s_0\b$
and $\FO'(\bar\sa^\sharp)=0.$
This means that $\bar\sa^\sharp=\bar\sa_j$
for some $j=1,\ldots,2N.$
Then by \eqref{octoberx}
$$
|\sa^\sharp-\sa_j(\hat\act )|
\leq
|\sa^\sharp-\bar\sa^\sharp|+
|\bar\sa_j-\sa_j(\hat\act )|
\leq
4\checco/s_0\b\,,
$$
which by \eqref{genesis} (recall \eqref{fusek})
contradicts \eqref{welcometothejungle}.
 \eproof
 
\subsection{Rescaled potentials}

Now let us introduce the affine functions
(considering $\hat\act $ as a parameter)
\begin{eqnarray}
\bar\g_{i}(\breve\sa)
&:=&
\frac{\bar\sa_i -\bar\sa_{i-1}}{2} \breve\sa + 
\frac{\bar\sa_i +\bar\sa_{i-1}}{2}
\,,\qquad
\g_{i}(\breve\sa,\hat\act )
:=
\frac{\sa_{i}(\hat\act )-\sa_{i-1}(\hat\act )}{2} \breve\sa + 
\frac{\sa_{i}(\hat\act )+\sa_{i-1}(\hat\act )}{2}
\nonumber
\\
\bar\l_i (E)
&:=&
(-1)^i
\frac{2E - \bar E_i - \bar E_{i-1} }{\bar E_i -\bar E_{i-1} }
\,,\qquad
\l_{i}(E,\hat\act )
:=
(-1)^i
\frac{2E - E_{i}(\hat\act )- E_{i-1}(\hat\act )}{E_{i}(\hat\act )-E_{i-1}(\hat\act )}
\,.
\label{acrobat}
\end{eqnarray}

\begin{lemma}\label{carciofo}
For every $0\leq i\leq 2N$ 
 \begin{equation}\label{carciofino}
|\partial_\xx \bar  \g_i|\leq \pi\,,
\qquad
{\modulo}\partial_\xx \g_i{\modulo}_{\hat D,r_0}\leq 2\pi\,,
\qquad
|\partial_{E} \bar  \l_i|\leq \frac{2}{\b}\,,
\qquad
{\modulo}\partial_{E} \l_i{\modulo}_{\hat D,r_0}\leq 
\frac{4}{\b}\,.
\end{equation}
Moreover, for $\s>0,$ we have
\begin{equation}\label{carciofino2}
\sup_{[-1,1]_\s\times \hat D_{r_0}} 
|\g_i-\bar  \g_i |
\leq 
\frac{2\checco}{\b s_0}(2+\s) \,,
\qquad
\sup_{\{|E|\leq 2\morse\}\times \hat D_{r_0}} 
|\l_i-\bar \l_i |\leq \frac{48 \morse \checco}{\b^2}
\end{equation}
and
\begin{equation}\label{carciofino3}
\sup_{[-1,1]_{\s}}\left| \Im \bar\g_{i}\right|\,,\ 
\sup_{[-1,1]_{\s}\times \hat D_{r_0}}\left| \Im \g_{i}\right|
\leq 
\frac{2\checco}{\b s_0}(2+\s)
+\pi\s\,.
\end{equation}

\end{lemma}
\proof
The first estimate follows by \eqref{recremisik};
then the second one follows by  \eqref{octoberx}, \eqref{genesis}
and \eqref{harlock2}.
The third and fourth estimates follow by
\eqref{ladispoli} and 
\eqref{ladispoli3}, respectively.
\\
The first estimate in \eqref{carciofino2}
follows by \eqref{octoberx}
and noting that $|\xx|\leq 1+\s$
for $\xx\in [-1,1]_\s.$
Regarding the second estimate in \eqref{carciofino2} we first note that
by \eqref{ladispoli},\eqref{ladispoli3}
and \eqref{october} 
$$
\left|\frac{1}{E_{i}(\hat\act )-E_{i-1}(\hat\act )}-
\frac{1}{\bar E_i -\bar E_{i-1} }\right|\leq
\frac{4\checco}{\b^2}\,.
$$
Then, by \eqref{genesis2}, \eqref{october}, \eqref{ladispoli} and the first estimate in 
\eqref{harlock}, we get
\begin{eqnarray*}
|\l_i-\bar \l_i |
&\leq&
2\big(|E|+|E_i|+|E_{i-1}|\big)
\left|\frac{1}{E_{i}(\hat\act )-E_{i-1}(\hat\act )}-
\frac{1}{\bar E_i -\bar E_{i-1} }\right|\\
&&\phantom{AAAAAAAAAAA} +
\frac{|E_i-\bar E_i|+|E_{i-1}-\bar E_{i-1} |}{|\bar E_i -\bar E_{i-1} |}
\\
&\leq&
\frac{40\morse\checco}{\b^2}+\frac{4\checco}{\b}
\leq \frac{48 \morse \checco}{\b^2}\,.
\end{eqnarray*}
\\
Finally, since $\bar\g_{i}(\Re \xx)\in \R,$
we have  $\Im \big(\g_{i}(\xx,\hat\act )\big)=
\Im \big(\g_{i}(\xx,\hat\act )-\bar\g_{i}(\Re \xx) \big)$ and
then, for $(\xx,\hat\act )\in [-1,1]_{\s}\times \hat D_{r_0},$
we get
$$
\left| \Im \big(\g_{i}(\xx,\hat\act )\big)\right|
\leq
\big|\g_{i}(\xx,\hat\act )-\bar\g_{i}(\Re \xx) \big|
\leq
\big|\g_{i}(\xx,\hat\act )-
\bar\g_{i}(\xx) \big|+
\big|\bar\g_{i}( \xx) -
\bar\g_{i}(\Re \xx) \big|
\,.
$$
Then \eqref{carciofino3}
follows by 
\eqref{carciofino} and \eqref{carciofino2}.
\eproof

We also have
\begin{eqnarray}\label{acrobat2}
\bar\l_i^{-1}(\breve E)
&=&
\frac12\left((-1)^i(\bar E_{i}-\bar E_{i-1}) \breve E
+\bar E_{i}+\bar E_{i-1}\right)\,,
\\
\l_{i}^{-1}(\breve E,\hat\act )
&=&
\frac12\left((-1)^i\big(E_{i}(\hat\act )-E_{i-1}(\hat\act )\big) \breve E
+E_{i}(\hat\act )+E_{i-1}(\hat\act )\right)\,,
\nonumber
\\
\l_{i}^{-1}(\bar\l_i(E),\hat\act )
&=&
\frac12\left(\big(E_{i}(\hat\act )-E_{i-1}(\hat\act )\big) 
\frac{2E - \bar E_i - \bar E_{i-1} }{\bar E_i -\bar E_{i-1} }
+E_{i}(\hat\act )+E_{i-1}(\hat\act )\right)
\nonumber
\\
&=&
\frac{E_{i}(\hat\act )-E_{i-1}(\hat\act )}{\bar E_i -\bar E_{i-1}} E
+ 
\frac{E_{i-1}(\hat\act )\bar E_i-E_{i}(\hat\act )\bar E_{i-1}}{\bar E_i -\bar E_{i-1} }
\,.
\nonumber
\end{eqnarray}
We  have
\begin{equation}\label{inverno}
|\l_{i}^{-1}(\bar\l_i(E),\hat\act )-E|
\leq 
\frac{4\checco}{\b}|E|+\frac{4\checco\morse}{\b}
\,,\qquad
\forall E\in\C\,, \ \hat\act \in \hat D_{r_0}\,,
\end{equation}
by \eqref{october},\eqref{goffredo} and \eqref{ladispoli}.

Let us introduced the ''rescaled''  functions
$\breve \FO_i(\xx)$ and
$\breve{\Gm}_i(\xx,\hat\act )$  by 
\begin{equation}\label{bastard}
\breve{\FO}_i:= \bar \l_i \circ \FO\circ \bar  \g_i \,,\qquad
\breve \Gm_i:= \l_i\circ \Gm\circ \g_i
\end{equation}
(recall \eqref{acrobat}).
Recalling \eqref{faciolata} and \eqref{faciolata2}
 For every real $\hat\act $, these functions are bijective from $[-1,1]$ to
 $[-1,1]$; in particular
 \begin{equation}\label{megadirettore}
 \breve{\FO}_i(\pm 1)=\breve{\Gm}_i(\pm 1,\hat\act )=\pm(-1)^i\,,
 \qquad\qquad
\partial_\xx \breve{\FO}_i(\pm 1)=
\partial_\xx\breve{\Gm}_i(\pm 1,\hat\act )=0 
 \,.
\end{equation}
 We also set
 \begin{equation}\label{megadirettore2}
\breve{\Gm}^*_i(\xx,\hat\act )=
\breve{\Gm}_i(\xx,\hat\act )
-\breve{\FO}_i(\xx)\,.
\end{equation}
Note that  by \eqref{megadirettore}
we have 
 \begin{equation}\label{megadirettore3}
\breve{\Gm}^*_i(\pm 1,\hat\act )=
\partial_\sa \breve{\Gm}^*_i(\pm 1,\hat\act )=0\,.
\end{equation}

\begin{lemma}
Let $0\leq i\leq 2N$ and $\checco\leq\checco_\diamond.$ Then
 $\breve{\FO}_i$, resp. $\breve{\Gm}_i,$
has analytic extension on
 $[-1,1]_{\breve s},$ resp.  $[-1,1]_{\breve s}\times \hat D_{r_0},$ where
\begin{equation}\label{exit}
 \breve s:=s_0/2\pi
\end{equation}
with
 \begin{equation}\label{gibbone}
\sup_{[-1,1]_{\breve s}}| \breve{\FO}_i|\,,\ 
\ \sup_{[-1,1]_{\breve s}\times \hat D_{r_0}}
| \breve{\Gm}_i|\leq 
\frac{8\morse}{\b}=:\breve \morse
\,.
\end{equation}
Moreover
\begin{equation}\label{borges}
\sup_{[-1,1]_{\breve s}\times \hat D_{r_0}}
| \breve{\Gm}^*_i|=
\sup_{[-1,1]_{\breve s}\times \hat D_{r_0}}
| \breve{\Gm}_i-\breve{\FO}_i|\leq 
\frac{70\morse}{\b^2 s_*^2}\checco
=:\breve \checco
\end{equation}
and
  \begin{equation}\label{ladispoli4}
\min_{\xx\in[-1-\breve s,1+\breve s]}
\Big( |\partial_{\xx}
\breve{\FO}_i(\xx)|+ 
|\partial_{\xx\xx}
\breve{\FO}_i(\xx)|\Big)\geq 
\frac{\b^2 s_0^3}{3\pi \morse^2}
=:\breve \b\,,
\end{equation}
  \begin{equation}\label{ladispoli5}
  \inf_{\hat\act \in \hat D_{r_0}}
\min_{\xx\in[-1-\breve s,1+\breve s]}
\Big( |\partial_{\xx}
\breve{\Gm}_i(\xx,\hat\act )|+ 
|\partial_{\xx\xx}
\breve{\Gm}_i(\xx,\hat\act )|\Big)\geq 
\breve \b/4\,.
\end{equation}
\end{lemma}
\proof
By \eqref{carciofino3} we get
\begin{equation}\label{suppli}
\sup_{[-1,1]_{\breve s}}\left| \Im \big(\bar\g_{i}(\xx)\big)\right|\,,\ 
\sup_{[-1,1]_{\breve s}\times \hat D_{r_0}}\left| \Im \big(\g_{i}(\xx,\hat\act )\big)\right|
\leq 
\frac{2\checco}{\b s_0}(2+\breve s)
+\pi\breve s
\stackrel{\eqref{harlock},\eqref{harlock2}}
\leq \frac{s_0}{2}\,,
\end{equation}
then (recall \eqref{ciccio}) $\breve{\Gm}_i$ in \eqref{bastard}
is well defined and holomorphic
in   $[-1,1]_{\breve s}\times \hat D_{r_0}.$
The case of $\breve{\FO}_i$ is similar.
 \\
\eqref{gibbone} follows by the definition of $\bar \l_i ,\l_i$
in \eqref{acrobat},
by
\eqref{ladispoli},
\eqref{ladispoli3}
and since here
$|E|,|\bar E_i|,|E_i(\hat\act )|
\leq \morse+\checco$
(recall \eqref{goffredo} and\eqref{ciccio2}).
\\
Moreover by \eqref{carciofino2}, \eqref{carciofino}, \eqref{ciccio}
\begin{eqnarray*}
| \breve{\Gm}_i-\breve{\FO}_i|
&\leq&
|(\l_i-\bar \l_i )\circ \Gm\circ\g_i|+
|\bar \l_i \circ \Gm\circ\g_i-\bar \l_i \circ \FO\circ\bar  \g_i  |
\\
&\leq&
\frac{48 \morse \checco}{\b^2}
+\frac2\b | \Gm\circ\g_i- \FO\circ\bar  \g_i  |
\\
&\leq&
\frac{48 \morse \checco}{\b^2}
+\frac2\b\Big(
| (\Gm-\FO)\circ\g_i|+
 | \FO\circ\g_i- \FO\circ\bar  \g_i  |
\Big)
\\
&\stackrel{\eqref{suppli}}\leq&
\frac{48 \morse \checco}{\b^2}
+\frac{2\checco}{\b}
+\frac2\b {\modulo}\partial_\sa \FO{\modulo}_{s_0/2}
\sup_{[-1,1]_\s\times \hat D_{r_0}} 
|\g_i-\bar  \g_i |
\\
&\stackrel{\eqref{carciofino2}}\leq&
\frac{48 \morse \checco}{\b^2}
+\frac{2\checco}{\b}
+
\frac{8\morse\checco}{\b^2 s_0^2}(2+\frac{s_0}{2\pi})
\\
&\stackrel{\eqref{harlock}}\leq&
\frac{70\morse\checco}{\b^2 s_*^2}\,,
\end{eqnarray*}
recalling the definition of $s_*$ in \eqref{badedas}.
\\
Noting that
$$
\partial_\xx \breve{\FO}_i=
(-1)^i
\frac{\bar\sa_i -\bar\sa_{i-1}}{\bar E_i -\bar E_{i-1} } 
\partial_\sa \FO\circ\bar  \g_i \,,
\qquad
\partial_{\xx\xx} \breve{\FO}_i=
(-1)^i
\frac{(\bar\sa_i -\bar\sa_{i-1})^2}{2(\bar E_i -\bar E_{i-1} )} 
\partial_{\sa\sa} \FO\circ\bar  \g_i \,,
$$
by \eqref{ladispoli} and \eqref{alcafone}
\begin{eqnarray*}
|\partial_\xx \breve{\FO}_i|+
|\partial_{\xx\xx} \breve{\FO}_i|
&\geq&
\frac{\xx_*}{\morse}|\partial_\sa \FO\circ \bar  \g_i |
+
\frac{\xx_*^2}{\morse}|\partial_{\sa\sa} \FO\circ \bar  \g_i |
\geq \frac{\xx_*^2}{\pi \morse}\big(
|\partial_\sa \FO\circ \bar  \g_i |+
|\partial_{\sa\sa} \FO\circ \bar  \g_i |
	\big)
	\\
	&\geq&
\frac{\xx_*^2\b}{\pi \morse}	=
\frac{\b^2 s_0^3}{3\pi \morse^2}\,.
\end{eqnarray*}
Similarly
by
$$
\partial_\xx \breve{\Gm}_i=
(-1)^i
\frac{\sa_{i}-\sa_{i-1}}{E_{i}-E_{i-1}} 
\partial_\sa \Gm\circ\g_i\,,
\qquad
\partial_{\xx\xx} \breve{\Gm}_i=
(-1)^i
\frac{(\sa_{i}-\sa_{i-1})^2}{2(E_{i}-E_{i-1})} 
\partial_{\sa\sa} \Gm\circ\g_i\,,
$$
and by \eqref{ladispoli3}, \eqref{octobery}
we get
\begin{eqnarray*}
|\partial_\xx \breve{\Gm}_i|+
|\partial_{\xx\xx} \breve{\Gm}_i|
&\geq&
\frac{\xx_*}{2\morse}|\partial_\sa \Gm\circ \g_i|
+
\frac{\xx_*^2}{4\morse}|\partial_{\sa\sa} \Gm\circ \g_i|\\
&\geq& 
\frac{\xx_*^2}{2\pi \morse}\big(
|\partial_\sa \Gm\circ \g_i|+
|\partial_{\sa\sa} \Gm\circ \g_i|
	\big)
	\\
	&\geq&
\frac{\xx_*^2\b}{4\pi \morse}	=
\frac{\b^2 s_0^3}{12\pi \morse^2}\,. \qedeq\end{eqnarray*}
Note that by \eqref{harlock} and \eqref{harlock2}
we get
\begin{equation}\label{harlock5}
\breve \morse\geq 4\,,\qquad
\breve \b\leq\frac{8}{3\pi}<1
\end{equation}
and
\begin{equation}\label{harlock4}
\frac{12\pi\breve\b}{\breve\morse}
\leq 1
\,,\quad
\frac{3\pi^2\breve\b\breve s}{\breve\morse}
\leq 1
\,,\quad
\frac{3\pi^3\breve\b\breve s^2}{\breve\morse}
\leq 1
\,,\quad
\frac{3\pi^4\breve\b\breve s^3}{\breve\morse}
\leq 1
\,.
\end{equation}

\begin{lemma}\label{parcoavventurabar}
For every $1\leq i\leq 2N,$ 
\begin{eqnarray}
&&\label{harlock1bar}
|\partial_\xx\breve{\FO}_i(\pm 1+\xx)|\geq \frac{\breve\b}{8}|\xx|\,,
\qquad\qquad \forall\, |\xx|\leq \breve\xx_\sharp
:=\frac{\breve\b \breve s^3}{ 2^9\breve\morse}
<\min \Big\{\frac{\breve s}{8},
\frac{1}{2^{16}}
\Big\}
\,,
\\
&&\label{harlock7bar}
\min_{[-1+ \breve \xx_0,1-\breve \xx_0]}
|\partial_\xx\breve\FO_i|
\geq 
\frac{\breve\b \breve \xx_0}{8}
\,,
\qquad\qquad
\forall\, 0\leq \breve \xx_0\leq \breve\xx_\sharp
\\
&&\label{harlock8bar}
\inf_{[-1+\breve \xx_0,
1-\breve \xx_0]_{\breve \xx_1}}
|\partial_\xx\breve\FO_i|
\geq 
\frac{\breve\b \breve \xx_0}{16}
\,,
\qquad
\breve \xx_1:=\frac{\breve\b \breve s^2\breve \xx_0}
{2^{7} \breve \morse}
<\frac{1}{2^{11}}\breve \xx_0
\,.
\end{eqnarray}
In particular
\begin{equation}\label{pajata}
\frac{1}{|\partial_\xx\breve\FO_i(\breve\sa)|}
\leq
\frac{8}{\breve\b}\left(
\frac{2}{\breve\xx_\sharp}
+
\frac{1}{|1-\breve\sa|}
+
\frac{1}{|1+\breve\sa|}
\right)\,,
\quad
\forall\, \breve\sa\in [-1,1]_{\breve\xx_\star}\,,
\quad \breve\xx_\star:=
\frac{\breve\b^2 \breve s^5}
{2^{16} \breve \morse^2}
<\frac{1}{2^{11}}\breve\xx_\sharp
\,.
\end{equation}
Finally
\begin{equation}\label{borges2}
|\partial_\xx \breve{\Gm}^*_i(\xx,\hat\act)|
=
|\partial_\xx  \breve{\Gm}_i (\xx,\hat\act)
-
\partial_\xx \breve{\FO}_i(\xx)
|
\leq 
\frac{2^{11}\breve\morse}{\breve\b^2\breve s^4}
\breve \checco |\partial_\xx \breve{\FO}_i(\xx)|\,,
\qquad\forall\, \xx\in [-1,1]_{\breve \xx_\star}\,,\ \ 
\hat\act\in \hat D_{r_0}
\,.
\end{equation}
\end{lemma}
\proof
We will consider only the case\footnote{The case
$i=2j-1$ is analogous.} $i=2j$.
For\footnote{The case $-\xx_\sharp\leq \xx\leq 0$ is analogous.} 
$|\xx|\leq \xx_{\breve\sharp}<\breve s/2$
 we get by Cauchy estimates
$$
|\partial_\xx\breve{\FO}(\pm 1+\xx)|
\geq 
|\partial_{\xx\xx}\breve{\FO}(\pm 1)||\xx|-\frac{48\breve\morse}{\breve s^3}|\xx|^2
\stackrel{\eqref{ladispoli5}}
\geq
\frac{\breve\b}{4} |\xx| -\frac{48\breve\morse}{\breve s^3}|\xx|^2
\geq 
\frac{\breve\b}{8}|\xx|\,.
$$
Regarding the minimum of $\partial_\xx\breve{\FO}$
in the interval\footnote{Note that
on such interval $\partial_\xx\breve{\FO}>0$.}
$[-1+\breve \xx_0,1-\breve \xx_0]$
 if it is achieved in an inner point $\breve\sa_*$, then 
$\partial_{\xx\xx}\breve{\FO}(\breve\sa_*)=0$ and, by \eqref{ladispoli4},
$
\partial_\xx\breve{\FO}(\breve\sa_*)
\geq 
\breve\b;
$
otherwise
the minimum is achieved at the endpoints and 
\eqref{harlock7bar} follows from 
\eqref{harlock1bar}.
\\
By \eqref{harlock7bar}, \eqref{gibbone}
and Cauchy estimates
we have that for every
$\xx\in[-1+\breve \xx_0,1-\breve \xx_0]$
and $|\tilde \xx|\leq \breve \xx_1$
$$
|\partial_\xx\breve\FO_i(\xx+\tilde\xx)|
\geq 
\frac{\breve\b \breve \xx_0}{8}-\frac{8\breve\morse}{\breve s^2}\breve \xx_1
\geq 
\frac{\breve\b \breve \xx_0}{16}
\,,
$$
showing  \eqref{harlock8bar}.
\\
\eqref{pajata} directly follows from
\eqref{harlock1bar} and \eqref{harlock8bar}
taking $\breve \xx_0=\breve \xx_\sharp,$
namely defining
$$
\breve\xx_\star:=
\frac{\breve\b \breve s^2\breve \xx_\sharp}
{2^{7} \breve \morse}
=
\frac{\breve\b^2 \breve s^5}
{2^{16} \breve \morse^2}
\stackrel{\eqref{harlock4}}<\frac{1}{2^{11}}\breve\xx_\sharp
\,.
$$
\\
Let us finally prove \eqref{borges2}.
By \eqref{harlock8bar} it follows that
\begin{equation}
\label{harlock5bar}
\inf_{[-1+\breve\xx_\sharp/2,
1-\breve\xx_\sharp/2]_{\breve\xx_\star}}
|\partial_\xx\breve\FO_i|\geq \frac{\breve\b^2 \breve s^3}{2^{10}\breve\morse}
\,.
\end{equation}
Recalling \eqref{megadirettore3} by Cauchy estimates
and \eqref{borges}, \eqref{harlock1bar} we get
$$
|\partial_\xx \breve{\Gm}^*_i(\pm 1+\xx,\hat\act)|
\leq
\frac{8\breve\checco}{\breve s^2}|\xx|
\leq 
\frac{64\breve\checco}{\breve\b \breve s^2}
|\partial_\xx\breve{\FO}_i(\pm 1+\xx)|\,,
\qquad
\forall\, |\xx|\leq \breve\xx_\sharp\,,\ \ 
\hat\act \in \hat D_{r_0}\,.
$$
For every 
$\xx\in [-1+\breve\xx_\sharp/2,
1-\breve\xx_\sharp/2]_{\breve\xx_\star}$
and 
$\hat\act \in \hat D_{r_0}$
by Cauchy estimates and \eqref{borges}
we get
$$
|\partial_\xx \breve{\Gm}^*_i(\xx,\hat\act)|
\leq \frac{2\breve\checco}{\breve s}
\stackrel{\eqref{harlock5bar}}
\leq 
\frac{2^{11}\breve\morse}{\breve\b^2\breve s^4}
\breve \checco |\partial_\xx \breve{\FO}_i(\xx)|
$$
proving \eqref{borges2}.
\eproof

\subsection{Inverting the rescaled unperturbed potential}
 Define, for $\rho>0,$
 the complex sets\footnote{
 In particular 
 $(\C_*-1)\cap (1-\C_*)
 =
 \{ z\in\C \ | \ \Im z=0 \implies -1<\Re z<1 \}.
 $}
\begin{equation}\label{bufala}
\C_*:=\{ z\in\C \ | \ \Im z=0 \implies \Re z>0 \}\,,
\qquad
\O_\rho :=[-1,1]_\rho \cap 
 (\C_*-1)\cap (1-\C_*)\,.
\end{equation}
We define the square root and the (natural) logarithm
on $\C_*$ in order to have  
$\sqrt{1}=1$ and $\ln 1=0.$

We want to invert the function $\breve \FO_{i}$ solving
\begin{equation}\label{hovistolaluce}
\breve \FO_{i}(\sa)=\breve E\,.
\end{equation}
 
\begin{lemma}\label{cerveza}
 Let
 \begin{equation}\label{nera}
\breve r_0:= \frac{\breve \b^4 \breve s^6}{2^{33}\breve \morse^3}
\stackrel{\eqref{exit}-\eqref{ladispoli4}}=
\frac{s_0^{18}}{2^{48}3^4\pi^{10}}
\left(\frac{\b}{\morse}\right)^{11}
\stackrel{\eqref{harlock}}<\frac{1}{2^{44}}\,.
\end{equation}
There exists a holomorphic function
$
\breve{\bar\Sa}_i $ defined on $\O_{\breve r_0}$
such that  $x=\breve{\bar\Sa}_i(\breve E)$ solves 
\eqref{hovistolaluce}, namely
\begin{equation}\label{blu}
\breve \FO_i\big(\breve{\bar\Sa}_i(\breve E)\big)=\breve E
\end{equation}
and 
\begin{equation}\label{amber}
\sup_{\O_{\breve r_0}}
|\breve{\bar\Sa}_i|\leq 2\,.
\end{equation}
Moreover there exist two holomorphic functions 
$\breve{\bar\Sa}_{i,+},\breve{\bar\Sa}_{i,-},$ defined 
on\footnote{Actually they are defined on the larger ball of radius $r_*/2$, with $r_*$ defined in \eqref{acero}.} 
$B_{\sqrt{\breve r_0}}(0),$ with 
\begin{equation}\label{siracide}
\sup_{B_{\sqrt{\breve r_0}}(0)}|\breve{\bar\Sa}_{i,\pm}|\leq \sqrt{\frac{2}{\breve \b}}
\stackrel{\eqref{ladispoli4}}=
\sqrt{\frac{6\pi}{s_0^3}}\frac{\morse}{\b}\,,
\qquad
\inf_{B_{\sqrt{\breve r_0}/4}(0)}|\breve{\bar\Sa}_{i,\pm}|
\geq
\frac{\breve s}{2\sqrt{\breve\morse}}
\end{equation}
and
\begin{equation}\label{siracide2}
\breve{\bar\Sa}_{i,\pm}(0)
=
\sqrt{\frac{2}{\mp (-1)^i \partial_{\xx\xx}
\breve \FO_i(\pm1)}}
\geq \frac{\breve s}{\sqrt{\breve \morse}}
\,,
\end{equation}
such that
\begin{equation}\label{giappone}
\breve{\bar\Sa}_i(\breve E) 
=(-1)^i\left(
\pm 1\mp\sqrt{1\mp \breve E}\  \breve{\bar\Sa}_{i,\pm}( \sqrt{1\mp \breve E})
\right)\,,\quad
\text{ when} \ \ \ 
1\mp \breve E \in \C_*\cap B_{\breve r_0}(0)
\,.
\end{equation}
Moreover, on the real,
 $\breve{\bar\Sa}_i$ is bijective from
 $[-1,1]$ to itself and is strictly increasing, resp. 
 decreasing, if $i$ is even, resp. odd;
 in particular
 $\breve{\bar\Sa}_{2j}(\mp 1)=
 \breve{\bar\Sa}_{2j-1}(\pm 1)=\mp 1.$
\\
Finally 
\begin{equation}\label{ventresca}
|\breve{\bar\Sa}_i(\breve E)\pm
(-1)^i|
\geq 
\frac{\breve s}{64}
\sqrt{\frac{\breve r_0}{\breve \morse}}
\sqrt{|\breve E\pm 1|}\,,\qquad\quad
\forall\, \breve E\in \O_{\breve r_1}\,,
\quad \breve r_1:=\frac{\breve s}{2^{10}}
\sqrt{\frac{\breve r_0^3}{\breve \morse}}
\end{equation}
and, for $\rho\leq \breve r_0/2,$
\begin{equation}\label{ragnetto}
\breve{\bar\Sa}_i(\O_\rho)
\, \subseteq\,
[-1,1]_{s(\rho)}\,,\qquad
\text{where}\qquad
s(\rho):=\max\left\{
\frac{4\rho^{1/3}}{\sqrt{\breve\b}},\,
\frac{4\rho}{\breve r_0}
\right\}\,.
\end{equation}

\end{lemma}

 \proof
 We consider only the case
 $i=2j$, the case $i=2j-1$ being analogous.
 We start  inverting  equation  \eqref{hovistolaluce}
for $x$ and $\breve E$ close to $-1$.
 Set
$$
\breve {\bar c}_{2j}:=\partial_{\psi_n \psi_n} 
\breve \FO_{2j}(-1)/2>0\,.
$$
Note that by Cauchy estimates and \eqref{ladispoli4}
\begin{equation}\label{intervallo}
\frac{\breve \b}{2}\leq \breve {\bar c}_{2j}\leq \frac{\breve \morse}{\breve s^2}
\,.
\end{equation}

We have
\begin{equation}\label{panontina0}
\breve \FO_{2j}(-1+\xx)-
\breve \FO_{2j}(-1)
\stackrel{\eqref{megadirettore}}=
\breve \FO_{2j}(-1+\xx)+1
=:
\xx^2 \hat \FO_{2j}(\xx)\,,
\end{equation}
for a suitable  function
  $\hat \FO_{2j}$ analytic
  for  
  $
  |\xx|< \breve s\,.
  $
 By definition
$$
\hat \FO_{2j}(\xx)=
\frac{\breve \FO_{2j}(-1+\xx)+1}{\xx^2}
$$
and, by Cauchy estimates,
\begin{equation}\label{intervallo3}
|\hat \FO_{2j}(\xx) - \breve {\bar c}_{2j}|\leq 
\frac{8\breve \morse}{\breve s^3} |\xx|
\,,\qquad
\forall\, |\xx|\leq \breve s/2\,,
\end{equation}
then, recalling \eqref{intervallo},
\begin{equation}\label{intervallo2}
\sup_{|\xx|\leq \breve s/2}|\hat \FO_{2j}(\xx)|\leq
\frac{5\breve \morse}{\breve s^2}\,.
\end{equation}
Moreover
\begin{equation}\label{locanda}
\hat \FO_{2j}(0)=\breve {\bar c}_{2j}\,.
\end{equation}
By \eqref{panontina0},
\eqref{hovistolaluce} is equivalent, in the variable $\xx=x+1,$ to
\begin{equation}\label{levito}
\xx^2 \hat \FO_{2j}(\xx)=
\breve \FO_{2j}(-1+\xx)+1
=\breve E+1\,.
\end{equation}
We define the square root $\sqrt \cdot$  on 
$\C_*$ defined in \eqref{bufala} such that it coincide with the
positive suqre root on the positive reals, namely
if $z=r e^{\ii \theta},$ $r>0,$ $-\pi<\theta<\pi,$ then
$\sqrt z:=\sqrt r e^{\ii \theta/2},$ so that $\Re \sqrt z>0.$
Thus for $\breve E+1\in \C_*$ we can define $\sqrt{\breve E+1}.$
Set
\begin{equation}\label{sanpietrino}
\breve \rho:=\frac{\breve \b \breve s^3}{2^7\breve \morse}
\end{equation}
and note that, by \eqref{harlock4},
\begin{equation}\label{sanpietrino2}
\breve \rho\leq \min\left\{
\frac{\breve s}{3 \pi^3 2^7}\,,\ 
\frac{1}{3 \pi^4 2^7}
\right\}\,.
\end{equation}
Then
\beqa{avoja}
&& |\xx|\leq \breve \rho
\quad\Longrightarrow\quad
|\hat \FO_{2j}(\xx)-\breve {\bar c}_{2j}|\leq 
\frac{8\breve \morse}{\breve s^3} |\xx|
\stackrel{\eqref{intervallo}}\leq
\frac{\breve\b}{2^4}\leq \frac18 \breve {\bar c}_{2j}
\nonumber\\
&&\phantom{AAAAA}
\Longrightarrow\quad
\Re \hat \FO_{2j}(\xx)\geq \frac12 \breve {\bar c}_{2j}
\quad\Longrightarrow\quad
\hat \FO_{2j}(\xx) \in \C_* \,.
\eeqa
Thus for every $|\xx|\leq \breve\rho$ we can define the holomorphic function $\sqrt{\hat \FO_{2j}(\xx)}$.
Let us consider now the equation
\begin{equation}\label{fanale}
\mathtt f(\xx):=\xx \sqrt{\hat \FO_{2j}(\xx)}= w\,.
\end{equation}
Setting 
$$
T:=1/\partial_\xx \mathtt f(0)=
1/\sqrt{\hat \FO_{2j}(0)}=1/\sqrt{\breve {\bar c}_{2j}}\,.
$$
recalling \eqref{locanda}.
If the smallness condition
\begin{equation}\label{giggetto}
\sup_{|\xx|\leq \breve \rho}
|1-T \partial_\xx\mathtt f(\xx)|
=
\sup_{|\xx|\leq \breve \rho}
\left|
1-\frac{1}{\sqrt{\hat \FO_{2j}(0)}}
\left(
\sqrt{\hat \FO_{2j}(\xx)}+
\frac{\xx \partial_\xx \hat \FO_{2j}(\xx)}{2\sqrt{\hat \FO_{2j}(\xx)}}
\right)
\right|
\leq \frac12
\end{equation}
is satisfied, then by a quantitative 
Inverse Function Theorem
there exists an analytic  function $g(w)$
defined for 
\begin{equation}\label{acero}
|w|\leq r_*:=\frac{\breve \rho}{2|T|}=\frac12 \breve \rho \sqrt{\breve {\bar c}_{2j}}\,.
\end{equation}
such that $g(0)=0$ and $\xx=g(w)$ satisfies equation 
\eqref{fanale}.
In order to prove \eqref{giggetto}
we first note that
 by  \eqref{intervallo} and \eqref{avoja}
we have
\begin{equation}\label{cornacchia}
|\hat \FO_{2j}(\xx)| \geq 
\frac{\breve \b}{4}
\qquad {\rm and}\qquad
|\sqrt{\hat \FO_{2j}(\xx)}-\sqrt{\hat \FO_{2j}(0)}|
\leq \frac{1}{\sqrt{\breve \b}}
|\hat \FO_{2j}(\xx)-\hat \FO_{2j}(0)|
\leq \frac{\sqrt{\breve \b}}{2^4}\,,
\qquad
\forall\, |\xx|\leq \breve \rho\,.
\end{equation}
Thus we get, recalling \eqref{intervallo},
$$
\left|1-\frac{\sqrt{\hat \FO_{2j}(\xx)}}{\sqrt{\hat \FO_{2j}(0)}}
\right|
=\frac{|\sqrt{\hat \FO_{2j}(\xx)}-\sqrt{\hat \FO_{2j}(0)}|}
{\sqrt{\breve {\bar c}_{2j}}}
\leq
\frac{1}{2^3\sqrt 2}
$$
and by Cauchy estimates 
\eqref{intervallo2}, \eqref{sanpietrino2}, \eqref{avoja}
and \eqref{sanpietrino}
$$
\left|
\frac{\xx \partial_\xx \hat \FO_{2j}(\xx)}{2
\sqrt{\breve {\bar c}_{2j}}\sqrt{\hat \FO_{2j}(\xx)}}
\right|
\leq \frac{20\sqrt 2 \breve \rho \breve \morse}{\breve\b \breve s^3}
=
\frac{5\sqrt 2}{2^5}\,,
$$
thus condition \eqref{giggetto} is satisfied.
Then we find 
an analytic  function $g(w)$
defined for 
$
|w|\leq r_*
$
with
\begin{equation}\label{briciola}
\sup_{|w|\leq r_*} |g(w)|\leq \breve \rho\,,
\end{equation}
such that $\xx=g(w)$ solves \eqref{fanale}.
Then we can define
\begin{equation}\label{baux}
\xx(E):=g(\sqrt{\breve E+1})\,,
\qquad \forall\, |\breve E+1|\leq r_*^2\,,\ \ \ 
\breve E+1\in \C_*\,,
\end{equation}
solving
$$
\xx(E) \sqrt{\hat \FO_{2j}(\xx(E))}=\sqrt{\breve E+1} 
$$
and, squaring,
$$
\xx^2(E) \hat \FO_{2j}(\xx(E))=\breve E+1
$$
We note that we can write 
$$
g(w)=w\tilde g(w) 
$$
for a suitable analytic function $\tilde g$
defined for
$
|w| \leq r_*/2
$
with
\begin{equation}\label{briciola2}
\sup_{|w|\leq r_*/2} |\tilde g(w)|\leq 
\frac{\breve \rho}{2r_*}\,.
\end{equation}
Indeed, by \eqref{briciola} and Cauchy estimate, 
$$
\sup_{|w|\leq r_*/2} |g'(w)|\leq 
\frac{\breve \rho}{2r_*}
$$
and, recalling that $g(0)=0,$
we get
$$
|g(w)|\leq \frac{\breve \rho}{2r_*} |w|\,,
\qquad \forall\, |w|\leq r_*/2\,,
$$
proving \eqref{briciola2}.

\nl
Moreover we claim that
\begin{equation}\label{ciaccona}
g([0,r_*/8])\supseteq [0,\breve \rho/24]\,.
\end{equation}
Indeed we first prove that
\begin{equation}\label{crotone}
\tilde g(r_*/8)\geq 
\frac{2}{3\sqrt{\breve {\bar c}_{2j}}}\,,
\end{equation}
from which we get
$$
g(r_*/8)\geq 
\frac{ r_*}{12\sqrt{\breve {\bar c}_{2j}}}\
\stackrel{\eqref{acero}}=\frac{\breve \rho}{24}
$$
and \eqref{ciaccona} follows.
Let us prove \eqref{crotone}.
By \eqref{briciola2} and Cauchy estimates
we get
$$
\sup_{|w|\leq r_*/8}|\tilde g'(w)|
\leq \frac{4\breve \rho}{3 r_*^2}\,,
$$
then, noting that 
$$
\tilde g(0)=g'(0)=
1/\sqrt{\hat \FO_{2j}(0)}=1/\sqrt{\breve {\bar c}_{2j}}\,,
$$
we get
$$
\tilde g(r_*/8)
\geq
\frac{1}{\sqrt{\breve {\bar c}_{2j}}}
-
\sup_{|w|\leq r_*/8}|\tilde g'(w)| \frac{r_*}{8}
=
\frac{1}{\sqrt{\breve {\bar c}_{2j}}}
-\frac{\breve \rho}{6 r_*}
\stackrel{\eqref{acero}}=
\frac{2}{3\sqrt{\breve {\bar c}_{2j}}}\,,
$$
proving \eqref{crotone} (and \eqref{ciaccona}).

\nl
Choosing  
$\breve{\bar\Sa}_{2j,-}:=\tilde g$ and recalling
\eqref{acero} and \eqref{intervallo}
 we get the first estimate in \eqref{siracide} in the case $E$ close to -1,
 the case close to +1 is analogous.
 The second estimate in 
 \eqref{siracide} follows from 
\eqref{siracide2}, \eqref{amber}, \eqref{nera},
\eqref{harlock5}, \eqref{harlock4} and Cauchy estimates.
 
 Now we consider the case of $E$ far away from $\pm 1.$
 Thanks to \eqref{ciaccona}
 (and the analogous estimate in the case
 we are close to +1),
 it remains to invert\footnote{Skipping for brevity the 
  $2j$ subscript from now on .} 
 $\breve \FO(\sa)=\breve \FO_{2j}(\sa)$ 
 for\footnote{Recall \eqref{sanpietrino2}.} 
 $$
 x\in
 [-1+\breve\rho/24,1-\breve\rho/24]\,.
 $$
 First of all we claim that 
\begin{equation}\label{suite}
\partial_\sa \breve \FO(\sa)\geq m:=\frac1{2^7} \breve \b \breve \rho\,,
\qquad \forall\,\sa\in
 [-1+\breve\rho/24,1-\breve\rho/24]\,.
\end{equation}
In order to prove \eqref{suite}
we note that, by \eqref{sanpietrino2}, $m\leq \breve \b/8.$
Then if, by contradiction, there exists
$\bar\sa\in [-1+\breve\rho/24,1-\breve\rho/24]$ with
$\partial_\sa \breve \FO(\bar\sa)< m,$
then by  \eqref{ladispoli4} we have
$
|\partial_{\sa\sa} \breve \FO(\bar\sa)|\geq \breve \b/2\,.
$
To fix ideas we consider the case 
$\partial_{\sa\sa} \breve \FO(\bar\sa) \leq - \breve \b/2$.
Then, recalling that $\breve \FO$ is strictly increasing,
we get for $x\geq \bar\sa,$
$$
0<\partial_\sa \breve \FO(\sa)< m-\frac{\breve \b}{2}(\sa-\bar\sa)
+3\frac{\breve \morse}{\breve s^3} (\sa-\bar\sa)^2
$$
and, choosing $x:=\bar\sa +\breve \rho/2^5,$
we have
$$
0< m- \frac{\breve \b \breve \rho}{2^6}+\frac{3\breve \morse \breve \rho^2}
{2^{10}\breve s^3}
\stackrel{\eqref{sanpietrino}}=
\left(-\frac1{2^7} +\frac{3}{2^{17}}\right) \breve \b \breve \rho<0\,,
$$
which is a contradiction, proving \eqref{suite}.
 
 \nl
 Then, fixing  
 $\tilde x\in [-1+\breve\rho/24,1-\breve\rho/24],$
 we want to apply the Inverse Function Theorem
 in the ball $|x-\tilde x|\leq \rho_1$
 where
 \begin{equation}\label{gianfranco}
 \rho_1:=\frac{\breve \b \breve s^2 \breve \rho}
 {2^{11}\breve \morse}
 \stackrel{\eqref{sanpietrino}}
 =
\frac{\breve \b^2 \breve s^5}{2^{18}\breve \morse^2} 
 \leq 
 \min\Big\{
 \frac{\breve \rho}{2^{11}\cdot3\pi^3}\,,\ 
 \frac{\breve s}{9 \pi^6 2^{18}}\,,\ 
\frac{1}{9 \pi^7 2^{18}}
 \Big\}\,, 
 \end{equation}
 recalling \eqref{harlock4}.
 We have, by Cauchy estimates
 $$
 \sup_{|x-\tilde x|\leq \rho_1}
 \left|
 1-\frac{\partial_\sa \breve \FO(\sa)}{\partial_\sa \breve \FO(\tilde x)}
 \right|
 \leq \frac{8\breve \morse\rho_1}{\breve s^2 m}
 =\frac{\breve \b \breve\rho}{2^8 m}
 \stackrel{\eqref{suite}}=
 \frac12
 $$
 We incidentally note that, by \eqref{suite}
 and Cauchy estimates, we get
 \begin{equation}\label{ghirlanda}
|\partial_\sa \breve \FO(\sa)|\geq \frac{m}{2}
=\frac{\breve \b \breve \rho}{2^8}\,,
\qquad
\forall\, 
x\in [-1+\breve\rho/24,1-\breve\rho/24]_{\rho_1}
\,.
\end{equation}
  Set\footnote{Actually we could choose a larger 
 $\breve r_0$, namely 
 $ \breve r_0:=\frac{\rho_1}{2|T|},$ where
 $T=1/\partial_\sa \breve \FO(\tilde x).$ Then recall \eqref{suite}.}
 $$
 \breve r_0=\frac{m^2 \breve s^2}{2^5 \breve \morse}
 =\frac{\breve \rho^2 \breve \b^2}{2^{19}\breve \morse}
 \stackrel{\eqref{sanpietrino}}=
 \frac{\breve \b^4 \breve s^6}{2^{33}\breve \morse^3}
 \,,\qquad
 \breve E_*:=\breve \FO(\tilde x)\,.
 $$
 Then by the Inverse Function Theorem
 there exists a holomorphic  function
 $$
 \breve{\bar\Sa} \ :\ 
\{ |\breve E-\breve E_*|< \breve r_0\}\ \to\ 
\{|x-\tilde x|\leq \rho_1\}
 $$
 inverting $\breve \FO.$
 We have by \eqref{baux}, \eqref{briciola}
 (and recalling that $\rho_1\leq \breve \rho$)
 $$
 \sup_{\O_{\breve r_0}}
|\breve{\bar\Sa}|
\leq 1+ \breve \rho\stackrel{\eqref{sanpietrino2}}
\leq 2\,,
 $$
 proving \eqref{amber}.
Anyway we note that
 $\breve{\bar\Sa}_+,\breve{\bar\Sa}_-$ are actually defined on 
 the larger domain 
 $$
 B_{r_*/2}(0)\supset B_{\sqrt{\breve r_0}}(0)\,.
 $$
 
 We now prove \eqref{ventresca}, showing
 only the $+$ case, the $-$ case being analogous.
 We first note that
by \eqref{giappone}, the second estimate in \eqref{siracide} (and \eqref{nera}),
 we have that 
 \eqref{ventresca} holds when
 $1\mp \breve E \in \C_*
 \cap B_{\breve r_1}(0).$
 Consider now a point
 $\breve E\in \O_{\breve r_0/32}$
with 
$|\breve E\pm 1|\geq \breve r_0/16$;
 then\footnote{Indeed the following result holds:
 if, for $0<r\leq1/2,$ $z\in \O_{r/2}$ but
 $|z\pm 1|\geq r$ then there exists a real
 $z_0\in[-1+r,1-r]$ such $|z-z_0|\leq
 \sqrt{2-\sqrt 3} \, r\leq 2r/3.$}
 there exists a real point 
 $\breve E_0\in [-1+\breve r_0/16,1-\breve r_0/16]$
 such that $|\breve E-\breve E_0|
 	\geq \breve r_0/24.$
Since the function $\breve{\bar\Sa}_{2j}$
	is increasing on the real we have that
	$$
\breve{\bar\Sa}_{2j}(\breve E_0)	+1
\geq 
\breve{\bar\Sa}_{2j}(-1+\breve r_0/16)+1
\geq
\sqrt{\frac{\breve r_0}{16}}
\frac{\breve s}{2\sqrt
{\breve \morse}}\,,
 $$
 by \eqref{giappone} and the second estimate in
 \eqref{siracide}.
  Therefore, by \eqref{amber} and Cauchy estimates, we get
  $$
|\breve{\bar\Sa}_{2j}(\breve E)	+1|
\geq
\frac{\breve s}{8}
\sqrt{\frac{\breve r_0}{\breve \morse}}
-
\frac{64 \breve r_1}{\breve r_0}
\geq
\frac{\breve s}{16}
\sqrt{\frac{\breve r_0}{\breve \morse}}\,.
  $$
  Then \eqref{ventresca} follows noting that
  $\sqrt{|\breve E\pm 1|}\leq 4.$
 
 We finally prove \eqref{ragnetto}.
 By \eqref{siracide}, \eqref{giappone}
 and noting that $\rho<1$ by \eqref{nera},
 we have that
 \begin{equation}\label{guardrail}
1\pm \breve E \in \C_*\cap B_{\rho^{2/3}}(0) 
\supset 
\C_*\cap B_{\rho}(0)
\qquad
\Longrightarrow
\qquad
1\pm \breve{\bar\Sa}_{2j}(\breve E) \in  B_{\rho^{1/3}\sqrt{2/\breve \b}}(0) \,.
 \end{equation}
 Take now
 $\breve E\in \O_\rho$ but
 $|\breve E\pm 1|\geq \rho^{2/3};$
 then 
 $$d:={\rm dist}(\breve E,\partial \O_{\breve r_0})
 \geq \min\{ \rho^{2/3},\breve r_0/2 \}
 $$
 and by \eqref{amber} and Cauchy estimates
 $$
 |\breve{\bar\Sa}_{2j}(\breve E)-
 \breve{\bar\Sa}_{2j}(\Re\breve E)|
 \leq \frac{2}{d}|\Im \breve E|
 \leq \frac{2\rho}{d}\leq s(\rho)\,,
 $$
 since $\breve \b< 1$ by \eqref{harlock5}.
 The above estimate and \eqref{guardrail}
 conclude the proof  of \eqref{ragnetto}.
 \eproof
 
 \subsection{Inverting the rescaled perturbed potential}
 
 We now find $\breve \Sa_i,$ the inverse of
 $\breve{\Gm}_i.$
 
\begin{lemma}\label{cerveza2}
Assume that $\breve\checco$ in \eqref{borges}
satisfies 
\begin{equation}\label{rochefort8}
\breve\checco 
\leq
\frac{\breve \b^2 \breve s^2 \rho_1^2}
{2^{11} \breve \morse}
\stackrel{\eqref{gianfranco}}=
\frac{\breve \b^6 \breve s^{12}}{2^{47}\breve \morse^5} 
\,.
\end{equation}
Then, for every $1\leq i\leq 2N,$ there exists an 
analytic function $\chi_i(\sa,\hat\act ),$
with\footnote{$\rho_1$ and $\breve{\Gm}^*_i$
where defined in \eqref{gianfranco}
and \eqref{megadirettore2}, respectively.
Recall also \eqref{borges}.}
\begin{equation}\label{rochefort10}
\sup_{[-1,+1]_{\rho_1}\times \hat D_{r_0}}
|\chi_i|\leq
\frac{16}{\breve \b \rho_1}
\sup_{[-1,1]_{\breve s}\times \hat D_{r_0}}
| \breve{\Gm}^*_i|
\leq \frac{16\breve \checco}{\breve \b \rho_1}
\stackrel{\eqref{gianfranco}}=
\frac{2^{22}\breve \morse^2}
{\breve \b^3 \breve s^5} \breve\checco
\qquad
\text{and}\qquad
\chi_i(\pm 1,\hat\act )=0\,,
\end{equation}
such that the analytic function\footnote{
$\breve{\bar \Sa}_i$ was defined in Lemma 
\eqref{cerveza}.
$\O_{\breve r_0}$
was defined in \eqref{bufala} and 
\eqref{nera}.}
 \begin{equation}\label{biretta}
\breve \Sa_i(\breve E,\hat\act ):= \breve{\bar \Sa}_i(\breve E)+ \chi_i (\breve{\bar \Sa}_i(\breve E),\hat\act )\,,\qquad
\breve E\in\O_{\breve r_0}\,,\ \ 
\hat\act \in \hat D_{r_0}\,,
\end{equation}
solves
\begin{equation}\label{rochefort6}
\breve{\Gm}_i\big(
\breve \Sa_i(\breve E,\hat\act ),\hat\act 
\big)=\breve E\,.
\end{equation}
Moreover\footnote{Recalling the definition of
$\breve r_0$ in \eqref{nera}.} 
\begin{equation}\label{giappone2}
\breve \Sa_i(\breve E,\hat\act ) 
=(-1)^i\left(
\pm 1\mp\sqrt{1\mp \breve E}\  \breve \Sa_{i,\pm}( \sqrt{1\mp \breve E},\hat\act )
\right)\,,\quad
\text{ when} \ \ \ 
1\mp \breve E \in \C_*\cap B_{\breve r_0}(0)
\,,
\end{equation}
 where
 \begin{equation}\label{biretta2}
\breve \Sa_{i,\pm}(y,\hat\act )
:=
\breve{\bar\Sa}_{i,\pm}(y)
\Big(
1+
\tilde \chi_{i,\pm}(y,\hat\act )
\Big)
\end{equation}
with
\begin{equation}\label{biretta3}
\tilde \chi_{i,\pm}(y,\hat\act )
:=
\hat \chi_{i,\pm}(\mp y \breve{\bar\Sa}_{i,\pm}(y),\hat\act )
\end{equation}
and\footnote{Note
that $\hat \chi_{i,\pm}$ is analytic
for $z$ close to zero, since 
$\chi_i(\pm 1,\hat\act )=0.$}
\begin{equation}\label{biretta4}
\hat \chi_{i,\pm}(z,\hat\act )
:=
\frac{(-1)^i}{z} 
\chi_i\big(
(-1)^i(\pm 1+z),\hat\act 
\big)\,.
\end{equation}
In particular
 $\breve \Sa_{2j}(\mp 1,\hat\act )=
 \breve \Sa_{2j-1}(\pm 1,\hat\act )=\mp 1.$
 Finally the following estimates hold
 \begin{equation}\label{oxiana}
\sup_{|z|<\rho_1,\hat\act \in\hat D_{r_0}}
|\hat \chi_{i,\pm}(z,\hat\act )|
\leq
\frac{32\breve \checco}{\breve \b \rho_1^2}
\leq \frac{1}{6^4}\,,
\end{equation}
and
\begin{equation}\label{oxiana2}
\sup_{|y|<\sqrt{\breve r_0},
\hat\act \in \hat D_{r_0}}
|\breve \Sa_{i,\pm}(y,\hat\act )|
\leq
\frac{2}{\sqrt{\breve \b}}\,.
\end{equation}

\end{lemma}
\proof 
Let us introduce some notations.
For $s>0,$ let $\breve B^{s}$ the Banach space of analytic function $\chi(\sa,\hat\act )$ on the complex neighborhood
$\chi:[-1,+1]_{s}\times \hat D_{r_0}\to\C,$  with bounded sup-norm
$$
{\modulo}\chi{\modulo}_{s}
:=\sup_{[-1,+1]_{s}\times \hat D_{r_0}}|\chi|
$$
 and such that
$\chi(\pm 1,\hat\act )=0.$  Let $\tilde B^{s}$ be the (closed) subspace 
of $\breve B^{s}$ with
$\partial_\sa \chi(\pm 1,\hat\act )=0.$ 
For $\rho>0,$
let $\breve B^{s}_\rho$, resp. $\tilde B^s_\rho$ be the closed ball  of $\breve B^s$, resp. 
$\tilde B^s,$ with center 0 and radius $\rho$.
Recalling the definition of
$\rho_1$ in \eqref{gianfranco}, let us consider
the two parameters $\rho_2,r_2\geq 0$ with
\begin{equation}\label{tacchino}
\rho_2\leq\frac{\breve \b \breve s^2 \rho_1}
{2^8 \breve \morse}
\leq \frac18 \rho_1^2
\leq\frac{1}{2^{15}}\breve\rho^2
\stackrel{\eqref{sanpietrino2}}
\leq
\frac{1}{2^{37}}\breve s
\,,
\qquad
r_2:=\frac{\breve \b \rho_1}{16}\rho_2
\leq
\frac{\breve \b^2 \breve s^2 \rho_1^2}
{2^{11} \breve \morse}
\,.
\end{equation}
Let us define the function
\beqa{sbs}
&&\mathcal F_i:
\breve B^{\rho_1}_{\rho_2}
 \times \tilde B^{ \breve s}_{r_2}
 \ \longrightarrow\ 
\tilde B^{\rho_1}\,,\nonumber\\
&& 
\big(\mathcal F_i(\chi,\breve G)\big)(\sa,\hat\act ):=
\breve \FO_i(\sa+\chi(\sa,\hat\act ))-\breve \FO_i(\sa)+\breve G(\sa+\chi(\sa,\hat\act ),\hat\act )\,,
\eeqa
where
$\breve \FO_i$
was defined in \eqref{bastard}.
In the following we will often
 drop the index $i$ for brevity,
 writing $\mathcal F, \breve \FO_i,$ 
instead of $\mathcal F_i,\breve \FO_i.$
Note that $\mathcal F$ is well defined
since
\begin{equation}\label{folaga}
\sup_{[-1,+1]_{\rho_1}\times \hat D_{r_0}}
|\Im (\sa+\chi(\sa,\hat\act ))|
\leq \rho_1+\rho_2\leq \frac{\breve s}{4}
\end{equation}
by \eqref{gianfranco}, \eqref{sanpietrino2}
and \eqref{tacchino}.
We want to find $\chi=\chi(\breve G)$
that
 solves the implicit function  $\mathcal F(\chi,\breve G)=0,$
for $\chi$ and $\breve G$ small since
 $\mathcal F(0,0)=0.$ Consider the functional 
$$
\chi\to\check\chi:=\partial_\chi \mathcal F(0,0)[\chi] 
$$
 defined as
$
\check\chi(\sa)=\partial_\sa \breve \FO(\sa)\chi(\sa)\,.
$
Then for a given $\check\chi\in\tilde B^{\rho_1}$ we have that 
\begin{equation}\label{struzzo}
\chi=\big(\partial_\chi \mathcal F(0,0)\big)^{-1}[\check\chi]
=:T\check\chi
\end{equation}
 is given by
\begin{equation}\label{airone}
\chi(\sa):=
\frac{\check\chi(\sa)}{\partial_\sa \breve \FO(\sa)}\,.
\end{equation}
We have to show that the above expression
is well defined for 
$\check\chi\in\tilde B^{\rho_1}$.
We first claim that
\begin{equation}\label{poiana}
\sup_{[-1,1]_{\rho_1}\cap\{|x+ 1|<\breve \rho\}}
|\chi(\sa)|\leq 
\frac{8{\modulo}\check\chi{\modulo}_{\rho_1}}{\breve\b\rho_1}\,.
\end{equation}
We will prove \eqref{poiana}
 only in the case 
$\breve \FO=\breve \FO_{2j}$,
the case $\breve \FO=\breve \FO_{2j-1}$ being analogous.
  Recalling the definition of
$\hat \FO=\hat \FO_{2j}$ in
\eqref{panontina0} and setting
$\xx:=x+1,$ 
we have that
$\breve \FO(\sa)=\breve \FO(\xx-1)=\xx^2 \hat \FO(\xx)-1$ and
$$
\partial_\sa \breve \FO(\sa)
=\partial_\sa \breve \FO(\xx-1)
=2\xx \hat \FO(\xx)+\xx^2\partial_\xx 
\hat \FO(\xx)\,.
$$
By \eqref{cornacchia},
\eqref{sanpietrino2},
Cauchy estimates and \eqref{sanpietrino}, 
we get for $|\xx|<\breve\rho$
\begin{equation*}
|2 \hat \FO(\xx)+\xx\partial_\xx \hat \FO(\xx)|
\geq
2| \hat \FO(\xx)|-|\xx||\partial_\xx \hat \FO(\xx)|
\geq
\frac{\breve \b}{2}-\breve\rho
\frac{20\breve \morse}{\breve s^3}
\geq \frac{\breve \b}{4}\,.
\end{equation*}
Then for every $ |x+ 1|=|\xx|<\breve \rho,$
\begin{equation}\label{fringuello}
|\partial_\sa \breve \FO(\sa)|
=
|\partial_\sa \breve \FO(\xx-1)|
\geq |\xx|\frac{\breve \b}{4}\,.
\end{equation}
Therefore, recalling \eqref{airone},
 we get
\begin{equation}\label{fregata}
\sup_{[-1,1]_{\rho_1}\cap\{\frac{\rho_1}{2}
\leq|x+ 1|<\breve \rho\}}
|\chi(\sa)|\leq 
\frac{{\modulo}\check\chi{\modulo}_{\rho_1}}{
(\rho_1/2)(\breve\b/4)}
=
\frac{8{\modulo}\check\chi{\modulo}_{\rho_1}}{\breve\b\rho_1}\,.
\end{equation}
Recalling that
$\check \chi(-1)=0$, by
 Cauchy estimates we get
$$
|\check \chi(\xx-1)|
\leq |\xx| \sup_{|\xx|\leq \rho_1/2}
|\check \chi(\xx-1)|
\leq |\xx|
\frac{2{\modulo}\check\chi{\modulo}_{\rho_1}}{\rho_1}\,,
\qquad \forall\, |\xx|\leq \frac{\rho_1}{2}\,.
$$
Then, by \eqref{fringuello}
we get
\begin{equation}\label{fregata2}
\sup_{[-1,1]_{\rho_1}\cap\{
|x+ 1|\leq\frac{\rho_1}{2}\}
}
|\chi(\sa)|\leq 
\frac{8{\modulo}\check\chi{\modulo}_{\rho_1}}{\breve\b\rho_1}\,.
\end{equation}
By \eqref{fregata} and 
 \eqref{fregata2} we finally get \eqref{poiana}.
 Moreover, since 
 $\check \chi(-1)=\partial_\sa
 \check \chi(-1)=0$ by \eqref{fringuello}
(recall \eqref{airone}) we have that
 $\chi(-1)=0.$ 
Analogously we can show that
\begin{equation}\label{poiana2}
\sup_{[-1,1]_{\rho_1}\cap\{|x- 1|<\breve \rho\}}
|\chi(\sa)|\leq 
\frac{8{\modulo}\check\chi{\modulo}_{\rho_1}}{\breve\b\rho_1}
\end{equation}
and $\chi(1)=0.$
Moreover, by \eqref{ghirlanda}
\begin{equation}\label{nibbio}
\sup_{[-1+\breve\rho/24,1-\breve\rho/24]_{\rho_1}}
|\chi|
\leq
\frac{2^8}{\breve \b \breve \rho}
{\modulo}\check\chi{\modulo}_{\rho_1}\,.
\end{equation}
Noting that by\footnote{In particular
$|\breve\rho/24+\ii \rho_1|<\breve \rho.$}
 \eqref{gianfranco}
$$
[-1,1]_{\rho_1}\subset
[-1+\breve\rho/24,1-\breve\rho/24]_{\rho_1}
\cup
\{|x+1|<\breve \rho\}
\cup
\{|x-1|<\breve \rho\}\,,
$$ 
we have that, by \eqref{nibbio}, \eqref{poiana}, \eqref{poiana2}
(and \eqref{gianfranco}),
 the expression in \eqref{airone} is well defined, $\chi\in \breve B^{\rho_1}$ and, moreover, we have
\beqa{falco}
&&\sup_{[-1,1]_{\rho_1}}|\chi|=
{\modulo}\chi{\modulo}_{\rho_1}\leq 
\frac{8{\modulo}\check\chi{\modulo}_{\rho_1}}{\breve\b\rho_1}\,,
\quad{\rm i.e.}\nonumber\\
&& {\modulo}\big(\partial_\chi \mathcal F(0,0)\big)^{-1}{\modulo}_{\mathcal L(\tilde B^{\rho_1},
\breve B^{\rho_1} )}
=
{\modulo}T{\modulo}_{\mathcal L(\tilde B^{\rho_1},
\breve B^{\rho_1} )}
\leq 
\frac{8}{\breve\b\rho_1}
\,.
\eeqa
We now apply the Implicit Function Theorem in Banach spaces.
First we note that, since $\rho_1\leq \breve s/4,$ we have
\begin{equation}\label{pinguino}
\sup_{{\modulo}\breve G{\modulo}_{\breve s}\leq r_2}
{\modulo}\mathcal F(0,\breve G){\modulo}_{\rho_1}
=
\sup_{{\modulo}\breve G{\modulo}_{\breve s}\leq r_2}
{\modulo}\breve G(\sa,\hat\act ){\modulo}_{\rho_1}
\leq r_2 
\stackrel{\eqref{tacchino},\eqref{falco}}\leq
\frac{\rho_2}{2{\modulo}\big(\partial_\chi \mathcal F(0,0)\big)^{-1}{\modulo}_{\mathcal L(\tilde B^{\rho_1},
\breve B^{\rho_1} )}}\,.
\end{equation}
For $\chi\in \breve B^{\rho_1}_{\rho_2}$ set
 $$
 R:= \partial_\sa \breve \FO(\sa+\chi)
 - \partial_\sa \breve \FO(\sa)
 $$
 and note that, by Cauchy estimates,
 \eqref{folaga} and \eqref{tacchino}, we get
 \begin{equation*}
 {\modulo}R{\modulo}_{\rho_1}
 \leq \frac{8\breve \morse}{\breve s^2}{\modulo}\chi{\modulo}_{\rho_1}
 \leq
\frac{8\breve \morse\rho_2}{\breve s^2}
\leq
\frac{\breve \b  \rho_1}
{2^5 }
 \end{equation*}
 and, for $\breve G\in \tilde B^{ \breve s}_{r_2} $, also
 $$
 {\modulo} \partial_\sa \breve G(\sa+\chi){\modulo}_{\rho_1}
 \leq \frac{2 r_2}{\breve s}=
 \frac{\breve \b \rho_1}{8\breve s}\rho_2
 \leq 
 \frac{\breve \b \rho_1}{2^{40}}
 \,.
 $$
 Then, by \eqref{falco}, we obtain
\begin{eqnarray}
&&
 \sup_{\breve B^{\rho_1}_{\rho_2}
 \times \tilde B^{ \breve s}_{r_2}
}
 {\modulo} Id-T \partial_\chi\mathcal  F
 (\chi,\breve G){\modulo}_{\mathcal L
 ( \breve B^{\rho_1},\breve B^{\rho_1}
 )}
 \nonumber
 \\
 &&
 =\sup_{\chi'\in \breve B^{\rho_1},\,{\modulo}\chi'{\modulo}_{\rho_1}=1}
 {\modulo}\chi' -T
 \Big(
 \partial_\sa \breve \FO(\sa+\chi)
 +
  \partial_\sa \breve G(\sa+\chi)
 \Big)\chi'{\modulo}_{\rho_1}
 \nonumber
 \\
 &&
 \leq 
 {\modulo}T{\modulo}
 \big(
{\modulo}R{\modulo}_{\rho_1}+{\modulo}\partial_\sa\breve G(\sa+\chi)
{\modulo}_{\rho_1}
\big)
\nonumber
\\
&&
\leq 
\frac{8}{\breve\b\rho_1}\left(
\frac{\breve \b  \rho_1}
{2^5 }
+
 \frac{\breve \b \rho_1}{2^{40}}
\right)\leq \frac12\,.
 \label{quaglia}
  \end{eqnarray}
 Since estimates 
 \eqref{pinguino} and \eqref{quaglia}
 are satisfied we can apply the Implicit 
 Function Theorem finding, for every
 $\breve G\in \tilde B^{ \breve s}_{r_2} $
 a function $\chi=\chi_i\in \breve B^{\rho_1}_{\rho_2}$ solving \eqref{sbs}, namely
 $\mathcal F_i(\chi_i,\breve G)=0.$
 In particular
the theorem can be applied with
$\breve G:=\breve{\Gm}^*_i$ defined in 
\eqref{megadirettore2}.
Indeed, by  \eqref{megadirettore3},
$\breve{\Gm}^*_i\in \tilde B^{ \breve s}.$
Moreover choosing
$r_2$ in \eqref{tacchino} as 
$$
r_2:=
\sup_{[-1,1]_{\breve s}\times \hat D_{r_0}}
| \breve{\Gm}^*_i|\,,
$$ 
we have that the conditions in \eqref{tacchino}
are satisfied by assumption \eqref{rochefort8}
(recall \eqref{borges}).
Then the (first) estimate in  \eqref{rochefort10}
directly follows by \eqref{tacchino}.

Then we have $\mathcal F_i(\chi_i,\breve{\Gm}^*_i)=0,$ which is equivalent  to
$$
\breve{\Gm}_i(\sa+\chi_i(\sa,\hat\act ),\hat\act )=\breve{\FO}_i(\sa)\,.
$$
Evaluating at $x=\breve{\bar \Sa}_i(\breve E)$
we get (recall \eqref{blu})
 \eqref{rochefort6}.

We now prove \eqref{oxiana}.
Recalling \eqref{biretta4}, 
 we distinguish the case 
$\rho_1/2\leq |z|<\rho_1$, where we directly
obtain \eqref{oxiana}, and the case
  $|z|<\rho_1/2$ where we have,
  uniformly for $\hat\act \in\hat D_{r_0},$
  $$
|\hat \chi_{i,\pm}(z,\hat\act )| \leq 
\sup_{[-1,1]_{\rho_1/2}\times \hat D_{r_0}}
|\partial_\sa \chi_i|
\stackrel{\eqref{rochefort10}}\leq
\frac{16\breve \checco}{\breve\b \rho_1}
\frac{1}{\rho_1/2}
  $$
  recalling that 
$\hat \chi_{i,\pm}(\pm 1,\hat\act )=0$
and by Cauchy estimates. The first inequality
in \eqref{oxiana} is proved.
The second one follows noting that
$$
\frac{32\breve \checco}{\breve\b \rho_1^2}
\stackrel{\eqref{rochefort8}}\leq
\frac{\breve \b \breve s^2}{2^6 \breve \morse}
\stackrel{\eqref{sanpietrino}}=
\frac{2\breve \rho}{\breve s}
\stackrel{\eqref{sanpietrino2}}
\leq \frac{1}{6^4}\,.
$$
Finally \eqref{oxiana2} follows by
\eqref{biretta2},
\eqref{oxiana} and \eqref{siracide}.
  \eproof

In the following we will assume the condition\footnote{Recall \eqref{harlock4}.} 
\begin{equation}\label{vascodegama2}
\breve\checco
\leq
\frac{\breve \b^9 \breve s^{18}}
{2^{73}\breve \morse^8}\,.
\end{equation}
\begin{remark}
  Note that, recalling
\eqref{exit}-\eqref{ladispoli4}, condition
\eqref{vascodegama2} is implied by\footnote{
The two conditions are equivalent w.r.t. the dependence
on the parameter, we choose the smaller constant
$2^{-197}$ in \eqref{vascodegama3} for brevity.}
\begin{equation}\label{vascodegama3}
\checco
\leq
\frac{\b^{28}s_0^{45}s_*^2}{2^{197}\morse^{27}}\,,
\end{equation}
which follows by \eqref{genesis}.
\end{remark}

  \begin{lemma}\label{lem:coppe}
Assume \eqref{vascodegama2}. Set\footnote{Recall
\eqref{exit}-\eqref{ladispoli4}}
\begin{equation}\label{coppe}
\breve{\mathtt r}:=
\frac{\breve\b^{10} \breve s^{19}}
{2^{88} \breve \morse^9}
=
\frac{ s_0^{49}\b^{29}}
{2^{134} 3^{10}\pi^{29}  \morse^{29}}
\,,
\qquad
\rho_\star:=
\frac{\breve \b^6 \breve s^{13}}{2^{49}
\breve \morse^6}
=
\frac{s_0^{31}\b^{18}}
{2^{70} 3^{6}\pi^{19}  \morse^{18}}
\,.
\end{equation}
Then 
\begin{equation}\label{fiodena}
\breve E \in \O_{\breve{\mathtt r}} 
 \qquad\Longrightarrow\qquad
\breve{\bar\Sa}_i(\breve E)\,,\ 
\breve{\Sa}_i(\breve E,\hat\act)\ \in \ 
[-1,1]_{\rho_\star/2}\,,\ \ 
\forall\, \hat\act\in \hat D_{r_0}\,.
\end{equation}
\end{lemma}
\proof
First 
we note that, recalling \eqref{nera},
$$
\frac{4\breve{\mathtt r}}{\breve r_0}=
\frac{\rho_\star}{2^7}
$$
(which implies $\breve{\mathtt r}\leq \breve r_0/2$)
and
$$
\frac{4\breve{\mathtt r}^{1/3}}{\sqrt{\breve\b}}
\geq 
\frac{4\breve{\mathtt r}}{\breve r_0}
$$
by \eqref{harlock5}.
Then by \eqref{ragnetto} 
we prove that, actually,
\begin{equation*}
\breve{\bar\Sa}_i(\breve E)\ \in \ 
[-1,1]_{\rho_\star/2^7}\,.
\end{equation*}
Then by \eqref{rochefort8},\eqref{biretta} and \eqref{rochefort10} we get
$\breve{\Sa}_i(\breve E,\hat\act )\ \in \ 
[-1,1]_{\breve \xx_\star/2}$
since
$$ 
\frac{2^{22}\breve \morse^2}
{\breve \b^3 \breve s^5} \breve\checco
\leq
\frac{\rho_\star}{4}\,,
$$
by \eqref{vascodegama2} and \eqref{harlock4}.
\eproof

\begin{lemma}
Assume \eqref{vascodegama2}.
Then
\begin{equation}\label{animelle}
|\partial_\xx \breve{\FO}_i(\sa+\chi_i(\sa))
-\partial_\xx \breve{\FO}_i(\sa)|
\leq
\frac{2^{71}\breve \morse^8}
{\breve \b^9 \breve s^{18}}
\breve\checco 
|\partial_\xx \breve{\FO}_i(\sa)|\,,
\qquad
\forall\, \sa\in[-1,+1]_{\rho_\star}
\end{equation} 
\end{lemma}
\proof
We start considering a neighborhood of $\pm 1.$
We apply Lemma \ref{fettuccine3}  with
$f \rightsquigarrow\partial_\xx \breve{\FO}_i$, 
$\chi\rightsquigarrow\chi_i$, 
$0\rightsquigarrow\pm1$, $r\rightsquigarrow\rho_1$
(defined in \eqref{gianfranco}),
$M_2\rightsquigarrow \frac{2\breve \morse}{\breve s}$
(recall \eqref{gibbone}),
$\checco\rightsquigarrow
\frac{16\breve \checco}{\breve \b \rho_1}$
(recall \eqref{rochefort10}), 
$\rho\rightsquigarrow \rho_3,$
where, recalling \eqref{ladispoli4},
\begin{equation}\label{magellano}
\rho_3:=\frac{\rho_1^2\breve\b\breve s}{32\breve\morse}
=
\frac{\breve \b^5 \breve s^{11}}{2^{41}
\breve \morse^5}
\,.
\end{equation}
By assumption\footnote{Which implies
condition $\checco\leq r/8$ of Lemma 
\ref{fettuccine3}.
} \eqref{vascodegama2}
we can apply Lemma \ref{fettuccine3}
obtaining
\begin{equation}\label{foco}
|\partial_\xx \breve{\FO}_i(\sa+\chi_i(\sa))
-\partial_\xx \breve{\FO}_i(\sa)|
\leq
\frac{32}{\breve \b \rho_1 \rho_3}
 \breve\checco 
|\partial_\xx \breve{\FO}_i(\sa)|\,,
\qquad
\forall\, |\sa\pm 1|<\rho_3\,.
\end{equation}
Applying
\eqref{harlock7bar}
with\footnote{Note that the condition
$\breve \xx_0\leq \breve\xx_\sharp$ is satisfied by
\eqref{harlock4}.
Note also that
$$
\rho_\star:=\frac{\breve\b \breve s^2\rho_3}
{2^{8} \breve M}
<\frac{1}{2^{12}}\rho_3\,.
$$
} $\breve \xx_0\rightsquigarrow \rho_3/2$
$$
\inf_{[-1+\rho_3/2,
1-\rho_3/2]_{\rho_\star}}
|\partial_\xx\breve\FO_i|
\geq 
\frac{\breve\b \breve \xx_0}{16}
\,.
$$
Then by \eqref{gibbone}, Cauchy estimates
\eqref{rochefort10}
\begin{eqnarray*}
|\partial_\xx \breve{\FO}_i(\sa+\chi_i(\sa))
-\partial_\xx \breve{\FO}_i(\sa)|
\leq
\frac{2^7\breve\morse}{\breve \b \breve s^2 \rho_1}
\breve\checco
\leq
\frac{2^{12}\breve\morse}{\breve \b^2 \breve s^2 \rho_1
\rho_3}
\breve\checco
|\partial_\xx \breve{\FO}_i(\sa)|
\end{eqnarray*}
\eproof

 Note that by \eqref{harlock4}, \eqref{rochefort10} and \eqref{vascodegama2}
 $$
 \rho_\star\leq
 \breve\xx_\star/2^{57}\qquad
 \text{and}
 \qquad
 \sup_{[-1,+1]_{\rho_1}\times \hat D_{r_0}}
|\chi_i|\leq 
\frac{\breve\xx_\star}{32}\,,
 $$
 then (recall \eqref{borges2}) 
\begin{equation}\label{confortably}
 \sa\in[-1,1]_{\rho_\star}\ \ \Longrightarrow\ \ 
 \sa+\chi_i(\sa)
 \in[-1,1]_{\breve\xx_\star/16}\,.
\end{equation}
 We claim that
\begin{equation}\label{salmone}
|
\partial_\xx  \breve{\Gm}_i (\sa+\chi_i(\sa),\hat\act)
-\partial_\xx \breve{\FO}_i(\sa)
|
\leq
\frac{2^{72}\breve \morse^8}
{\breve \b^9 \breve s^{18}}
\breve\checco 
|\partial_\xx \breve{\FO}_i(\sa)| \,,
\qquad
\forall\,
 \sa\in[-1,1]_{\rho_\star}\,, \ 
 \hat \act \in \hat D_{r_0}\,.
\end{equation}
Indeed 
by \eqref{confortably}, \eqref{borges2} and \eqref{animelle}
we have
 \begin{eqnarray*}
&&|
\partial_\xx  \breve{\Gm}_i (\sa+\chi_i(\sa),\hat\act)
-\partial_\xx \breve{\FO}_i(\sa)
|
\\
&\leq&
|
\partial_\xx  \breve{\Gm}_i (\sa+\chi_i(\sa),\hat\act)
-\partial_\xx \breve{\FO}_i(\sa+\chi_i(\sa))
|
+
|
\partial_\xx \breve{\FO}_i(\sa+\chi_i(\sa))
-
\partial_\xx \breve{\FO}_i(\sa)
|
\\
&\leq&
\frac{2^{11}\breve\morse}{\breve\b^2\breve s^4}
\breve \checco 
|\partial_\xx \breve{\FO}_i(\sa+\chi_i(\sa))|
+
\frac{2^{71}\breve \morse^8}
{\breve \b^9 \breve s^{18}}
\breve\checco 
|\partial_\xx \breve{\FO}_i(\sa)| 
\\
&\leq&
\left(
1+\frac{\breve \b^{7} \breve s^{14}}{2^{60}\breve \morse^7}
+
\frac{2^{11}\breve\morse}{\breve\b^2\breve s^4}
\breve \checco 
\right)
\frac{2^{71}\breve \morse^8}
{\breve \b^9 \breve s^{18}}
\breve\checco 
|\partial_\xx \breve{\FO}_i(\sa)| 
\\
&\stackrel{\eqref{vascodegama2}}\leq&
\frac{2^{72}\breve \morse^8}
{\breve \b^9 \breve s^{18}}
\breve\checco
|\partial_\xx \breve{\FO}_i(\sa)| 
\,.
\end{eqnarray*}
By \eqref{vascodegama2}
and \eqref{salmone}
we also have
\begin{equation}\label{salmone2}
\frac12 
|\partial_\xx \breve{\FO}_i(\sa)
|
\leq
|
\partial_\xx  \breve{\Gm}_i (\sa+\chi_i(\sa),\hat\act)
|
\leq 2
|\partial_\xx \breve{\FO}_i(\sa)
|
 \,,
\quad
\forall\,
 \sa\in[-1,1]_{\rho_\star}\,, \ 
 \hat \act \in \hat D_{r_0}\,.
\end{equation}


\subsection{Inverting the original potentials}

Recalling \eqref{bastard} we set
\begin{equation}\label{blindness}
\bar\Sa_i:= \bar  \g_i\circ \breve{\bar\Sa}_i\circ \bar  \l_i\,,\qquad
\Sa_i:= \g_i\circ \breve \Sa_i\circ \l_i\,,
\end{equation}
solving 
\begin{equation}\label{premiata}
\FO(\bar \Sa_i(E))=E\,,\qquad
\Gm \big(\Sa_i(E,\hat\act  ), \hat\act  \big)=E\,,
\qquad
\forall\, 1\leq i\leq 2N
\,, 
\end{equation}
respectively.
Note that for $\hat\act \in \hat D$ (namely $\hat\act $ real)
\begin{eqnarray}
\Sa_{2j}(\cdot,\hat\act ) 
&:& 
\big[E_{2j-1}(\hat\act  ), E_{2j}(\hat\act  )\big]
\to \big[\sa_{2j-1}(\hat\act  ),\sa_{2j}(\hat\act  )\big]\,,
\nonumber
\\
\Sa_{2j-1}(\cdot,\hat\act  ) 
&:& 
\big[E_{2j-1}(\hat\act  ), E_{2j-2}(\hat\act  )\big]
\to \big[\sa_{2j-2}(\hat\act  ),\sa_{2j-1}(\hat\act  )\big]\,,
\label{foggia}
\end{eqnarray}
moreover $\Sa_i$ is increasing, resp. decreasing (as a function of $E$), if 
$i$ is even, resp. odd.
Note also that
\begin{equation}\label{9till5}
\partial_E \Sa_i(E,\hat\act  )=
1/\partial_\sa \Gm \big(\Sa_i(E,\hat\act  ), \hat\act  \big)
\end{equation}
and
\begin{eqnarray}\nonumber
&&\Sa_{2j-1}(E_{2j-2}(\hat\act  ),\hat\act  )
=\sa_{2j-2}(\hat\act  )\,,
\qquad
\Sa_{2j-1}(E_{2j-1}(\hat\act  ),\hat\act  )
=\sa_{2j-1}(\hat\act  )\,,
\\
&&\Sa_{2j}(E_{2j-1}(\hat\act  ),\hat\act  )
=\sa_{2j-1}(\hat\act  )\,,
\qquad
\Sa_{2j}(E_{2j}(\hat\act  ),\hat\act  )
=\sa_{2j}(\hat\act  )\,.
\label{cook}
\end{eqnarray}
Regarding the derivatives we have
by \eqref{faciolata} and \eqref{faciolata2}
\begin{equation}\label{faciolata3}
\partial_{\hat\act } \sa_i(\hat\act )=-
\frac{\partial_{\sa\hat\act } G(\sa_i(\hat\act ), \hat\act )}
{\partial_{\sa\sa} \Gm(\sa_i(\hat\act ), \hat\act )}
\end{equation}
and
\begin{equation}\label{faciolata4}
\partial_{\hat\act } E_i(\hat\act )=
\partial_{\hat\act } G(\sa_i(\hat\act ), \hat\act )-
\frac{\partial_\sa \Gm(\sa_i(\hat\act ), \hat\act )}
{\partial_{\sa\sa} \Gm(\sa_i(\hat\act ), \hat\act )}
\partial_{x\hat\act } G(\sa_i(\hat\act ), \hat\act )\,.
\end{equation}

Now, recalling 
that $\breve\Sa_{i,\pm}$ (defined in 
\eqref{biretta2})
are holomorphic in $B_{\sqrt{\breve r_0}}(0)\times
\hat D_{r_0}$ (see \eqref{oxiana2})
and noting that
$$
|y|<\frac12 \sqrt{\breve r_0  \b }
\stackrel{\eqref{ladispoli3}}\leq
 \sqrt{\frac{\breve r_0}{2}}
\sqrt{E_{2j}(\hat\act )-E_{2j-1}(\hat\act )}\,,
$$
we can define, for
\begin{equation}\label{burrata}
|y|<\frac12 \sqrt{\breve r_0  \b }
\stackrel{\eqref{nera}}=
\frac{s_0^{9}}{2^{25}9\pi^{5}}
\frac{\b^6}{\morse^{11/2}}
=:
r_\diamond
\,,\qquad
\hat\act \in \hat D_{r_0}\,,
\end{equation}
  the holomorphic 
functions
\begin{equation}\label{BWV1013}
\Sa_{i,\pm}(y,\hat\act )
:=
\frac{\sa_{i}(\hat\act )-\sa_{i-1}(\hat\act )}{
\sqrt 2
\sqrt{(-1)^i\big(
E_{i}(\hat\act )-E_{i-1}(\hat\act )
\big)}}
\breve\Sa_{i,\pm}
\bigg(
\frac{\sqrt 2}{
\sqrt{(-1)^i\big(
E_{i}(\hat\act )-E_{i-1}(\hat\act )
\big)}
}
y,\hat\act 
\bigg)
\end{equation}

By \eqref{oxiana2} and \eqref{ladispoli3}
 we get
\begin{equation}\label{oxiana3}
\sup_{|y|<r_\diamond,
\hat\act \in \hat D_{r_0}}
|\Sa_{i,\pm}(y,\hat\act )|
\leq
\frac{4\pi}{\sqrt{\breve \b \b}}
\stackrel{\eqref{ladispoli4}}=
\frac{4\sqrt 3 \pi \morse}{s_0^{3/2} \b^{3/2}}\,.
\end{equation}

\begin{lemma}\label{cumuli}
The following equalities hold\footnote{Where they are meaningful,
namely, e.g., in the third equality $E-E_{2j-1}(\hat\act )\in\C_*,$
$|E-E_{2j-1}(\hat\act )|< r_\diamond^2$
and $\hat\act \in \hat D_{r_0}$.
}
\begin{eqnarray}
\Sa_{2j-1}(E,\hat\act ) &=& \sa_{2j-1}(\hat\act ) -\sqrt{E-E_{2j-1}(\hat\act )}
\ \Sa_{2j-1,-}\Big(\sqrt{E-E_{2j-1}(\hat\act )},\hat\act \Big)\,,
\nonumber
\\
 \Sa_{2j-1}(E,\hat\act ) &=& \sa_{2j-2}(\hat\act ) +
 \sqrt{E_{2j-2}(\hat\act )-E}
\  \Sa_{2j-1,+}\Big(\sqrt{E_{2j-2}(\hat\act )-E},\hat\act \Big)\,,
\nonumber
\\
\Sa_{2j}(E,\hat\act ) &=& \sa_{2j-1}(\hat\act ) +\sqrt{E-E_{2j-1}(\hat\act )}
\ \Sa_{2j,-}\Big(\sqrt{E-E_{2j-1}(\hat\act )},\hat\act \Big)\,,
\nonumber
\\
 \Sa_{2j}(E,\hat\act ) &=& \sa_{2j} (\hat\act )-\sqrt{E_{2j}(\hat\act )-E}
\  \Sa_{2j,+}\Big(\sqrt{E_{2j}(\hat\act )-E},\hat\act \Big)\,.
 \label{gargano}
\end{eqnarray}
Moreover
\begin{eqnarray}\label{stambecco}
&&\Sa_{2j,-}(0,\hat\act )=\sqrt{\frac{2}{\partial_{\sa\sa}\Gm(
\sa_{2j-1}(\hat\act ),\hat\act )}}
=\Sa_{2j-1,-}(0,\hat\act )\,,
\nonumber
\\
&&\Sa_{2j,+}(0,\hat\act )=\sqrt{\frac{2}{-\partial_{\sa\sa}\Gm(\sa_{2j}(\hat\act ),
\hat\act )}}
=\Sa_{2j+1,+}(0,\hat\act )
\end{eqnarray}
and
\begin{equation}\label{lapulcedacqua}
\Sa_{2j,-}(y,\hat\act )=\Sa_{2j-1,-}(-y,\hat\act )\,,
\qquad
\Sa_{2j,+}(y,\hat\act )=\Sa_{2j+1,+}(-y,\hat\act )\,.
\end{equation}
Finally, for every $\hat\act \in \hat D_{r_0}$
\begin{equation}\label{spadadefoco}
|\Sa_{i,\pm}(0,\hat\act )|\geq \frac{s_0}{2\sqrt\morse}
\end{equation}
and
\begin{equation}\label{spadadefoco2}
\sup_{|y|\leq r_\dag,\, \hat\act \in \hat D_{r_0}}
\frac{1}{|\Sa_{i,\pm}(y,\hat\act )|}
\leq \frac{4\sqrt\morse}{s_0}\,,
\end{equation}
where
\begin{equation}\label{spadadefoco3}
r_\dag:=
\frac{s_0^{23/2} \b^{15/2}}{2^{47} \morse^7}
\stackrel{\eqref{harlock}}
<
\frac{r_\diamond}{2^{30}}
\,.
\end{equation}
\end{lemma}
\proof
The equalities in \eqref{gargano}
follows by \eqref{blindness}, \eqref{giappone2} and \eqref{BWV1013}.
\\
\eqref{stambecco} follows by \eqref{gargano} and \eqref{premiata}.
For example, in order to prove the first\footnote{The other being analogous.}
 equality in \eqref{stambecco},
we insert the third expression in \eqref{gargano} into the second equality
 in \eqref{premiata} obtaining\footnote{Developing in Taylor expansion w.r.t. $y$.}
 $$
0= \frac{\Gm\big(\sa_{2j-1}(\hat\act )
+y\,\Sa_{2j,-}(y,\hat\act )\big)-E}{y^2}
=
-1+\frac12 \partial_{\sa\sa} \Gm(\sa_{2j-1}(\hat\act ),\hat\act )
\big(\Sa_{2j,-}(y,\hat\act )\big)^2 + O(y)
 $$
where $y$ is short for $\sqrt{E-E_{2j-1}(\hat\act )}.$
Then we conclude taking $y\to 0.$
\\
We finally prove the first equality in \eqref{lapulcedacqua}.
Let us fix $\hat\act $ and, for brevity, let us omit to write it.
Again let $y$ be short for $\sqrt{E-E_{2j-1}}.$
Inserting the first and the third equalities in \eqref{gargano} into
\eqref{premiata}
we get
$$
\Gm\big(\sa_{2j-1}
-y\,\Sa_{2j-1,-}(y)\big)
=
y^2+E_{2j-1}
=
\Gm\big(\sa_{2j-1}
+y\,\Sa_{2j,-}(y)\big)\,.
$$
In particular we have that  $\Gm\big(\sa_{2j-1}
-y\,\Sa_{2j-1,-}(y)\big)$ is an even function of $y$.
Then
\begin{equation}\label{albinoni}
\Gm\big(\sa_{2j-1}+
y\,\Sa_{2j-1,-}(-y)\big)
=
y^2+E_{2j-1}
=
\Gm\big(\sa_{2j-1}
+y\,\Sa_{2j,-}(y)\big)\,.
\end{equation}
Now we claim that for $y$ close to $0$
there exists a {\sl unique} function $g(y)$ with $\Re g(0)>0$
such that
\begin{equation}\label{galuppi}
\Gm\big(\sa_{2j-1}+
y\,g(y)\big)
=
y^2+E_{2j-1}\,.
\end{equation}
In order to prove the claim we  define the analytic function
 $$
 \tilde \Gm(\xx):=(\Gm(\sa_{2j-1}+\xx)-E_{2j-1})/\xx^2\,,
 $$
 and note that 
 \begin{equation}\label{purcello}
\Re \tilde \Gm(\xx)=\Re \frac12\partial_{\sa\sa}f(\sa_{2j-1})>0\,.
\end{equation}
 Then, for $y$ close to 0, it is well defined the function
 $\sqrt{\tilde \Gm \big(yg(y)\big)}.$
We can rewrite \eqref{galuppi} as
 $$
 g^2(y)\, \tilde \Gm \big(yg(y)\big)=y^2
 $$
 and the above equation has
 a {\sl unique} solution $g(y)$ with $\Re g(0)>0$.
Indeed, taking the square root, it is equivalent to
$$
 g(y)\, \sqrt{\tilde \Gm \big(yg(y)\big)}=y\,,
 $$
 which,
 by the implicit function theorem  and since
 $\sqrt{\tilde \Gm(0)}\neq 0$ (by \eqref{purcello}),
 has
 a {\sl unique} solution $g(y)$ with $\Re g(0)>0$. 
So we have proved that there exists 
a {\sl unique} function $g(y)$ with $\Re g(0)>0$ solving \eqref{galuppi}.
Then by \eqref{albinoni} we conclude that $\Sa_{2j,-}(y)=\Sa_{2j-1,-}(-y)$
noting that by \eqref{stambecco} (and our definition of the square root)
$$
\Re \Sa_{2j,-}(0)=\Re \Sa_{2j-1,-}(0) >0\,.
$$
Finally recalling that by \eqref{octoberx}
$|\sa_i(\hat\act )-\bar \sa_i|\leq s_0/8$
and that $\bar \sa_i\in\R$, by Cauchy estimates
and \eqref{ciccio3}
we have
$$
|\partial_{\sa\sa}\Gm(
\sa_i(\hat\act ),\hat\act )|
\leq
\frac{4\morse}{(7s_0/8)^2}
\leq 
\frac{8\morse}{s_0^2}\,.
$$
Then by \eqref{stambecco}
we get \eqref{spadadefoco}.
Finally by \eqref{oxiana3}, Cauchy estimates,
\eqref{spadadefoco}
and \eqref{spadadefoco3}
we get \eqref{spadadefoco2}.
\eproof

\begin{remark}\label{ledimensionicontano}
 Note that $\morse,\b,\checco,E, \Gm,$ etc., respectively
  $r_0,y,\hat\act $, have the same
 ``physical'' dimension of an energy,
 respectively of an action (the square root of an energy).
 The parameters with ``$\ \breve{}\ $'' are adimensional.
\end{remark}

\section{The action variables}

\rem\label{cleopatra} 
As well known, in the proof of Arnold--Liouville's theorem, 
one first introduces the actions  through  line integrals
\beq{4-1}
\act_n^{(i)}(E,\hat\act):= \frac1{2\pi}\oint_{\Hpend^{-1}(E;\hat\act)} \ p_ndq_n\ 
\eeq 
as 
functions of energy $E$ and then defines the integrated Hamiltonian  ${\mathtt E}^{(i)} $ inverting such functions; in particular, one has 
\begin{equation}
\label{enrico}
\boxed{\act_n^{(i)}\Big({\mathtt E}^{(i)} (\hat \act, \act_n), \hat \act \Big)=\act_n}
\end{equation}
Indeed, all the fine analytic properties of ${\mathtt E}^{(i)}$   will be described through  the functions $E\to \act_n^{(i)}(E,\hat\act)$.
\erem

\subsection{Actions for the unperturbed Hamiltonian}

Recall \eqref{maieli}.
For $1\leq j\leq N$ set
\begin{equation}\label{padula}
j_\diamond:=j-1 \ \ \ \text{if}\ \ \ \bar E_{2j-2}<
\bar E_{2j}\quad \text{and}\quad
 j_\diamond:=j\ \ \ \text{otherwise}\,.
\end{equation}
Note that
$$
\bar E_{2j_\diamond}=\min\{ \bar E_{2j-2}, 
\bar E_{2j}\}\,.
$$
For $1\leq j< N$ set
\begin{eqnarray}\label{21stcentury}
&&
j_-:=\max\{ i<j\ \ {\rm s.t.}\ \ \bar E_{2i}>
 \bar E_{2j} \}\,,\qquad
 j_+:=\min\{ i>j\ \ {\rm s.t.}\ \ \bar E_{2i}>
 \bar E_{2j} \}
 \nonumber
 \\
 &&
  j_*=j_- \  \ {\rm if}\ \ \  \bar E_{2j_-} < \bar E_{2j_+}\ \ \ 
 {\rm and}  \ \ \ j_*=j_+ \ \ \ {\rm otherwise.}
\end{eqnarray}
Note that
\begin{equation}\label{anvedi}
\bar E_{2j}< \bar E_{2j_*}=\min\{ \bar E_{2j_-}, \bar E_{2j_+}\}\,.
\end{equation}
Let us set
\begin{eqnarray}\label{gruffalo}
&&\bar E^{(i)}_-:=\bar E_i\,,
\qquad   E^{(i)}_-(\hat\act):= E_i(\hat\act)\,,\qquad{\rm for}\ \ 1\leq i\leq 2N\,,
\nonumber
\\
&&\bar E^{(2j-1)}_+:=\bar E_{2j_\diamond}\,,
\quad
 E^{(2j-1)}_+(\hat\act):= E_{2j_\diamond}(\hat\act)\,,
\qquad {\rm for}\ \ 1\leq j\leq N\,,
\nonumber
\\
&&\bar E^{(2j)}_+:=\bar E_{2j_*}\,,
\quad
E^{(2j)}_+(\hat\act):= E_{2j_*}(\hat\act)\,,
\qquad {\rm for}\ \ 1\leq j< N\,.
\end{eqnarray}
 Let us define the functions
\begin{eqnarray}
\bar \act_n^{(2j-1)}(E)&:=&
\frac{1}{\pi}
\int_{\bar \Sa_{2j-1}(E )}^{\bar \Sa_{2j}(E )}
\sqrt{E-\FO(\sa)}\, d\sa\,,
\qquad 
\bar E^{(2j-1)}_- \leq E\leq \bar E^{(2j-1)}_+\,,
\quad 1\leq j\leq N\,,
\nonumber
\\
\bar \act_n^{(2j)}(E)&:=& 
\frac{1}{\pi}
\int_{\bar \Sa_{2j_-+1}(E )}^{\bar \Sa_{2j_+}(E )}
\sqrt{E-\FO(\sa)}\, d\sa\,,
\qquad 
\bar E^{(2j)}_- \leq E\leq \bar E^{(2j)}_+\,,
\quad 1\leq j\leq N-1\,,
\nonumber
\\
\bar \act_n^{(2N)}(E)&:=& 
\frac{1}{2\pi}
\int_{-\pi}^{\pi}
\sqrt{E-\FO(\sa)} \, d\sa\,,
\qquad
E\geq \bar E^{(2N)}_-
\nonumber
\\
\bar \act_n^{(0)}(E)&:=&
-\bar \act_n^{(2N)}(E)=
-\frac{1}{2\pi}
\int_{-\pi}^{\pi}
\sqrt{E-\FO(\sa)} \, d\sa\,,
\qquad
E\geq \bar E^{(2N)}_-\,.
\label{sunday}
\end{eqnarray}

\begin{remark}\label{telemann1}
 We want to show that the functions in \eqref{sunday}
 have an analytic extension for complex $E.$
 Moreover, while it is immediate to evaluate the derivative
 of the last two, namely
 \begin{eqnarray}
\partial_E \bar \act_n^{(2N)}(E)&=& 
\frac{1}{4\pi}
\int_{-\pi}^{\pi}
\frac{1}{\sqrt{E-\FO(\sa)}} \, d\sa\,,
\nonumber
\\
\partial_E \bar \act_n^{(0)}(E)&=&
-\frac{1}{4\pi}
\int_{-\pi}^{\pi}
\frac{1}{\sqrt{E-\FO(\sa)}} \, d\sa\,,
\label{sundayD}
\end{eqnarray}
it is not obvious to justify the formal derivation\footnote{
See Remark \ref{telemann2} below.}
 \begin{eqnarray}
\partial_E \bar \act_n^{(2j-1)}(E)&=&
\frac{1}{2\pi}
\int_{\bar \Sa_{2j-1}(E )}^{\bar \Sa_{2j}(E )}
\frac{1}{\sqrt{E-\FO(\sa)}}\, d\sa\,,
\nonumber
\\
\partial_E \bar \act_n^{(2j)}(E)&=& 
\frac{1}{2\pi}
\int_{\bar \Sa_{2j_-+1}(E )}^{\bar \Sa_{2j_+}(E )}
\frac{1}{\sqrt{E-\FO(\sa)}}\, d\sa
\,.
\label{cippone}
\end{eqnarray}
\end{remark}

Set
\begin{eqnarray}
\bar \act_n^{(2j-1),-}(E)
&:=&
\frac{1}{\pi}
\int_{\bar \Sa_{2j-1}(E )}^{\bar\sa_{2j-1}}
\sqrt{E-\FO(\sa)}\, d\sa\,,
\qquad 
\bar E_{2j-1}\leq E\leq \bar E_{2j-2}\,,
\ 1\leq j\leq N\,,
\nonumber
\\
\bar \act_n^{(2j-1),+}(E)
&:=&
\frac{1}{\pi}
\int_{\bar\sa_{2j-1}}^{\bar \Sa_{2j}(E )}
\sqrt{E-\FO(\sa)}\, d\sa\,,
\qquad 
\bar E_{2j-1}\leq E\leq \bar E_{2j}\,,
\ 1\leq j\leq N\,,
\nonumber
\\
\bar \act_n^{(2j),-}(E)
&:=&
\frac{1}{\pi}
\int_{\bar\sa_{2j-1}}^{\bar \sa_{2j}}
\sqrt{E-\FO(\sa)}\, d\sa\,,
\qquad 
E\geq \bar E_{2j}\,,
\ 1\leq j\leq N\,,
		\nonumber
\\
\bar \act_n^{(2j),+}(E)
&:=&
\frac{1}{\pi}
\int_{\bar\sa_{2j}}^{\bar \sa_{2j+1}}
\sqrt{E-\FO(\sa)}\, d\sa
\,,
\qquad 
E\geq \bar E_{2j}\,,
\ 0\leq j< N\,,\,.
\label{ancomarzio}
\end{eqnarray}
\begin{remark}\label{calamaro}
 Note that in $\bar \act_n^{(2j-1),\pm}(E)$
 the interval of integration does depend on 
 $E$, while this does not happen 
 in $\bar \act_n^{(2j),\pm}(E).$
 This fact makes this last case  simpler. 
 In particular $\bar \act_n^{(2j),\pm}(E)$
 are defined for every $E\in \C_* +\bar E_{2j}.$
\end{remark}

\noindent
In view of \eqref{sunday}
and \eqref{ancomarzio} we get
\begin{eqnarray}\label{citola}
\bar \act_n^{(2j-1)}(E)
&=&
\bar \act_n^{(2j-1),-}(E)+
\bar \act_n^{(2j-1),+}(E)\,,
\qquad {\rm for }\ \ \ 1\leq j\leq N
\\
\bar \act_n^{(2j)}(E)
&=&
\bar \act_n^{(2j_- +1),-}(E)+
\bar \act_n^{(2j_+ -1),+}(E)
+
\sum_{j_- < i < j_+}
\bar \act_n^{(2i),-}(E)+
\bar \act_n^{(2i),+}(E)
\,,
\nonumber
\\
&&\qquad\qquad\qquad\qquad
\qquad\qquad\qquad\qquad
\qquad\qquad{\rm for }\ \ \ 1\leq j< N\,,
\nonumber
\\
2\bar \act_n^{(2N)}(E)
&=&
\bar \act_n^{(0),+}(E)+
\bar \act_n^{(2N),-}(E)
+
\sum_{1\leq i\leq N-1}
\bar \act_n^{(2i),-}(E)+
\bar \act_n^{(2i),+}(E)
\,.
\nonumber
\end{eqnarray}
We want to express the functions
$\bar \act_n^{(i)}(E)$ 
in terms of the rescaled potentials
$\breve \FO_{i}$
(defined in \eqref{bastard}).
Given $\FO$ as in \eqref{goffredo}
we set, for $1\leq i\leq 2N,$ 
\begin{equation}\label{remo}
\FO_i:=\FO_{\big| [\sa_{i-1},\sa_i]}\,.
\end{equation}
Recalling \eqref{bastard} we get
\begin{equation}\label{bastard2}
\breve{\FO}_i= \bar \l_i \circ \FO_i\circ \bar  \g_i \,.
\end{equation}
Recalling \eqref{acrobat}, \eqref{bastard},
\eqref{blindness} we set
\begin{eqnarray}
&&
\breve{\bar \act}_n^{(2j-1),-}(\breve E):=
\frac{1}{\pi}
\int_{\breve{\bar \Sa}_{2j-1}(\breve E )}^1
\sqrt{\breve E-\breve{\FO}_{2j-1}(\breve\sa)}\, d\breve\sa\,,
\nonumber
\\
&&
\breve{\bar \act}_n^{(2j-1),+}(\breve E):=
\frac{1}{\pi}
\int_{-1}^{\breve{\bar \Sa}_{2j}(\breve E )}
\sqrt{\breve E-\breve{\FO}_{2j}(\breve\sa)}\, d\breve\sa\,,
\nonumber
\\
&&
\breve{\bar \act}_n^{(2j),-}(\breve E):=
\frac{1}{\pi}
\int_{-1}^{1}
\sqrt{\breve E-\breve{\FO}_{2j}(\breve\sa)}\, d\breve\sa\,,
\nonumber
\\
&&
\breve{\bar \act}_n^{(2j),+}(\breve E):=
\frac{1}{\pi}
\int_{-1}^{1}
\sqrt{\breve E-\breve{\FO}_{2j+1}(\breve\sa)}\, d\breve\sa\,.
\label{corneliabar}
\end{eqnarray}

Then
\begin{eqnarray}
\bar \act_n^{(2j-1),-}(E)
&=&
\frac{\bar\sa_{2j-1}-\bar\sa_{2j-2}}{2}
\sqrt{\frac{\bar E_{2j-2}-\bar E_{2j-1}}{2}}
\breve{\bar \act}_n^{(2j-1),-}
\big(\bar\l_{2j-1}(E)\big)\,,
\nonumber
\\
\bar \act_n^{(2j-1),+}(E)
&=&
\frac{\bar\sa_{2j}-\bar\sa_{2j-1}}{2}
\sqrt{\frac{\bar E_{2j}-\bar E_{2j-1}}{2}}
\breve{\bar \act}_n^{(2j-1),+}
\big(\bar\l_{2j}(E)\big)\,,
\nonumber
\\
\bar \act_n^{(2j),-}(E)
&=&
\frac{\bar\sa_{2j}-\bar\sa_{2j-1}}{2}
\sqrt{\frac{\bar E_{2j}-\bar E_{2j-1}}{2}}
\breve{\bar \act}_n^{(2j),-}
\big(\bar\l_{2j}(E)\big)\,,
\nonumber
\\
\bar \act_n^{(2j),+}(E)
&=&
\frac{\bar\sa_{2j+1}-\bar\sa_{2j}}{2}
\sqrt{\frac{\bar E_{2j}-\bar E_{2j+1}}{2}}
\breve{\bar \act}_n^{(2j),+}
\big(\bar\l_{2j+1}(E)\big)
\,.
\label{valeriabar}
\end{eqnarray}

\begin{lemma}\label{mardibering}
 The functions $\breve{\bar \act}_n^{(2j-1),\pm}(\breve E)$ in\footnote{$\breve{\mathtt r}$ was defined in \eqref{coppe}
 and recall \eqref{bufala}.} \eqref{corneliabar} has holomorphic extension on $ \O_{\breve{\mathtt r}}$. Moreover,
 setting 
\begin{equation}\label{spade}
\breve E(t):=\breve E-(\breve E+1 )t
=-t  + (1-t)\breve E\,, \qquad
\text{for}\ \ 0\leq t\leq 1\,,
\end{equation} 
the following formulas hold:
\begin{eqnarray}
&&
\breve{\bar \act}_n^{(2j-1),-}(\breve E)=
-\frac{(\breve E+1)^{3/2}}{\pi}
\int_{0}^{1}
\sqrt{t}
\frac{
1
}
{\partial_{\breve \sa}
\breve{\FO}_{2j-1}\Big(\breve  {\bar\Sa}_{2j-1}\big(\breve E(t) \big)\Big)}
\, dt\,,
\nonumber
\\
&&
\breve{\bar\act}_n^{(2j-1),+}(\breve E)=
\frac{(\breve E+1)^{3/2}}{\pi}
\int_{0}^{1}
\sqrt{t}
\frac{
1
}
{\partial_{\breve \sa}
\breve{\FO}_{2j}\Big(\breve{\bar\Sa}_{2j}\big(\breve E(t) \big)\Big)}
\, dt
\label{cornelia2bar}
\end{eqnarray}
and, for the derivatives,
\begin{eqnarray}
\partial_{\breve E}\breve{\bar \act}_n^{(2j-1),-}(\breve E)
&=&
-\frac{\sqrt{\breve E+1}}{2\pi}
\int_{0}^{1}
\frac{1}{\sqrt{t}}
\frac{
1
}
{\partial_{\breve \sa}
\breve{\FO}_{2j-1}\Big(\breve  
{\bar\Sa}_{2j-1}\big(\breve E(t) \big)\Big)}
\, dt\,,
\nonumber
\\
\partial_{\breve E}\breve{\bar \act}_n^{(2j-1),+}(\breve E)
&=&
\frac{\sqrt{\breve E+1}}{2\pi}
\int_{0}^{1}
\frac{1}{\sqrt{t}}
\frac{
1
}
{\partial_{\breve \sa}
\breve{\FO}_{2j}\Big
(\breve{\bar \Sa}_{2j}\big(\breve E(t) \big)\Big)}
\, dt\,.
\label{Dcornelia2bar}
\end{eqnarray}
\end{lemma}
\proof
Note that 
$\breve E(t)$ defined in \eqref{spade} describes the segment joining
$\breve E$ and $-1.$
Moreover
\begin{equation}\label{tronodidenari}
\breve E\in \O_\rho
\quad
\Longrightarrow
\quad
\breve E(t)\in \O_\rho\,,
\quad \forall\, \rho>0\,,\ \  \forall\, 0\leq t \leq 1\,.
\end{equation}
Then, for  $\breve E\in  \O_{\breve{\mathtt r}}$
we have that 
$\breve E(t)\in  \O_{\breve{\mathtt r}}$
and, recalling Lemma \ref{lem:coppe}, we also get
that   the following  changes of variables
are well defined:
\begin{eqnarray}
&&\breve\sa = \breve{\bar\Sa}_{2j-1}\big(
\breve E(t) \big)\in [-1,1]_{\rho_\star/2}
\,, \qquad  0\leq t\leq 1\,, 
\nonumber
\\
&&\breve\sa = \breve{\bar\Sa}_{2j}\big(
\breve E(t) \big)\in [-1,1]_{\rho_\star/2}\,,
\qquad
0\leq t\leq 1\,,
\label{mammolabar}
\end{eqnarray}
with $\rho_\star$ defined in \eqref{coppe}.
By \eqref{blu} we get
$$
t=\frac{\breve E-\breve{\FO}_{2j-1}(\breve \sa)}{\breve E+1 }\,,
\qquad
t=\frac{\breve E-\breve{\FO}_{2j}(\breve \sa)}{\breve E+1 }
\,,
$$
and
$$
d\breve \sa=-\frac{\breve E+1}{\partial_{\breve \sa}
\breve{\FO}_{2j-1}
\Big(\breve  {\bar\Sa}_{2j-1}\big(\breve E(t) \big)\Big)} \, dt \qquad
\text{or} \qquad
d\breve \sa=-\frac{\breve E+1}{\partial_{\breve \sa}
\breve{\FO}_{2j}\Big(\breve{\bar\Sa}_{2j}\big(\breve E(t) \big)\Big)} \, dt \,,
$$
(deriving \eqref{blu}).
Making the changes of variables 
\eqref{mammolabar}
 in the first and second integral
in \eqref{corneliabar} respectively, we get \eqref{cornelia2bar}.

We now show  that the functions in \eqref{cornelia2bar}  are holomorphic 
in $\O_{\breve{\mathtt r}}$
and that the expressions in \eqref{Dcornelia2bar}
hold.
In order to prove the claim we first note that, for $t\sim 1$ (namely $\breve E(t)\sim -1$)
\begin{equation}\label{notte}
\breve{\bar\Sa}_{2j}\big(\breve E(t) \big)=
- 1+\sqrt{(1-t)(1+\breve E)}\  \breve{\bar\Sa}_{2j,-}\big( \sqrt{(1-t)(1+\breve E)}\big)
\end{equation}
by \eqref{giappone}. 
This implies that the function
\begin{equation}\label{nottatona}
\breve E \ \mapsto\ 
\int_{0}^{1}
\frac{
\sqrt{t}
}
{\partial_{\breve \sa}
\breve{\FO}_{2j}\Big(\breve{\bar\Sa}_{2j}\big(\breve E(t) \big)\Big)}
\, dt
\end{equation}
is holomorphic in $\O_{\breve{\mathtt r}}$, since, for every $\breve E_0\in\O_{\breve{\mathtt r}}$
and every closed ball $B$ centered in $\breve E_0$ and   contained in 
$\O_{\breve{\mathtt r}},$
the function
$$
\frac{d}{d\breve E}
\frac{
\sqrt{t}
}
{\partial_{\breve \sa}
\breve{\FO}_{2j}\Big(\breve{\bar\Sa}_{2j}\big(\breve E(t) \big)\Big)}\, =\,
\frac{
\sqrt{t}(t-1)
\partial_{\breve \sa\breve \sa}
\breve{\FO}_{2j}\Big(\breve{\bar\Sa}_{2j}\big(\breve E(t) \big)\Big)
}
{\Big(\partial_{\breve \sa}
\breve{\FO}_{2j}\Big(\breve{\bar\Sa}_{2j}\big(\breve E(t) \big)\Big)
\Big)^3}
$$
is dominated\footnote{This follows by \eqref{fiodena} and \eqref{pajata}. In particular note that the denominator
in the r.h.s. vanishes only for $t=1$ where it behaves as
$(1-t)^{3/2}$ (recall \eqref{notte}); then the whole function
behaves as $(1-t)^{-1/2}$, which is in $\mathcal L^1$.} 
by an $\mathcal L^1$-function
uniformly on $B.$
Then, by the Lebesgue's theorem, the function in \eqref{nottatona}
is holomorphic in $B$ and we can exchange the derivative with the integral 
obtaining that
$$
\frac{d}{d\breve E}\,
\int_0^1 \,
\frac{
\sqrt{t}
}
{\partial_{\breve \sa}
\breve{\FO}_{2j}\Big(\breve{\bar\Sa}_{2j}\big(\breve E(t) \big)\Big)}\, dt\, =\,
\int_0^1
\frac{
\sqrt{t}(t-1)
\partial_{\breve \sa\breve \sa}
\breve{\FO}_{2j}\Big(\breve{\bar\Sa}_{2j}\big(\breve E(t) \big)\Big)
}
{\Big(\partial_{\breve \sa}
\breve{\FO}_{2j}\Big(\breve{\bar\Sa}_{2j}\big(\breve E(t) \big)\Big)
\Big)^3}\, dt\,.
$$
As a direct consequence we have that
 $\breve{\bar\act}_n^{(2j-1),+}$ in \eqref{cornelia2bar}
 is holomorphic in $\O_{\breve{\mathtt r}}$ 
and
\begin{eqnarray*}
&&\partial_{\breve E}\breve{\bar\act}_n^{(2j-1),+}(\breve E)
\\
&=&
\frac{1}{2\pi}
\int_{0}^{1}
\left(
\frac{
3\sqrt{t}
\sqrt{\breve E+1}
}
{\partial_{\breve \sa}
\breve{\FO}_{2j}\Big(\breve{\bar\Sa}_{2j}\big(\breve E(t) \big)\Big)}
+
\frac{
2\sqrt{t}(t-1)
(\breve E+1)^{3/2}
\partial_{\breve \sa\breve \sa}
\breve{\FO}_{2j}\Big(\breve{\bar\Sa}_{2j}\big(\breve E(t) \big)\Big)
}
{\Big(\partial_{\breve \sa}
\breve{\FO}_{2j}\Big(\breve{\bar\Sa}_{2j}\big(\breve E(t) \big)\Big)
\Big)^3}
\right)
\, dt
\\
&=&
\frac{\sqrt{\breve E+1}}{2\pi}
\int_{0}^{1}
\frac{1}{\sqrt{t}}
\frac{
1
}
{\partial_{\breve \sa}
\breve{\FO}_{2j}\Big
(\breve{\bar \Sa}_{2j}\big(\breve E(t) \big)\Big)}
\, dt
\nonumber
\\
&&
+\
\frac{\sqrt{\breve E+1}}{\pi}
\int_{0}^{1}
\frac{d}{dt}
\left(
\frac{
\sqrt{t}(t-1)
}
{\partial_{\breve \sa}
\breve{\FO}_{2j}\Big(\breve{\bar\Sa}_{2j}\big(\breve E(t) \big)\Big)}
\right)
\, dt\,,
\end{eqnarray*}
where the last integral vanishes by the fundamental theorem of calculus\footnote{
By the above considerations 
$\partial_{\breve \sa}
\breve{\FO}_{2j}\Big(\breve{\bar\Sa}_{2j}\big(\breve E(t) \big)\Big)\sim \sqrt{1-t}$
for $t\sim 1.$}
proving (the second\footnote{The proof of the first one is analogous.} expression in) \eqref{Dcornelia2bar}.
\eproof

\begin{definition}
 Given  $A\subset\C$ and $r>0$ we define\footnote{$\ii=\sqrt{-1}.$}
\begin{equation}\label{antille}
A_{(r)}:\{ 
z\in\C\ \ :\ \ z=z_0+\ii\, t\,, \ z_0\in A\,,\  |t|<r
\}
\end{equation}

\end{definition}

\begin{lemma}\label{sargassibar}
 $\bar{\act}_n^{(2j-1),+}(E),$ respectively 
 $\bar{\act}_n^{(2j-1),-}(E),$
  in \eqref{ancomarzio} has holomorphic extension on
  $\bar \l_{2j}^{-1}(\O_{\breve{\mathtt r}}),$
respectively
  $\bar \l_{2j-1}^{-1}(\O_{\breve{\mathtt r}}).$
  Moreover
  \begin{eqnarray}
\bar \l_{2j}^{-1}(\O_{\breve{\mathtt r}/2})
&\supset&
(\bar E_{2j-1},\bar E_{2j})_{(\mathtt r_1)}\,,
\nonumber
\\
\bar \l_{2j-1}^{-1}(\O_{\breve{\mathtt r}/2})
&\supset&
(\bar E_{2j-1},\bar E_{2j-2})_{(\mathtt r_1)}\,,
\label{tortuga}
\end{eqnarray}
with
\begin{equation}\label{caraibi}
\mathtt r_1:=\frac{\b \breve{\mathtt r}}{4}\,.
\end{equation}
\end{lemma}
\proof
The first part is a direct consequence of  \eqref{valeriabar} and of Lemma \ref{mardibering}. The inclusions in \eqref{tortuga}
follow by \eqref{acrobat} and \eqref{carciofino}.
\eproof

\begin{remark}\label{telemann2}
 Note that making the inverses  of the change of 
 variables  \eqref{mammolabar}
 in the expression \eqref{Dcornelia2bar}
  we get
 \begin{eqnarray}
\partial_{\breve E}\breve{\bar \act}_n^{(2j-1),-}(\breve E)
&=&
\frac{1}{2\pi}
\int_{\breve{ \bar\Sa}_{2j-1}(\breve E )}^1
\frac{
1
}{\sqrt{\breve E-\breve{\FO}_{2j-1}(\breve\sa)}}
\, d\breve\sa\,,
\nonumber
\\
\partial_{\breve E}\breve{\bar \act}_n^{(2j-1),+}(\breve E)
&=&
\frac{1}{2\pi}
\int_{-1}^{\breve{\bar \Sa}_{2j}(\breve E )}
\frac{
1
}
{\sqrt{\breve E-\breve{\FO}_{2j}(\breve\sa)}}
\, d\breve\sa
\,.
\label{saval}
\end{eqnarray}
 showing also that the formal derivation in \eqref{cippone}
 is actually correct.
 Indeed, recalling \eqref{blindness},\eqref{bastard2},\eqref{acrobat}, and making the change of variables
 $\breve\sa=\bar\g_{2j-1}^{-1}(\sa)$ and
$\breve\sa=\bar\g_{2j}^{-1}(\sa)$ 
 in the first and in the second integral in \eqref{saval}, respectively,
 and summing the results, we get  \eqref{cippone}.
\end{remark}

\begin{lemma}
The functions 
$\bar \act_n^{(2j),\pm}(E)$
are holomorphic for
$E\in \C_*+\bar E_{2j}.$
Their derivatives are 
 \begin{eqnarray}
\partial_E\bar \act_n^{(2j),-}(E)
&:=&
\frac{1}{2\pi}
\int_{\bar\sa_{2j-1}}^{\bar \sa_{2j}}
\frac{1}{\sqrt{E-\FO(\sa)}}\, d\sa\,,
		\nonumber
\\
\partial_E\bar \act_n^{(2j),+}(E)
&:=&
\frac{1}{2\pi}
\int_{\bar\sa_{2j}}^{\bar \sa_{2j+1}}
\frac{1}{\sqrt{E-\FO(\sa)}}\, d\sa
\,.
\label{pigiamini}
\end{eqnarray}
Analogously
the functions 
$\breve{\bar \act}_n^{(2j),\pm}(\breve E)$
are holomorphic for
$\breve E\in \C_*+1.$
Their derivatives are 
 \begin{eqnarray}
&&
\partial_{\breve E}\breve{\bar \act}_n^{(2j),-}(\breve E):=
\frac{1}{2\pi}
\int_{-1}^{1}
\frac{1}{\sqrt{\breve E-\breve{\FO}_{2j}(\breve\sa)}}\, d\breve\sa\,,
\nonumber
\\
&&
\partial_{\breve E}\breve{\bar \act}_n^{(2j),+}(\breve E):=
\frac{1}{2\pi}
\int_{-1}^{1}
\frac{1}{\sqrt{\breve E-\breve{\FO}_{2j+1}(\breve\sa)}}\, d\breve\sa\,.
\label{pigiamini2}
\end{eqnarray}
\end{lemma}
\proof
We simply note that for $E\in \C_*+\bar E_{2j}$
we have that
$$
E-\FO(\sa)\in \C_*\,, \qquad \forall \bar\sa_{2j-1}\leq \sa\leq
\bar\sa_{2j}\quad {\rm and}\quad
\forall \bar\sa_{2j+1}\leq \sa\leq
\bar\sa_{2j}
$$
in particular $\min_{\bar\sa_{2j-1}\leq \sa\leq
\bar\sa_{2j}} |E-\FO(\sa)|>0$ and 
$\min_{\bar\sa_{2j+1}\leq \sa\leq
\bar\sa_{2j}} |E-\FO(\sa)|>0.$
So we can derive inside the integral obtaining \eqref{pigiamini}.
The case of $\breve{\bar \act}_n^{(2j),\pm}(\breve E)$ 
is analogous.
\eproof


\begin{lemma}\label{mardibering+}
For $\breve E\in \C_*+1$
 the following formulas hold:
\begin{eqnarray}
&&
\breve{\bar \act}_n^{(2j),-}(\breve E)=
\frac{1}{\pi}
\int_{-1}^{1}
\frac{
\sqrt{\breve E-y}
}
{\partial_{\breve \sa}
\breve{\FO}_{2j}\Big(\breve  {\bar\Sa}_{2j}(y)\Big)}
\, dy\,,
\nonumber
\\
&&
\breve{\bar\act}_n^{(2j),+}(\breve E)
=-
\frac{1}{\pi}
\int_{-1}^{1}
\frac{
\sqrt{\breve E-y}
}
{\partial_{\breve \sa}
\breve{\FO}_{2j+1}\Big(\breve  {\bar\Sa}_{2j+1}(y)\Big)}
\, dy
\label{cornelia2bar+}
\end{eqnarray}
and, for the derivatives,
\begin{eqnarray}
&&
\partial_{\breve E}
\breve{\bar \act}_n^{(2j),-}(\breve E)=
\frac{1}{2\pi}
\int_{-1}^{1}
\frac{
1
}
{
\sqrt{\breve E-y}\,
\partial_{\breve \sa}
\breve{\FO}_{2j}\Big(\breve  {\bar\Sa}_{2j}(y)\Big)}
\, dy\,,
\nonumber
\\
&&
\partial_{\breve E}
\breve{\bar\act}_n^{(2j),+}(\breve E)
=-
\frac{1}{2\pi}
\int_{-1}^{1}
\frac{
1
}
{
\sqrt{\breve E-y}\, 
\partial_{\breve \sa}
\breve{\FO}_{2j+1}\Big(\breve  {\bar\Sa}_{2j+1}(y)\Big)}
\, dy
\,.
\label{Dcornelia2bar+}
\end{eqnarray}
\end{lemma}
\proof
Recalling Lemma \ref{lem:coppe},
for  $y\in  \O_{\breve{\mathtt r}}$
   the following  changes of variables
are well defined:
\begin{eqnarray}
&&\breve\sa = \breve{\bar\Sa}_{2j}\big(
 y \big)\in [-1,1]_{\rho_\star/2}
\,,  
\nonumber
\\
&&\breve\sa = \breve{\bar\Sa}_{2j+1}\big(
 y \big)\in [-1,1]_{\rho_\star/2}\,,
\label{mammolabar+}
\end{eqnarray}
with $\rho_\star$ defined in \eqref{coppe}.
Note that  $ \breve{\bar\Sa}_{2j}(\pm1)=
 \breve{\bar\Sa}_{2j+1}(\mp 1)=\pm 1$
(recall Lemma \ref{cerveza}).
Recalling \eqref{blu} we get
$$
d\breve \sa=\frac{1}{\partial_{\breve \sa}
\breve{\FO}_{2j}
\Big(\breve  {\bar\Sa}_{2j}(y )\Big)} \, dy \qquad
\text{or} \qquad
d\breve \sa=\frac{1}{\partial_{\breve \sa}
\breve{\FO}_{2j+1}\Big(\breve{\bar\Sa}_{2j+1}(y )\Big)} \, dt \,,
$$
(deriving \eqref{blu}).
Making the changes of variables 
\eqref{mammolabar+}
 in the third and fourth integral
in \eqref{corneliabar} respectively, we get\footnote{One can easily control that the integrals
in \eqref{cornelia2bar+} absolutely converge.} \eqref{cornelia2bar+}.
Note that if $\breve E\in \C_*+1$ and 
$-1\leq y\leq 1$ then $\breve E-y \in \C_*,
$ so that $\sqrt{\breve E-y}$ is well defined.
Deriving \eqref{cornelia2bar+} we get
\eqref{Dcornelia2bar+}.
\eproof

\medskip

Set
  \begin{equation}\label{vana0}
c_{\FO}:=
\min_{1\leq i\leq 2N-1 } \inf_{E\in (\bar E^{(i)}_-,\bar E^{(i)}_+)}
\partial_E \bar \act_n^{(i)}(E)
\,.
\end{equation}

\begin{lemma}\label{lorien}
 We have that 
 \begin{equation}\label{november}
\partial_E \bar \act_n^{(2j-1),\pm}(E)
	\geq \frac{s_0}{2\pi\sqrt\morse}\,,
	\qquad 
	\forall\,
\bar E_{2j-1}<E<\bar E_{2j-1\pm 1}
\end{equation}
and
\begin{equation}\label{rain}
\partial_E \bar \act_n^{(2j),\pm}(E)
\geq \frac{1}{2\pi\sqrt{E-\bar E_{2j-1}}}
\sqrt{\frac{\b s_0^3}{3 \morse}}\,,
\qquad
\forall\,
E>\bar E_{2j}\,.
\end{equation}
As a consequence
\begin{equation}\label{galadriel}
c_{\FO}\geq \frac{\sqrt\b s_0^{3/2}}{32\morse}\,.
\end{equation}
Moreover
\begin{equation}\label{galadriel2}
\partial_E \bar \act_n^{(2N)}(E)
\geq \frac{1}{2\sqrt{E+\morse}}
\,,
\qquad
\forall\,
E>\bar E_{2N}\,.
\end{equation}
\end{lemma}
\proof
First we note that, since
$$
\FO(\bar\sa_{2j-1}+\sa)-\FO(\bar\sa_{2j-1})
=
\FO(\bar\sa_{2j-1}+\sa)-\bar E_{2j-1}
\leq \frac{\morse}{s_0^2}\sa^2
$$
for every $\sa$, then
\begin{equation}\label{lothlorien}
\bar\Sa_{2j}(E)-\bar\sa_{2j-1}\,,\ 
\bar\sa_{2j-1}-\bar\Sa_{2j-1}(E)\,\geq
\frac{s_0}{\sqrt\morse}\sqrt{E-\bar E_{2j-1}}\,.
\end{equation}
Therefore\footnote{Recall
\eqref{cippone} and \eqref{ancomarzio}.}
$$
\partial_E \bar \act_n^{(2j-1),+}(E)
=
\frac{1}{2\pi}
\int_{\bar\sa_{2j-1}}^{\bar \Sa_{2j}(E )}
\frac{1}{\sqrt{E-\FO(\sa)}}\, d\sa
\geq
\frac{1}{2\pi}
\int_{\bar\sa_{2j-1}}^{\bar \Sa_{2j}(E )}
\frac{1}{\sqrt{E-\bar E_{2j-1}}}\, d\sa
\geq \frac{s_0}{2\pi\sqrt\morse}
$$
for $\bar E_{2j-1}<E<\bar E_{2j}.$
The estimates for 
$\partial_E \bar \act_n^{(2j-1),-}$
is analogous.
\\
For $\bar\sa_{2j-1}\leq \sa\leq \bar\sa_{2j}$ 
we have
$$
E-\FO(\sa)\leq E-\bar E_{2j-1}\,,
$$
then
 we get
 $$
\partial_E \bar \act_n^{(2j),-}(E)
=
\frac{1}{2\pi}
\int_{\bar\sa_{2j-1}}^{\bar \sa_{2j}}
\frac{1}{\sqrt{E-\FO(\sa)}}\, d\sa
\geq
\frac{1}{2\pi}
\int_{\bar\sa_{2j-1}}^{\bar \sa_{2j}}
\frac{1}{\sqrt{E-\bar E_{2j-1}}}\, d\sa
$$
and \eqref{rain}
follows\footnote{The estimates for 
$\partial_E \bar \act_n^{(2j),+}$
is analogous.}   by \eqref{draghetto}.
\\
 \eqref{galadriel} follows by \eqref{november},
\eqref{rain}, \eqref{harlock} and \eqref{citola}.
\\
Finally \eqref{galadriel2} directly follows by
\eqref{sundayD} and \eqref{goffredo}.
\eproof

Recallig \eqref{gruffalo},\eqref{vana0} and \eqref{galadriel}
 we have that, for $1\leq i< 2N,$
the function
$$
\bar \act_n^{(i)}: E\in  [\bar E^{(i)}_-,\bar E^{(i)}_+]\mapsto  \bar \act_n^{(i)}(E)
$$
is strictly monotone increasing and, therefore, invertible.
Its inverse
\beq{mortazza}
\bar E^{(i)}: 
\act_n\in  [\bar\acci^{(i)}_-,\bar\acci^{(i)}_+]
\mapsto  \bar E^{(i)}(\act_n)\in 
[\bar E^{(i)}_-,\bar E^{(i)}_+]
\,,
\eeq
where
\begin{equation}\label{torpignattara}
\bar\acci^{(i)}_\pm:=\bar \act_n^{(i)}(\bar E^{(i)}_\pm)\,,
\ \ {\rm for}\ \ 1\leq i<2N\,,\qquad
\bar\acci^{(2N)}_-:=\bar \act_n^{(2N)}(\bar E^{(2N)}_-)
\,.
\end{equation}
Note that actually
\begin{equation}\label{torpignattara2}
\bar\acci^{(2j-1)}_-:=0\,,\qquad
\forall\, 1\leq j\leq N\,.
\end{equation}
Moreover, by \eqref{galadriel2}, also 
the function
$$
\bar \act_n^{(2N)}=
-\bar \act_n^{(0)}
\ :\  E\in  [\bar E^{(2N)}_-,\infty)\mapsto  
[\bar\acci^{(2N)}_-,\infty)
$$
is strictly monotone increasing and, therefore, invertible.
We have the corresponding inverse functions
\beq{mortazzabis}
\bar E^{(2N)}: 
\act_n\in   [\bar\acci^{(2N)}_-,\infty)\mapsto  
[\bar E^{(2N)}_-,\infty)
\,,\quad
\bar E^{(0)}: 
\act_n\in   (-\infty,-\bar\acci^{(2N)}_-]\mapsto  
[\bar E^{(2N)}_-,\infty)\,.
\eeq

\subsection{The domains of definitions of the action functions for the perturbed Hamiltonian}

Let us consider now a real analytic Hamiltonian
\beq{Pasqua}
\boxed{\Hpend(p,q_n)=\big(1+ b (p,q_n)\big) \big(p_n-{\mathtt P}^*_n(\hat p)\big)^2  
  + \Gm (\hat p,q_n )}
\eeq
with holomorphic extension on
$$
(p,q_n)=(\hat p,P_n,q_n)\in 
\hat D_{r_0}\times(-R_0,R_0)_{r_0/2}\times
\T_{s_0}\,.
$$
Assume also that for some $\checco\geq 0,$
\begin{equation}\label{Urina3}
{\modulo}\Gm-\FO {\modulo}_{\hat D,r_0,s_0}
\leq \checco
\,,
\qquad
{\modulo}{\mathtt P}^*_n{\modulo}_{\hat D,r_0}\leq \checco/r_0
\end{equation}
with $\FO$ as in \eqref{goffredo}-\eqref{maieli}
and
\begin{equation}\label{Straussbis}
{\modulo}b{\modulo}_*\leq
\frac{\checco}{r_0^2}
\,,
\qquad
{\modulo}\partial_{p_n}
b{\modulo}_*\leq
\frac{\checco}{r_0^3}
\,,
\qquad
 {\modulo} p_n \cdot  b(p,q_n){\modulo}_*\leq
\frac{\checco}{r_0}\,,
\qquad
 {\modulo} p_n \cdot  \partial_{p_n}b(p,q_n){\modulo}_*\leq
\frac{\checco}{r_0^2}
\,,
\end{equation}
where\footnote{Recall \eqref{wagner}.} 
\begin{equation}\label{Wagner}
|\cdot|_*:=\sup_{\hat D_{r_0}\times(-R_0,R_0)_{r_0/2}\times
\T_{s_0}}|\cdot|\,. 
\end{equation}
\begin{remark}
Note that the Hamiltonian $\Hpend$ in \eqref{pasqua}
is of the form \eqref{Pasqua}-\eqref{Straussbis}
with $\checco=100\checco_0$
by \eqref{urina} and \eqref{straussbis}.
\end{remark}

Assume\footnote{See \eqref{adcazzum} below.}  
\begin{equation}\label{legna}
 R_0\geq 2\sqrt \morse\,,
\end{equation} 
with $\morse$ defined in \eqref{goffredo}.
Note that 
\begin{equation}\label{pizzoli}
\checco\leq \frac{r_0^2}{32}\,,
\end{equation}
as it is implied by \eqref{genesis}.

\medskip

Note that $\hat p$ is a parameter with no dynamical meaning since its conjugated variable $\hat q$
does not appear in the Hamiltonian.
Note that
$$
\nabla_{(p_n,q_n)}\Hpend (p,q_n )=0
\qquad \Longrightarrow\qquad
p_n={\mathtt P}^*_n(\hat p)
\,,\qquad \partial_{q_n}\Gm(p,q_n )=0\,.
$$

\bigskip

Let us consider the
  equation
\begin{equation}\label{meaux}
 {\mathtt P}^*_n(\hat p)+\frac{z}{\sqrt{1+ b(p,q_n )}}-p_n=0\,.
\end{equation}
Since
${\modulo}b{\modulo}_*\leq
\frac{\checco}{r_0^2}\leq \frac12$ by \eqref{Straussbis}
and \eqref{pizzoli}, we have
\begin{equation}\label{aterno}
\Re(1+ b(p,q_n ))\geq \frac12\,,\qquad
\forall\ 
\hat p\in \hat D_{r_0}\,,\ \ 
p_n\in (-R_0,R_0)_{r_0/2}\,,\ \ 
q_n\in \T_{s_0}
\end{equation}
and, therefore, 
$\sqrt{1+ b(p,q_n )}$ is well defined 
on $\hat D_{r_0}\times
(-R_0,R_0)_{r_0/2}
 \times \T_{s_0}.$

\begin{lemma}\label{amatriciana}
Assume \eqref{pizzoli}.
Then there exists a (unique) real analytic function
  $\tilde{\mathcal P}:(-R_0,R_0)_{r_0/4}\times\mathbb T_{s_0}\times \hat D_{r_0}\, \to\, \C$
with
 \begin{equation}\label{3holes}
 |\tilde{\mathcal P}|_\dagger:=
\sup_{z\in (-R_0,R_0)_{r_0/4}} 
{\modulo}\tilde{\mathcal P}(z,	\cdot,\cdot){\modulo}_{\hat D,r_0,s_0}
\leq  \frac{4}{r_0}\checco\leq \frac{r_0}{8}\,,
\end{equation} 
   such that
\begin{equation}\label{spiderman}
p_n=\mathcal P(z,q_n,\hat p ):=
z +\tilde{\mathcal P}(z,q_n,\hat p )
\end{equation}
solves \eqref{meaux}.
Moreover
\begin{equation}\label{piove}
\mathcal P: 
(-R_0,R_0)_{r_0/4}\times\mathbb T_{s_0}\times \hat D_{r_0}
\to (-R_0,R_0)_{r_0/2}
\end{equation}
\end{lemma}
\proof
We first note that 
if $\tilde{\mathcal P}$ satisfies the first inequality in \eqref{3holes},
then, by \eqref{pizzoli}, it also satisfies the second one.
Therefore,
if $z\in (-R_0,R_0)_{r_0/4}$
then
  $z+\tilde{\mathcal P}\in (-R_0,R_0)_{r_0/2}$
  and \eqref{piove} holds.
Let us define 
 $\tilde{\mathcal P}=
 \tilde{\mathcal P}(z,q_n,\hat p )$ as the solution of 
the fixed point equation
\begin{equation}\label{meauxbis}
\tilde{\mathcal P}=\Phi(\tilde{\mathcal P}):=
 {\mathtt P}^*_n(\hat p)+\left(\frac{1}{\sqrt{1+ 
 b(\hat p,z+\tilde{\mathcal P},q_n )}}-1\right)z
\end{equation}
in the closed ball $\mathtt B$ of $\tilde{\mathcal P}$ satisfying
\eqref{3holes}.
We immediately see that $\Phi(\mathtt B)\subseteq
\mathtt B,$ since, by\footnote{Recall  the notation in\eqref{Wagner}.} \eqref{Straussbis} 
and \eqref{aterno},
$$
|\Phi(\tilde{\mathcal P})|_\dagger
\leq \frac{\checco}{r_0}+
2|b z|_* \leq \frac{3\checco}{r_0}\,.
$$
Moreover $\Phi$ is a contraction 
since by  the fifth and the third equation in\footnote{Omitting for brevity to write
$\hat p, q_n$ we have, setting
$\xx=\xx(t)=(1-t)\tilde{\mathcal P}' +t\tilde{\mathcal P}$
\begin{eqnarray*}
&&\big(b(z+\tilde{\mathcal P})-
b(z+\tilde{\mathcal P}' )\big)z
=\big(\tilde{\mathcal P}-\tilde{\mathcal P}'\big)\int_0^1 \partial_{p_n} 
 b(z+\xx )\cdot z \,dt
 \\
 &&=
 \big(\tilde{\mathcal P}-\tilde{\mathcal P}'\big)
 \left[
 \int_0^1 \partial_{p_n} 
b(z+\xx )\cdot (z+	\xx) \,dt
\, -\, 
 \int_0^1 \partial_{p_n} 
b(z+\xx )\cdot \xx \,dt
 \right]\,.
\end{eqnarray*}
Finally note that, for every $0\leq t\leq 1$ 
we have $|\xx(t)|_\dagger \leq 8\checco_0/r_0
\leq r_0/4$ by \eqref{3holes} and \eqref{pizzoli}.
}
\eqref{straussbis}
$$
|\Phi(\tilde{\mathcal P})-\Phi(\tilde{\mathcal P}')|_\dagger
\leq 2 \big|\big(b(\hat p,z+\tilde{\mathcal P},q_n )-
b(\hat p,z+\tilde{\mathcal P}',q_n )\big)z\big|_\dagger
\leq \frac{4 \checco}{r_0^2}
|\tilde{\mathcal P}-\tilde{\mathcal P}'|_\dagger
\stackrel{\eqref{pizzoli}}\leq \frac12 
|\tilde{\mathcal P}-\tilde{\mathcal P}'|_\dagger
\,.
$$
Then equation \eqref{meauxbis} is solved  by the Contraction Theorem.
\eproof

Obviously\footnote{For real values of $\hat p,$ $q_n,$ $E$
such that $E\geq\Gm(\hat p,q_n )$.}
\begin{equation}\label{ummagamma}
p_n=\mathcal P\Big(\pm\sqrt{E-\Gm(\hat p,q_n )},q_n,\hat p \Big)\quad
{\rm solves \ (w.r.t.} \ p_n) \quad
\Hpend (\hat p,p_n,q_n )=E\,,
\end{equation}
according to $\pm \big( p_n- {\mathtt P}^*_n(\hat p) \big)\geq 0.$
Note that, for real $q_n,\hat p,$ 
$z\mapsto \mathcal P(z,q_n,\hat p)$ is an increasing function of (real) $z\in (-R_0,R_0)$,
since 
\begin{equation}\label{fragile}
\partial_z \mathcal P\geq 1-\frac{16\checco}{r_0^2}
\geq \frac12
\end{equation}
by \eqref{pizzoli}-\eqref{spiderman} and Cauchy estimates.
Note also that
$$
\mathcal P(0,q_n,\hat p)= {\mathtt P}^*_n(\hat p)\,.
$$
\Giu

\begin{remark}\label{adpenem}
We will assume that\footnote{Recall \eqref{legna}.}, for real $E$, 
\begin{equation}\label{adcazzum}
E+\morse< R_0^2\,,
\end{equation} 
so that $\mathcal P\Big(\pm\sqrt{E-\Gm(\hat p,q_n )},q_n,\hat p \Big)$
is well defined.
\end{remark}

Outside a zero measure set contained in the set of  critical energies $\{ \Hpend =E_i \},$
$1\leq i\leq 2N,$ 
the phase space
$\hat D\times(-R_0,R_0)\times \T^n$ 
is decomposed in $2N+1$ open connected components
$
\mathcal C^i,$ $0\leq i\leq 2N,$
on which we will define the action-angle transformation.
The
$
\mathcal C^i$ are
normal sets with respect to the variable
$p_n$ and are
defined as follows.
\\
For $i=2j-1$ odd, $1\leq j\leq N,$
\begin{eqnarray}\label{skys}
\mathcal C^{2j-1}   
&:=&
\Big\{ (p,q)\in \hat D\times(-R_0,R_0)\times \T^n \ \  {\rm s.t.}
\nonumber
\\
&&\ \ 
\mathcal P
\Big(-\sqrt{E_{2j_\diamond}(\hat p) -\Gm(\hat p,q_n)},q_n,\hat p \Big)
< p_n<
\mathcal P
\Big(\sqrt{E_{2j_\diamond}(\hat p) -\Gm(\hat p,q_n)},q_n, \hat p \Big)\,,
\nonumber
\\
&&\ \ \Sa_{2j-1}\big(E_{2j_\diamond}(\hat p),   \hat p\big)
< q_n <
\Sa_{2j}\big(E_{2j_\diamond}(\hat p),\hat p   \big)
\Big\}
\nonumber
\\ &&\setminus \ 
\Big\{  
p_n=\mathtt P_n^*(\hat p) \,, q_n=
\sa_{2j-1}(\hat p) 
\Big\}
\end{eqnarray}
where  $j_\diamond$ was defined in \eqref{padula}.
\\
For $i=2j$ even, $1\leq j\leq N-1,$
 \begin{eqnarray}\label{skys2}
\mathcal C^{2j}   
&:=&
\Big\{ (p,q)\in \hat D\times(-R_0,R_0)\times \T^n \ \  {\rm s.t.}
\nonumber
\\
&&\ \ 
\mathcal P
\Big(-\sqrt{E_{2j_*}(\hat p) -\Gm(\hat p,q_n)},q_n,\hat p \Big)
< p_n<
\mathcal P
\Big(\sqrt{E_{2j_*}(\hat p) -\Gm(\hat p,q_n)},q_n, \hat p \Big)\,,
\nonumber
\\
&&\ \ \Sa_{2j_- +1}\big(E_{2j_*}(\hat p),   \hat p\big)
< q_n <
\Sa_{2j_+}\big(E_{2j_*}(\hat p),\hat p   \big)
\Big\}
\nonumber
\\
&&\setminus \Big\{
\mathcal P
\Big(-\sqrt{E_{2j}(\hat p) -\Gm(\hat p,q_n)},q_n,\hat p \Big)
\leq p_n\leq
\mathcal P
\Big(\sqrt{E_{2j}(\hat p) -\Gm(\hat p,q_n)},q_n, \hat p \Big)\,,
\nonumber
\\
&&\ \ \ \  \Sa_{2j_- +1}\big(E_{2j}(\hat p),   \hat p\big)
\leq q_n \leq
\Sa_{2j_+}\big(E_{2j}(\hat p),\hat p   \big)
\Big\}\,,
\end{eqnarray}
where   $j_-,j_+,j_*$ were defined in \eqref{21stcentury}.

Finally
\begin{eqnarray}\label{skys3}
\mathcal C^{2N}   
&:=&
\Big\{
(p,q)\in \hat D\times(-R_0,R_0)\times \T^n \ \  {\rm s.t.}
\nonumber
\\
&&\ \
R_0>p_n>
\mathcal P
\Big(\sqrt{E_{2N}(\hat p) -\Gm(\hat p,q_n)},q_n,\hat p \Big)
\Big\}
\\
\label{skys4}
\mathcal C^{0}   
&:=&
\Big\{
(p,q)\in \hat D\times(-R_0,R_0)\times \T^n \ \  {\rm s.t.}
\nonumber
\\
&&\ \ 
-R_0<p_n<
\mathcal P
\Big(-\sqrt{E_{2N}(\hat p) -\Gm(\hat p,q_n)},q_n,\hat p \Big)
\Big\}
\end{eqnarray}
Note that actually in $\mathcal C^i  $ with $1\leq i<2N,$
$q_n$ is not an angle!

\Giu

\vglue0.5truecm

\centerline{\boxed{\includegraphics[width=6cm]{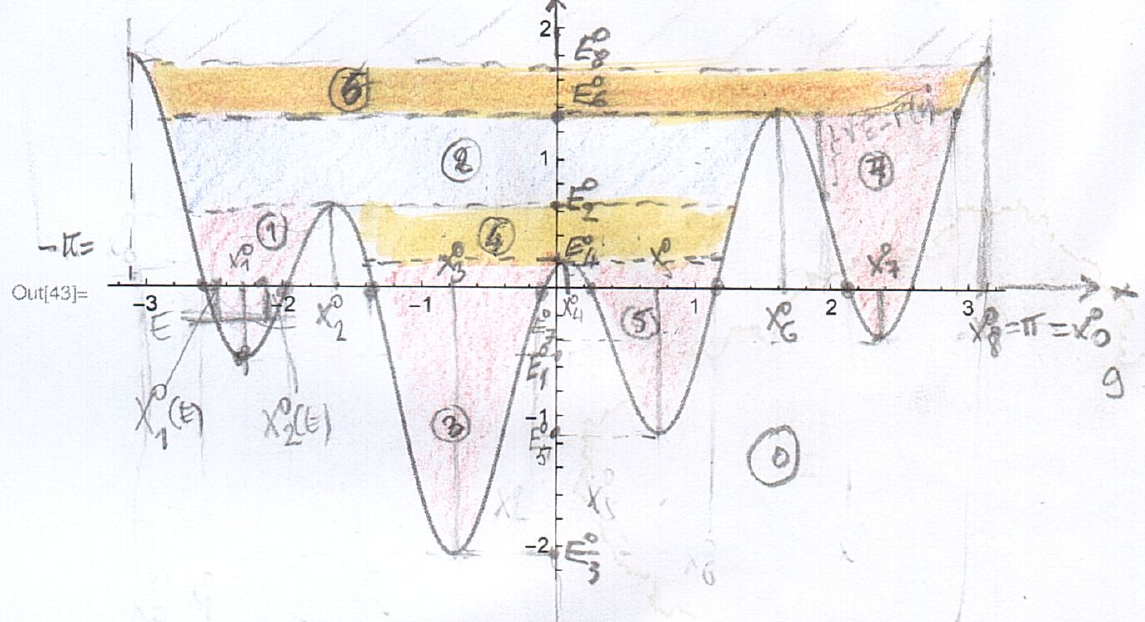}}}\nopagebreak
\centerline{\footnotesize \bf  Figure 1. \it  A Morse potential $\FO$ with 8 critical points }

\vglue0.5truecm

\centerline{\boxed{\includegraphics[width=6cm]{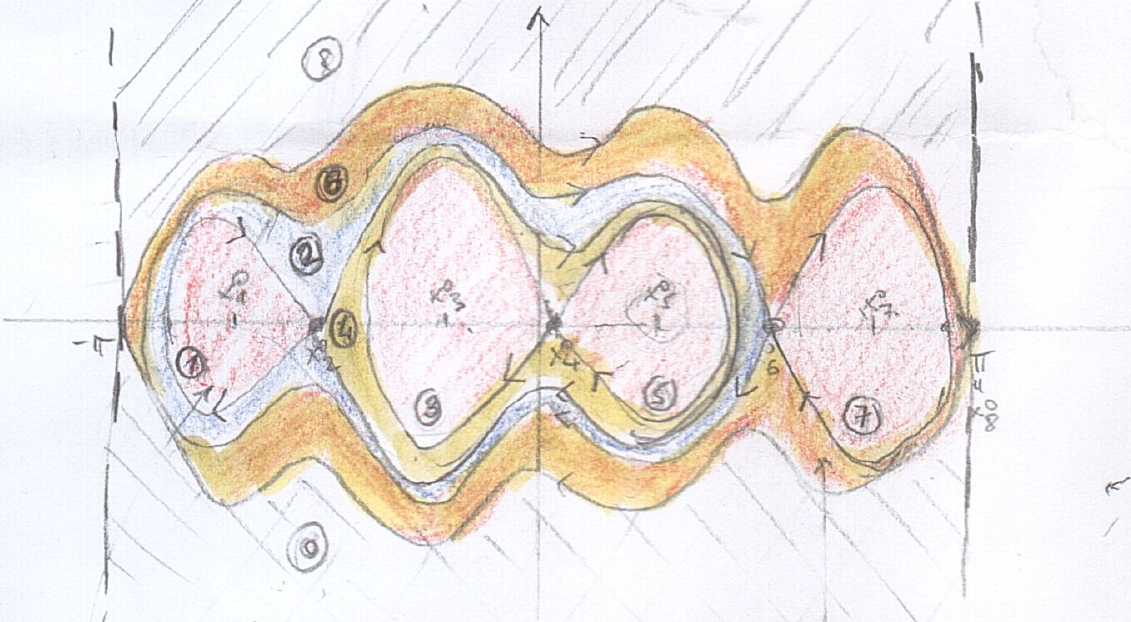}}}\nopagebreak
\centerline{\footnotesize \bf Figure 2. \it  The phase portrait of the above Morse potential}

\vglue0.5truecm

\rem\label{icm2018}
 The phase space regions $\mathcal C^i $ are closely related to the phase portrait of the Morse potential
$\FO$, which is the limiting potential, as $\checco$ goes to zero, of the potential $\GG$ of $\Hpend$. 
Roughly speaking, $\dst\bigcup_{i=0}^{2N}\mathcal C^i $ is a bounded region around $p_n=0$ of the phase space of $\Hpend$ minus the separatrices (stable/unstable manifolds) issuing from critical hyperbolic lower dimensional manifolds (in the 1D projection: critical hyperbolic points). The suffix  $i$, for  $1\le i \le 2N-1$, labels the internal regions trapped by separatrices, where the actual motion
of the 1D projection $(p_n,q_n)$--system is oscillatory,  while $i=0,2N$ labels the two external regions, where the $q_n$--coordinate rotates (positive velocity for $i=2N$ and negative velocity for $i=0$).
\erem

\subsection{Action variables and their adimensional version}

On the above connected components $\mathcal C^i,$
$0\leq i\leq 2N,$
we want to define action angle variables integrating
$\Hpend .$
We first define the action variables 
as a function of the energy $E$ and of the dummy
variable $\hat p$.

\smallskip

\noindent
Define the analytic function $b_\sharp=
b_\sharp(z,\sa,\hat p)$ as follows:
\begin{equation}\label{4holes}
 b_\sharp(z,\sa,\hat p ):=
\frac{1}{2\sqrt{1+ b\big(\hat p,\mathcal P(z,\sa,\hat p ),\sa \big)}}
+
\frac{1}{2\sqrt{1+ b\big(\hat p,\mathcal P(-z,\sa,\hat p ),\sa \big)}}
-1
\end{equation}
($b$ defined in \eqref{Pasqua}).
Note that $b_\sharp$ is {\sl even} w.r.t. $z$ and 
that\footnote{Note that, if $\Re (1+b)\geq 1/2$ (recall \eqref{aterno}), then $|(1+b)^{-1/2}-1|\leq |b|.$}
\begin{equation}\label{rodimento}
\sup_{z\in (-R_0,R_0)_{r_0/4}}  {\modulo}b_\sharp(z,\hat p,\sa){\modulo}_{\hat D,r_0,s_0}
\leq |b|_*\leq \frac{\checco}{r_0^2}
\,,
\end{equation}
by \eqref{Straussbis}.

\begin{remark}\label{chimay1}
 In the following, for brevity, we will often
 omit to write the immaterial
dependence on $\hat p$ and/or on $\sa$.
\end{remark}

Moreover we have
$$
\partial_z b_\sharp(z)
=\frac{-\partial_{p_n}b(\mathcal P(z))
\partial_z \mathcal P(z)
}{4(1+b(\mathcal P(z)))^{3/2}}
+
\frac{\partial_{p_n}b(\mathcal P(-z))
\partial_z \mathcal P(-z)
}{4(1+b(\mathcal P(-z)))^{3/2}}
$$
and, since  
$z \partial_{p_n} 
b(\mathcal P(z))=\mathcal P(z)\partial_{p_n} 
b(\mathcal P(z))-\tilde{\mathcal P}(z)
\partial_{p_n} 
b(\mathcal P(z)),$ we get
$$
\sup_{z\in (-R_0,R_0)_{r_0/4}}  {\modulo}
z \partial_{p_n} 
b(z){\modulo}_{\hat D,r_0,s_0}
\leq \frac{\checco}{r_0^2}+
\frac{r_0}{8}
\frac{\checco}{r_0^3}
\leq 
\frac{2\checco}{r_0^2}\,.
$$
by \eqref{Straussbis} and \eqref{3holes}.
Therefore, since 
$$
\sup_{z\in (-R_0,R_0)_{r_0/8}}
|\partial_z \mathcal P|_{\hat D,r_0,s_0}
\leq 2
$$
by Cauchy estimates and \eqref{3holes},\eqref{spiderman},
we finally obtain that
\begin{equation}\label{prurito}
\sup_{z\in (-R_0,R_0)_{r_0/8}}
|z\partial_z b_\sharp|_{\hat D,r_0,s_0}
\leq
\frac{8\checco}{r_0^2}\,.
\end{equation}

\begin{lemma}\label{geronimo}
 Assume that $g(z)$ is holomorphic 
 on $[0,R]_r$. If $g$ is even\footnote{Namely
 $g(z)=g(-z)$ for $z$ close to zero.} then
 one can define $G(v)$ holomorphic on  $[0,R^2]_{r^2}$
 such that $G(z^2)=g(z).$
\end{lemma}
\proof
We incidentally note that, since $g$ is even, it is actually 
holomorphic 
 on $[-R,R]_r$. 
Denoting by $D_r(0):=\{ |z|<r\},$
we have that, since
$g$ is holomorphic and even on $D_r(0),$
 $g(z)=\sum_{j\geq 0} a_{2j} z^{2j},$
 where the power series has a radius of convergence
 $\geq r.$ Then $G(v):=\sum_{j\geq 0} a_{2j}v^j$
 has radius of convergence $\geq r^2.$ 
 It remains to define $G$ in the set
 $\O:=[0,R^2]_{r^2}\setminus D_{r^2}(0).$
 Note that $\O\subset \C_*;$
 then we set $G(v):=g(\sqrt v)$ for $v	\in \O,$
 noting that $z:=\sqrt v\in [0,R]_r.$
 This follows noting that if $v\in D_{r^2}(v_0^2),$
 with $v_0\in\R,$ $v_0>r,$ then $\sqrt v\in
 D_r (v_0).$
 This last fact is equivalent to proving that\footnote{This inclusion follows 
 noting that, for every angle $\theta,$
 we have 
 $
 |(v_0 + re^{\ii\theta})^2-v_0^2|\geq r^2.
 $
 The last inequality follows noting that
 it is equivalent to
 $|re^{\ii 2\theta}+2 v_0 e^{\ii\theta}|=
 |re^{\ii \theta}+2 v_0|\geq r$,
 that follows from $v_0>r.$
 }
 $D_{r^2}(v_0^2)\subseteq s(D_r (v_0)),$
 where $s(v):=v^2.$
\eproof

\noindent
Since 
$b_\sharp$ is  even w.r.t. $z,$
by Lemma \ref{geronimo}
we can define  the analytic function
$b_\dag(v,\hat p,\sa)$ such that
\begin{equation}\label{pachelbel}
b_\dag(z^2,\hat p,\sa):=
b_\sharp(z,\hat p,\sa)
\end{equation}
with (recall \eqref{rodimento})
\begin{equation}\label{3holesbis}
\sup_{v\in (0,R_0^2)_{r_0^2/16}}  {\modulo}b_\dag(v,\hat p,\sa){\modulo}_{\hat D,r_0,s_0}\leq \frac{\checco}{r_0^2}\,.
\end{equation}
Moreover, since 
$v\partial_v b_\dag (v)=\frac12 \sqrt v\partial_z
b_\sharp (\sqrt v),$
\begin{equation}\label{3holester}
\sup_{v\in (0,R_0^2)_{r_0^2/64}}  {\modulo}v\partial_v b_\dag (v){\modulo}_{\hat D,r_0,s_0}
\leq \frac{4\checco}{r_0^2}\,,
\end{equation}
by \eqref{prurito}.
\medskip

\noindent
For $i=2j-1$ odd, $1\leq j\leq N,$ 
and $E_{2j-1} <E<E_{2j_\diamond
} $, we set
\begin{eqnarray}\label{musicalbox}
&&
\act_n^{(2j-1)}(E)=\act_n^{(2j-1)}(E, \hat p )
\\
\nonumber
&&
:=
\frac{1}{2\pi}
\int_{\Sa_{2j-1}(E )}^{\Sa_{2j}(E )}
\Big[
\mathcal P
\Big(\sqrt{E-\Gm(\sa)},\sa \Big) 
-
\mathcal P
\Big(-\sqrt{E-\Gm(\sa)},\sa \Big)
\Big]
\, d\sa
\\
\nonumber
&&
\stackrel{\eqref{meaux},\eqref{4holes}}=
\frac{1}{\pi}
\int_{\Sa_{2j-1}(E )}^{\Sa_{2j}(E )}
\sqrt{E-\Gm(\sa)}\Big(
1
+ b_\sharp\big(\sqrt{E-\Gm(\sa)},\sa \big)
\Big)\, d\sa
\\
\nonumber
&&
\stackrel{\eqref{pachelbel}}=
\frac{1}{\pi}
\int_{\Sa_{2j-1}(E )}^{\Sa_{2j}(E )}
\sqrt{E-\Gm(\sa)}\Big(
1
+ b_\dag\big(E-\Gm(\sa),\sa \big)
\Big)\, d\sa\,,
\end{eqnarray}
Analogously, for $i=2j$ even, $1\leq j\leq N-1$
and $E_{2j} <E<E_{2j_*} $, 
we set (recall \eqref{21stcentury})
\begin{eqnarray}\label{musicalbox2}
&&
\act_n^{(2j)}(E)=
\act_n^{(2j)}(E, \hat p )
\\
\nonumber
&&
:=
\frac{1}{2\pi}
\int_{\Sa_{2j_-+1}(E )}^{\Sa_{2j_+}(E )}
\Big[
\mathcal P
\Big(\sqrt{E-\Gm(\sa)},\sa \Big) 
-
\mathcal P
\Big(-\sqrt{E-\Gm(\sa)},\sa \Big)
\Big]
\, d\sa
\\
\nonumber
&&
=\frac{1}{\pi}
\int_{\Sa_{2j_-+1}(E )}^{\Sa_{2j_+}(E )}
\sqrt{E-\Gm(\sa)}\Big(
1
+ b_\sharp\big(\sqrt{E-\Gm(\sa)},\sa \big)
\Big)\, d\sa
\\
\nonumber
&&
=\frac{1}{\pi}
\int_{\Sa_{2j_-+1}(E )}^{\Sa_{2j_+}(E )}
\sqrt{E-\Gm(\sa)}\Big(
1
+ b_\dag\big(E-\Gm(\sa),\sa \big)
\Big)\, d\sa\,.
\end{eqnarray}
Finally for $E>E_{2N} $ we set
\begin{eqnarray}\label{musicalbox3}
\act_n^{(2N)}(E)=
\act_n^{(2N)}(E, \hat p )
&:=&
\frac{1}{2\pi}
\int_{-\pi}^{\pi}
\mathcal P
\Big(\sqrt{E-\Gm(\sa)},\sa \Big) \, d\sa\,,
\\
\label{musicalbox4}
\act_n^{(0)}(E)=
\act_n^{(0)}(E, \hat p )
&:=&
\frac{1}{2\pi}
\int_{-\pi}^{\pi}
\mathcal P
\Big(-\sqrt{E-\Gm(\sa)},\sa \Big) \, d\sa\,.
\end{eqnarray}

Recalling the definition of 
$b_\dag$ in \eqref{pachelbel} we set
\begin{equation}\label{palidoro}
\tilde b(v)=\tilde b(v,\hat p,\sa)
:=b_\dag(v)+2 v \partial_v b_\dag(v)\,,
\end{equation}
with
\begin{equation}\label{3holesquater}
\sup_{v\in (0,R_0^2)_{r_0^2/64}}  {\modulo}
\tilde b (v){\modulo}_{\hat D,r_0,s_0}
\leq \frac{9\checco}{r_0^2}
\leq 
\frac{9}{32}\,,
\end{equation}
by \eqref{3holesbis}, \eqref{3holester} and 
\eqref{pizzoli}.

Recalling Remark \ref{telemann1} 
we have the following

\begin{remark}\label{telemann3}
 We want to show that the functions in
  \eqref{musicalbox}-\eqref{musicalbox4}
 have an analytic extension for complex $E.$
 Moreover, while it is immediate to evaluate the derivative
 of the last two\footnote{And also the second derivatives:
 \begin{eqnarray}
\partial_{EE} \act_n^{(2N)}(E)
&=&
\frac{1}{8\pi}
\int_{-\pi}^{\pi}
\frac{\sqrt{E-\Gm(\sa)}
\partial_{zz}\mathcal P
\Big(\sqrt{E-\Gm(\sa)},\sa \Big)
-
\partial_z\mathcal P
\Big(\sqrt{E-\Gm(\sa)},\sa \Big)}
{(E-\Gm(\sa))^{3/2}}
 \, d\sa\,,\ \ 
\label{lamentino3}
\\
\partial_{EE} \act_n^{(0)}(E)
&=&
\frac{1}{8\pi}
\int_{-\pi}^{\pi}
\frac{\sqrt{E-\Gm(\sa)}
\partial_{zz}\mathcal P
\Big(-\sqrt{E-\Gm(\sa)},\sa \Big)
+
\partial_z\mathcal P
\Big(-\sqrt{E-\Gm(\sa)},\sa \Big)}
{(E-\Gm(\sa))^{3/2}}
 \, d\sa\,. 
 \nonumber
\end{eqnarray}}, 
namely
 \begin{eqnarray}
\partial_E \act_n^{(2N)}(E)
&=&
\frac{1}{4\pi}
\int_{-\pi}^{\pi}
\frac1{\sqrt{E-\Gm(\sa)}}
\partial_z\mathcal P
\Big(\sqrt{E-\Gm(\sa)},\sa \Big) \, d\sa\,,
\nonumber
\\
\partial_E \act_n^{(0)}(E)
&=&
-\frac{1}{4\pi}
\int_{-\pi}^{\pi}
\frac1{\sqrt{E-\Gm(\sa)}}
\partial_z \mathcal P
\Big(-\sqrt{E-\Gm(\sa)},\sa \Big) \, d\sa\,.
\label{lamento3}
\end{eqnarray}
 it is not obvious to justify the formal derivation\footnote{
See Remark \ref{telemann4} below.}
\begin{eqnarray}
\partial_E \act_n^{(2j-1)}(E)
&=&
\frac{1}{2\pi}
\int_{\Sa_{2j-1}(E )}^{\Sa_{2j}(E )}
\frac1{\sqrt{E-\Gm(\sa)}}\Big(
1
+ \tilde b\big(E-\Gm(\sa),\sa \big)
\Big)\, d\sa\,,
\nonumber
\\
\partial_E \act_n^{(2j)}(E)
&=&
\frac{1}{2\pi}
\int_{\Sa_{2j_-+1}(E )}^{\Sa_{2j_+}(E )}
\frac1{\sqrt{E-\Gm(\sa)}}\Big(
1
+ \tilde b\big(E-\Gm(\sa),\sa \big)
\Big)\, d\sa\,,
\label{lamento2}
\end{eqnarray}
\end{remark}


Recalling \eqref{musicalbox} and \eqref{citola} we split
the integral
$
\int_{ \Sa_{2j-1}(E )}^{\Sa_{2j}(E )}
=
\int_{ \Sa_{2j-1}(E )}^{\sa_{2j-1}}+
\int_{\sa_{2j-1}}^{\Sa_{2j}(E )}
$
obtaining
\begin{equation}\label{noce}
 \act_n^{(2j-1)}(E)=
 \act_n^{(2j-1),+}(E)+
 \act_n^{(2j-1),-}(E)\,,
\end{equation}
where, for $1\leq j\leq N,$
\begin{eqnarray}
 \act_n^{(2j-1),-}(E)
 =
  \act_n^{(2j-1),-}(E, \hat\act)
 &:=&
 \frac{1}{\pi}
\int_{ \Sa_{2j-1}(E )}^{\sa_{2j-1}}
\sqrt{E-\Gm(\sa)}\Big(
1
+ b_\dag\big(E-\Gm(\sa),\sa \big)
\Big)\, d\sa\,,
\nonumber
\\
 \act_n^{(2j-1),+}(E)
 =
  \act_n^{(2j-1),+}(E, \hat\act)
 &:=&
 \frac{1}{\pi}
\int_{\sa_{2j-1}}^{\Sa_{2j}(E )}
\sqrt{E-\Gm(\sa)}\Big(
1
+ b_\dag\big(E-\Gm(\sa),\sa \big)
\Big)\, d\sa\,.
\nonumber
\\
\label{ancomarzio2}
\end{eqnarray}
Analogously, recalling \eqref{musicalbox2} and \eqref{citola} 
we split the integral 
\begin{eqnarray}\label{noce2}
\int_{\Sa_{2j_-+1}(E )}^{\Sa_{2j_+}(E )}
&=&
\int_{\Sa_{2j_-+1}(E )}^{\sa_{2j_-+1}}
+
\sum_{j_- < i < j_+}\left(
\int_{\sa_{2i-1}}^{\sa_{2i}}+
\int_{\sa_{2i}}^{\sa_{2i+1}}
\right)
+
\int_{\sa_{2j_+-1}}^{\Sa_{2j_+}(E )}\,,
\\
\nonumber
\act_n^{(2j)}(E)
&=&
\act_n^{(2j_- +1),-}(E)
+
\sum_{j_- < i < j_+}
\Big( \act_n^{(2i),-}(E)+
 \act_n^{(2i),+}(E)\Big)
 +
\act_n^{(2j_+ -1),+}(E)
 \,,
\end{eqnarray}
for $1\leq j<N,$
where
\begin{eqnarray}
 \act_n^{(2j),-}(E)
 =
  \act_n^{(2j),-}(E, \hat\act)
 :=
 \frac{1}{\pi}
\int_{ \sa_{2j-1}}^{\sa_{2j}}
\sqrt{E-\Gm(\sa)}\Big(
1
+ b_\dag\big(E-\Gm(\sa),\sa \big)
\Big)\, d\sa\,,
\ 1\leq j\leq N\,,
\nonumber
\\
 \act_n^{(2j),+}(E)
 =
  \act_n^{(2j),+}(E, \hat\act)
 :=
 \frac{1}{\pi}
\int_{\sa_{2j}}^{\sa_{2j+1}}
\sqrt{E-\Gm(\sa)}\Big(
1
+ b_\dag\big(E-\Gm(\sa),\sa \big)
\Big)\, d\sa\,, \
0\leq j<N
\,.
\nonumber
\\
\label{ancomarzio2bis}
\end{eqnarray}
Finally
\begin{equation}\label{ancomarzio4}
2 \act_n^{(2N)}(E)
=
 \act_n^{(0),+}(E)+
\act_n^{(2N),-}(E)
+
\sum_{1\leq i\leq N-1}
 \act_n^{(2i),-}(E)+
 \act_n^{(2i),+}(E)\,.
\end{equation}

\subsubsection*{The adimensional action}

Recalling the definition of $\breve{\Gm}_i$
given in \eqref{bastard} we set
\begin{eqnarray}
\breve{ \act}_n^{(2j-1),-}(\breve E)
&=&
\breve{ \act}_n^{(2j-1),-}(\breve E,\hat\act)
\nonumber
\\
&:=&
\frac{1}{\pi}
\int_{\breve{ \Sa}_{2j-1}(\breve E )}^1
\sqrt{\breve E-\breve{\Gm}_{2j-1}(\breve\sa)}
\Big(
1
+ \breve b_{\dag,2j-1}\big(\breve E-\breve\Gm_{2j-1}(\breve\sa),\breve\sa \big)
\Big)
\, d\breve\sa\,,
\nonumber
\\
\breve{ \act}_n^{(2j-1),+}(\breve E)
&=&
\breve{ \act}_n^{(2j-1),+}(\breve E,\hat\act)
\nonumber
\\
&:=&
\frac{1}{\pi}
\int_{-1}^{\breve{ \Sa}_{2j}(\breve E )}
\sqrt{\breve E-\breve{\Gm}_{2j}(\breve\sa)}
\Big(
1
+ \breve b_{\dag,2j}\big(\breve E-\breve\Gm_{2j}(\breve\sa),\breve\sa \big)
\Big)
\, d\breve\sa\,,
\nonumber
\\
\breve{ \act}_n^{(2j),-}(\breve E)
&=&
\breve{ \act}_n^{(2j),-}(\breve E,\hat\act)
\nonumber
\\
&:=&
\frac{1}{\pi}
\int_{-1}^1
\sqrt{\breve E-\breve{\Gm}_{2j}(\breve\sa)}
\Big(
1
+ \breve b_{\dag,2j}\big(\breve E-\breve\Gm_{2j}(\breve\sa),\breve\sa \big)
\Big)
\, d\breve\sa\,,
\nonumber
\\
\breve{ \act}_n^{(2j),+}(\breve E)
&=&
\breve{ \act}_n^{(2j),+}(\breve E,\hat\act)
\nonumber
\\
&:=&
\frac{1}{\pi}
\int_{-1}^{1}
\sqrt{\breve E-\breve{\Gm}_{2j+1}(\breve\sa)}
\Big(
1
+ \breve b_{\dag,2j+1}\big(\breve E-\breve\Gm_{2j+1}(\breve\sa),\breve\sa \big)
\Big)
\, d\breve\sa\,,
\label{cornelia}
\end{eqnarray}
where
\begin{equation}\label{borotalco}
\breve b_{\dag,i}(\breve v,\breve \sa)=
\breve b_{\dag,i}(\breve v,\hat p,\breve \sa):=
b_\dag\Big(
(-1)^i \frac{E_i-E_{i-1}}{2}\breve v,\g_i(\breve\sa)
\Big)\,,
\end{equation}
($b_\dag$ defined in  \eqref{pachelbel}).
Note that, by \eqref{3holesbis} we have
(recalling \eqref{october})
\begin{equation}\label{japan}
\sup_{\breve v\in (0,R_0^2/\D_i)_{r_0^2/2^6\morse}}  
{\modulo}\breve b_{\dag,i}(\breve v,\hat p,\breve\sa){\modulo}_{\hat D,r_0,s_0}\leq \frac{\checco}{r_0^2}\,,\qquad
{\rm where}\qquad
\D_i:=(-1)^i (\bar E_i-\bar E_{i-1})/2\,.
\end{equation}

\noindent
Recalling \eqref{acrobat}, \eqref{bastard},
\eqref{blindness}, \eqref{ancomarzio2}
 we have\footnote{Note that
$$
E-\Gm(\sa)=
(-1)^i \frac{E_i-E_{i-1}}{2}
\Big(
\l_i(E)-\breve\Gm_i(\g_i^{-1}(\sa))
\Big)\,.
$$
}
\begin{eqnarray}
 \act_n^{(2j-1),-}(E)
&=&
\frac{\sa_{2j-1}-\sa_{2j-2}}{2}
\sqrt{\frac{ E_{2j-2}- E_{2j-1}}{2}}
\breve{ \act}_n^{(2j-1),-}
\big(\l_{2j-1}(E)\big)\,,
\nonumber
\\
 \act_n^{(2j-1),+}(E)
&=&
\frac{\sa_{2j}-\sa_{2j-1}}{2}
\sqrt{\frac{ E_{2j}- E_{2j-1}}{2}}
\breve{ \act}_n^{(2j-1),+}
\big(\l_{2j}(E)\big)\,,
\nonumber
\\
 \act_n^{(2j),-}(E)
&=&
\frac{\sa_{2j}-\sa_{2j-1}}{2}
\sqrt{\frac{ E_{2j}- E_{2j-1}}{2}}
\breve{ \act}_n^{(2j),-}
\big(\l_{2j}(E)\big)\,,
\nonumber
\\
 \act_n^{(2j),+}(E)
&=&
\frac{\sa_{2j+1}-\sa_{2j}}{2}
\sqrt{\frac{ E_{2j}- E_{2j+1}}{2}}
\breve{ \act}_n^{(2j),+}
\big(\l_{2j+1}(E)\big)\,.
\label{valeria}
\end{eqnarray}

Let set\footnote{Recall \eqref{borotalco}.}
\begin{eqnarray}\label{borotalco2}
&&\breve{\tilde b}_i(\breve v,\breve\sa)
=
\breve{\tilde b}_i(\breve v,\hat p,\breve\sa)
=
\breve b_{\dag,i}(\breve v,\breve\sa)
+2 \breve v \partial_{\breve v}
\breve b_{\dag,i}(\breve v,\breve\sa)
=\tilde b
\Big(
(-1)^i \frac{E_i-E_{i-1}}{2}\breve v,\g_i(\breve\sa)
\Big)\,,
\nonumber
\\
&&
\sup_{\breve v\in (0,R_0^2/\D_i)_{r_0^2/2^8\morse}} {\modulo}
\breve{\tilde b}_i (\breve v){\modulo}_{\hat D,r_0,s_0}
\leq \frac{9\checco}{r_0^2}
\,,
\end{eqnarray}
recalling \eqref{3holesquater} and \eqref{japan}.

\begin{lemma}\label{mardibering2}
 The functions $\breve{\act}_n^{(2j-1),\pm}(\breve E,\hat\act)$ in\footnote{$\breve{\mathtt r}$ was defined in \eqref{coppe}
 and recall \eqref{bufala}.} \eqref{cornelia} are holomorphic in 
 $ \O_{\breve{\mathtt r}}\times \hat D_{r_0}$. Moreover 
the following formulas hold:
\begin{eqnarray}
&&
\breve{ \act}_n^{(2j-1),-}(\breve E)=
-\frac{(\breve E+1)^{3/2}}{\pi}
\int_{0}^{1}
\sqrt{t}
\frac{
1
+  \breve b_{\dag,2j-1}\big(t(1+\breve E),
\breve\Sa_{2j-1}\big(
\breve E(t) \big)
 \big)
}
{\partial_{\breve \sa}
\breve\Gm_{2j-1}\Big(\breve  \Sa_{2j-1}\big(\breve E(t) \big)\Big)}
\, dt\,,
\nonumber
\\
&&
\breve{ \act}_n^{(2j-1),+}(\breve E)=
\frac{(\breve E+1)^{3/2}}{\pi}
\int_{0}^{1}
\sqrt{t}
\frac{
1
+ \breve b_{\dag,2j}\big(t(1+\breve E),
\breve\Sa_{2j}\big(
\breve E(t) \big)
 \big)
}
{\partial_{\breve \sa}
\breve\Gm_{2j}\Big(\breve  \Sa_{2j}\big(\breve E(t) \big)\Big)}
\, dt\,.
\label{cornelia2}
\end{eqnarray}
and, for the derivatives,
\begin{eqnarray}
&&
\partial_{\breve E}\breve{ \act}_n^{(2j-1),-}(\breve E):=
-\frac{\sqrt{\breve E+1}}{2\pi}
\int_{0}^{1}
\frac{1}{\sqrt{t}}
\frac{
1
+ \breve{\tilde b}_{2j-1}\big(t(1+\breve E),
\breve\Sa_{2j-1}\big(
\breve E(t) \big)
 \big)
}
{\partial_{\breve \sa}
\breve\Gm_{2j-1}\Big(\breve  \Sa_{2j-1}\big(\breve E(t) \big)\Big)}
\, dt\,,
\nonumber
\\
&&
\partial_{\breve E}\breve{ \act}_n^{(2j-1),+}(\breve E):=
\frac{\sqrt{\breve E+1}}{2\pi}
\int_{0}^{1}
\frac{1}{\sqrt{t}}
\frac{
1
+\breve{\tilde b}_{2j}\big(t(1+\breve E),
\breve\Sa_{2j}\big(
\breve E(t) \big)
 \big)
}
{\partial_{\breve \sa}
\breve\Gm_{2j}\Big(\breve  \Sa_{2j}\big(\breve E(t) \big)\Big)}
\, dt\,.
\label{Dcornelia2}
\end{eqnarray}
\end{lemma}
\proof
We proceed as in the proof of Lemma \ref{mardibering}.
Recalling \eqref{tronodidenari}
we have that
if  $\breve E\in  \O_{\breve{\mathtt r}}$
then 
$\breve E(t)\in  \O_{\breve{\mathtt r}}$
(with $\breve E(t)$ defined in \eqref{spade})
and, recalling Lemma \ref{lem:coppe}, we also get
that   the following  changes of variables
are well defined:
\begin{eqnarray}
&&\breve\sa = \breve\Sa_{2j-1}\big(
\breve E(t),\hat\act  \big)
\in [-1,1]_{\breve \rho_\star/2}
\,, \qquad \forall\, \   0\leq t\leq 1\,, \
\hat\act \in \hat D_{r_0}\,,
\nonumber
\\
&&\breve\sa = \breve \Sa_{2j}\big(
\breve E(t),\hat\act  \big)
\in [-1,1]_{\rho_\star/2}\,,
\qquad \forall\, \   0\leq t\leq 1\,, \
\hat\act \in \hat D_{r_0}\,,
\label{mammola}
\end{eqnarray}
with $\rho_\star$ defined in \eqref{coppe}.
By \eqref{rochefort6} we get
$$
t=\frac{\breve E-\breve \Gm_{2j-1}(\breve \sa)}{\breve E+1 }\,,
\qquad
t=\frac{\breve E-\breve \Gm_{2j}(\breve \sa)}{\breve E+1 }
$$
and
$$
d\breve \sa=-\frac{\breve E+1}{\partial_{\breve \sa}
\breve\Gm_{2j-1}\Big(\breve  \Sa_{2j-1}\big(\breve E(t) \big)\Big)} \, dt \qquad
\text{or} \qquad
d\breve \sa=-\frac{\breve E+1}{\partial_{\breve \sa}
\breve\Gm_{2j}\Big(\breve  \Sa_{2j}\big(\breve E(t) \big)\Big)} \, dt \,.
$$
Making the change of variables
\eqref{mammola}  in the first and second integral
in \eqref{cornelia}, respectively, we get
\eqref{cornelia2}.
The proof of \eqref{Dcornelia2} is analogous 
to the one of \eqref{Dcornelia2bar}.
\eproof


\begin{definition}\label{segment}
 For $z,w\in\C$ we define the  segment
 \begin{equation}\label{havana}
(z,w):=\{ (1-t)z+tw\,,\ \  0< t< 1\}\,.
\end{equation}
\end{definition}

\begin{lemma}\label{sargassi}
 $\act_n^{(2j-1),+}(E, \hat\act),$ respectively 
 $\act_n^{(2j-1),-}(E, \hat\act),$
  in \eqref{ancomarzio2} is holomorphic 
  on\footnote{Where, obviously,
  $(\l_{2j}^{-1},id)(\O_{\breve{\mathtt r}}\times \hat D_{r_0}):=
  \{ (E,\hat\act)=(\l_{2j}^{-1}(\breve E,\hat\act),\hat\act)\ |\ 
  (\breve E,\hat\act)\in \O_{\breve{\mathtt r}}\times \hat D_{r_0}
  \}.$
  } 
  $ (\l_{2j}^{-1},id)(\O_{\breve{\mathtt r}}\times \hat D_{r_0}),$
respectively
  $ (\l_{2j-1}^{-1},id)(\O_{\breve{\mathtt r}}\times \hat D_{r_0}).$
  In particular
  \begin{eqnarray}
\partial_E \act_n^{(2j-1),-}(E)
&=&
\frac{\sa_{2j-1}-\sa_{2j-2}}{2}
\sqrt{\frac{2}{ E_{2j-2}- E_{2j-1}}}
\partial_{\breve E}\breve{ \act}_n^{(2j-1),-}
\big(\l_{2j-1}(E)\big)\,,
\nonumber
\\
\partial_E  \act_n^{(2j-1),+}(E)
&=&
\frac{\sa_{2j}-\sa_{2j-1}}{2}
\sqrt{\frac{2}{ E_{2j}- E_{2j-1}}}
\partial_{\breve E}\breve{ \act}_n^{(2j-1),+}
\big(\l_{2j}(E)\big)\,.
\label{valerione}
\end{eqnarray}
  Moreover for every fixed $\hat\act\in \hat D_{r_0}$
  the function $E\mapsto \act_n^{(2j-1),+}(E, \hat\act)$,
  respectively $E\mapsto \act_n^{(2j-1),-}(E, \hat\act)$,
  is holomorphic on 
  \begin{equation}\label{tortuga2}
\big( E_{2j-1}(\hat\act), E_{2j}(\hat\act)\big)_{\mathtt r_1}\,,\ 
\text{respectively}\ 
\big( E_{2j-1}(\hat\act), E_{2j-2}(\hat\act)\big)_{\mathtt r_1}\,,
\end{equation}
with
$
\mathtt r_1=\b \breve{\mathtt r}/4$
define in \eqref{caraibi}.
\end{lemma}
\proof
The first part is a direct consequence of  \eqref{valeria} and of Lemma \ref{mardibering2}.  
\eqref{valerione} follows by \eqref{valeria} 
and \eqref{acrobat}.
\eqref{tortuga2}
follows by \eqref{acrobat} and \eqref{carciofino}.
\eproof

Analogously to Remark  \ref{telemann2} we have
\begin{remark}\label{telemann4}
 Note that making the inverses  of the change of 
 variables  \eqref{mammola}
  in the expression \eqref{Dcornelia2}
  we get
 \begin{eqnarray}
&&
\partial_{\breve E}\breve{ \act}_n^{(2j-1),-}(\breve E):=
\frac{1}{2\pi}
\int_{\breve{ \Sa}_{2j-1}(\breve E )}^1
\frac{
1
+ \breve{\tilde b}_{2j-1}\big(\breve E-\breve\Gm_{2j-1}(\breve\sa),\breve\sa \big)
}{\sqrt{\breve E-\breve{\Gm}_{2j-1}(\breve\sa)}}
\, d\breve\sa\,,
\nonumber
\\
&&
\partial_{\breve E}\breve{ \act}_n^{(2j-1),+}(\breve E):=
\frac{1}{2\pi}
\int_{-1}^{\breve{ \Sa}_{2j}(\breve E )}
\frac{
1
+ \breve{\tilde b}_{2j}\big(\breve E-\breve\Gm_{2j}(\breve\sa),\breve\sa \big)
}
{\sqrt{\breve E-\breve{\Gm}_{2j}(\breve\sa)}}
\, d\breve\sa\,,
\label{Dcornelia}
\end{eqnarray}
 showing also that the formal derivation in \eqref{lamento2}
 is actually correct.
 Indeed, recalling \eqref{blindness},\eqref{bastard2},\eqref{acrobat}, and making the change of variables
 $\breve\sa=\g_{2j-1}^{-1}(\sa)$ and
$\breve\sa=\g_{2j}^{-1}(\sa)$ 
 in the first and in the second integral in \eqref{Dcornelia}, respectively,
 and summing the results, we get  \eqref{lamento2}.
\end{remark}

Let us define
\begin{equation}\label{margherita}
E(t):=E-(E-E_{2j-1} )t
=t E_{2j-1} + (1-t)E\,.
\end{equation}
Recalling \eqref{acrobat}, 
\eqref{bastard}, \eqref{blindness}, 
\eqref{borotalco2}
we have that
\begin{eqnarray}
\breve E=\l_i(E) &&
\Longrightarrow
\nonumber
\\
&&
\breve E(t)=\l_i\big(E(t)\big)\,,
\nonumber
\\
&& 
(-1)^i \frac{E_i-E_{i-1}}{2}(1+\breve E)=E-E_{2j-1}\,,\qquad
{\rm for}\ \ i=2j, 2j-1
\nonumber
\\
&&
\Sa_i(E(t))=\g_i\circ \breve \Sa_i(\breve E(t))\,,
\nonumber
\\
&&
\tilde b\big(t(E-E_{2j-1}),\Sa_{i}(E(t))\big)=
\breve{\tilde b}_{i}\big(t(1+\breve E),
\breve\Sa_{i}\big(
\breve E(t) \big)
 \big)\,,
 \nonumber
 \\
&&
\partial_{\sa}
\Gm_{i}\Big( \Sa_{i}\big( E(t) \big)\Big)
=
(-1)^i\frac{E_i-E_{i-1}}{\sa_i-\sa_{i-1}}
\partial_{\breve \sa}
\breve\Gm_{i}\Big(\breve  \Sa_{i}\big(\breve E(t) \big)\Big)
\,,
\label{allemande}
\end{eqnarray}
where $\breve E(t)$ and $E(t)$
were defined in \eqref{spade} and \eqref{margherita}
respectively.
By \eqref{allemande} 
and recalling \eqref{Dcornelia2}, \eqref{valerione}
we get the following
\begin{lemma}
We have
\begin{eqnarray}\label{marconi}
&&\partial_E \act_n^{(2j-1)}(E)=
\frac{\sqrt{E-E_{2j-1} }}{2\pi}\cdot
\\
&&
\cdot\int_0^1
\frac{1}{\sqrt t}
\left(
\frac{1+\tilde b\big(t(E-E_{2j-1}),\Sa_{2j}(E(t))\big)}{\partial_{\sa} \Gm\big(
\Sa_{2j}(E(t))  \big)}
-
\frac{1+\tilde b\big(t(E-E_{2j-1}),\Sa_{2j-1}(E(t))\big)}{\partial_{\sa} \Gm\big(
\Sa_{2j-1}(E(t))  \big)}
\right)
\, dt\,.
 \nonumber
\end{eqnarray}

\end{lemma}
\proof
We also give a more direct proof.
it is convenient to split the first integral
in \eqref{lamento2} as
$$
\int_{\Sa_{2j-1}(E )}^{\Sa_{2j}(E )}
=\int_{\Sa_{2j-1}(E )}^{\sa_{2j-1} }
+\int_{\sa_{2j-1} }^{\Sa_{2j}(E )}
$$
and make the changes of variables
\begin{eqnarray}
&&\sa = \Sa_{2j-1}\big(E(t) \big)
\,, \qquad  0\leq t\leq 1\,, 
\nonumber
\\
&&\sa = \Sa_{2j}\big(E(t) \big)\,,
\qquad
0\leq t\leq 1\,,
\nonumber
\end{eqnarray}
respectively in the first and second integral.
Note that in both cases, by \eqref{premiata}, 
$$
t=\frac{E-\Gm(\sa)}{E-E_{2j-1} }\,,
$$
while
$$
d\sa=-\frac{E-E_{2j-1}}{\partial_{\sa}
\Gm\Big( \Sa_{2j-1}\big(E(t) \big)\Big)} \, dt\,, \qquad
d\sa=-\frac{E-E_{2j-1}}{\partial_{\sa}
\Gm\Big( \Sa_{2j}\big(E(t) \big)\Big)} \, dt\,,
$$
respectively
(recalling \eqref{9till5}).
Recalling \eqref{cook} we get \eqref{marconi}. 
\eproof

\medskip

\begin{lemma}\label{mardibering3}
 The functions $\breve\act_n^{(2j),\pm}(\breve E)$ in \eqref{cornelia} are holomorphic 
 in\footnote{The intervals in \eqref{vltava} are non empty
 by \eqref{legna}.} 
\begin{eqnarray}\label{vltava}
&&
\nonumber
\O^{(2j),-}_{r_0^2/2^8\morse}
:=\Big(\C_*+1\Big) \cap \left(1,\frac{R_0^2}{\D_{2j}}-1\right)_{r_0^2/2^8\morse}\,,
\\
&&\O^{(2j),+}_{r_0^2/2^8\morse}
:=\Big(\C_*+1\Big) \cap \left(1,\frac{R_0^2}{\D_{2j+1}}-1\right)_{r_0^2/2^8\morse}
\end{eqnarray}
 (recall \eqref{japan} and \eqref{borotalco2}).
  Moreover
 the following formulas hold:
\begin{eqnarray}
&&
\breve{ \act}_n^{(2j),-}(\breve E)=
\frac{1}{\pi}
\int_{-1}^{1}
\frac{
\sqrt{\breve E-y}
\Big(
1
+ \breve b_{\dag,2j}\big(\breve E-y,
\breve  {\Sa}_{2j}(y) \big)
\Big)
}
{\partial_{\breve \sa}
\breve{\Gm}_{2j}\Big(\breve  {\Sa}_{2j}(y)\Big)}
\, dy\,,
\nonumber
\\
&&
\breve{\act}_n^{(2j),+}(\breve E)
=-
\frac{1}{\pi}
\int_{-1}^{1}
\frac{
\sqrt{\breve E-y}
\Big(
1
+ \breve b_{\dag,2j+1}\big(\breve E-y,
\breve  {\Sa}_{2j+1}(y) \big)
\Big)
}
{\partial_{\breve \sa}
\breve{\Gm}_{2j+1}\Big(\breve  {\Sa}_{2j+1}(y)\Big)}
\, dy
\label{cornelia2+}
\end{eqnarray}
(where $\breve b_{\dag,i}$ were defined in \eqref{borotalco})
and, for the derivatives,
\begin{eqnarray}
&&
\partial_{\breve E}
\breve{ \act}_n^{(2j),-}(\breve E)=
\frac{1}{2\pi}
\int_{-1}^{1}
\frac{
\Big(
1
+ \breve{\tilde b}_{2j}\big(\breve E-y,
\breve  {\Sa}_{2j}(y) \big)
\Big)
}
{
\sqrt{\breve E-y}\, 
\partial_{\breve \sa}
\breve{\Gm}_{2j}\Big(\breve  {\Sa}_{2j}(y)\Big)}
\, dy\,,
\nonumber
\\
&&
\partial_{\breve E}
\breve{\act}_n^{(2j),+}(\breve E)
=-
\frac{1}{2\pi}
\int_{-1}^{1}
\frac{
\Big(
1
+ \breve{\tilde b}_{2j+1}\big(\breve E-y,
\breve  {\Sa}_{2j+1}(y) \big)
\Big)
}
{
\sqrt{\breve E-y}\, 
\partial_{\breve \sa}
\breve{\Gm}_{2j+1}\Big(\breve  {\Sa}_{2j+1}(y)\Big)}
\, dy
\,,
\label{Dcornelia2+}
\end{eqnarray}
(where $\breve{\tilde b}_{i}$ were defined in \eqref{borotalco2}).
\end{lemma}
\proof
Recalling Lemma \ref{lem:coppe},
for  $y\in  \O_{\breve{\mathtt r}}$
   the following  changes of variables
are well defined:
\begin{eqnarray}
&&\breve\sa = \breve{\Sa}_{2j}\big(
 y \big)\in [-1,1]_{\rho_\star/2}
\,,  
\nonumber
\\
&&\breve\sa = \breve{\Sa}_{2j+1}\big(
 y \big)\in [-1,1]_{\rho_\star/2}\,,
\label{mammolabar+2}
\end{eqnarray}
with $\rho_\star$ defined in \eqref{coppe}.
Note that  $ \breve{\Sa}_{2j}(\pm1)=
 \breve{\Sa}_{2j+1}(\mp 1)=\pm 1$
(recall Lemma \ref{cerveza2}).
Recalling \eqref{rochefort6} we get
$$
d\breve \sa=\frac{1}{\partial_{\breve \sa}
\breve{\Gm}_{2j}
\Big(\breve  {\Sa}_{2j}(y )\Big)} \, dy \qquad
\text{or} \qquad
d\breve \sa=\frac{1}{\partial_{\breve \sa}
\breve{\Gm}_{2j+1}\Big(\breve{\Sa}_{2j+1}(y )\Big)} \, dt \,,
$$
(deriving \eqref{rochefort6}).
Making the changes of variables 
\eqref{mammolabar+2}
 in the third and fourth integral
in \eqref{cornelia} respectively, we get\footnote{One can easily control that the integrals
in  absolutely converge.} .
Note that if $\breve E\in \C_*+1$ and 
$-1\leq y\leq 1$ then $\breve E-y \in \C_*,
$ so that $\sqrt{\breve E-y}$ is well defined.
Moreover 
\begin{equation}\label{starlight}
\breve E\in 
(1,R_0^2/\D_i-1)_{r_0^2/2^8\morse}
\quad\Longrightarrow\quad
\breve E-y\in 
(0,R_0^2/\D_i)_{r_0^2/2^8\morse}
\,,
\qquad \forall\, -1\leq y\leq 1\,.
\end{equation}
Deriving   (and recalling
\eqref{borotalco2}) we get
\eqref{Dcornelia2+}.
\eproof

\subsection{Closeness of  the rescaled unperturbed and perturbed actions}

By 
 \eqref{salmone}
and \eqref{salmone2}
we get, for every 
$-1<y<1$ and
$\hat\act\in \hat D_{r_0},$
\begin{eqnarray}
&&\left|
\frac{
1
}
{\partial_{\breve \sa}
\breve\Gm_{i}\Big(\breve  \Sa_{i}(y)\Big)}
-
\frac{
1
}
{\partial_{\breve \sa}
\breve{\FO}_{i}\Big(\breve  
{\bar\Sa}_{i}(y)\Big)}
\right|
\ = \ 
\left|
\frac{
\partial_{\breve \sa}
\breve{\FO}_{i}\Big(\breve  
{\bar\Sa}_{i}(y)\Big)
-
\partial_{\breve \sa}
\breve\Gm_{i}\Big(\breve  \Sa_{i}(y)\Big)
}
{\partial_{\breve \sa}
\breve\Gm_{i}\Big(\breve  \Sa_{i}(y)\Big)
\ 
\partial_{\breve \sa}
\breve{\FO}_{i}\Big(\breve  
{\bar\Sa}_{i}(y)\Big)}
\right|
\nonumber
\\
&&\leq
\frac{2^{73}\breve \morse^8}
{\breve \b^9 \breve s^{18}}
\breve\checco 
\frac{
1
}
{\Big|\partial_{\breve \sa}
\breve{\FO}_{i}\Big(\breve  
{\bar\Sa}_{i}(y)\Big)\Big|}
\ \leq\ 
\frac{2^{193} \morse^{27}}
{ \b^{28}  s_0^{45}s_*^2} 
\frac{
\checco
}
{\Big|\partial_{\breve \sa}
\breve{\FO}_{i}\Big(\breve  
{\bar\Sa}_{i}(y)\Big)\Big|}
\label{scogli3}
\\
&&\leq
\frac{
1
}
{\Big|\partial_{\breve \sa}
\breve{\FO}_{i}\Big(\breve  
{\bar\Sa}_{i}(y)\Big)\Big|}\,,
\label{scogli4}
\end{eqnarray}
where the second inequality holds by
\eqref{exit}-\eqref{ladispoli4} and the last inequality
follows from \eqref{genesis}.
Moreover for every $\breve E \in \O^{(2j),\pm}_{r_0^2/2^8\morse},$
$-1<y<1$ and
$\hat\act\in \hat D_{r_0},$
we have (recalling \eqref{starlight})
\begin{eqnarray}
&&\left|
\frac{
1
+ \breve{\tilde b}_{i}\big(\breve E-y,
\breve\Sa_{i}(
y)
 \big)
}
{\partial_{\breve \sa}
\breve\Gm_{i}\Big(\breve  \Sa_{i}(y)\Big)}
-
\frac{
1
}
{\partial_{\breve \sa}
\breve{\FO}_{i}\Big(\breve  
{\bar\Sa}_{i}(y)\Big)}
\right|
\nonumber
\\
&&\leq
\left|
\frac{
 \breve{\tilde b}_{i}\big(\breve E-y,
\breve\Sa_{i}(
y)
 \big)
}
{\partial_{\breve \sa}
\breve\Gm_{i}\Big(\breve  \Sa_{i}(y)\Big)}
\right|
+
\left|
\frac{
1
}
{\partial_{\breve \sa}
\breve\Gm_{i}\Big(\breve  \Sa_{i}(y)\Big)}
-
\frac{
1
}
{\partial_{\breve \sa}
\breve{\FO}_{i}\Big(\breve  
{\bar\Sa}_{i}(y)\Big)}
\right|
\nonumber
\\
&&\leq
\left(
\frac{18}{r_0^2}
+
\frac{2^{193} \morse^{27}}
{ \b^{28}  s_0^{45}s_*^2} 
\right)
\frac{
\checco
}
{\Big|\partial_{\breve \sa}
\breve{\FO}_{i}\Big(\breve  
{\bar\Sa}_{i}(y)\Big)\Big|}
\label{scogli}
\\
&&\leq
\frac{
1
}
{\Big|\partial_{\breve \sa}
\breve{\FO}_{i}\Big(\breve  
{\bar\Sa}_{i}(y)\Big)\Big|}
\label{scogli2}
\end{eqnarray}
where the second  inequality
follows by \eqref{3holesquater}, \eqref{borotalco2},
\eqref{cogli3}
 and the last inequality
follows from \eqref{genesis}.

By \eqref{fiodena} and \eqref{tronodidenari}
we have that
\begin{eqnarray}\label{correggia}
&&\breve E \in \O_{\breve{\mathtt r}} 
 \qquad\Longrightarrow\qquad
\\
&&\breve{\bar\Sa}_i(\breve E)\,,\ 
\breve{\Sa}_i(\breve E,\hat\act)\,,\ 
\breve{\bar\Sa}_i(\breve E(t))\,,\ 
\breve{\Sa}_i(\breve E(t),\hat\act)\ \in \ 
[-1,1]_{\rho_\star/2}\,,\qquad
\forall\, t\in[0,1]\,, 
\forall\, \hat\act\in \hat D_{r_0} \,,
\nonumber
\end{eqnarray}
where $\breve{\mathtt r}$ and $\rho_\star$ were defined in \eqref{coppe}.
Then, for every 
$\breve E \in \O_{\breve{\mathtt r}},$
$\hat\act\in \hat D_{r_0},$
$0\leq t\leq 1,$
\begin{eqnarray}
&&\left|
\frac{
1
}
{\partial_{\breve \sa}
\breve\Gm_{i}\Big(\breve  \Sa_{i}\big(\breve E(t) \big)\Big)}
-
\frac{
1
}
{\partial_{\breve \sa}
\breve{\FO}_{i}\Big(\breve  
{\bar\Sa}_{i}\big(\breve E(t) \big)\Big)}
\right|
\ \leq\ 
\frac{2^{193} \morse^{27}}
{ \b^{28}  s_0^{45}s_*^2} 
\frac{
\checco
}
{\Big|\partial_{\breve \sa}
\breve{\FO}_{i}\Big(\breve  
{\bar\Sa}_{i}\big(\breve E(t) \big)\Big)\Big|}
\label{cogli3}
\\
&&\leq
\frac{
1
}
{\Big|\partial_{\breve \sa}
\breve{\FO}_{i}\Big(\breve  
{\bar\Sa}_{i}\big(\breve E(t) \big)\Big)\Big|}\,,
\label{cogli4}
\end{eqnarray}
where the second inequality holds by
\eqref{exit}-\eqref{ladispoli4} and the last inequality
follows from \eqref{genesis}
and 
\begin{eqnarray}
&&\left|
\frac{
1
+ \breve{\tilde b}_{i}\big(t(1+\breve E),
\breve\Sa_{i}\big(
\breve E(t) \big)
 \big)
}
{\partial_{\breve \sa}
\breve\Gm_{i}\Big(\breve  \Sa_{i}\big(\breve E(t) \big)\Big)}
-
\frac{
1
}
{\partial_{\breve \sa}
\breve{\FO}_{i}\Big(\breve  
{\bar\Sa}_{i}\big(\breve E(t) \big)\Big)}
\right|
\nonumber
\\
&&\leq
\left(
\frac{18}{r_0^2}
+
\frac{2^{193} \morse^{27}}
{ \b^{28}  s_0^{45}s_*^2} 
\right)
\frac{
\checco
}
{\Big|\partial_{\breve \sa}
\breve{\FO}_{i}\Big(\breve  
{\bar\Sa}_{i}\big(\breve E(t) \big)\Big)\Big|}
\label{cogli}
\\
&&\leq
\frac{
1
}
{\Big|\partial_{\breve \sa}
\breve{\FO}_{i}\Big(\breve  
{\bar\Sa}_{i}\big(\breve E(t) \big)\Big)\Big|}
\label{cogli2}
\end{eqnarray} 
 
 $\bullet${\sl The odd case}

\begin{lemma}\label{giga3}
We have that 
$\breve{\bar\act}_n^{(2j-1),\pm} (\breve E)$
and $\breve\act_n^{(2j-1),\pm}(\breve E,\hat\act) $
are holomorphic functions in the sets\footnote{Recall \eqref{bufala}
and \eqref{coppe}.}
$
\O_{\breve{\mathtt r}}
$
and
$\O_{\breve{\mathtt r}}
\times
\hat D_{r_0}
$
respectively.
Moreover, for
$\hat\act\in \hat D_{r_0}$ and
  $\breve E\in \O_{\breve{\mathtt r}}$ 
  with $\Re \breve E<1$,
  the following estimate holds
\begin{equation}
\label{collemaggio}
|\partial_{\breve E}\breve\act_n^{(2j-1),\pm}(\breve E,\hat\act)
-
\partial_{\breve E}\breve{\bar\act}_n^{(2j-1),\pm}(\breve E)
|
\leq
\checco
\left(
\frac{36}{r_0^2}
+
\frac{2^{194} \morse^{27}}
{ \b^{28}  s_0^{45}s_*^2} 
\right)
\frac{ 2^{55}\morse^7}{\b^7  s_0^{13}}
\left( 
1
+
\ln \frac{1}{1-\Re \breve E}
\right)
\,.
\end{equation}
\end{lemma}
\proof
The   ``-'' term in \eqref{collemaggio} is bounded by
\begin{equation}\label{morrone}
\frac{\sqrt{|\breve E+1
|}}{2\pi}
\int_{0}^{1}
\frac{1}{\sqrt{t}}
\left|
\frac{
1
+ \breve{\tilde b}_{2j-1}\big(t(1+\breve E),
\breve\Sa_{2j-1}\big(
\breve E(t) \big)
 \big)
}
{\partial_{\breve \sa}
\breve\Gm_{2j-1}\Big(\breve  \Sa_{2j-1}\big(\breve E(t) \big)\Big)}
-
\frac{
1
}
{\partial_{\breve \sa}
\breve{\FO}_{2j-1}\Big(\breve  
{\bar\Sa}_{2j-1}\big(\breve E(t) \big)\Big)}
\right|
\, dt\,.
 \end{equation}
By \eqref{cogli} 
it
is bounded by
$$
\checco
\left(
\frac{18}{r_0^2}
+
\frac{2^{193} \morse^{27}}
{ \b^{28}  s_0^{45}s_*^2} 
\right)
\frac{\sqrt{|\breve E+1
|}}{2\pi}
\int_{0}^{1}
\frac{1}{\sqrt{t}}
\frac{
1
}
{\Big|\partial_{\breve \sa}
\breve{\FO}_{2j-1}\Big(\breve  
{\bar\Sa}_{2j-1}\big(\breve E(t) \big)\Big)\Big|}
\, dt\,.
$$
The ``+'' term is estimated analogously\footnote{With $2j$ instead of $2j-1$.}.
Then
we get
\begin{eqnarray}\label{tartufo}
&&
|\partial_{\breve E}\breve\act_n^{(2j-1),\pm}(\breve E,\hat\act)
-
\partial_{\breve E}\breve{\bar\act}_n^{(2j-1),\pm}(\breve E)|
\nonumber
\\
&&
\leq
\checco
\left(
\frac{36}{r_0^2}
+
\frac{2^{194} \morse^{27}}
{ \b^{28}  s_0^{45}s_*^2} 
\right)
\frac{\sqrt{|\breve E+1
|}}{2\pi}
\int_{0}^{1}
\frac{1}{\sqrt{t}}
\frac{
1
}
{\Big|\partial_{\breve \sa}
\breve{\FO}_{i}\Big(\breve  
{\bar\Sa}_{i}\big(\breve E(t) \big)\Big)\Big|}
\, dt\,,
\end{eqnarray}
where $i=2j$ in the $+$ case and 
$i=2j-1$ in the $-$ case.
\\
We estimate the integral in \eqref{tartufo}
in the following

\begin{lemma}
 For $\breve E \in \O_{\breve{\mathtt r}} $
 with $\Re \breve E<1$ we have
\begin{equation}\label{dave}
 \int_{0}^{1}
\frac{
1
}
{\sqrt{t}\Big|\partial_{\breve \sa}
\breve{\FO}_{i}\Big(\breve  
{\bar\Sa}_{i}\big(\breve E(t) \big)\Big)\Big|}
\, dt
\leq 
\frac{ 2^{31}\breve\morse}{\breve\b^3 \breve s^4}
\left( 
\frac{1}{\sqrt{|1+\breve E|}}
+
\ln \frac{1}{1-\Re \breve E}
\right)
\,.
 \end{equation}
\end{lemma}
\proof
  Recalling \eqref{correggia}, \eqref{pajata} and 
 \eqref{harlock1bar}
 we get
$$
\frac{
1
}
{\Big|\partial_{\breve \sa}
\breve{\FO}_{i}\Big(\breve  
{\bar\Sa}_{i}\big(\breve E(t) \big)\Big)\Big|}
\leq
\frac{8}{\breve\b}\left(
\frac{ 2^{10}\breve\morse}{\breve\b \breve s^3} 
+
\frac{1}{\big|1-\breve  
{\bar\Sa}_{i}\big(\breve E(t) \big)\big|}
+
\frac{1}{\big|1+\breve  
{\bar\Sa}_{i}\big(\breve E(t) \big)\big|}
\right)\,,
\,.
$$ 
Recalling
\eqref{ventresca}, \eqref{coppe}, \eqref{nera}
and noting that 
$\breve{\mathtt r}\leq \breve r_0,$
we get
\begin{eqnarray*}
\frac{
1
}
{\Big|\partial_{\breve \sa}
\breve{\FO}_{i}\Big(\breve  
{\bar\Sa}_{i}\big(\breve E(t) \big)\Big)\Big|}
\leq
\frac{ 2^{26}\breve\morse}{\breve\b^3 \breve s^4}
\left( 
1
+
\frac{1}{\sqrt{|\breve E(t)+ 1|}}
+
\frac{1}{\sqrt{|\breve E(t)- 1|}|}
\right)
\,.
\end{eqnarray*}
We claim that the integral in \eqref{dave}
is bounded by	
$$
\frac{ 2^{26}\breve\morse}{\breve\b^3 \breve s^4}
\left( 
2
+
\frac{\pi}{\sqrt{|1+\breve E|}}
+
4+2\ln \frac{1}{1-\Re \breve E}
\right)
$$
and \eqref{dave} follows\footnote{Note that
for $\breve E \in \O_{\breve{\mathtt r}} $
we have $1/\sqrt{|1+\breve E|}\geq 1/2.$ }.
We now  prove the above claim. Note that
$
1+\breve E(t)=(1-t)(1+\breve E)$
and
$|\breve E(t)- 1|
\geq 1-\Re\breve E(t)
$ 
$=(1-\Re\breve E)+(1+\Re\breve E)t,$ for every $t\in[0,1].$
Then $\int_0^1 dt/\sqrt{t|\breve E(t)+ 1|}=\pi/\sqrt{|\breve E+ 1|}.$
Moreover in estimating $\mathfrak I:=\int_0^1 dt/\sqrt{t|\breve E(t)- 1|}$
we have two cases: if $\Re \breve E\leq 0$ then $|\breve E(t)- 1|\geq 1,$
$\forall t\in[0,1]$ and $\mathfrak I\leq 2;$ otherwise, when 
$\Re \breve E\geq 0,$
we have $|\breve E(t)- 1|\geq 1-\Re \breve E+t,$ $\forall t\in[0,1],$ then,
setting $\xi:=1-\Re \breve E$ (note that $0<\xi\leq 1$), we get
$$
\mathfrak I\leq \int_0^1\frac{dt}{\sqrt t \sqrt{1-\Re \breve E+t}}
=2\int_0^{1/\sqrt{\xi}} \frac{ds}{\sqrt{1+s^2}}
\leq 
4\int_0^{1/\sqrt{\xi}} \frac{ds}{1+s}
=4\ln (1+\frac{1}{\sqrt\xi})
\leq 4+2\ln\frac{1}{\xi}\,,
$$
proving the claim.
  \eproof

Inserting the estimate \eqref{dave} in \eqref{tartufo} we get
$$
|\partial_{\breve E}\breve\act_n^{(2j-1),\pm}(\breve E,\hat\act)
-
\partial_{\breve E}\breve{\bar\act}_n^{(2j-1),\pm}(\breve E)
|
\leq
\checco
\left(
\frac{36}{r_0^2}
+
\frac{2^{194} \morse^{27}}
{ \b^{28}  s_0^{45}s_*^2} 
\right)
\frac{ 2^{31}\breve\morse}{\breve\b^3 \breve s^4}
\left( 
1
+
\ln \frac{1}{1-\Re \breve E}
\right)
\,.
$$
By \eqref{exit}-\eqref{ladispoli4} we get \eqref{collemaggio}. \eproof

$\bullet${\sl The even case}

\begin{lemma}
 For   
$\Re \breve E>1$ we have
\begin{equation}\label{dave+}
\int_{-1}^{1}
\frac{
1
}
{\left|
\sqrt{\breve E-y}\,
\partial_{\breve \sa}
\breve{\FO}_{i}\big(\breve  {\bar\Sa}_{i}(y)\big)\right|}
\, dy
\leq 
\frac{ 2^{31}\breve\morse}{\breve\b^3 \breve s^4}
\frac{1}{\sqrt{\Re \breve E}}
\ln\left(
4+\frac{1}{\Re\breve E-1}
\right)
\,.
 \end{equation}
\end{lemma}
\proof
  Recalling \eqref{correggia}, \eqref{pajata} and 
 \eqref{harlock1bar}
 we get, for every $-1<y<1,$
$$
\frac{
1
}
{\Big|\partial_{\breve \sa}
\breve{\FO}_{i}\Big(\breve  
{\bar\Sa}_{i}(y)\Big)\Big|}
\leq
\frac{8}{\breve\b}\left(
\frac{ 2^{10}\breve\morse}{\breve\b \breve s^3} 
+
\frac{1}{\big|1-\breve  
{\bar\Sa}_{i}(y)\big|}
+
\frac{1}{\big|1+\breve  
{\bar\Sa}_{i}(y)\big|}
\right)\,,
\,.
$$ 
Recalling
\eqref{ventresca}, \eqref{coppe}, \eqref{nera}
we get
$$
\frac{
1
}
{\Big|\partial_{\breve \sa}
\breve{\FO}_{i}\Big(\breve  
{\bar\Sa}_{i}(y)\Big)\Big|}
\leq
\frac{ 2^{26}\breve\morse}{\breve\b^3 \breve s^4}
\left( 
1
+
\frac{1}{\sqrt{|y+ 1|}}
+
\frac{1}{\sqrt{|y- 1|}}
\right)
\leq
\frac{ 2^{26}\breve\morse}{\breve\b^3 \breve s^4}
\left( 
2
+
\frac{1}{\sqrt{1-y}}
\right)
\,.
$$
Since $1/|\sqrt{\breve E-y}|\leq 1/\sqrt{\Re \breve E-y}$
we have that the integral in \eqref{dave+}
is bounded by	
$$
\frac{ 2^{26}\breve\morse}{\breve\b^3 \breve s^4}
\int_{-1}^1
\left( 
2
+
\frac{1}{\sqrt{1-y}}
\right)
\frac{1}{\sqrt{\Re \breve E-y}}\, dy
\leq
\frac{ 2^{31}\breve\morse}{\breve\b^3 \breve s^4}
\frac{1}{\sqrt{\Re \breve E}}
\ln\left(
4+\frac{1}{\Re\breve E-1}
\right)
$$
and \eqref{dave+} follows.
  \eproof
\begin{lemma}\label{giga4}
For every
$  
\hat\act\in D_{r_0}$ and 
  $\breve E\in \O^{(2j),\pm}_{r_0^2/2^8\morse}$ with
$\Re \breve E>1$, we have
\begin{eqnarray}
&&|\partial_{\breve E}\breve\act_n^{(2j),\pm}(\breve E,\hat\act)
-
\partial_{\breve E}\breve{\bar\act}_n^{(2j),\pm}(\breve E)
|
\nonumber
\\
&&
\leq
\checco
\left(
\frac{18}{r_0^2}
+
\frac{2^{193} \morse^{27}}
{ \b^{28}  s_0^{45}s_*^2} 
\right)
\frac{ 2^{55}\morse^7}{\b^7  s_0^{13}}
\frac{1}{\sqrt{\Re \breve E}}
\ln\left(
4+\frac{1}{\Re\breve E-1}
\right)
\,.
\label{collemaggio+}
\end{eqnarray}
\end{lemma}
\proof
Let us consider the ``-'' case, the ``+'' one is analogous.
Then the quantity on the l.h.s. of \eqref{collemaggio+}
is bounded by 
\begin{eqnarray*}
&&\frac{1}{2\pi}
\int_{-1}^{1}
\frac{1}{|\sqrt{\breve E-y}|}
\left|
\frac{
1
+ \breve{\tilde b}_{2j}\big(\breve E-y,
\breve  {\Sa}_{2j}(y) \big)
}
{ 
\partial_{\breve \sa}
\breve{\Gm}_{2j}\Big(\breve  {\Sa}_{2j}(y)\Big)}
-
\frac{1}{\partial_{\breve \sa}
\breve{\FO}_{2j}\Big(\breve  {\bar\Sa}_{2j}(y)\Big)}
\right|
\, dy
\\
&&
\stackrel{\eqref{scogli}}\leq
\frac{1}{2\pi}
\int_{-1}^{1}
\frac{1}{|\sqrt{\breve E-y}|}
\left(
\frac{18}{r_0^2}
+
\frac{2^{193} \morse^{27}}
{ \b^{28}  s_0^{45}s_*^2} 
\right)
\frac{
\checco
}
{\Big|\partial_{\breve \sa}
\breve{\FO}_{i}\Big(\breve  
{\bar\Sa}_{i}(\breve E)\Big)\Big|}
\, dy
\\
&&\stackrel{\eqref{dave+}}\leq
\checco
\left(
\frac{18}{r_0^2}
+
\frac{2^{193} \morse^{27}}
{ \b^{28}  s_0^{45}s_*^2} 
\right)
\frac{ 2^{31}\breve\morse}{\breve\b^3 \breve s^4}
\frac{1}{\sqrt{\Re \breve E}}
\ln\left(
4+\frac{1}{\Re\breve E-1}
\right)\,.
\end{eqnarray*}
By \eqref{exit}-\eqref{ladispoli4} we get \eqref{collemaggio+}.
\eproof

\subsection{Estimates on  $\partial_E \act_n^{(2j-1)}$}

\begin{remark}\label{stiffe}
 We will assume, to fix ideas, that 
 $\bar E_{2j-2}<\bar E_{2j}$
 so that $j_\diamond =j-1$ in \eqref{padula},\eqref{gruffalo}
 and $\bar E_{2j_\diamond}=\bar E_{2j-2}=\bar E^{(2j-1)}_+;$
 analogously 
 $ E_{2j_\diamond}= E_{2j-2}= E^{(2j-1)}_+.$
\end{remark}

\begin{lemma}\label{giga}
 $\partial_E \act_n^{(2j-1)}$ is a holomorphic
 function of the complex variable
 $\z=E-E_{2j-1}$ for
\begin{equation}\label{epistassi}
  |\z|<\mathtt r_\star
 := \min\left\{
\frac{s_0^{18}}{2^{77}}
\frac{\b^{12}}{\morse^{11}},
\frac{r_0^2}{2^8}
\right\}\,.
\end{equation}
 In particular
 \begin{equation}\label{sarabanda}
\partial_E \act_n^{(2j-1)}(E,\hat\act)=
\frac{1}{2\pi}
\int_0^1
\frac{G\big(E-E_{2j-1},
(1-t)(E-E_{2j-1}),\hat\act\big)}
{\sqrt t\sqrt{1-t}}
\, dt\,,
\end{equation}
 for a 
suitable holomorphic function 
$G(v,y,\hat\act)$ satisfying
\begin{equation}\label{tritone}
\sup_{v\in (0,R_0^2)_{r_0^2/64}}
\sup_{|y|<\mathtt r_\star}|G|_{\hat D,r_0}
\leq \frac{32\sqrt\morse}{s_0\b}\,.
\end{equation}
Then
\begin{equation}\label{gattoboy}
\sup_{|E-E_{2j-1}|<\mathtt r_\star}|\partial_E \act_n^{(2j-1)}(E,\hat\act)|_{\hat D,r_0}
\leq \frac{16\sqrt\morse}{s_0\b}\,.
\end{equation}
 \end{lemma}
\proof
Set 
\begin{equation}\label{capricciosa}
w=w(t):=\sqrt{E-E_{2j-1}} \sqrt{1-t}
\stackrel{\eqref{margherita}}=\sqrt{E(t)-E_{2j-1}}
\,,
\end{equation}
then
by \eqref{gargano}
\begin{eqnarray}\label{tonino}
\Sa_{2j-1}(E(t))=
\sa_{2j-1}-w(t) \Sa_{2j-1,-}(w(t))\,,
\quad
\Sa_{2j}(E(t))=
\sa_{2j-1}+w(t) \Sa_{2j,-}(w(t))\,,\!\!\!\!\!\!\!\!\!&&
\nonumber
\\
\text{for}\ |w(t)|<r_\diamond\,,&&
\end{eqnarray}
 ($r_\diamond$ 
defined in \eqref{burrata}),
which is implied by 
\begin{equation}\label{romeo}
|E-E_{2j-1}|<\mathtt r_\star
<\min\{r_\diamond^2,\frac{r_0^2}{2^8}\}
=
\min\left\{
\frac{s_0^{18}}{2^{50}3^4\pi^{10}}
\frac{\b^{12}}{\morse^{11}},
\frac{r_0^2}{2^8}
\right\}
\,,
\qquad E-E_{2j-1}\in \C_*\,.
\end{equation}
Since $\partial_{\sa} \Gm(\sa_{2j-1})=0$ we have that the function
$$
\hat \Gm(\xx):=\partial_{\sa} \Gm(\sa_{2j-1}+\xx)
/\xx
=\int_0^1 \partial_{\sa\sa} \Gm(\sa_{2j-1}+\xx y)dy\,. 
$$
 is 
holomorphic.
By Taylor expansion, Cauchy estimates, \eqref{ciccio3}
and \eqref{ladispoli3}
\begin{equation}\label{maria}
|\hat \Gm(\xx)|=|\frac{\partial_{\sa} \Gm(\sa_{2j-1}+\xx)}{\xx}|
\geq 
|\partial_{\sa \sa} \Gm(\sa_{2j-1})|
-\frac{8\morse}{s_0^3}|\xx|
\geq \frac{\b}{2}-\frac{\b}{4}=
\frac{\b}{4}
\end{equation}
for 
$$
|\xx|
\leq
\frac{\b s_0^3}{32 \morse}=:\xx_\diamond\,.
$$
By \eqref{tonino} we have that
\begin{eqnarray}\label{ocio}
\partial_{\sa} \Gm\big(
\Sa_{2j}(E(t))  \big)
&=&
w(t) \Sa_{2j,-}(w(t))\hat\Gm\big(
w(t) \Sa_{2j,-}(w(t))
\big)
\nonumber
\\
\partial_{\sa} \Gm\big(
\Sa_{2j-1}(E(t))  \big)
&=&
-w(t) \Sa_{2j-1,-}(w(t))\hat\Gm\big(
-w(t) \Sa_{2j-1,-}(w(t))
\big)
\end{eqnarray}
for $|w(t)|<r_\diamond$.

Let us define the holomorphic function
\begin{eqnarray}\label{gufetta}
g(v,w,\hat \act)&:=&
\frac{1+\tilde b\big(v-w^2,
\sa_{2j-1}+w \Sa_{2j,-}(w)\big)}
{ \Sa_{2j,-}(w) \hat \Gm\big(
w \Sa_{2j,-}(w)\big)}
\nonumber
\\
&&\ \ \ 
+
\frac{1+\tilde b\big(v-w^2,
\sa_{2j-1}-w \Sa_{2j-1,-}(w)\big)}
{ \Sa_{2j-1,-}(w) \hat \Gm\big(
-w \Sa_{2j-1,-}(w)\big)}\,.
\end{eqnarray}
By \eqref{lapulcedacqua}
it follows that $g$ is an {\sl even}
function w.r.t. $w$. Then 
\begin{equation}\label{grattaebasta}
g(v,w,\hat\act)=:G(v,w^2,\hat\act)
\end{equation}
 for a suitable function
$G$.
Recalling \eqref{tonino} 
and noting that, for $|w|<r_\diamond,$
\begin{eqnarray}
|w \Sa_{2j,\pm}(w)|
&\stackrel{\eqref{oxiana3}}\leq&
 r_\diamond
\frac{4\sqrt 3 \pi \morse}{s_0^{3/2} \b^{3/2}}
\stackrel{\eqref{burrata}}=
\frac{s_0^{9}}{2^{25}9\pi^{5}}
\frac{\b^6}{\morse^{11/2}}
\frac{4\sqrt 3 \pi \morse}{s_0^{3/2} \b^{3/2}}
=
\frac{s_0^{15/2}}{2^{23}3\sqrt 3 \pi^{4}}
\frac{\b^{9/2}}{\morse^{9/2}}
\nonumber
\\
&\stackrel{\eqref{harlock}}\leq&
\min\left\{ \sa_\diamond,
 \frac{s_0}{8}\right\}\,,
 \label{grattaevinci}
\end{eqnarray}
we have that
\begin{equation}\label{dolores}
\sup_{v\in (0,R_0^2)_{r_0^2/64}}
\sup_{|w|<\sqrt{\mathtt r_\star}}|g(v,w,\hat\act)|_{\hat D,r_0}
=
\sup_{v\in (0,R_0^2)_{r_0^2/64}}
\sup_{|y|<\mathtt r_\star}|G(v,y,\hat\act)|_{\hat D,r_0}
\leq \frac{32\sqrt\morse}{s_0\b}
\end{equation}
by \eqref{3holesquater}, \eqref{grattaevinci}, 
\eqref{maria} and \eqref{spadadefoco}.

By \eqref{ocio} and \eqref{grattaebasta} we get
\begin{eqnarray}\label{gufetta2}
&&
\frac{1+\tilde b\big(t(E-E_{2j-1}),\Sa_{2j}(E(t))\big)}{\partial_{\sa} \Gm\big(
\Sa_{2j}(E(t))  \big)}
-
\frac{1+\tilde b\big(t(E-E_{2j-1}),\Sa_{2j-1}(E(t))\big)}{\partial_{\sa} \Gm\big(
\Sa_{2j-1}(E(t))  \big)}
\nonumber
\\
&&=
\frac{g\big(E-E_{2j-1},w(t),\hat \act	\big)}{w(t)}
=
\frac{G\big(E-E_{2j-1},w^2(t),\hat \act	\big)}{w(t)}\,.
\end{eqnarray}
By \eqref{marconi}
\begin{eqnarray}\label{edison}
\partial_E \act_n^{(2j-1)}(E)=
\frac{1}{2\pi}
\int_0^1
\frac{g\big(E-E_{2j-1},w(t),\hat \act	\big)}{\sqrt t\sqrt{1-t}}
\, dt
=
\frac{1}{2\pi}
\int_0^1
\frac{G\big(E-E_{2j-1},w^2(t),\hat \act	\big)}{\sqrt t\sqrt{1-t}}
\, dt
\nonumber
\\
\end{eqnarray}
and\eqref{sarabanda} follows.
Then \eqref{gattoboy} follows from \eqref{tritone}
and since $\int_0^1 1/\sqrt t \sqrt{1-t}\, dt=\pi.$
\eproof

\begin{remark}\label{figaro}
By Lemma \ref{giga} the function
 $\partial_E \act_n^{(2j-1)}(E,\hat\act)$ is a holomorphic
 function of the complex variable
 $\z=E-E_{2j-1}$.
  This is not the case of the functions 
 $\partial_E \act_n^{(2j-1),\pm}(E,\hat\act)$ and 
 $\partial_E \bar\act_n^{(2j-1),\pm}(E)$ 
 that are holomorphic
 functions of $\sqrt\z$ only.
 The holomorphicity of
 $\partial_E \act_n^{(2j-1)}=\partial_E \act_n^{(2j-1),+}
 +\partial_E \act_n^{(2j-1),-}$
 is due to parity cancellations. 
\end{remark}

\begin{lemma}\label{giga2}
Set\footnote{$r_\dag$ was defined in \eqref{spadadefoco3}.
Use \eqref{harlock} to prove the inequality.}
\begin{equation}\label{tropicana}
\mathtt r_2:=\min\left\{
\frac{s_0^{49}\b^{30}}{2^{214} \morse^{30}},\
\frac{r_0^2}{2^{10}\morse}\right\}
<
\frac{\b r_\dag^2}{2^{10}\morse^2}
\,.
\end{equation}
The function 
$\partial_E \act_n^{(2j-1)}(E_{2j-2}(\hat\act)-z\morse,\hat\act),$
initially defined for $0<z<\mathtt r_2$ and
 $\hat\act\in \hat D,$
 has holomorphic extension to the complex set
 $\{z\in\C_*\ \ \text{s.t.}\ \  |z|<\mathtt r_2\}\times \hat D_{r_0}.$
In particular
\begin{equation}\label{maracaibo}
\partial_E \act_n^{(2j-1)}(E_{2j-2}(\hat\act)-z \morse,\hat\act)
=
\f^{(2j-1)}(z,\hat\act)+\psi^{(2j-1)}(z,\hat\act)\ln z\,,
\end{equation}
where $\f^{(2j-1)}( z,\hat\act)$ and $\psi^{(2j-1)}(z,\hat\act)$
are holomorphic function in the set $\{|z|<\mathtt r_2\}\times \hat D_{r_0}$
with
\begin{equation}\label{trocadero}
\sup_{\{|z|<\mathtt r_2\}\times \hat D_{r_0}} |\f^{(2j-1)}(z,\hat\act)|
\leq
\frac{ 2^{84}\morse^{8}}{s_* s_0^{13}\b^{17/2} }
\,,\qquad 
\sup_{\{|z|<\mathtt r_2\}\times \hat D_{r_0}} |\psi^{(2j-1)}(z,\hat\act)|
\leq 
\frac{2^7\sqrt\morse}{\b s_0}\,,
\end{equation}
and
\begin{equation}\label{otello}
\inf_{\hat\act\in \hat D_{r_0}}
|\psi^{(2j-1)}(0,\hat\act)|\geq \frac{s_0}{32 \sqrt{\morse}}\,.
\end{equation}
Moreover the functions 
$\partial_E \act_n^{(2j-1),-}(E_{2j-2}(\hat\act)-z \morse,\hat\act)$
and
$\partial_E \act_n^{(2j-1),+}(E_{2j-2}(\hat\act)-z \morse,\hat\act)$
have holomorphic extension to the complex sets
$\{z\in\C_*\ \ \text{s.t.}\ \  |z|<\mathtt r_2\}\times \hat D_{r_0}$
and
 $\{  |z|<\mathtt r_2\}\times \hat D_{r_0},$ respectively,
 with
 \begin{eqnarray}\label{maracaibo2}
\partial_E \act_n^{(2j-1),-}(E_{2j-2}(\hat\act)-z \morse,\hat\act)
&=&
\f^{(2j-1),-}(z,\hat\act)+\psi^{(2j-1)}(z,\hat\act)\ln z\,,
\nonumber
\\
\partial_E \act_n^{(2j-1),+}(E_{2j-2}(\hat\act)-z \morse,\hat\act)
&=&
\f^{(2j-1),+}(z,\hat\act)
\end{eqnarray}
and
\begin{equation}\label{trocadero2}
\sup_{\{|z|<\mathtt r_2\}\times \hat D_{r_0}} |\f^{(2j-1),\pm}(z,\hat\act)|
\leq
\frac{ 2^{84}\morse^{8}}{s_* s_0^{13}\b^{17/2} }\,.
\end{equation}
Finally 
\begin{equation}\label{trocadero*}
\sup_{\{|z|<\mathtt r_2\}\times \hat D_{r_0/2}} 
|\partial_{\hat\act}\f^{(2j-1)}(z,\hat\act)|
\leq
M_\f \checco
\,,\qquad 
\sup_{\{|z|<\mathtt r_2\}\times \hat D_{r_0/2}} 
|\partial_{\hat\act}\psi^{(2j-1)}(z,\hat\act)|
\leq 
M_\psi\checco\,,
\end{equation}
where
$M_\f,M_\psi$ are suitable large constants.
\end{lemma}
\begin{remark}
 The constants $M_\f,M_\psi$ can be explicitly evaluated
 but we do not do this here!
\end{remark}
\proof
First set
$$
\z:=z\morse\,.
$$
Note that
the function 
\begin{eqnarray}\label{scarborough}
&&\sqrt{E-E_{2j-1} }=
\sqrt{\Delta E -\zeta}=\sqrt{u(\z)}\,,
\nonumber
\\
&& \text{with}\ \ \ 
\Delta E:=E_{2j-2} -E_{2j-1}\,,
\ \ u=u(\zeta):=\Delta E -\zeta
\end{eqnarray}
is holomorphic for $|\zeta|<
\b/2.$ 
Noting that by \eqref{goffredo} 
we have $\b\leq \bar E_{2j-2} -\bar E_{2j-1}\leq 2\morse,$
recalling \eqref{october}, \eqref{legna} (and \eqref{genesis})
we have 
\begin{equation}\label{mandrake}
u(\z)\, \in\, (\b,2\morse)_{r_0^2/2^7} 
\, \subset\, 
(0,R_0^2/2)_{r_0^2/2^7}
\qquad
{\rm and}
\qquad
|u(\z)|< 5\morse
\end{equation}
for $|\z|<\mathtt r_2\morse$.
Indeed by \eqref{ladispoli}
we have $\Re (\bar E_{2j-2} -
\bar E_{2j-1})=
\bar E_{2j-2} -\bar E_{2j-1}\geq \b$
and
$$
\Re \Delta E=
\Re(E_{2j-2} -E_{2j-1})
\stackrel{\eqref{october}}\geq 
\b -4\checco 
\stackrel{\eqref{genesis}}\geq 
\b -\frac{\b^2}{64\morse}
\stackrel{\eqref{harlock}}\geq 
\frac{\b}{2}\,.
$$
In particular note that
\begin{equation}\label{cipollotto}
|\z|\leq \b/4 \qquad \Longrightarrow\qquad
|\Delta E-\z|\geq \b/4\,.
\end{equation}
Recalling that
$$
E=E_{2j-2}(\hat\act)-z\morse\,,
$$
we set
\begin{eqnarray}
\mathcal I_1 &:=&
\frac{\sqrt u}{2\pi \sqrt t}
\frac{1+\tilde b\big(t(E-E_{2j-1}),\Sa_{2j}(E(t))\big)}{\partial_{\sa} \Gm\big(
\Sa_{2j}(E(t))  \big)}
\nonumber
\\
&\stackrel{\eqref{allemande}}=&
\frac{\sa_{2j}-\sa_{2j-1}}{E_{2j}-E_{2j-1}}
\frac{\sqrt u}{2\pi \sqrt t}
\frac{1+\breve{\tilde b}_{2j}\big(t(1+\breve E),
\breve\Sa_{2j}\big(
\breve E(t) \big)
 \big)}
{
\partial_{\breve \sa}
\breve\Gm_{2j}\Big(\breve  \Sa_{2j}\big(\breve E(t) \big)\Big)
}
\,,
\nonumber
\\
&{\rm for}&\breve E=\l_{2j}\big(E_{2j-2}(\hat\act)-z\morse\big)
\,,
\nonumber
\\
&&\quad
\nonumber
\\
\mathcal I_2 &:=&
-\frac{\sqrt u}{2\pi \sqrt t}
\frac{1+\tilde b\big(t(E-E_{2j-1}),\Sa_{2j-1}(E(t))\big)}{\partial_{\sa} \Gm\big(
\Sa_{2j-1}(E(t))  \big)}
\nonumber
\\
&\stackrel{\eqref{allemande}}=&
\frac{\sa_{2j-1}-\sa_{2j-2}}{E_{2j-1}-E_{2j-2}}
\frac{\sqrt u}{2\pi \sqrt t}
\frac{1+\breve{\tilde b}_{2j-1}\big(t(1+\breve E),
\breve\Sa_{2j-1}\big(
\breve E(t) \big)
 \big)}
{
\partial_{\breve \sa}
\breve\Gm_{2j-1}\Big(\breve  \Sa_{2j-1}\big(\breve E(t) \big)\Big)
}\,,
\nonumber
\\
&{\rm for}&\breve E=\l_{2j-1}\big(E_{2j-2}(\hat\act)-z\morse\big)
\,,
\label{mela}
\end{eqnarray} 
(recalling that 
$\breve E(t)
=-t  + (1-t)\breve E$
was defined in \eqref{spade}).
Then we split
$$
\partial_E \act_n^{(2j-1)}=I_1+I_2+I_3+I_4\,,
$$
where
\begin{eqnarray}
I_1 &:=&\int_{1-t_1}^1 (\mathcal I_1+\mathcal I_2)dt\,,
\nonumber
\\
I_2 &:=& \int_0^{t_1} \mathcal I_2\, dt\,,
\nonumber
\\
I_3 &:=& \int_0^{1-t_1} \mathcal I_1\, dt\,,
\nonumber
\\
I_4 &:=& \int_{t_1}^{1-t_1} \mathcal I_2\, dt\,,
\label{prima}
\end{eqnarray}
and
\begin{equation}\label{arancia}
 t_1:=\frac{r_\dag^2}{2^6\morse} 
 =
 \frac{s_0^{23} \b^{15}}{2^{100} \morse^{15}}
 <\frac{r_\diamond^2}{2^{66}\morse}
 \,,
 \end{equation}
recalling  \eqref{spadadefoco3}.

\smallskip

{\sl In the beginning we consider the real case, namely
 $0<\z<\mathtt r_2\morse$ (and, therefore,
$0<z<\mathtt r_2$) and $\hat\act\in D$. Note that in this case
we have $u\geq \b/4>0.$ Then we will
rewrite the functions $I_i$ in a different way
such that it clearly appears that $I_1,I_3,I_4$ actually have
a holomorphic extension for 
$|\z|<\mathtt r_2\morse$ (and $\hat\act\in D_{r_0}$), while $I_2$
has a  holomorphic extension for 
$|\z|<\mathtt r_2\morse$, $\z\in\C_*,$
(and $\hat\act\in D_{r_0}$)
due to the presence of a logarithmic term.
Note that we will omit to explicitly write the dependence on
the dummy variable $\hat\act$, with respect to which
all the estimates are uniform.}

\smallskip

$\bullet$ {\sl Study of } $I_1.$

\noindent
Recalling the definition of 
$w(t)$ in \eqref{capricciosa}, we have that
$$
1-t_1\leq t\leq 1 \quad\Longrightarrow\quad
|w(t)|\leq 2\sqrt\morse \sqrt{t_1}
\stackrel{\eqref{arancia}}<r_\diamond
\quad\Longrightarrow\quad
\eqref{tonino}\ \ \text{holds}\,.
$$
Then, recalling the definition of $G$ in \eqref{grattaebasta}
and  \eqref{gufetta2}, we get
$$
I_1=
\frac{1}{2\pi}
\int_{1-t_1}^1
\frac{G\big(E-E_{2j-1},
(1-t)(E-E_{2j-1}),\hat\act\big)}
{\sqrt t\sqrt{1-t}}
\, dt\,,
$$
which by \eqref{tritone}  and \eqref{legna}
gives
that $I_1$ is actually a holomorphic function of
$z=(E-E_{2j-1})/\morse$ in the ball $\{|z|<\mathtt r_2\}$ 
with estimate
\begin{equation}\label{strunz}
\sup_{\{|z|<\mathtt r_2\}\times \hat D_{r_0}}
|I_1|\leq 
\frac{32\sqrt\morse}{s_0\b}\,.
\end{equation}

$\bullet$ {\sl Study of } $I_2.$

\noindent

Set 
$$
\mathtt v=\mathtt v(t):=\sqrt{E_{2j-2}-E(t)}
=\sqrt{t\Delta E+(1-t)\zeta}
=\sqrt{\zeta+ u t}\,,\ \ 
\text{where}\ \ 
u=u(\zeta):=\Delta E -\zeta
\,.
$$
By \eqref{gargano}
\begin{equation}\label{salieri}
\Sa_{2j-1}(E(t))
=\sa_{2j-2}+\mathtt v \Sa_{2j-1,+}(\mathtt v)\,,
\end{equation}
when, recalling \eqref{burrata},
$
|\mathtt v|<
r_\diamond
.$
Since (recall \eqref{spadadefoco3}) 
$$
\mathtt r_2\leq r_\dag^2/8\morse
$$ by \eqref{tropicana} and \eqref{harlock}
we have that 
\begin{equation}\label{arancia2}
0\leq t\leq t_1=\frac{r_\dag^2}{2^6\morse} \quad
\Longrightarrow
\quad
|\mathtt v(t)|\leq \frac{r_\dag}{2}\leq \frac{r_\diamond}{2^{31}}\,.
\end{equation}
Then, by
\eqref{oxiana3}
$$
|\mathtt v(t) \Sa_{2j-1,+}(\mathtt v(t))|
\leq
\frac{16 r_\dag \morse}{s_0^{3/2} \b^{3/2}}
\stackrel{\eqref{spadadefoco3}}\leq
\frac{s_0^{10}}{2^{43}}
\frac{\b^{6}}{\morse^{6}}
\stackrel{\eqref{fusek},\eqref{harlock}}\leq
\frac{\sa_\diamond}{2^{30}}
\,,\qquad
\forall\, 0\leq t\leq t_1\,.
$$
In conclusion, by Lemma \ref{anna}
we have that, for every $0\leq t\leq t_1,$
\begin{equation}\label{mozart}
\partial_{\sa} \Gm\big(
\Sa_{2j-1}(E(t))  \big)
=
\mathtt v(t) \Sa_{2j-1,+}(\mathtt v(t))
\, \check\Gm\big(
\mathtt v(t) \Sa_{2j-1,+}(\mathtt v(t))
\big)
=
\mathtt v(t) \Gm_\diamond(\mathtt v(t))\,,
\end{equation}
where
\begin{equation}\label{guardian}
\Gm_\diamond(\mathtt v):=
\Sa_{2j-1,+}(\mathtt v)
\, \check\Gm\big(
\mathtt v \Sa_{2j-1,+}(\mathtt v)
\big)\,,\qquad
{\rm with}\quad
\check \Gm(0)=\partial_{\sa\sa}\Gm(\sa_{2j-2})
\end{equation}
is a holomorphic function in the ball $|\mathtt v|<r_\diamond$
and, again by Lemma \ref{anna} and  \eqref{spadadefoco2}, we get
\begin{equation}\label{corelli} 
\sup_{\{|\mathtt v|<r_\dag\}\times \hat D_{r_0}}
 \frac{1}{|\Gm_\diamond|}
 \leq
\frac{16\sqrt\morse}{\b s_0}\,.
 \end{equation}
We have
\begin{eqnarray*}
I_2
&=&
\frac{\sqrt u}{2\pi}
\int_0^{ t_1}
\frac{1+\tilde b\big(t(\Delta E-\zeta),\Sa_{2j-1}(E(t))\big)}{\sqrt t\ \partial_{\sa} \Gm\big(
\Sa_{2j-1}(E(t))  \big)}
\, dt
\\
&=&
\frac{\sqrt u}{2\pi}
\int_0^{ t_1}
\frac{1+\tilde b\big(t(\Delta E-\zeta),
\sa_{2j-2}+\mathtt v(t) \Sa_{2j-1,+}(\mathtt v(t))
\big)}
{\sqrt t\ \mathtt v(t) \Gm_\diamond(\mathtt v(t))}
\, dt
\,.
\end{eqnarray*}
Let us consider the holomorphic function
\begin{equation}\label{blind}
\frak b(\mathtt v)=
\frak b(y,\mathtt v)
=
\frak b(y,\mathtt v,\hat\act):= 2
\frac{1+\tilde b\big(y,
\sa_{2j-2}+\mathtt v \Sa_{2j-1,+}(\mathtt v)
\big)}
{ \Gm_\diamond(\mathtt v)}
\end{equation}
with\footnote{Recall \eqref{corelli} and
\eqref{3holesquater}.}
\begin{equation}\label{vivaldi}
\sup_{(0,R_0^2)_{r_0^2/64}\times\{|\mathtt v|<r_\dag\}\times \hat D_{r_0}}
|\frak b(v,\mathtt v,\hat\act)|
 \leq
\frac{64\sqrt\morse}{\b s_0}
\end{equation}
Split it in its  even and odd part w.r.t. $\mathtt v,$
namely\footnote{Omitting for brevity the dependence on
$v$ and $\hat\act.$}
$\frak b(\mathtt v)=\frak b_{\rm e}(\mathtt v)+\frak b_{\rm o}(\mathtt v),$
where  $\frak b_{\rm e}(\mathtt v):=(\frak b(\mathtt v)+
\frak b(-\mathtt v))/2$ and
$\frak b_{\rm o}(\mathtt v):=(\frak b(\mathtt v)-
\frak b(-\mathtt v))/2,$ for which the same estimate as \eqref{vivaldi}
holds.
Since $\frak b_{\rm e}(\mathtt v)$ is an even function 
there exists a holomorphic function
$b_{\rm e}(w)$ such that $\frak b_{\rm e}(\mathtt v)=
b_{\rm e}(\mathtt v^2)$ with estimate
\begin{equation}\label{vivaldi2}
\sup_{(0,R_0^2)_{r_0^2/64}\times\{|w|<r_\dag^2\}\times \hat D_{r_0}}
|b_{\rm e}(v,w,\hat\act)|
 \leq
\frac{64\sqrt\morse}{\b s_0}\,.
\end{equation}
On the other hand, since $\frak b_{\rm o}(\mathtt v)$ is odd,
the function $\frak b_{\rm o}(\mathtt v)/\mathtt v$ is also holomorphic
with estimate
\begin{equation}\label{vivaldi3}
\sup_{(0,R_0^2)_{r_0^2/64}\times\{|\mathtt v|<r_\dag/2\}\times \hat D_{r_0}}
|\frak b_{\rm o}(v,\mathtt v,\hat\act)/\mathtt v|
 \leq
\frac{2^7\sqrt\morse}{r_\dag\b s_0}\,.
\end{equation}
by Cauchy estimates.
Since $\frak b_{\rm o}(\mathtt v)/\mathtt v$ is an even function 
there exists a holomorphic function
$b_{\rm o}(w)$ such that $\frak b_{\rm o}(\mathtt v)/\mathtt v=
b_{\rm o}(\mathtt v^2)$ with estimate
\begin{equation}\label{vivaldi4}
\sup_{(0,R_0^2)_{r_0^2/64}\times\{|w|<r_\dag^2/4\}\times \hat D_{r_0}}
|b_{\rm o}(v,w,\hat\act)|
 \leq
\frac{2^7\sqrt\morse}{r_\dag\b s_0}\,.
\end{equation}
Recollecting 
$$
\frak b(v,\mathtt v)=\mathtt v  b_{\rm o}(v,\mathtt v^2)+   b_{\rm e}(v,\mathtt v^2)\,.
$$ 
Noting that 
$$
\mathtt v^2(t)=\z + u t
$$
and setting
\begin{eqnarray*}
I_{2,\rm e}(\z)
&:=&
\frac{\sqrt u}{4\pi}
\int_0^{ t_1}
\frac{b_{\rm e}\big(u t, \z + u t\big)}
{\sqrt t \ \sqrt{\z + u t}}
\ dt\,,
\\
I_{2,\rm o}(\z)
&:=&
\frac{\sqrt u}{4\pi}
\int_0^{t_1}
\frac{b_{\rm o}\big(u t, \z + u t\big)}{\sqrt t}
\ dt
=
\frac{\sqrt u}{2\pi}
\int_0^{\sqrt t_1}
b_{\rm o}\big(u s^2, \z + u s^2\big)
\ ds\,,
\end{eqnarray*}
we get
$$
I_2=I_{2,\rm e} + I_{2,\rm o}\,.
$$
Recalling that $|u|< 5\morse$ and 
$\mathtt r_2\leq r_\dag^2/8\morse$ we get, for every
$0\leq s\leq \sqrt{t_1},$
\begin{equation}\label{olivia}
|us^2|\leq |u t_1|< r_\dag^2/8\,,
\qquad
|\z +us^2|< r_\dag^2/4
\,,
\end{equation}
for $|\z|\leq \mathtt r_2\morse.$
It is obvious that $I_{2,\rm o}$ is holomorphic,
moreover, by \eqref{vivaldi4}
\begin{equation}\label{gigante}
\sup_{|\z|\leq \mathtt r_2\morse}
|I_{2,\rm o}(\z)|
\leq
\frac{2^7\morse}{r_\dag\b s_0}\sqrt t_1
\stackrel{\eqref{arancia}}=
\frac{2^4\sqrt\morse}{\b s_0}\,.
\end{equation}
Regarding $I_{2,\rm e}$ we split it as
$$
I_{2,\rm e}=I_{2,1}+I_{2,2}:=\frac{\sqrt u}{4\pi}
\int_0^{16\z/u} +
\frac{\sqrt u}{4\pi}\int_{16\z/u}^{t_1}\,.
$$
We claim that $I_{2,1}$ is a holomorphic function of
$\z$ in the ball $\{|\z|<\mathtt r_2\morse\}.$
Indeed changing variable
$t=y^2\z/u$ we get
$$
I_{2,1}
=
\frac{1}{2\pi}
\int_0^4
\frac{b_{\rm e}\big(\z y^2, \z(1+y^2)\big)}
{\sqrt{1+y^2}}
\ dy\,,
$$
which is obviously holomorphic on $\{|\z|<\mathtt r_2\morse\}$
with estimate
\begin{equation}\label{melanzane}
\sup_{|\z|<\mathtt r_2\morse} |I_{2,1}(\z)|\leq 
\frac{2^{10}\sqrt\morse}{\b s_0}
\end{equation}
by  \eqref{vivaldi2}.
On the other hand
$$
I_{2,2}
=
\frac{\sqrt u}{4\pi}
\int_{16\z/u}^{ t_1}
\frac{b_{\rm e}\big(u t, \z + u t\big)}
{\sqrt t \ \sqrt{\z + u t}}
\ dt
\,,
$$
substituting $w=u t$
becomes
\begin{equation}\label{king}
I_{2,2}
=
\frac{1}{4\pi}
\int_{16 \z}^{u t_1}
\frac{b_{\rm e}(w, \z + w)}
{\sqrt{w}\sqrt{\z + w}}
\ dw
=
\frac{1}{4\pi}
\int_{16 \z}^{u t_1}
\frac{b_{\rm e}(w, \z + w)}
{w\sqrt{1 + \frac{\z}{w}}}
\ dw
\,.
\end{equation}
By \eqref{vivaldi2} we write, for $|\z|,|w|<r_\dag^2/2$
\begin{equation}\label{king2}
b_{\rm e}(w, \z + w)=\sum_{h\geq 0} b_h(\z)w^h\,,
\end{equation}
for suitable holomorphic functions $b_h(\z)$
satisfying
\begin{equation}\label{helmut}
\sup_{|\z|<r_\sharp}
|b_h(\z)|
 \leq
 M_\sharp r_\sharp^{-h}\,,
 \qquad
 {\rm with}\qquad
 M_\sharp:=
\frac{64\sqrt\morse}{\b s_0}\,,\quad
r_\sharp:=r_\dag^2/2\,.
\end{equation}
Let us develop, for $|y|<1,$
\begin{equation}\label{king4}
\frac{1}{\sqrt{1+y}}=\sum_{k\geq 0} c_k y^k\,,
\qquad c_k:=\binom{-1/2}{k}
\end{equation}
and note that $|c_k|<1.$
For $|\z|<|w|<r_\sharp,$ 
we have
\begin{equation}\label{king5}
\frac{b_{\rm e}(w, \z + w)}
{\sqrt{1 + \frac{\z}{w}}}
=
\sum_{n\in\Z} d_n(\z)w^n\,,\qquad
{\rm where}\qquad
d_n(\z):=\sum_{k\geq \max\{0,-n\}} c_k b_{k+n}(\z)\z^k\,,
\end{equation}
in particular, for $n\geq 0$
\begin{equation}\label{loosing}
d_n(\z)=\sum_{k\geq 0} c_k b_{k+n}(\z)\z^k\,,
\qquad
d_{-n}(\z)=\sum_{k\geq n} c_k b_{k-n}(\z)\z^k
=\z^n\sum_{k\geq 0} c_{k+n} b_{k}(\z)\z^k
\,.
\end{equation}
By \eqref{helmut} we have, for $n\geq 0,$
\begin{equation}\label{walcha}
\sup_{|\z|<\frac{r_\sharp}{2}} |d_n(\z)|\leq 2M_\sharp r_\sharp^{-n}\,,
\qquad
\sup_{|\z|<\frac{r_\sharp}{2}} |d_{-n}(\z)|\leq 2M_\sharp |\z|^n\,.
\end{equation}
Recalling \eqref{olivia} and we note  that,
in the real case, for every  $0<\z<\mathtt r_2\morse$
(and $\hat\act\in D$)
we have  $16\z\leq w\leq ut_1<r_\dag^2/8=r_\sharp/4$
(recall \eqref{helmut}) and then
$|\z/w|\leq 1/16$. Therefore, for every
$0<\z<\mathtt r_2\morse$,
 the first series in \eqref{king5}
totally converges in the interval $16\z\leq w\leq ut_1$
and we get
\begin{equation}
I_{2,2}=\frac{1}{4\pi}\int_{16 \z}^{u t_1}
\frac{b_{\rm e}(w, \z + w)}
{w\sqrt{1 + \frac{\z}{w}}}
\ dw
=
\frac{1}{4\pi}
\sum_{n\in\Z} d_n(\z)
\int_{16 \z}^{u t_1} w^{n-1}\, dw
=\psi(z)\ln z+ I_{2,3}
\,,
\label{crimson}
\end{equation}
where
\begin{equation}\label{king3}
\psi(z):=-\frac{d_0(z\morse)}{4\pi}
\end{equation}
and
\begin{equation}\label{recremisi}
I_{2,3}:=
\frac{1}{4\pi}d_0(\z)\ln(u t_1/16\morse)
+\frac{1}{4\pi}\sum_{n\neq 0} \frac{d_n(\z)}{n}\big( (ut_1)^n-(16\z)^n \big)\,.
\end{equation}
Note that, except for the first one, all the other addenda
in the last line are holomorphic functions of $\z$ in the ball
$\{|\z|<\mathtt r_2\morse\};$
for example  by \eqref{loosing}
$$
\sum_{n>0} \frac{d_{-n}(\z)}{n(16\z)^n}
=
\sum_{n>0} \frac{1}{n 16^n}
\sum_{k\geq 0} c_{k+n} b_{k}(\z)\z^k\,.
$$
Recalling \eqref{walcha} and \eqref{cipollotto}
\begin{equation}\label{girotondo}
\sup_{|z|<\frac{r_\sharp}{2\morse}} |\psi(z)|\leq 
2M_\sharp\,,
\end{equation}
which, recalling \eqref{helmut}, implies \eqref{trocadero}.
By \eqref{king2},\eqref{king4},\eqref{king5},
\eqref{blind},\eqref{guardian}, \eqref{stambecco}
we have
\begin{eqnarray}\label{king6}
d(0)
&=&
c_0 b_0(0)=b_{\rm e}(0,0)=\frak b(0,0)
=
2
\frac{1+\tilde b(0,
\sa_{2j-2}
)}
{ \Gm_\diamond(0)}
=
2
\frac{1+\tilde b(0,
\sa_{2j-2}
)}
{\Sa_{2j-1,+}(0)
 \check\Gm(0)}
 \nonumber
 \\
&=&
2
\frac{1+\tilde b(0,
\sa_{2j-2}
)}
{\sqrt{-2/\partial_{\sa\sa}\Gm(\sa_{2j-2})}
\partial_{\sa\sa}\Gm(\sa_{2j-2})
}
=
-\sqrt 2
\frac{1+\tilde b(0,
\sa_{2j-2}
)}
{\sqrt{-\partial_{\sa\sa}\Gm(\sa_{2j-2})}
}
 \end{eqnarray}
 so that
 \begin{equation}\label{jenny}
\psi(0)=
\frac{1+\tilde b(0,
\sa_{2j-2}
)}
{
\sqrt{8\pi}
\sqrt{-\partial_{\sa\sa}\Gm(\sa_{2j-2})}
}\,.
\end{equation}
Then, by \eqref{ciccio3} and Cauchy estimates 
we get \eqref{otello}.
\\
From \eqref{recremisi}, \eqref{king},\eqref{crimson},\eqref{walcha},\eqref{olivia},
\eqref{cipollotto},\eqref{arancia},\eqref{tropicana} we get
\begin{eqnarray}
\sup_{|z|< \mathtt r_2}
\left|
I_{2,3}(z\morse)
\right|
&\leq&
 2M_\sharp \ln(2^7\morse/r_\dag^2)
+ 2M_\sharp
\sum_{n>0}\left(
\frac{1}{4^n} + \frac{1}{2^n}
+\Big(\frac{4\mathtt r_2}{\b t_1}\Big)^n+\frac{1}{16^n}
\right)
\nonumber
\\
&\leq&
 2M_\sharp\ln(2^7\morse/r_\dag^2)+ 16 M_\sharp
 \stackrel{\eqref{spadadefoco3}}=
 2M_\sharp\ln\left(
 \frac{2^{101} \morse^{15}}{s_0^{23} \b^{15}}
 \right)+ 16 M_\sharp
 \nonumber
\\
&\leq&
2^{14}M_\sharp\sqrt{\frac{\morse}{s_0 \b}}
\stackrel{\eqref{helmut}}=
2^{20} \frac{\morse}{s_0^{3/2} \b^{3/2}}
\label{kingcrimson}
\,,
\end{eqnarray}
by \eqref{ocarina} and also\footnote{Noting that $\ln x\leq \sqrt x.$}
\begin{equation}\label{gomme}
\ln \frac{2^{101} \morse^{15}}
{s_0^{23} \b^{15}}
\leq
\ln \frac{2^{109} \morse^{23}}
{s_0^{23} \b^{23}}
\leq 
109+23\ln\frac{\morse}{s_0 \b}
\leq
109+23\sqrt{\frac{\morse}{s_0 \b}}
\leq
2^{12}\sqrt{\frac{\morse}{s_0 \b}}\,.
\end{equation}
Recollecting we have
\begin{equation}\label{something}
I_2(z)=I_{2,{\rm o}}(z\morse)+I_{2,1}(z\morse)
+I_{2,3}(z\morse)+\psi(z)\ln z\,.
\end{equation}

\medskip

$\bullet$ {\sl Study of } $I_3.$

\noindent
We claim that
\begin{equation}\label{numb}
|z|<\mathtt r_2
\qquad
\Longrightarrow
\qquad
\breve E=\l_{2j}(E_{2j-2}-z\morse)\,
\in\, \O_{\breve{\mathtt r}}\,.
\end{equation}
In order to prove \eqref{numb}
we note that, since $\bar\l_{2j}$ (defined in \eqref{acrobat})
is an increasing function\footnote{Note that
$\bar\l_i(\bar E_{i'})$ is real.},
\begin{equation}\label{patato}
-1=\bar\l_{2j}(\bar E_{2j-1})<
\bar\l_{2j}(\bar E_{2j-2})<
\bar\l_{2j}(\bar E_{2j})=1
\,,
\end{equation}
recalling that
\begin{equation}\label{torrone}
\bar E_{2j-1}<\bar E_{2j-2}<\bar E_{2j}
\end{equation}
 by Remark \ref{stiffe}.
 By \eqref{carciofino} we have
 $$
 | \l_{2j}(E_{2j-2}-z\morse)-\l_{2j}(\bar E_{2j-2})|
 \leq \frac{4}{\b}\big(
 |E_{2j-2}-\bar E_{2j-2}|+\morse \mathtt r_2
 \big)
 \stackrel{\eqref{october}}\leq
 \frac{4}{\b}\big(
 2\checco+\morse \mathtt r_2
 \big)
 $$
and
$$
|\l_{2j}(\bar E_{2j-2})-\bar\l_{2j}(\bar E_{2j-2})|
\leq \frac{48\morse\checco}{\b^2}
\stackrel{\eqref{harlock}}\leq
\frac{96\checco}{\b}
$$
by \eqref{carciofino2}.
Then
\begin{equation}\label{erpomata}
 | \l_{2j}(E_{2j-2}-z\morse)-\bar\l_{2j}(\bar E_{2j-2})|
 \leq \frac{2^7\checco}{\b}+
 \frac{4\morse \mathtt r_2}{\b}
  \stackrel{\eqref{tropicana}}=
 \frac{2^7\checco}{\b}+
 \frac{s_0^{49}\b^{29}}{2^{212} \morse^{29}}
\stackrel{\eqref{genesis},\eqref{harlock}}<
 \frac{\breve{\mathtt r}}{2}
\end{equation}
($\breve{\mathtt r}$ defined in \eqref{coppe}).
\\
Moreover, since
$$
1-\bar\l_{2j}(\bar E_{2j-2})
=2\frac{\bar E_{2j}-\bar E_{2j-2}}{\bar E_{2j}-\bar E_{2j-1}}
\stackrel{\eqref{goffredo},\eqref{ladispoli}}\geq
\frac{\b}{\morse}\,,
$$
by \eqref{erpomata} we get
\begin{equation}\label{bellicapelli}
1-\Re \breve E \geq \frac{\b}{2\morse}\,.
\end{equation}
Analogously we have
\begin{equation}\label{bellicapelli2}
|1+ \breve E|\geq
1+\Re \breve E \geq \frac{\b}{2\morse}\,.
\end{equation}
Finally   \eqref{patato}, \eqref{erpomata},
\eqref{bellicapelli}, 	\eqref{bellicapelli2} imply \eqref{numb}.

\noindent
By \eqref{mela} and recalling \eqref{prima} we get
\begin{equation}\label{roscio}
 I_3
=
\frac{\sa_{2j}-\sa_{2j-1}}{E_{2j}-E_{2j-1}}
\frac{\sqrt u}{2\pi }
\tilde{I}_3\,,\qquad
{\rm with}\qquad
\tilde{I}_3:=
\int_0^{1-t_1}
\frac{1+\breve{\tilde b}_{2j}\big(t(1+\breve E),
\breve\Sa_{2j}\big(
\breve E(t) \big)
 \big)}
{\sqrt t\,
\partial_{\breve \sa}
\breve\Gm_{2j}\Big(\breve  \Sa_{2j}\big(\breve E(t) \big)\Big)
}\, dt\,,
\end{equation}
where 
$
\breve E=\l_{2j}(E_{2j-2}-z\morse)\,.
$
Since $\breve E \in \O_{\breve{\mathtt r}}$
by \eqref{numb}, 
we can apply estimate \eqref{cogli2} obtaining
$$
|\tilde{I}_3|
\leq
 \int_{0}^{1}
\frac{
2
}
{\sqrt{t}\Big|\partial_{\breve \sa}
\breve{\FO}_{2j}\Big(\breve  
{\bar\Sa}_{2j}\big(\breve E(t) \big)\Big)\Big|}
\, dt\,.
$$
Since by \eqref{bellicapelli} we also have 
  $\Re \breve E<1$, by \eqref{dave}, 
  \eqref{bellicapelli}, \eqref{bellicapelli2}, \eqref{ocarina} 
  we get\footnote{Noting that $\ln x\leq \sqrt x.$}
\begin{equation*}
 |\tilde{I}_3|
\leq 
\frac{ 2^{32}\breve\morse}{\breve\b^3 \breve s^4}
\left( 
\frac{1}{\sqrt{|1+\breve E|}}
+
\ln \frac{1}{1-\Re \breve E}
\right)
\leq 
\frac{ 2^{38}\breve\morse}{\breve\b^3 \breve s^4}
\frac{\sqrt\morse}{\sqrt\b}
\leq \frac{ 2^{68}\morse^{15/2}}{s_0^{13}\b^{15/2} }
\,,
 \end{equation*}
by \eqref{exit}-\eqref{ladispoli4}.
Then by \eqref{roscio}, \eqref{mandrake}, \eqref{ladispoli3}
we get
\begin{equation}\label{soldatino}
|I_3|\leq 
\frac{ 2^{71}\morse^{8}}{s_0^{13}\b^{17/2} }\,.
\end{equation}

\medskip

$\bullet$ {\sl Study of } $I_4.$

It is similar to the case $I_3.$
By \eqref{mela} and recalling \eqref{prima} we get
\begin{equation}\label{rosciobis}
 I_4
=
\frac{\sa_{2j-1}-\sa_{2j-2}}{E_{2j-1}-E_{2j-2}}
\frac{\sqrt u}{2\pi }
\tilde{I}_4\,,\qquad
{\rm with}\qquad
\tilde{I}_4:=
\int_{t_1}^{1-t_1}
\frac{1+\breve{\tilde b}_{2j-1}\big(t(1+\breve E),
\breve\Sa_{2j-1}\big(
\breve E(t) \big)
 \big)}
{\sqrt t\,
\partial_{\breve \sa}
\breve\Gm_{2j-1}\Big(\breve  \Sa_{2j-1}\big(\breve E(t) \big)\Big)
}\, dt\,,
\end{equation}
where 
$$
\breve E=\l_{2j-1}(E_{2j-2}-z\morse)\,.
$$
In the integral in \eqref{rosciobis} we make the change of variable
$$
t=(1-t_1)\tilde t +t_1
$$
such that
\begin{equation}\label{caccola}
\breve E(t)=\tilde E(\tilde t):=\tilde E-(\tilde E+1)\tilde t\,,
\qquad
{\rm where}\quad
\tilde E:=\breve E-(\breve E+1)t_1
\end{equation}
and 
\begin{equation}\label{noia}
\tilde I_4=\sqrt{1-t_1}
\int_{0}^{1-t_2}
\frac{1+\breve{\tilde b}_{2j-1}\Big((1+\tilde E)\big(\tilde t+t_2\big),
\breve\Sa_{2j-1}\big(
\tilde E(\tilde t) \big)
 \Big)}
{\sqrt{\tilde t +t_2}\,
\partial_{\breve \sa}
\breve\Gm_{2j-1}\Big(\breve  \Sa_{2j-1}\big(\tilde E(\tilde t) \big)\Big)
}\, d\tilde t\,,
\end{equation}
where
$$
t_2:=\frac{t_1}{1-t_1}\,.
$$
Since $\bar\l_{2j-1}$ (defined in \eqref{acrobat})
is an increasing function\footnote{Note that
$\bar\l_i(\bar E_{i'})$ is real.},
\begin{equation}\label{patatobis}
-1=\bar\l_{2j-1}(\bar E_{2j-1})<
\bar\l_{2j-1}(\bar E_{2j-2})=1
\,,
\end{equation}
recalling \eqref{torrone}.
 By \eqref{carciofino} we have
 $$
 | \l_{2j-1}(E_{2j-2}-z\morse)-\l_{2j-1}(\bar E_{2j-2})|
 \leq \frac{4}{\b}\big(
 |E_{2j-2}-\bar E_{2j-2}|+\morse \mathtt r_2
 \big)
 \stackrel{\eqref{october}}\leq
 \frac{4}{\b}\big(
 2\checco+\morse \mathtt r_2
 \big)
 $$
and
$$
|\l_{2j-1}(\bar E_{2j-2})-\bar\l_{2j-1}(\bar E_{2j-2})|
\leq \frac{48\morse\checco}{\b^2}
\stackrel{\eqref{harlock}}\leq
\frac{96\checco}{\b}
$$
by \eqref{carciofino2}.
Then
\begin{equation}\label{erpomatabis}
 | \breve E-1|
 \leq \frac{2^7\checco}{\b}+
 \frac{4\morse \mathtt r_2}{\b}
  \stackrel{\eqref{tropicana}}=
 \frac{2^7\checco}{\b}+
 \frac{s_0^{49}\b^{29}}{2^{212} \morse^{29}}
\stackrel{\eqref{genesis},\eqref{harlock}}<
 \frac{\breve{\mathtt r}}{2}
\end{equation}
($\breve{\mathtt r}$ defined in \eqref{coppe}).
By \eqref{caccola} and  \eqref{arancia}
we have 
$$
|\tilde E-\breve E|\leq 3 t_1
$$
and, therefore, by \eqref{erpomatabis}
\begin{equation}\label{dartagnan}
|\tilde E-1|\leq \breve{\mathtt r}+3 t_1\,.
\end{equation}
Set\footnote{With $x_i\in\R.$}
\begin{equation}\label{gajardo}
 \breve E-1=:x_1+\ii x_2\,,\qquad
 |x_i|\leq  \frac{\breve{\mathtt r}}{2}\,,
\end{equation}
by \eqref{erpomatabis}.
We note that
\begin{equation}\label{cornetti}
-\frac12<\Re \tilde E-1=
x_1-(1+x_1)t_1\leq -\frac{t_1}{2}<0\,,
\end{equation}
indeed, the first inequality is immediate by \eqref{dartagnan};
 regarding the second one 
we note that, if $x_1\leq 0$ it is obvious, otherwise, when
$x_1> 0,$ we have
$x_1-(1+x_1)t_1<x_1-t_1\leq\frac{\breve{\mathtt r}}{2}-t_1
\leq -\frac{t_1}{2}<0$
since
$$
\breve{\mathtt r}<t_1
$$
(recalling \eqref{coppe}, \eqref{arancia}, \eqref{harlock}).
Moreover
\begin{equation}\label{asta}
|\Im \tilde E|=
|\Im \breve E|(1-t_1)<|\Im \breve E|
\leq \frac{\breve{\mathtt r}}{2}\,,
\end{equation}
by \eqref{gajardo}.
Recollecting by \eqref{cornetti} and \eqref{asta}
we get 
\begin{equation}\label{numbbis}
\tilde E\,
\in\, \O_{\breve{\mathtt r}}\,,
\qquad
\Re \tilde E<1
\,.
\end{equation}
By
\eqref{3holesquater}, \eqref{borotalco2},
we have that for every $0\leq t\leq 1-t_1$
$$
\left|
\breve{\tilde b}_{2j-1}\Big((1+\tilde E)\big(\tilde t+t_2\big),
\breve\Sa_{2j-1}\big(
\tilde E(\tilde t) \big)
 \Big)
\right|
\leq 
\frac{9\checco}{r_0^2}\leq 1
$$
by \eqref{genesis}.
Then by \eqref{noia} we get
$$
|\tilde{I}_4|
\leq
\int_{0}^{1}
\frac{2}
{\sqrt{\tilde t }\,\left|
\partial_{\breve \sa}
\breve\Gm_{2j-1}\Big(\breve  \Sa_{2j-1}\big(\tilde E(\tilde t) \big)\Big)
\right|}\, d\tilde t\,,
$$
Since $\tilde E \in \O_{\breve{\mathtt r}}$
by\footnote{With $\tilde E$ instead of $\breve E.$} \eqref{cogli3} we get
$$
|\tilde{I}_4|
\leq
\int_{0}^{1}
\frac{
4
}
{\sqrt{\tilde t}\Big|\partial_{\breve \sa}
\breve{\FO}_{2j-1}\Big(\breve  
{\bar\Sa}_{2j-1}\big(\tilde E(\tilde t) \big)\Big)\Big|}
\, d\tilde t\,.
$$
Since \eqref{numbbis}  holds, we can apply 
estimate\footnote{Again with $\tilde E$ instead of $\breve E.$} 
 \eqref{dave} obtaining 
\begin{equation*}
 |\tilde{I}_4|
\leq 
\frac{ 2^{33}\breve\morse}{\breve\b^3 \breve s^4}
\left( 
\frac{1}{\sqrt{|1+\tilde E|}}
+
\ln \frac{1}{1-\Re \tilde E}
\right)
\,.
 \end{equation*}
Then by \eqref{dartagnan}, \eqref{cornetti} 
we get
\begin{eqnarray*}
|\tilde{I}_4|
&\leq& 
\frac{ 2^{35}\breve\morse}{\breve\b^3 \breve s^4}
\ln \frac{1}{t_1}
\leq
\frac{2^{68}\morse^7}{\b^7 s_0^{13}}
\ln \frac{2^{100} \morse^{15}}
{s_0^{23} \b^{15}}
\leq
\frac{2^{80}\morse^{15/2}}{\b^{15/2} s_0^{27/2}}
\end{eqnarray*}
by \eqref{exit}-\eqref{ladispoli4}, \eqref{arancia} and \eqref{gomme}.
Then by \eqref{rosciobis}, \eqref{mandrake}, \eqref{ladispoli3}
we get
\begin{equation}\label{soldatinobis}
|I_4|\leq 
\frac{2^{83}\morse^{8}}{\b^{17/2} s_0^{27/2}}
\,.
\end{equation}

$\bullet$ {\sl Proof of \eqref{trocadero}}

\noindent
Recalling \eqref{prima} and \eqref{something}
we set
\begin{equation}\label{intheway}
\f(z)=I_1(z\morse)+ I_{2,{\rm o}}(z\morse)+I_{2,1}(z\morse)
+I_{2,3}(z\morse)+I_3(z\morse)+I_4(z\morse)\,.
\end{equation}
Then by 
\eqref{strunz},
\eqref{gigante},
\eqref{melanzane},
\eqref{kingcrimson},
\eqref{soldatino},
\eqref{soldatinobis}
we get
\begin{equation*}
\sup_{\{|z|<\mathtt r_2\}\times \hat D_{r_0}} |\f(z,\hat\act)|
\leq
\frac{32\sqrt\morse}{s_0\b}
+
\frac{2^4\sqrt\morse}{\b s_0}
+
\frac{2^{10}\sqrt\morse}{\b s_0}
+
2^{20} \frac{\morse}{s_0^{3/2} \b^{3/2}}
+
\frac{ 2^{71}\morse^{8}}{s_0^{13}\b^{17/2} }
+
\frac{2^{83}\morse^{8}}{\b^{17/2} s_0^{27/2}}\,.
\end{equation*}
Then by \eqref{ocarina} and \eqref{harlock} also 
 the first estimate in \eqref{trocadero} follows.
 
 $\bullet$ {\sl Concerning the functions  $\partial_E \act_n^{(2j-1),\pm}$}

\noindent
Recalling \eqref{mela} and \eqref{prima} we have that 
 $$
 \partial_E \act_n^{(2j-1),+}=
 I_3+I_5\,,
 \qquad
  \partial_E \act_n^{(2j-1),-}=
 I_2+I_4+I_6\,,
 $$ 
 with
 \begin{equation}\label{melamela}
I_5:= \int_{1-t_1}^1 \mathcal I_1 dt\,,\qquad
I_6:= \int_{1-t_1}^1 \mathcal I_2 dt\,.
\end{equation}

We have to consider only the term
$$ 
I_5
=
\frac{\sqrt u}{2\pi }
\int_{1-t_1}^1
\frac{1+\tilde b\big(t(E-E_{2j-1}),\Sa_{2j}(E(t))\big)}
{\sqrt t\, \partial_{\sa} \Gm\big(
\Sa_{2j}(E(t))  \big)}\, dt\,,
$$
 the term  
 $I_6$ being analogous.
Noting that
 $\mathtt r_2\morse<\mathtt r_\star$ we can argue as in
  Lemma \ref{giga}, obtaining 
\begin{eqnarray*}
\left|I_5\right|
&=&
\left|
\frac{1}{2\pi }
\int_{1-t_1}^1
\frac{1+\tilde b\big(t(E-E_{2j-1})-w^2(t),
\sa_{2j-1}+w(t) \Sa_{2j,-}(w(t))\big)}
{\sqrt t \sqrt{1-t}\, \Sa_{2j,-}(w(t)) \hat \Gm\big(
w(t) \Sa_{2j,-}(w(t))\big)}
\, dt
\right|
\\
&\leq&
 \frac{32\sqrt\morse}{s_0\b}
\end{eqnarray*}
(arguing as in \eqref{dolores} and recall the definition of $w(t)$ in \eqref{capricciosa}).
This proves \eqref{trocadero2}.
\eproof

 \medskip

 \begin{lemma}\label{bratislava}
If 
\begin{equation}\label{vienna}
|z|\leq \frac{s_0^4\b^2\mathtt r_2^2}{2^{32}\morse^2}	\,,
\qquad z\in\C_*
\end{equation}
then
\begin{equation}\label{budapest}
\left|
\frac{\partial_{EE} \act_n^{(2j-1)}(E_{2j-2}(\hat\act)-z \morse,\hat\act)}
{\big(\partial_{E} \act_n^{(2j-1)}(E_{2j-2}(\hat\act)-z \morse,\hat\act)\big)^3}
\right|
\geq
\frac{s_0^4\b^3}{2^{30} \morse^{3}}\frac{1}{\sqrt{|z|}}
\,.
\end{equation}
\end{lemma}
\proof
Deriving \eqref{maracaibo} w.r.t. $z$ 
we get, for $|z|\leq \mathtt r_2/2,z\in\C_*$
\begin{eqnarray}\nonumber
&&\morse |\partial_{EE} \act_n^{(2j-1)}(E_{2j-2}(\hat\act)-z \morse,\hat\act)|
\\
&&\geq
\frac{1}{|z|}|\psi^{(2j-1)}(0,\hat\act)|
-\frac{2}{\mathtt r_2} \sup_{|\tilde z|< \mathtt r_2}
|\psi^{(2j-1)}(\tilde z,\hat\act)| (1+|\ln z|)
-
\frac{2}{\mathtt r_2} \sup_{|\tilde z|< \mathtt r_2}
|\f^{(2j-1)}(\tilde z,\hat\act)|
\nonumber
\\
&&
\geq
\frac{1}{|z|}\frac{s_0}{32 \sqrt{\morse}}
-
\frac{2^8\sqrt\morse}{\b s_0\mathtt r_2}(1+|\ln z|)
-
\frac{ 2^{85}\morse^{8}}{s_* s_0^{13}\b^{17/2} \mathtt r_2}
\nonumber
\end{eqnarray}
by Cauchy estimates, \eqref{trocadero} and \eqref{otello}.
Then, using that for $|z|\leq 1/e^2,$ $z\in\C_*,$
$|\ln z|\leq 1/\sqrt{|z|},$ we get,
for $z$ as in \eqref{vienna}, 
that the following estimate holds:
\begin{eqnarray*}
&&|\partial_{EE} \act_n^{(2j-1)}(E_{2j-2}(\hat\act)-z \morse,\hat\act)|
\\
&&\geq
\frac{1}{|z|}\frac{s_0}{2^5 \morse^{3/2}}
-
\frac{2^9}{\b s_0\mathtt r_2\sqrt\morse\sqrt{|z|}}
-
\frac{ 2^{85}\morse^{7}}{s_* s_0^{13}\b^{17/2} \mathtt r_2}
\geq
\frac{1}{|z|}\frac{s_0}{2^6 \morse^{3/2}}\,.
\end{eqnarray*}
Since  $|\ln z|\leq 1/ |z|^{1/6}$ for  $|z|\leq 2^{36}$, $z\in\C_*,$
and recalling \eqref{maracaibo} and \eqref{trocadero}
 we get for $z$ as in \eqref{vienna}
  $$
 |\partial_{E} \act_n^{(2j-1)}(E_{2j-2}(\hat\act)-z \morse,\hat\act)|
 \leq
 \frac{ 2^{84}\morse^{8}}{s_* s_0^{13}\b^{17/2} }
 +
 \frac{2^7\sqrt\morse}{\b s_0|z|^{1/6}}
 \leq 
 \frac{2^8\sqrt\morse}{\b s_0|z|^{1/6}}\,.
 $$
\eqref{budapest} follows.
\eproof

\begin{lemma}\label{giga2bis}
The functions 
$\partial_E \act_n^{(2j-1),\pm}(E,\hat\act)$
 have holomorphic extension 
 to\footnote{$\mathtt r_2$ was defined in \eqref{tropicana}.}
 $$
E\in (\bar E^{(2j-1)}_-,
\bar E^{(2j-1)}_+ -3\mathtt r_2\morse/2)_{\mathtt r_2\morse}
\cap
\{ \Re E > \bar E^{(2j-1)}_- + \mathtt r_2 \morse/2^9\}
\qquad\hat\act\in \hat D_{r_0}
 $$
 with uniform estimate
\begin{equation}\label{trocaderobis}
|\partial_E \act_n^{(2j-1),\pm}(E,\hat\act)|
\leq
\frac{ 2^{79}\morse^{15/2}}{s_0^{13}\b^{8} \sqrt{\mathtt r_2}}
\,.
\end{equation}
\end{lemma}
\proof
We proceed in a way similar to Lemma \ref{giga2}.
First we define $\mathcal I_1$, respectively $\mathcal I_2$
as in \eqref{mela} but with
\begin{equation}\label{melabis}
\breve E=\l_{2j}(E)\,,\qquad
{\rm respectively}\qquad
\breve E=\l_{2j-1}(E)\,.
\end{equation}
Then we define $I_1$ and $I_3$ as in \eqref{prima}, while
\begin{equation}
I_7 := \int_{0}^{1-t_1} \mathcal I_2\, dt\,,
\label{primabis}
\end{equation}
so that 
$$
\partial_E \act_n^{(2j-1)}=I_1+I_3+I_7\,.
$$
The estimate of the term $I_1$ and $I_3$
are as in Lemma \ref{giga2}.
The estimates of the term in $I_7$ is similar 
to the one of $I_4$ in Lemma \ref{giga2}.
More precisely,
since $E=E_1+E_2$
with $E_1\in (\bar E_{2j-1}+ \mathtt r_2 \morse/2^5,
\bar E_{2j-2}-3\mathtt r\morse/2)$ 
 and $|E_2|<\mathtt r_2\morse,$ $\Re E_2\geq 0,$   we have, recalling \eqref{acrobat} and the
 definition of $\breve{\mathtt r}$ in \eqref{coppe},
\begin{eqnarray}
&&|\l_{2j-1}(E)-\bar \l_{2j-1}(E)|
\stackrel{\eqref{carciofino2}}\leq
\frac{48 \morse \checco}{\b^2}
\stackrel{\eqref{genesis}}
\leq
\frac{\mathtt r_2}{2^6}
\leq
\frac{\breve{\mathtt r}}{4}
\,,
\nonumber
\\
&&|\bar\l_{2j-1}(E)-\bar \l_{2j-1}(E_1+\Re E_2)|
=
\frac{2|\Im E_2|}{\bar E_{2j-2}-\bar E_{2j-1}}
\leq
|\bar\l_{2j-1}(E)-\bar \l_{2j-1}(E_1)|
\nonumber
\\
&&
\quad=
\frac{2|E_2|}{\bar E_{2j-2}-\bar E_{2j-1}}
<
\frac{2\mathtt r_2\morse}{\bar E_{2j-2}-\bar E_{2j-1}}
\stackrel{\eqref{ladispoli}}
\leq
\frac{2\mathtt r_2\morse}{\b}
\leq
\frac{\breve{\mathtt r}}{4}\,,
\nonumber
\\
&&
-1+\frac{\mathtt r_2}{2^9}\leq
-1+\frac{\mathtt r_2 \morse}{2^8(\bar E_{2j-2}-\bar E_{2j-1})}
=\bar \l_{2j-1}(\bar E_{2j-1}+ \mathtt r_2 \morse/2^9)
<
\bar \l_{2j-1}(E_1)
\nonumber
\\
&&
\quad=\Re\big(\bar \l_{2j-1}(E) \big)
\leq \bar\l_{2j-1}(E_1+\Re E_2)
<\bar \l_{2j-1}(\bar E_{2j-2}-3\mathtt r\morse/2+\mathtt r_2 \morse)
\nonumber
\\
&&\quad
=
1-\frac{\mathtt r_2\morse}{\bar E_{2j-2}-\bar E_{2j-1}}
\leq
1-\frac{\mathtt r_2}{2}\,,
\label{torroncino}
\end{eqnarray}
(recalling  \eqref{goffredo}) and
noting that
$$
\Re\big(\bar \l_{2j-1}(E) \big)=\bar\l_{2j-1}(E_1+\Re E_2)\,,
\qquad
\Im\big(\bar \l_{2j-1}(E) \big)=\frac{2\Im E_2}{\bar E_{2j-2}-\bar E_{2j-1}}\,.
$$
Recalling \eqref{torroncino}  we get
\begin{equation}\label{puledro}
\breve E=\l_{2j-1}(E)\,
\in\, \O_{\breve{\mathtt r}}\,, 
\qquad
-1+\frac{\mathtt r_2}{2^{10}}\leq
\Re \breve E
\leq
1-\frac{\mathtt r_2}{4}
\,.
\end{equation}
By \eqref{mela},\eqref{melabis} and recalling \eqref{prima} we get\footnote{Recalling that
$u=E-E_{2j-1}$ by \eqref{scarborough}.}
\begin{eqnarray}\label{roscioter}
 &&I_7
=
\frac{\sa_{2j-1}-\sa_{2j-2}}{E_{2j-1}-E_{2j-2}}
\frac{\sqrt{E-E_{2j-1}}}{2\pi }
\tilde{I}_7\,,\qquad 
{\rm with}
\nonumber
\\ 
&&\tilde{I}_7:=
\int_{0}^{1-t_1}
\frac{1+\breve{\tilde b}_{2j-1}\big(t(1+\breve E),
\breve\Sa_{2j-1}\big(
\breve E(t) \big)
 \big)}
{\sqrt t\,
\partial_{\breve \sa}
\breve\Gm_{2j-1}\Big(\breve  \Sa_{2j-1}\big(\breve E(t) \big)\Big)
}\, dt\,.
\end{eqnarray}
By \eqref{puledro}, 
we can apply estimate \eqref{cogli2} obtaining
$$
|\tilde{I}_7|
\leq
 \int_{0}^{1}
\frac{
2
}
{\sqrt{t}\Big|\partial_{\breve \sa}
\breve{\FO}_{2j-1}\Big(\breve  
{\bar\Sa}_{2j-1}\big(\breve E(t) \big)\Big)\Big|}
\, dt\,.
$$
Since by \eqref{puledro} we also have 
  $\Re \breve E<1$, by \eqref{dave}
  we get
\begin{equation*}
 |\tilde{I}_7|
\leq 
\frac{ 2^{32}\breve\morse}{\breve\b^3 \breve s^4}
\left( 
\frac{1}{\sqrt{|1+\breve E|}}
+
\ln \frac{1}{1-\Re \breve E}
\right)
\leq \frac{ 2^{68}\morse^{7}}{s_0^{13}\b^{7} }
\left( 
\frac{2^5}{\sqrt{\mathtt r_2}}
+
\ln \frac{4}{\mathtt r_2}
\right)
\leq \frac{ 2^{74}\morse^{7}}{s_0^{13}\b^{7} \sqrt{\mathtt r_2}}
\,,
 \end{equation*}
by \eqref{exit}-\eqref{ladispoli4} and \eqref{puledro}.
Then by \eqref{roscioter}, \eqref{puledro}, \eqref{ladispoli3}
we get
\begin{equation*}
|I_7|\leq 
\frac{ 2^{78}\morse^{15/2}}{s_0^{13}\b^{8} \sqrt{\mathtt r_2}}
\,.
\end{equation*}
Then \eqref{trocaderobis} follows recalling \eqref{strunz} and \eqref{soldatino}.
\eproof

For every $\mathtt r,\tilde\checco\geq 0$
set\footnote{Recall the definition of
$\bar E^{(i)}_\pm$
in \eqref{gruffalo}.}
\begin{eqnarray}
&&\mathcal E^{2j-1}(\mathtt r,\tilde\checco)
:=
(\bar E^{(2j-1)}_- ,
\bar E^{(2j-1)}_+)_{\mathtt r}
\cap \{\Re  E<\bar E^{(2j-1)}_+ - \tilde \checco\}
\,,
\nonumber
\\
&&\mathcal E^{2j-1}_*(\mathtt r,\tilde\checco)
:=
\mathcal E^{2j-1}(\mathtt r,\tilde\checco)
\cap \{\Re  E>\bar E^{(2j-1)}_- + \mathtt r_2\morse/2^9\}
\,,
\nonumber
\\
&&\mathcal E^{2j-1}_{**}(\mathtt r,\tilde\checco)
:=
\mathcal E^{2j-1}(\mathtt r,\tilde\checco)
\cap \{\Re  E>\bar E^{(2j-1)}_- + \mathtt r_2\morse/2^8\}
\,.
\label{autunno2}
\end{eqnarray}
Note that the above families of sets satisfy
$\mathcal E^{2j-1}(\mathtt r,\tilde\checco)\supset \mathcal E^{2j-1}_*(\mathtt r,\tilde\checco)
\supset \mathcal E^{2j-1}_{**}(\mathtt r,\tilde\checco)$ and
 are increasing with  
$\mathtt r$ and decreasing with
$\tilde \checco.$

Note that
\begin{equation}\label{primavera}
\mathtt r_2\leq \min\left\{ 
\frac{\b\breve{\mathtt r}}{8\morse},
\frac{\mathtt r_\star}{\morse},\mathtt r_3
\right\}\,,
\end{equation}
where
$\breve{\mathtt r},$ 
$\mathtt r_\star,$
$\mathtt r_2$
and $\mathtt r_3$ were defined in
\eqref{coppe},\eqref{epistassi},\eqref{tropicana} and \eqref{tropicana+},
respectively.
Set
\begin{equation}\label{acri}
\mathtt r_4:=\frac{\mathtt r_2\morse}{2^5}
=
\min\left\{
\frac{s_0^{49}\b^{30}}{2^{219} \morse^{29}},\
\frac{r_0^2}{2^{15}}\right\}\,.
\end{equation}
The  following lemma is justified in view of Remark \ref{figaro}.

\begin{lemma}\label{estate}
 i) The functions 
 $\partial_E \act_n^{(2j-1),\pm}(E,\hat\act)$
 and 
 $\partial_E \bar \act_n^{(2j-1),\pm}(E)$
 are holomorphic on
$\mathcal E^{2j-1}_*(2\mathtt r_4,2\checco)
 \times  \hat D_{r_0}$ 
 and 
 $\mathcal E^{2j-1}_*(2\mathtt r_4,2\checco),
$
respectively.
\\
ii) Moreover, for 
 $
 2\checco
 \leq
  \tilde\checco \leq \mathtt r_4$
 we have
 \begin{eqnarray}\label{damasco}
&&\sup_{\mathcal E^{2j-1}_*(2\mathtt r_4,2\tilde\checco)
\times\hat D_{r_0}}
\Big|
\partial_E \act_n^{(2j-1),\pm}(E,\hat\act)
\Big|\,,\ \ 
    \sup_{\mathcal E^{2j-1}_*(2\mathtt r_4,2\tilde\checco)}
\Big|\partial_E \bar \act_n^{(2j-1),\pm}(E)
\Big|
\nonumber
\\
&&
\leq
\frac{ 2^{80}\morse^{15/2}}{s_0^{13}\b^{8} \sqrt{\mathtt r_2}}
+ 
\frac{2^8\sqrt\morse}{\b s_0}
\ln \frac{\morse}{\tilde\checco}
\end{eqnarray} 
and
 \begin{eqnarray}\label{damasco2}
&&\sup_{\mathcal E^{2j-1}_{**}(3\mathtt r_4/2,2\tilde\checco)
\times\hat D_{r_0}}
\Big|
\partial_{EE} \act_n^{(2j-1),\pm}(E,\hat\act)
\Big|\,,\ \ 
    \sup_{\mathcal E^{2j-1}_{**}(3\mathtt r_4/2,2\tilde\checco)}
\Big|\partial_{EE} \bar \act_n^{(2j-1),\pm}(E)
\Big|
\nonumber
\\
&&
\leq
\frac{ 2^{90}\morse^{13/2}}{s_0^{13}\b^{8} \mathtt r_2^{3/2}}
+\frac{2^9\sqrt\morse}{\b s_0\tilde\checco}\,.
\end{eqnarray}
\end{lemma}
\proof
Here we prove only the minus case,
namely the case of
$\partial_E \act_n^{(2j-1),-}$ and
 $\partial_E \bar \act_n^{(2j-1),-};$
 the case of 
 $\partial_E \act_n^{(2j-1),+}$ and
 $\partial_E \bar \act_n^{(2j-1),+}$
 is completely analogous.

i) 
By Lemma \ref{giga2bis}
we have that $\partial_E \act_n^{(2j-1),-}
(E,\hat\act)$
is holomorphic on\footnote{Recall \eqref{gruffalo} and Remark \ref{stiffe}.}
\begin{equation}\label{medio4}
E\in (\bar E_{2j-1},
\bar E_{2j-2} -3\mathtt r_2\morse/2)_{\mathtt r_2\morse}
\cap
\{ \Re E > \bar E_{2j-1} + \mathtt r_2 \morse/2^9\}
\qquad\hat\act\in \hat D_{r_0}
\end{equation}
with uniform estimate
\begin{equation}\label{anulare}
|\partial_E \act_n^{(2j-1),\pm}(E,\hat\act)|
\leq
\frac{ 2^{79}\morse^{15/2}}{s_0^{13}\b^{8} \sqrt{\mathtt r_2}}
\end{equation}
and the same also holds 
for\footnote{\label{aqualung} Indeed
this is a particular case of Lemma \ref{giga2bis} with $\checco=0$ in \eqref{ciccio}.}
$\partial_E \bar\act_n^{(2j-1),-}(E).$
Moreover by Lemma \ref{giga2}
$\partial_E \act_n^{(2j-1),-}(E,\hat\act)$
is holomorphic for $\hat\act \in \hat D_{r_0}$
and $E$ belonging to the set   
\begin{equation}\label{medio3}
B_{\mathtt r_2 \morse}\big(E_{2j-2}(\hat\act)\big)
\,\cap\, (E_{2j-2}(\hat\act)-\C_*)
\end{equation}
and, as above (recall footnote \ref{aqualung}), the same also holds for
$\partial_E \bar\act_n^{(2j-1),-}(E)$
on the set
$$
B_{\mathtt r_2 \morse}\big(\bar E_{2j-2}\big)
\,\cap\, (\bar E_{2j-2}-\C_*)\,.
$$
By \eqref{october} we have that
\begin{eqnarray}
&&B_{\mathtt r_2 \morse}\big(E_{2j-2}(\hat\act)\big)
\,\cap\, (E_{2j-2}(\hat\act)-\C_*)
\
\supset
\
B_{\mathtt r_2 \morse-2\checco}
\big(\bar E_{2j-2}\big)\,\cap\,
\{\Re E<\bar E_{2j-2}-2\checco\}
\nonumber
\\
&&\supset\ 
B_{7\mathtt r_2 \morse/8}
\big(\bar E_{2j-2}\big)\,\cap\,
\{\Re E<\bar E_{2j-2}-2\checco\}
\,,
\label{medio}
\end{eqnarray}
since $\checco\leq \mathtt r_2 \morse/16$
by \eqref{genesis}.
Moreover we have that
\begin{equation}\label{medio2}
(\bar E_{2j-1},
\bar E_{2j-2}-\frac32\mathtt r _2\morse)_{\mathtt r_2\morse}
\ \cup\ 
B_{7\mathtt r_2 \morse/8}
\big(\bar E_{2j-2}\big)
\ \supset \ 
(\bar E_{2j-1}-\mathtt r_2\morse/2,
\bar E_{2j-2})_{\mathtt r_2\morse/4}\,.
\end{equation}
Then by \eqref{medio4}-\eqref{medio2}
we get
that
$\partial_E \bar\act_n^{(2j-1),-}(E)$
and 
$\partial_E \act_n^{(2j-1),-}(E,\hat\act)$
are holomorphic for
$$
E\in 
\mathcal E^{2j-1}(2\mathtt r_4,2\checco)
\cap
\{ \Re E > \bar E_{2j-1} + \mathtt r_2 \morse/2^9\}
=
\mathcal E^{2j-1}_*(2\mathtt r_4,2\checco)
$$
and $\hat\act\in\hat D_{r_0}.$
This proves the part {\it i)} of the lemma.

\medskip

{\it ii)} Assume now $
 2\checco
 \leq
  \tilde\checco \leq \mathtt r_4$.
By \eqref{maracaibo} and \eqref{trocadero}
we have that for
$E\in B_{\mathtt r_2 \morse}\big(E_{2j-2}(\hat\act)\big)
\,\cap\, (E_{2j-2}(\hat\act)-\C_*),$
$\hat\act\in \hat D_{r_0}$
it results
\begin{equation}\label{pollice}
|\partial_E \act_n^{(2j-1),-}
(E,\hat\act)|
\leq
\frac{ 2^{84}\morse^{8}}{s_* s_0^{13}\b^{17/2} }
+ 
\frac{2^7\sqrt\morse}{\b s_0}
\left|\ln
\frac{E_{2j-2}(\hat\act)-E}{\morse}
\right|
\,.
\end{equation}
Take
$$
E\in B_{7\mathtt r_2 \morse/8}
\big(\bar E_{2j-2}\big)\,\cap\,
\{\Re E<\bar E_{2j-2}-2\tilde\checco\}\,.
$$
Recalling \eqref{medio} and \eqref{pollice}
we get
\begin{equation}\label{pollice2}
|\partial_E \act_n^{(2j-1),-}
(E,\hat\act)|
\leq
\frac{ 2^{84}\morse^{8}}{s_* s_0^{13}\b^{17/2} }
+ 
\frac{2^8\sqrt\morse}{\b s_0}
\ln \frac{\morse}{\tilde\checco}
\,,
\end{equation}
since 
\begin{equation}\label{jaffa}
|E_{2j-2}(\hat\act)-E|\geq |\bar E_{2j-2}-E|
-|E_{2j-2}(\hat\act)-\bar E_{2j-2}|
\geq 2\tilde\checco- 2\checco
\geq \tilde\checco
\end{equation}
and $\morse/\tilde\checco> \morse/\mathtt r_4\geq 2^{10}.$
By \eqref{medio2}, \eqref{anulare} and \eqref{pollice2}
we get, recalling the definition of $\mathtt r_2$ in \eqref{tropicana}
(and  \eqref{badedas}),
\begin{eqnarray*}
\sup_{\mathcal E^{2j-1}_*(2\mathtt r_4,2\tilde\checco)
\times\hat D_{r_0}}
\Big|
\partial_E \act_n^{(2j-1),-}(E,\hat\act)
\Big|
&\leq&
\frac{ 2^{79}\morse^{15/2}}{s_0^{13}\b^{8} \sqrt{\mathtt r_2}}
+
\frac{ 2^{84}\morse^{8}}{s_* s_0^{13}\b^{17/2} }
+ 
\frac{2^8\sqrt\morse}{\b s_0}
\ln \frac{\morse}{\tilde\checco}
\\
&\leq&
\frac{ 2^{80}\morse^{15/2}}{s_0^{13}\b^{8} \sqrt{\mathtt r_2}}
+ 
\frac{2^8\sqrt\morse}{\b s_0}
\ln \frac{\morse}{\tilde\checco}
\end{eqnarray*}
proving \eqref{damasco}
(the estimate on $\partial_E \bar\act_n^{(2j-1),-}(E)$ is analogous).

Let us now prove \eqref{damasco2}.
We use  \eqref{damasco} with 
$\tilde\checco=\mathtt r_4,$ namely (recall \eqref{acri}, \eqref{tropicana}, \eqref{harlock})
\begin{eqnarray*}
\sup_{\mathcal E^{2j-1}_*(2\mathtt r_4,2\mathtt r_4)
\times\hat D_{r_0}}
\Big|
\partial_E \act_n^{(2j-1),-}(E,\hat\act)
\Big|
&\leq&
\frac{ 2^{80}\morse^{15/2}}{s_0^{13}\b^{8} \sqrt{\mathtt r_2}}
+ 
\frac{2^8\sqrt\morse}{\b s_0}
\ln \frac{2^5}{\mathtt r_2}\nonumber
\\
&\leq&
\frac{ 2^{81}\morse^{15/2}}{s_0^{13}\b^{8} \sqrt{\mathtt r_2}}\,.
\end{eqnarray*}
Then by Cauchy estimates
\begin{equation}\label{aleppo}
\sup_{\mathcal E^{2j-1}_{**}(3\mathtt r_4/2,\mathtt r_4)
\times\hat D_{r_0}}
\Big|
\partial_{EE} \act_n^{(2j-1),-}(E,\hat\act)
\Big|
\leq
\frac{ 2^{90}\morse^{13/2}}{s_0^{13}\b^{8} \mathtt r_2^{3/2}}\,.
\end{equation}
By \eqref{maracaibo}
we have, for $|z|<\mathtt r_2,$
\begin{equation*}
-\morse
\partial_{EE} \act_n^{(2j-1)}(E_{2j-2}(\hat\act)-z \morse,\hat\act)
=
\partial_z\f(z,\hat\act)+\partial_z\psi(z,\hat\act)\ln z
+\frac{1}{z}\psi(z,\hat\act)\,,
\end{equation*}
with, by \eqref{trocadero} and Cauchy estimates,
\begin{equation*}
\sup_{\{|z|<\mathtt r_2/2\}\times \hat D_{r_0}} |\partial_z\f(z,\hat\act)|
\leq
\frac{ 2^{85}\morse^{8}}{s_* s_0^{13}\b^{17/2} \mathtt r_2}
\,,\qquad 
\sup_{\{|z|<\mathtt r_2/2\}\times \hat D_{r_0}} |\partial_z\psi(z,\hat\act)|
\leq 
\frac{2^8\sqrt\morse}{\b s_0\mathtt r_2}\,;
\end{equation*}
then, using again \eqref{trocadero} and \eqref{jaffa}, we get, for 
$$
E\in B_{\mathtt r_2 \morse/2}\big(E_{2j-2}(\hat\act)\big)
\,\cap\, (E_{2j-2}(\hat\act)-\C_*)\,,\qquad
\hat\act\in \hat D_{r_0}\,,$$
that\footnote{\label{tripoli}
 Noting that, for $z=\frac{E_{2j-2}(\hat\act)-E}{\morse}=r e^{\ii \theta},$  $r>0,$
$-\pi<\theta<\pi,$  with $|z|=r\leq \mathtt r_2/2<1/16,$
we have $$
|\ln z|\leq |\ln r|+\pi\leq 2\ln \frac1r
\leq 2\ln\frac{1}{\mathtt r_2}+\frac{\mathtt r_2}{r}\,.
$$
}
\begin{eqnarray*}
|\partial_{EE} \act_n^{(2j-1)}(E,\hat\act)|
&\leq&
\frac{ 2^{85}\morse^{7}}{s_* s_0^{13}\b^{17/2} \mathtt r_2}
+
\frac{2^8}{s_0\b \sqrt\morse \mathtt r_2}
\left|\ln
\frac{E_{2j-2}(\hat\act)-E}{\morse}
\right|
+\frac{2^7\sqrt\morse}{\b s_0|E_{2j-2}(\hat\act)-E|}
\\
&\leq&
\frac{ 2^{85}\morse^{7}}{s_* s_0^{13}\b^{17/2} \mathtt r_2}
+
\frac{2^9}{s_0\b \sqrt\morse \mathtt r_2}
\ln
\frac{1}{\mathtt r_2}
+\frac{2^9\sqrt\morse}{\b s_0|E_{2j-2}(\hat\act)-E|}
\,.
\end{eqnarray*}
By \eqref{jaffa}
we have that,
for
$$
E\in B_{\mathtt r_2 \morse/2}\big(E_{2j-2}(\hat\act)\big)
\,\cap\, \{ \Re E\leq \bar E_{2j-2}-2\tilde\checco\}\,,\qquad
\hat\act\in \hat D_{r_0}\,,
$$ 
\begin{equation}\label{montsegur}
|\partial_{EE} \act_n^{(2j-1)}(E,\hat\act)|
\leq
\frac{ 2^{85}\morse^{7}}{s_* s_0^{13}\b^{17/2} \mathtt r_2}
+
\frac{2^9}{s_0\b \sqrt\morse \mathtt r_2}
\ln
\frac{1}{\mathtt r_2}
+\frac{2^9\sqrt\morse}{\b s_0\tilde\checco}\,.
\end{equation}
Recalling \eqref{acri}, by \eqref{aleppo} and \eqref{montsegur} 
(and \eqref{harlock} and  \eqref{badedas})
the estimate \eqref{damasco2} follows
(the estimate on $\partial_E \bar\act_n^{(2j-1),-}(E)$ is analogous).
\eproof

By \eqref{noce} and lemmata \ref{giga},\ref{giga2},\ref{estate} 
we get the following
\begin{corollary}
 The function $\act_n^{(2j-1)}(E,\hat\act)$
 is holomorphic for $\hat\act \in \hat D_{r_0}$ and 
$$
E\in (\bar E^{(2j-1)}_-,\bar E^{(2j-1)}_+)_{2\mathtt r_4}
\cap \{ E^{(2j-1)}_+(\hat\act)-\C_*\}\,.
$$ 
\end{corollary}

\subsection{Estimates on  $\partial_E \act_n^{(2j)}$}

Deriving \eqref{ancomarzio2bis} w.r.t. $E$ we get

\begin{eqnarray}
 \partial_E\act_n^{(2j),-}(E)
 =
 \frac{1}{2\pi}
\int_{ \sa_{2j-1}}^{\sa_{2j}}
\frac{
1
+ \tilde b\big(E-\Gm(\sa),\sa \big)
}
{\sqrt{E-\Gm(\sa)}}\, d\sa\,,
\nonumber
\\
\partial_E \act_n^{(2j),+}(E)
 =
 \frac{1}{2\pi}
\int_{\sa_{2j}}^{\sa_{2j+1}}
\frac{
1
+ \tilde b\big(E-\Gm(\sa),\sa \big)
}
{\sqrt{E-\Gm(\sa)}}
\, d\sa\,,
\label{ancomarzio3}
\end{eqnarray}
with $\tilde b$ defined in \eqref{palidoro}.
Note that by \eqref{3holesquater}
the function
$ \tilde b\big(E-\Gm(\sa),\sa \big)$
is well defined for
\begin{equation}\label{elba}
E\in (E_{2j}(\hat\act),R_0^2-2\morse)_{r_0^2/2^7}\,,
\qquad
\hat\act\in \hat D_{r_0}\,,\qquad
\sa\in \T^1_{s_0/2}\,.
\end{equation}

\begin{lemma}\label{giga2+}
Set
\begin{equation}\label{tropicana+}
\mathtt r_3:=
\min\left\{
\frac{s_0^{6}\b^{3}}{2^{25} \morse^{3 }},\
\frac{r_0^2}{2^{8}\morse}\right\}
\,.
\end{equation}
The function 
$\partial_E \act_n^{(2j),-}(E_{2j}(\hat\act)+z\morse,\hat\act)+
\partial_E \act_n^{(2j),+}(E_{2j}(\hat\act)+z\morse,\hat\act),$
initially defined for $0<z<\mathtt r_3$ and
 $\hat\act\in \hat D,$
 has holomorphic extension to the complex set
 $\{z\in\C_*\ \ \text{s.t.}\ \  |z|<\mathtt r_3\}\times \hat D_{r_0}.$
In particular
\begin{equation}\label{maracaibo+}
\partial_E \act_n^{(2j),-}(E_{2j}(\hat\act)+z\morse,\hat\act)+
\partial_E \act_n^{(2j),+}(E_{2j}(\hat\act)+z\morse,\hat\act)
=
\f(z,\hat\act)+\psi(z,\hat\act)\ln z\,,
\end{equation}
where $\f( z,\hat\act)$ and $\psi(z,\hat\act)$
are holomorphic function in the set $\{|z|<\mathtt r_3\}\times \hat D_{r_0}$
with
\begin{equation}\label{trocadero+}
\sup_{\{|z|<\mathtt r_3\}\times \hat D_{r_0}} |\f(z,\hat\act)|
\leq \frac{2^{20} \morse}{\b^{3/2} s_0^2 s_*}
\,,\qquad 
\sup_{\{|z|<\mathtt r_3\}\times \hat D_{r_0}} |\psi(z,\hat\act)|
\leq  \frac{16}{\sqrt\b}\,,
\end{equation}
and
\begin{equation}\label{otello+}
\inf_{\hat\act\in \hat D_{r_0}}
|\psi(0,\hat\act)|\geq \frac{s_0}{4 \sqrt{\morse}}\,.
\end{equation}
Finally
\begin{equation}\label{trocadero+*}
\sup_{\{|z|<\mathtt r_3\}\times \hat D_{r_0/2}} 
|\partial_{\hat\act}\f(z,\hat\act)|
\leq M_\f\checco
\,,\qquad 
\sup_{\{|z|<\mathtt r_3\}\times \hat D_{r_0/2}} 
|\partial_{\hat\act}\psi(z,\hat\act)|
\leq  M_\psi\checco\,.
\end{equation}
$M_\f,M_\psi$ definite in \eqref{trocadero*}
\end{lemma}
\proof
We first note that, recalling \eqref{elba},
$$
\tilde b\big(E_{2j}(\hat\act)+z\morse-\Gm(\sa),\sa \big)
$$
is well defined for 
$$
|z|\leq \frac{r_0^2}{2^{8}\morse}\,,
\qquad
\hat\act\in \hat D_{r_0}\,,\qquad
\sa\in \T^1_{s_0/2}\,.
$$
Setting
\begin{equation}\label{vaccino}
\sa_\bullet:=\frac{\b s_0^3}{2^{10} \morse}
\stackrel{\eqref{harlock}}\leq s_0/4\,,
\end{equation}
we split
\begin{eqnarray}
&&2\pi\Big(\partial_E \act_n^{(2j),-}(E)+
\partial_E \act_n^{(2j),+}(E)\Big)
=I_1+I_2+I_3
\nonumber
\\
&&:=
\int_{ \sa_{2j-1}}^{\sa_{2j}-\sa_\bullet}
 +
\int_{ \sa_{2j}-\sa_\bullet}^{\sa_{2j}+\sa_\bullet}
 +
\int_{\sa_{2j}+\sa_\bullet}^{\sa_{2j+1}}
\ \frac{
1
+ \tilde b\big(E-\Gm(\sa),\sa \big)
}
{\sqrt{E-\Gm(\sa)}}
\, d\sa\,.
\label{bisanzio}
\end{eqnarray}
We first consider the more relevant integral which
is 
$\int_{ \sa_{2j}-\sa_\bullet}^{\sa_{2j}+\sa_\bullet}.$ Since the interval of integration 
is symmetric w.r.t. $\sa_{2j}$ we can consider the ``even part'' of the integrand, namely,
changing variable $\sa=\sa_{2j}+\vartheta$
\begin{eqnarray*}
&&
I_2=\int_{ \sa_{2j}-\sa_\bullet}^{\sa_{2j}+\sa_\bullet}
\frac{
1
+ \tilde b\big(E-\Gm(\sa),\sa \big)
}
{\sqrt{E-\Gm(\sa)}}
\, d\sa
\\
&&=
\int_{0}^{\sa_\bullet}
\left(
\frac{
1
+ \tilde b\big(E-\Gm(\sa_{2j}+\vartheta),
\sa_{2j}+\vartheta \big)
}
{\sqrt{E-\Gm(\sa_{2j}+\vartheta)}}
+
\frac{
1
+ \tilde b\big(E-\Gm(\sa_{2j}-\vartheta),
\sa_{2j}-\vartheta \big)
}
{\sqrt{E-\Gm(\sa_{2j}-\vartheta)}}
\right)
\, d\vartheta\,.
\end{eqnarray*}

\noindent
Since $\Gm$ has a maximum at $\sa_{2j}$
we have that
$$
\Gm(\sa_{2j}+\vartheta)
=E_{2j}- \b_0 \vartheta^2 - \vartheta^3\Gm_*(\vartheta)\,,
$$
where $\b_0=-\partial_{\sa\sa}\Gm(\sa_{2j}(\hat	\act),\hat\act)/2$
with 
\begin{equation}\label{confortablynumb}
\b/4\leq |\b_0|\leq \morse/s_0^2
\end{equation}
 by \eqref{ladispoli3}, \eqref{ciccio3} 
and
\begin{equation}\label{yoga}
\sup_{|\vartheta|\leq 2\sa_\bullet}|\Gm_*(\vartheta)|
\leq 8 \morse /s_0^3\,.
\end{equation}

We set
$$
\z=w^2\morse:=z\morse\,,\qquad
E=E_{2j}(\hat\act)+\z=
E_{2j}(\hat\act)+w^2\morse=
E_{2j}(\hat\act)+z\morse\,.
$$
{\sl For the moment being, we consider only real $\hat\act\in\hat D,$ so that 
$\Gm$ is real on real. In particular we think
$\z>0,w>0.$ We also have $\b_0>0.$}

Then we have
$$
I_2
=
\int_{0}^{\sa_\bullet}
\left(
\frac{
1
+ \tilde b\big(\z
+\b_0 \vartheta^2 + \vartheta^3\Gm_*(\vartheta),
\sa_{2j}+\vartheta \big)
}
{\sqrt{\z
+\b_0 \vartheta^2 + \vartheta^3\Gm_*(\vartheta)}}
+
\frac{
1
+ \tilde b\big(\z
+\b_0 \vartheta^2 - \vartheta^3\Gm_*(-\vartheta),
\sa_{2j}-\vartheta \big)
}
{\sqrt{\z
+\b_0 \vartheta^2 -\vartheta^3\Gm_*(-\vartheta)}}
\right)
\, d\vartheta\,.
$$
Then we split
$$
I_2=I_4+I_5:=\int_0^{4w \sqrt{\morse /\b_0}} +\int_{4 w \sqrt{\morse /\b_0}}^{\sa_\bullet}\,.
$$
Changing  variable $\vartheta=w y$ we get
\begin{eqnarray*}
&&I_4
=
\int_{0}^{4\sqrt{\morse /\b_0}}
\Bigg(
\frac{
1
+ \tilde b\big(w^2\morse
+\b_0 w^2 y^2 + w^3 y^3\Gm_*(wy),
\sa_{2j}+wy \big)
}
{\sqrt{\morse
+\b_0 y^2 + w y^3\Gm_*(w y)}}
\\
&&\qquad\qquad\quad
+
\frac{
1
+ \tilde b\big(w^2\morse
+\b_0 w^2 y^2 - w^3 y^3\Gm_*(-wy),
\sa_{2j}-wy \big)
}
{\sqrt{\morse
+\b_0 y^2 - w y^3\Gm_*(-w y)}}
\Bigg)
\, d y\,.
\end{eqnarray*}
We note that for $|w|\leq \sqrt{\mathtt r_3}$
and $0\leq y\leq 4\sqrt{\morse /\b_0}$ we have
$$
|w y^3\Gm_*(wy)|
\stackrel{\eqref{yoga}}\leq 
2^9\sqrt{\mathtt r_3}\morse^{5/2}/\b_0^{3/2} s_0^3\leq \morse/8
$$
by \eqref{tropicana+}.
Regarding the term $\b_0 y^2$ it is positive
when $\hat\act\in\hat D;$
if $\hat\act\in\hat D_{r_0}$ then 
$\hat\act=\hat\act_1+\hat\act_2$
with $\hat\act_1\in\hat D$ and $|\hat\act_2|<r_0;$ then 
$\b_0=\b_1+\b_2$ with
$$
\b_1:=-\partial_{\sa\sa}\Gm(\sa_{2j}(\hat\act_1),\hat\act_1)/2\geq \b/4>0
$$ and 
$\b_2:=-\frac12(\partial_{\sa\sa}\Gm(\sa_{2j}(\hat\act_2),\hat\act_2)-\partial_{\sa\sa}\Gm(\sa_{2j}(\hat\act_1),\hat\act_1))$, with
(recall \eqref{ciccio},\eqref{octoberx},\eqref{goffredo})
$$
|\b_2|
\leq |\partial_{\sa\sa}\FO(\sa_{2j}(\hat\act_2))-
\partial_{\sa\sa}\FO(\sa_{2j}(\hat\act_1))|
+8\checco/s_0^2
\leq 16\checco\morse/\b s_0^3
+8\checco/s_0^2
\stackrel{\eqref{genesis}}\leq
\frac{\b}{2^{16}}
\stackrel{\eqref{confortablynumb}}	\leq
\frac{|\b_0|}{2^{12}}
\,. 
$$
Recollecting we have, also in the complex case,
$$
\Re(\morse
+\b_0 y^2 \pm w y^3\Gm_*(\pm w y))\geq \morse/2\,.
$$
Then the modulus of the integrand function in $I_4$
is, for every $0\leq y\leq 4\sqrt{\morse /\b_0},$ 
bounded 
(recall also \eqref{3holesquater})
by $8/\sqrt\morse$.
Then $I_4$ defines a even\footnote{Since the integrand is even w.r.t. $w.$} holomorphic function of $w$ in $|w|\leq \sqrt{\mathtt r_3}$
 with 
 \begin{equation}\label{praga}
 |I_4|\leq 16/\sqrt\morse
\end{equation}
 uniformly; equivalently $I_4$ is a 
 holomorphic function of $z$ in $|z|\leq \mathtt r_3$
 with the same bound.
 
 Let us consider now the term $I_5.$
 We rewrite it as
 $$
 I_5=\int_{4 w \sqrt{\morse /\b_0}}^{\sa_\bullet}
 \frac{1}{\vartheta}G(w^2/\vartheta^2,\vartheta)d\vartheta\,,
 $$
 where
$$
G(\xi,\vartheta):=
\frac{
1
+ \tilde b\big(\morse \xi \vartheta^2
+\b_0 \vartheta^2 + \vartheta^3\Gm_*(\vartheta),
\sa_{2j}+\vartheta \big)
}
{\sqrt{\b_0 + \morse\xi
 + \vartheta\Gm_*(\vartheta)}}
+
\frac{
1
+ \tilde b\big(\morse \xi \vartheta^2
+\b_0 \vartheta^2 - \vartheta^3\Gm_*(-\vartheta),
\sa_{2j}-\vartheta \big)
}
{\sqrt{\b_0 + \morse\xi
  -\vartheta\Gm_*(-\vartheta)}}\,.
$$
By \eqref{yoga} and \eqref{vaccino} we obtain that
$$
\sup_{|\vartheta|<2\sa_\bullet}|\vartheta\Gm_*(\vartheta)|
\leq 16\morse \sa_\bullet/s_0^3\leq |\b_0|/8\,.
$$ 
Then $G$ is holomorphic and bounded 
by
$$
\frac{8}{\sqrt{|\b_0|}}
\stackrel{\eqref{confortablynumb}}\leq\frac{16}{\sqrt\b}
$$
on the set
$$
\{|\xi|\leq|\b_0|/2\morse\}\times\{|\vartheta|\leq 2\sa_\bullet\}\,.
$$
and even w.r.t. $\vartheta.$ In particular 
$$
G(\xi,\vartheta)=\sum_{h,k\geq 0}
G_{hk}\xi^h \vartheta^{2k}
$$
for suitable coefficients $G_{hk}$ satisfying
\begin{equation}\label{relax}
|G_{hk}|\leq \frac{16}{\sqrt\b} (2\morse/|\b_0|)^h
(1/4 \sa_\bullet^2)^k\,.
\end{equation}
Then, using that the series totally converges
on the above set,
we get
\begin{eqnarray}
I_5
&=&
\sum_{h,k\geq 0}  G_{hk} w^{2h}
\int_{4 w \sqrt{\morse /\b_0}}^{\sa_\bullet}
 \vartheta^{2(k-h)-1}\, d\vartheta
 \nonumber
 \\
 &=&
 \sum_{h\neq k}  G_{hk} w^{2h}\left(
 \frac{{\sa_\bullet}^{2(k-h)}}{2(k-h)}
 -
  \frac{(4w)^{2(k-h)}(\morse /\b_0)^{k-h}}{2(k-h)}
 \right)+
 \sum_{h\geq 0}
 G_{hh} w^{2h}(\ln \frac{\sa_\bullet\sqrt\b_0}{4\sqrt\morse}-\ln w)
 \nonumber
 \\
 &=&
 \psi(z)\ln z+\f_1(z)\,,
 \label{praga2}
 \end{eqnarray}
where
\begin{equation}\label{kukident}
\psi(z):=- \sum_{h\geq 0}
\frac12  G_{hh} z^h 
\end{equation}
and
\begin{equation}\label{crema}
\f_1(z):=-2\ln \frac{\sa_\bullet\sqrt\b_0}{4\sqrt\morse}\psi(z)+
\sum_{h\geq 0}\f_{1,h}z^h\,,
\qquad
\f_{1,h}:=
\sum_{k\geq 0,k\neq h}
 \frac{G_{hk}{\sa_\bullet}^{2(k-h)}+G_{kh}(16\morse /\b_0)^{h-k}}{2(k-h)}
\end{equation}
{\sl Now we note that the above representation formula for $I_5$
also holds for complex value of $\hat\act$ and that the functions
$\psi(z)$ and $\f_1(z)$ are well defined (since their series totally converge) and 
holomorphic for $|z|< \mathtt r_3$.}
\\
Indeed by \eqref{relax},\eqref{confortablynumb} and \eqref{tropicana+} we get
$$
\sup_{|z|\leq \mathtt r_3}|\psi(z)|\leq \frac{16}{\sqrt\b}
$$
(namely the second estimate in \eqref{trocadero+} holds)
and
$$
|\f_{1,h}|\leq \frac{32}{\sqrt\b} \left(\frac{4\morse}{|\b_0|\sa_\bullet^2} \right)^h
$$
Then, recalling \eqref{crema} and \eqref{confortablynumb},
we also get\footnote{In the last inequality we use that
$\ln x<x$ and \eqref{harlock}.}
\begin{equation}\label{ravel}
\sup_{|z|\leq \mathtt r_3}|\f_1(z)|\leq
\frac{2^9}{\sqrt\b}\left( 1+\ln\frac{\morse}{\b \sa_\bullet^2}\right)
=
\frac{2^9}{\sqrt\b}\left( 1+\ln\frac{2^{20} \morse^3}{\b^3 s_0^6}\right)
\leq \frac{2^{16}\morse}{\b^{3/2}s_0^2}
\end{equation}
Note that by \eqref{kukident}
$$
\psi(0)=-\frac12 G_{00}
=-\frac12 G(0,0)=-\frac{1+\tilde b(0,\sa_{2j})}
{\sqrt{\b_0}}\,,
$$
then by \eqref{confortablynumb} and
\eqref{3holesquater} we get \eqref{otello+}.
\\
We finally consider the terms $I_1$ 
and $I_3$ (recall \eqref{bisanzio}), which are analogous.
First we note that 
\begin{equation}\label{fritter}
\bar E_{2j}-\FO(\bar\sa)\geq \frac{\b\sa_\bullet^2}{4}\,,\qquad
\forall\, \bar\sa_{2j-1}\leq \bar\sa\leq 
\bar\sa_{2j}-\sa_\bullet\,,
\ \ \forall\, 
\bar\sa_{2j}+\sa_\bullet
\leq \bar\sa\leq \bar\sa_{2j+1}\,.
\end{equation}
Indeed, considering the case
$\bar\sa_{2j}+\sa_\bullet
\leq \bar\sa\leq \bar\sa_{2j+1}$
(the other case being analogous),
we have, since in such interval $\FO$
is decreasing, 
$$
\bar E_{2j}-\FO(\bar\sa)\geq 
\bar E_{2j}-\FO(\bar\sa_{2j}+\sa_\bullet)
\geq
-\frac12\partial_{\sa\sa}\FO(\bar\sa_{2j})\sa_\bullet^2- \frac{\morse}{s_0^3}\sa_\bullet^3
\geq 
\frac{\b}{2}\sa_\bullet^2
-\frac{\morse}{s_0^3}\sa_\bullet^3
\geq
\frac{\b\sa_\bullet^2}{4}
$$
by \eqref{goffredo}, Cauchy estimates,
\eqref{ladispoli} and \eqref{vaccino}.
We now consider a point 
$\sa\in(\sa_{2j}+\sa_\bullet,\sa_{2j+1})$
(recall the definition in \ref{segment}).
This means that there exists $0\leq t\leq 1$
such that 
$$
\sa=t(\sa_{2j}+\sa_\bullet)+(1-t)\sa_{2j+1}.
$$
Then, set
$$
\bar\sa=t(\bar\sa_{2j}+\sa_\bullet)+(1-t)\bar\sa_{2j+1},
$$
with
$\bar\sa_{2j}+\sa_\bullet
\leq \bar\sa\leq \bar\sa_{2j+1}$.
By \eqref{octoberx} we get
\begin{equation}\label{europe}
|\sa-\bar\sa|\leq 4\checco/\b s_0\,.
\end{equation}
For
$$
E=E_{2j}(\hat\act)+z\morse\,,\qquad
|z|<\mathtt r_3\,,
$$
we have
\begin{eqnarray*}
\big|
(E-\Gm(\sa))-(\bar E_{2j}-\FO(\bar\sa))
\big|
&\stackrel{\eqref{ciccio}}\leq&
 |E_{2j}-\bar E_{2j}|+\mathtt r_3\morse
+\checco+|\FO(\sa)-\FO(\bar\sa)|
\\
&\leq&
3\checco+\mathtt r_3\morse
+\frac{4\checco\morse}{\b s_0^2}
\leq \frac{\b\sa_\bullet^2}{8}
\end{eqnarray*}
by \eqref{october},\eqref{goffredo},\eqref{europe},
\eqref{tropicana+} and \eqref{genesis}.
Then by \eqref{fritter} we get
$$
|E-\Gm(\sa)|\geq \frac{\b\sa_\bullet^2}{8}\,,
\qquad
{\rm for}\qquad
\sa\in(\sa_{2j}+\sa_\bullet,\sa_{2j+1})
$$
(analogously for $\sa\in(\sa_{2j-1},\sa_{2j}-\sa_\bullet)$).
Hence, recalling \eqref{bisanzio} (and \eqref{3holesquater})
$$
|I_1|\,,\ |I_3|\,\leq\, 
\frac{2^5}{\sqrt{\b}\,\sa_\bullet}
=
\frac{2^{15} \morse}{\b^{3/2} s_0^3}
\,.
$$
Recalling \eqref{praga},\eqref{praga2},
\eqref{ravel} (and \eqref{harlock}) we get \eqref{trocadero+}.

We omit the proof of \eqref{trocadero+*}.
\eproof

For every $\mathtt r,\tilde\checco\geq 0$
set\footnote{Recall the definition of
$\bar E^{(i)}_\pm$
in \eqref{gruffalo}.}
\begin{equation}\label{autunno2+}
\mathcal E^{2j}(\mathtt r,\tilde\checco)
:=
 (\bar E^{(2j)}_- ,R_0^2-2\morse)_{\mathtt r}
\cap \{\Re  E>\bar E^{(2j)}_- + \tilde \checco\}
\,.
\end{equation}
Note that this family of sets is
increasing with  
$\mathtt r$ and decreasing with
$\tilde \checco.$

\begin{lemma}\label{giga2bis+}
The functions $ \partial_E\act_n^{(2j),\pm}(E,\hat\act)$
 are holomorphic for
 \begin{equation}\label{suzanne}
 E\in \mathcal E^{2j}(r_0^2/2^7,\tilde\checco)\,,\qquad
 \forall\, \tilde\checco\geq \frac{4\morse}{\b s_0^2}\checco
\end{equation}
 and $\hat\act \in \hat D_{r_0}$
 with uniform estimate
\begin{equation}\label{suzanne2}
\sup_{\mathcal E^{2j}(r_0^2/2^7,\tilde\checco)
\times \hat D_{r_0}}
 |\partial_E\act_n^{(2j),\pm}(E,\hat\act)|\leq 
\frac{2^9}{\sqrt\b}\left(
\frac{\morse}{\b s_0^3}+
\ln \frac{\b^3 s_0^6}{\tilde\checco \morse^2}\right)\,.
\end{equation}
\end{lemma}
\proof
We consider only $\partial_E\act_n^{(2j),+}(E)$
since the argument for $\partial_E\act_n^{(2j),-}(E)$
is analogous.
Recalling \eqref{ancomarzio} and splitting the integral we have that
\begin{equation}\label{armenia}
2\pi\partial_E \act_n^{(2j),+}(E)
 =
 \int_{\sa_{2j}}^{\bar\sa_{2j}}
 +
 \int_{\bar\sa_{2j}}^{\bar\sa_{2j+1}}
 +
\int_{\bar\sa_{2j+1}}^{\sa_{2j+1}}
\frac{
1
+ \tilde b\big(E-\Gm(\sa),\sa \big)
}
{\sqrt{E-\Gm(\sa)}}
\, d\sa\,.
\end{equation}
For $|\Im \sa|<s_0$ by \eqref{ciccio} we have
\begin{equation}\label{cilicia}
|E-\Gm(\sa)|\geq |E-\FO(\sa)|-\checco
\geq |\Re(E-\FO(\sa))|-\checco
\end{equation}
and, by \eqref{harlock1} and recalling that 
$\FO(\sa)$ is decreasing for 
$\sa\in[\bar\sa_{2j},\bar\sa_{2j+1}]$
we have, for $E$ as in \eqref{suzanne},
\begin{equation}\label{cilicia2}
|E-\Gm(\sa)|\geq
\tilde\checco+
\bar E_{2j}-\FO(\sa)
-\checco
\geq
\left\{
\begin{array}{ll} 
\displaystyle
\frac{\tilde\checco}{2}+\frac{\b}{2}(\sa-\bar\sa_{2j})^2 
& {\rm if}\ \ \bar\sa_{2j}\leq \sa\leq \bar\sa_{2j}+\sa_\sharp\\ 
& \\
\displaystyle
\frac{\b}{2}\sa_\sharp^2 
& {\rm if}\ \  \bar\sa_{2j}+\sa_\sharp
\leq
\sa\leq \bar\sa_{2j+1}
\end{array}
\right.
\end{equation}
where
$
 \xx_\sharp
=\frac{\b s_0^3}{6\morse}
$
was defined in \eqref{harlock1}.
Moreover by \eqref{3holesquater}
we get
\begin{equation}\label{cilicia3}
|1
+ \tilde b\big(E-\Gm(\sa),\sa \big)|
\leq 2\,.
\end{equation}
Then by \eqref{cilicia2}, \eqref{cilicia3}
we get\footnote{Using that, for $x_0\geq 9$
we have
$\int_0^{x_0}(1+x^2)^{-1/2}dx=
\ln(\sqrt{1+x_0^2}+x_0)\leq 2\ln x_0.$}
\begin{equation}
\label{cilicia4}
\left| \int_{\bar\sa_{2j}}^{\bar\sa_{2j+1}}
\frac{
1
+ \tilde b\big(E-\Gm(\sa),\sa \big)
}
{\sqrt{E-\Gm(\sa)}}
\, d\sa
\right|
\leq
\int_{\bar\sa_{2j}}^{\bar\sa_{2j+1}}
\frac{
2}
{\sqrt{|E-\Gm(\sa)|}}
\, d\sa
\leq
\frac{2^8}{\sqrt\b}\left(
\frac{\morse}{\b s_0^3}+
\ln \frac{\b^3 s_0^6}{\tilde\checco \morse^2}\right)\,.
\end{equation}

For $|\sa-\bar\sa_{2j}|\leq 2\checco/\b s_0$
and
$E$ as in \eqref{suzanne}
we have by \eqref{cilicia}, \eqref{octoberx} and Cauchy estimates
\begin{eqnarray*}
&&|E-\Gm(\sa)|\geq |E-\FO(\sa)|-\checco
=
|E-\bar E_{2j}+\FO(\bar\sa_{2j})-\FO(\sa)|-\checco
\nonumber
\\
&&\geq 
\tilde\checco -
\frac{2\checco\morse}{\b s_0^2}-\checco
\geq 
\frac{\tilde\checco}{4}\,,
\end{eqnarray*}
by \eqref{harlock};
then, by \eqref{cilicia3} and \eqref{octoberx}, we get
\begin{equation}\label{cilicia5}
\left| \int_{\sa_{2j}}^{\bar\sa_{2j}}
\frac{
1
+ \tilde b\big(E-\Gm(\sa),\sa \big)
}
{\sqrt{E-\Gm(\sa)}}
\, d\sa
\right|
\leq
\frac{8\checco}{\sqrt{\tilde\checco}\b s_0}
\leq 
\frac{1}{\sqrt{\morse}}\,,
\end{equation}
where the last inequality follows by \eqref{genesis} and \eqref{suzanne}.
The estimate for the last integral in \eqref{armenia}
is analogous (even better).
Then by \eqref{cilicia4},\eqref{cilicia5} (and \eqref{harlock})
we get \eqref{suzanne2}.
\eproof

By \eqref{noce2} and lemmata \ref{giga2+}, \ref{giga2bis+} 
we get the following
\begin{corollary}
 The function $\act_n^{(2j)}(E,\hat\act)$
 is holomorphic for $\hat\act \in \hat D_{r_0}$ and 
$$
E\in (\bar E^{(2j)}_-,\bar E^{(2j)}_+)_{2\mathtt r_4}
\cap \{ E^{(2j)}_-(\hat\act)+\C_*\}
\cap \{ E^{(2j)}_+(\hat\act)-\C_*\}\,.
$$ 
\end{corollary}



\subsection{Closeness of  the  unperturbed and perturbed actions}

\begin{lemma}
Let $1\leq j\leq N.$
 For 
 \begin{equation}\label{bassora}
    \frac{2^5\morse\checco}{\b}
  \leq
  \tilde\checco
  \leq \mathtt r_4\stackrel{\eqref{acri}}=
  \frac{\mathtt r_2\morse}{2^5}
  \,,
\end{equation}
 we have
 \begin{equation}
\Big|
\partial_E \act_n^{(2j-1),\pm}(E,\hat\act)
-
\partial_E \bar \act_n^{(2j-1),\pm}(E)
\Big|
\leq
\left(
\frac{ 2^{78}\morse^{7}}{s_0^{12}\b^{7} \mathtt r_2^{3/2}}
+\frac{s_*^{12}\morse}{\tilde\checco}
\right)
\frac{2^{17}\morse^{1/2}}{\b^2 s_0}
\checco
\label{damietta}
\end{equation}
for\footnote{
$\mathcal E^{2j-1}(\mathtt r_4,\tilde\checco)$ was
defined in \eqref{autunno2}.} 
\begin{equation}\label{autunno}
E\in\mathcal E^{2j-1}(\mathtt r_4,\tilde\checco)
\cap \{\Re  E\geq\bar E^{(2j-1)}_- + \mathtt r_2\morse/2^7\}
\,,\qquad
\hat\act\in \hat D_{r_0}
\,.
\end{equation}
\end{lemma}
\proof
By \eqref{valeriabar}, \eqref{valerione}
and \eqref{collemaggio}
we have
\begin{eqnarray}
&&\Big|
\frac{\sqrt 2\sqrt{E_{2j-2}- E_{2j-1}}}{\sa_{2j-1}-\sa_{2j-2}}
\partial_E \act_n^{(2j-1),-}\big(\l_{2j-1}^{-1}(\breve E)\big)
\nonumber
\\
&&\qquad-
\frac{\sqrt 2\sqrt{\bar E_{2j-2}- \bar E_{2j-1}}}{\bar \sa_{2j-1}-\bar \sa_{2j-2}}
\partial_E \bar \act_n^{(2j-1),-}\big(\bar \l_{2j-1}^{-1}(\breve E)\big)
\Big|
\nonumber
\\
&&=
|\partial_{\breve E}\breve\act_n^{(2j-1),-}(\breve E,\hat\act)
-
\partial_{\breve E}\breve{\bar\act}_n^{(2j-1),-}(\breve E)
|
\leq
\nonumber
\\
&&
\checco
\left(
\frac{36}{r_0^2}
+
\frac{2^{194} \morse^{27}}
{ \b^{28}  s_0^{45}s_*^2} 
\right)
\frac{ 2^{55}\morse^7}{\b^7  s_0^{13}}
\left( 
1
+
\ln \frac{1}{1-\Re \breve E}
\right)
\,,
\label{afrodite}
\end{eqnarray}
for
$\hat\act\in \hat D_{r_0}$ and
  $\breve E\in \O_{\breve{\mathtt r}}$ 
  with $\Re \breve E<1$.
  Making the substitution 
$\breve E=\bar\l_{2j-1}(E)$, \eqref{afrodite} becomes
\begin{eqnarray}
&&\Big|
\frac{ \sqrt{E_{2j-2}- E_{2j-1}}}{\sa_{2j-1}-\sa_{2j-2}}
\partial_E \act_n^{(2j-1),-}\big(\l_{2j-1}^{-1}(\bar\l_{2j-1}(E))\big)
-
\frac{\sqrt{\bar E_{2j-2}- \bar E_{2j-1}}}{\bar \sa_{2j-1}-\bar \sa_{2j-2}}
\partial_E \bar \act_n^{(2j-1),-}(E)
\Big|
\nonumber
\\
&&\leq
\checco
\left(
\frac{36}{r_0^2}
+
\frac{2^{194} \morse^{27}}
{ \b^{28}  s_0^{45}s_*^2} 
\right)
\frac{ 2^{55}\morse^7}{\b^7  s_0^{13}}
\left( 
1
+
\ln \frac{\bar E_{2j-2}-\bar E_{2j-1}}{2(\bar E_{2j-2}-\Re E)}
\right)
\,,
\label{afrodite2}
\end{eqnarray}
for
 $\hat\act\in \hat D_{r_0}$
 and 
 $$
 E\in \bar\l_{2j-1}^{-1}\Big(\O_{\breve{\mathtt r}}
\cap \{\Re \breve E<1\}\Big)\,.
 $$
Recalling \eqref{acrobat2} and \eqref{ladispoli} we note that
\begin{eqnarray}\label{venere}
&&\bar\l_{2j-1}^{-1}\Big(\O_{\breve{\mathtt r}}
\cap \{\Re \breve E<1\}\Big)
\nonumber
\\
&&
\supseteq
\
(\bar E_{2j-1},\bar E_{2j-2})_{\b\breve{\mathtt r}/2}
\cap (\C_*+\bar E_{2j-1})
\cap \{\Re  E<\bar E_{2j-2}\}
\nonumber
\\
&&\supseteq
\mathcal E^{2j-1}(2\mathtt r_4,0)
\cap (\C_*+\bar E_{2j-1})
\,,
\end{eqnarray}
where the last inclusion holds
since $4\mathtt r_4\leq \b\breve{\mathtt r}/2$
 by \eqref{primavera}
(recall also \eqref{gruffalo}).
By \eqref{venere} we have that 
 \begin{equation}\label{afrodite3}
 \eqref{afrodite2}
\ \text{holds for}\   
 E\in 
 \mathcal E^{2j-1}(2\mathtt r_4,0)
\cap (\C_*+\bar E_{2j-1})
\end{equation}
and $\hat\act\in \hat D_{r_0}$.
By \eqref{inverno} we have, for 
$|E|\leq 2\morse$ and 
$\hat\act\in \hat D_{r_0},$
\begin{equation}\label{baghdad}
|\l_{i}^{-1}(\bar\l_i(E),\hat\act)-E|
\leq 
\frac{12\checco\morse}{\b}
\leq \min\left\{ \frac{\tilde\checco}{4}, \ \frac{\mathtt r_4}{2^5}\right\}
\leq\frac{\mathtt r_2\morse}{2^{10}}
\,,
\end{equation}
where the last inequality follows by \eqref{genesis}, \eqref{bassora} and \eqref{harlock}.
By \eqref{baghdad} and \eqref{harlock}
\begin{equation}\label{edessa}
\tilde\checco-\frac{12\checco\morse}{\b}
\geq \frac14\tilde\checco
\geq \frac{4\checco\morse}{\b}
\geq 2\checco\,.
\end{equation}
Then
\begin{eqnarray}\label{hormuz}
&&E\in\mathcal E^{2j-1}(\mathtt r_4,\tilde\checco)
\cap \{\Re  E\geq\bar E^{(2j-1)}_- + \mathtt r_2\morse/2^7\}
\nonumber
\\
&&
\Longrightarrow
\quad
\l_{2j-1}^{-1}(\bar\l_{2j-1}(E),\hat\act)
\in
\mathcal E^{2j-1}_{**}(3\mathtt r_4/2,\tilde\checco/4)
\subset
 \mathcal E^{2j-1}_*(2\mathtt r_4,2\tilde\checco)
\end{eqnarray}
for every $\hat\act\in \hat D_{r_0}.$
For
$$
E\in\mathcal E^{2j-1}(\mathtt r_4,\tilde\checco)
\cap \{\Re  E\geq\bar E^{(2j-1)}_- + \mathtt r_2\morse/2^7\}
\,
\subset
\,
\mathcal E^{2j-1}_{**}(3\mathtt r_4/2,\tilde\checco/4)
$$
by
 \eqref{hormuz} and \eqref{baghdad} we have 
that\footnote{Note that $\mathcal E^{2j-1}_{**}(3\mathtt r_4/2,\tilde\checco/4)$
is a convex set.} 
\begin{eqnarray}
&&\Big|
\partial_E \act_n^{(2j-1),-}\big(\l_{2j-1}^{-1}(\bar\l_{2j-1}(E))\big)
-
\partial_E \act_n^{(2j-1),-}(E)
\Big|
\nonumber
\\
&&\leq
\sup_{\mathcal E^{2j-1}_{**}(3\mathtt r_4/2,\tilde\checco/4)\times \hat D_{r_0}}
|\partial_{EE} \act_n^{(2j-1),-}| 
\frac{12\checco\morse}{\b}
\nonumber
\\
&&\leq
\left(
\frac{ 2^{94}\morse^{15/2}}{s_0^{13}\b^{9} \mathtt r_2^{3/2}}
+\frac{2^{16}\morse^{3/2}}{\b^2 s_0\tilde\checco}
\right)
\checco
\label{medina}
\end{eqnarray}
by \eqref{damasco2} (used with
$\tilde\checco\rightsquigarrow\tilde\checco/8$).

Let us estimate
\begin{eqnarray}\label{beirut}
& &
\left|\sqrt{E_{2j-2}- E_{2j-1}}-\sqrt{\bar E_{2j-2}- \bar E_{2j-1}}\right|
\nonumber
\\
&=&
\sqrt{|\bar E_{2j-2}- \bar E_{2j-1}|}
\left|
\sqrt{1+\frac{E_{2j-2}-\bar E_{2j-2}+ \bar E_{2j-1}- E_{2j-1}}{\bar E_{2j-2}- \bar E_{2j-1}}}
-1
\right|
\nonumber
\\
&\leq&
8\frac{\sqrt\morse}{\b}\checco\,, 
\end{eqnarray}
by \eqref{goffredo},\eqref{october}.
Then we have
\begin{eqnarray}
& &\left|
\frac{ \sqrt{E_{2j-2}- E_{2j-1}}}{\sa_{2j-1}-\sa_{2j-2}}
-
\frac{\sqrt{\bar E_{2j-2}- \bar E_{2j-1}}}{\bar \sa_{2j-1}-\bar \sa_{2j-2}}
\right|
\nonumber
\\
&\leq&
\frac{|\sqrt{E_{2j-2}- E_{2j-1}}-\sqrt{\bar E_{2j-2}- \bar E_{2j-1}}|}
{|\sa_{2j-1}-\sa_{2j-2}|}
\nonumber
\\
& &
+
\frac{|\sqrt{\bar E_{2j-2}- \bar E_{2j-1}}|
(|\sa_{2j-2}-\bar\sa_{2j-2}|+|\sa_{2j-1}-\bar\sa_{2j-1}|)}
{|\sa_{2j-1}-\sa_{2j-2}||\bar \sa_{2j-1}-\bar \sa_{2j-2}|}
\nonumber
\\
&\leq&
\sqrt{\frac{3 \morse}{\b s_0^3}}
\left(
8\frac{\sqrt\morse}{\b}\checco
+
\frac{8\morse\checco}{s_0^{5/2}\b^{3/2}}
\right)
\nonumber
\\
&\leq&
\frac{2^6\morse^{3/2}}{s_0^4\b^2}\checco
\,,
\label{palmira}
\end{eqnarray}
by \eqref{beirut},\eqref{octoberx},\eqref{alcafone},\eqref{octobery}
and \eqref{harlock}.
Note that by \eqref{ladispoli}
\begin{equation}\label{pasargade}
\left|
\frac{\bar \sa_{2j-1}-\bar \sa_{2j-2}}{\sqrt{\bar E_{2j-2}- \bar E_{2j-1}}}
\right|
\leq \frac{2\pi}{\sqrt\b}\,.
\end{equation}
For $E$ as in \eqref{autunno}
we have  by \eqref{hormuz},
 \eqref{palmira} and
 \eqref{damasco} that
\begin{eqnarray}
&&\left|
\left(\frac{ \sqrt{E_{2j-2}- E_{2j-1}}}{\sa_{2j-1}-\sa_{2j-2}}
-
\frac{\sqrt{\bar E_{2j-2}- \bar E_{2j-1}}}{\bar \sa_{2j-1}-\bar \sa_{2j-2}}
\right)
\partial_E \act_n^{(2j-1),-}\big(\l_{2j-1}^{-1}(\bar\l_{2j-1}(E))\big)
\right|
\nonumber
\\
&\leq&
\frac{2^6\morse^{3/2}}{s_0^4\b^2}\checco
\left(
\frac{ 2^{80}\morse^{15/2}}{s_0^{13}\b^{8} \sqrt{\mathtt r_2}}
+ 
\frac{2^8\sqrt\morse}{\b s_0}
\ln \frac{\morse}{\tilde\checco}
\right)\,.
\label{krak}
\end{eqnarray} 
Then
we have for $E$ as in \eqref{autunno}
\begin{eqnarray}
& &
\Big|
\partial_E \act_n^{(2j-1),-}\big(\l_{2j-1}^{-1}(\bar\l_{2j-1}(E))\big)
-
\partial_E \bar \act_n^{(2j-1),-}(E)
\Big|
\nonumber
\\
&\stackrel{\eqref{pasargade}}\leq&
\frac{2\pi}{\sqrt\b}
\left|\frac{\sqrt{\bar E_{2j-2}- \bar E_{2j-1}}}{\bar \sa_{2j-1}-\bar \sa_{2j-2}}
\Big(
\partial_E \act_n^{(2j-1),-}\big(\l_{2j-1}^{-1}(\bar\l_{2j-1}(E))\big)
-
\partial_E \bar \act_n^{(2j-1),-}(E)
\Big)
\right|
\nonumber
\\
&\stackrel{\eqref{krak},\eqref{afrodite3}}\leq&
\frac{2^9\morse^{3/2}}{s_0^4\b^{5/2}}\checco
\left(
\frac{ 2^{80}\morse^{15/2}}{s_0^{13}\b^{8} \sqrt{\mathtt r_2}}
+ 
\frac{2^8\sqrt\morse}{\b s_0}
\ln \frac{\morse}{\tilde\checco}
\right)
\nonumber
\\
& &+
\checco
\left(
\frac{36}{r_0^2}
+
\frac{2^{194} \morse^{27}}
{ \b^{28}  s_0^{45}s_*^2} 
\right)
\frac{ 2^{58}\morse^7}{\b^{15/2}  s_0^{13}}
\left( 
1
+
\ln \frac{\bar E_{2j-2}-\bar E_{2j-1}}{2\tilde\checco}
\right)
\nonumber
\\
&\leq&
\checco
\frac{ 2^{89}\morse^{9}}{s_0^{17}\b^{21/2} \sqrt{\mathtt r_2}}
+
\checco
\left(
\frac{1}{r_0^2}
+
\frac{2^{194} \morse^{27}}
{ \b^{28}  s_0^{45}s_*^2} 
\right)
\frac{ 2^{65}\morse^7}{\b^{15/2}  s_0^{13}}
\ln\frac{\morse}{\tilde\checco}
\nonumber
\\
&&
\leq
\checco
\left(
\frac{1}{r_0^2}
+
\frac{2^{194} \morse^{27}}
{ \b^{28}  s_0^{45}s_*^2} 
\right)
\frac{ 2^{66}\morse^7}{\b^{15/2}  s_0^{13}}
\ln\frac{\morse}{\tilde\checco}
\label{antiochia}
\end{eqnarray}
using (\eqref{harlock} and) that by \eqref{bassora} and \eqref{acri}
$$
\frac{\morse}{\tilde\checco}
\geq
\frac{\morse}{\mathtt r_4}=\frac{2^5}{\mathtt r_2}
$$
and, in the last inequality, \eqref{tropicana}.
\\
Recollecting for $E$ as in \eqref{autunno} we have
that \eqref{damietta} follows by
\begin{eqnarray*}
& &\Big|
\partial_E \act_n^{(2j-1),-}(E)
-
\partial_E \bar \act_n^{(2j-1),-}(E)
\Big|
\\
&\stackrel{\eqref{medina}}\leq&
\Big|
\partial_E \act_n^{(2j-1),-}\big(\l_{2j-1}^{-1}(\bar\l_{2j-1}(E))\big)
-
\partial_E \bar \act_n^{(2j-1),-}(E)
\Big|
\\
& &
+
\left(
\frac{ 2^{94}\morse^{15/2}}{s_0^{13}\b^{9} \mathtt r_2^{3/2}}
+\frac{2^{16}\morse^{3/2}}{\b^2 s_0\tilde\checco}
\right)
\checco
\\
&\stackrel{\eqref{antiochia}}\leq&
\checco
\left(
\frac{1}{r_0^2}
+
\frac{2^{194} \morse^{27}}
{ \b^{28}  s_0^{45}s_*^2} 
\right)
\frac{ 2^{66}\morse^7}{\b^{15/2}  s_0^{13}}
\ln\frac{\morse}{\tilde\checco}
+
\left(
\frac{ 2^{94}\morse^{15/2}}{s_0^{13}\b^{9} \mathtt r_2^{3/2}}
+\frac{2^{16}\morse^{3/2}}{\b^2 s_0\tilde\checco}
\right)
\checco
\\
&=&
\left[
\left(
\frac{2^{10}\morse}{r_0^2}
+
\frac{2^{204} \morse^{28}}
{ \b^{28}  s_0^{45}s_*^2} 
\right)
\frac{ 2^{40}\morse^{11/2}}{\b^{11/2}  s_0^{12}}
\ln \frac{\morse}{\tilde\checco}
+
\left(
\frac{ 2^{78}\morse^{7}}{s_0^{12}\b^{7} \mathtt r_2^{3/2}}
+\frac{\morse}{\tilde\checco}
\right)
\right]
\frac{2^{16}\morse^{1/2}}{\b^2 s_0}
\checco
\\
&\stackrel{\eqref{badedas},\eqref{tropicana}}\leq&
\left[
\frac{ 2^{41}\morse^{11/2}}{\b^{11/2}  s_0^{12}\mathtt r_2}
\ln \frac{\morse}{\tilde\checco}
+
\left(
\frac{ 2^{78}\morse^{7}}{s_0^{12}\b^{7} \mathtt r_2^{3/2}}
+\frac{\morse}{\tilde\checco}
\right)
\right]
\frac{2^{16}\morse^{1/2}}{\b^2 s_0 }
\checco
\\
&\leq&
\left(
\frac{ 2^{78}\morse^{7}}{s_0^{12}\b^{7} \mathtt r_2^{3/2}}
+\frac{s_*^{12}\morse}{\tilde\checco}
\right)
\frac{2^{17}\morse^{1/2}}{\b^2 s_0}
\checco\,,
\end{eqnarray*}
where in the last inequality we have used that
\begin{equation}\label{compleannoTommy}
\frac{ 2^{41}\morse^{11/2}s_*^{12}}{\b^{11/2}  s_0^{12}\mathtt r_2}
\ln \frac{\morse}{\tilde\checco}
\leq
\left(
\frac{ 2^{78}\morse^{7}s_*^{12}}{s_0^{12}\b^{7} \mathtt r_2^{3/2}}
+\frac{\morse}{\tilde\checco}
\right)\,.
\end{equation}
In order to prove
\eqref{compleannoTommy}
we first note that, for $a,b>0,$ 
$$
\min_{x>0} (x+b-a \ln x)=a+b-a\ln a\,.
$$
Using the above formula with
\begin{equation}\label{compleanno Betta}
a=\frac{ 2^{41}\morse^{11/2}s_*^{12}}{\b^{11/2}  s_0^{12}\mathtt r_2}\,,\qquad
b=\frac{ 2^{78}\morse^{7}s_*^{12}}{s_0^{12}\b^{7} \mathtt r_2^{3/2}}
\end{equation}
(and
$x=\frac{\morse}{\tilde\checco}$)
we have that \eqref{compleannoTommy}
follows if we show that
$$
b\geq  a\ln a\,.
$$
This last estimate follows by \eqref{tropicana} and \eqref{harlock}.
\eproof

\begin{corollary}\label{sarabanda3}
Let $1\leq j\leq N.$
If \eqref{bassora} holds
we have that the functions 
$
\partial_E \act_n^{(2j-1)}(E,\hat\act),$
$
\partial_E \bar \act_n^{(2j-1)}(E)
$ are holomorphic on $\mathcal E^{2j-1}(\mathtt r_4/4,\tilde\checco)
\times\hat D_{r_0}$ with 
 \begin{eqnarray}
&&
\sup_{\mathcal E^{2j-1}(\mathtt r_4/4,\tilde\checco)
\times\hat D_{r_0}}
\Big|
\partial_E \act_n^{(2j-1)}(E,\hat\act)
-
\partial_E \bar \act_n^{(2j-1)}(E)
\Big|
\nonumber
\\
&&
\qquad
\leq
\left(
\frac{ 2^{78}\morse^{7}}{s_0^{12}\b^{7} \mathtt r_2^{3/2}}
+\frac{s_*^{12}\morse}{\tilde\checco}
\right)
\frac{2^{22}\morse^{1/2}}{\b^2 s_0}
\checco
\,.
\label{rosetta}
\end{eqnarray}
Then\footnote{By Cauchy estimates.}
\begin{equation}
\sup_{\mathcal E^{2j-1}(\mathtt r_4/4,\tilde\checco)
\times\hat D_{r_0/2}}
\Big|
\partial_{E\hat\act} \act_n^{(2j-1)}(E,\hat\act)
\Big|
\leq
\left(
\frac{ 2^{78}\morse^{7}}{s_0^{12}\b^{7} \mathtt r_2^{3/2}}
+\frac{s_*^{12}\morse}{\tilde\checco}
\right)
\frac{2^{23}\morse^{1/2}}{\b^2 s_0 r_0}
\checco
\,.
\label{rosetta*}
\end{equation}
\end{corollary}
\proof
By \eqref{damietta}, \eqref{autunno}, \eqref{citola} and \eqref{noce} we have
 \begin{eqnarray}
&&\Big|
\partial_E \act_n^{(2j-1)}(E,\hat\act)
-
\partial_E \bar \act_n^{(2j-1)}(E)
\Big|
\leq
\left(
\frac{ 2^{78}\morse^{7}}{s_0^{12}\b^{7} \mathtt r_2^{3/2}}
+\frac{s_*^{12}\morse}{\tilde\checco}
\right)
\frac{2^{18}\morse^{1/2}}{\b^2 s_0}
\checco
\label{damietta2}
\\
&&\text{for}\ \ E\in\mathcal E^{2j-1}(\mathtt r_4,\tilde\checco)
\cap \{\Re  E\geq\bar E^{(2j-1)}_- + \mathtt r_4/4\}
\,,\qquad
\hat\act\in \hat D_{r_0}
\,.
\nonumber
\end{eqnarray}
By Lemma \ref{giga}, \eqref{epistassi} and \eqref{acri}
the function
$$
\partial_E \act_n^{(2j-1)}(E,\hat\act)
-
\partial_E \bar \act_n^{(2j-1)}(E)
$$
is holomorphic
on the closed ball
$$
\overline{B_{3\mathtt r_4/4}(\bar E^{(2j-1)}_- + \mathtt r_4/4)}
$$
and 
 \eqref{damietta2}
holds on one of its diameter, namely on
$$
\overline{B_{3\mathtt r_4/4}(\bar E^{(2j-1)}_- + \mathtt r_4/4)}
\cap \{\Re E=\bar E^{(2j-1)}_- + \mathtt r_4/4\}\,.
$$
Then, by the Borel-Caratheodory 
Theorem\footnote{
The Borel-Caratheodory 
Theorem says that {\it if $f$ is a holomorphic function on the closed disc of radius $R$
centered in $z_0$ 
and $|f(z)|$ is bounded by $M>0$ on a diameter of the above disc
then, for every $0<r<R,$}
$$
\sup_{|z-z_0|\leq r} |f(z)|\leq 
\frac{2r}{R-r}M+\frac{R+r}{R-r}|f(z_0)|
\leq \frac{R+3r}{R-r}M
\,.$$
We use the theorem with $z:=E$, 
$z_0:=\bar E^{(2j-1)}_- + \mathtt r_4/4,$
$f:=\partial_E \act_n^{(2j-1)}
-
\partial_E \bar \act_n^{(2j-1)},$ $R:=3\mathtt r_4/4,$
$r:=\mathtt r_4/2,$ 
$M:=$the right had side of \eqref{damietta2}.
},
we get, for every $\hat\act\in \hat D_{r_0}$,
 \begin{eqnarray*}
&&
\sup_{|E-\bar E^{(2j-1)}_- - \mathtt r_4/4|\leq \mathtt r_4/2}
\Big|
\partial_E \act_n^{(2j-1)}(E,\hat\act)
-
\partial_E \bar \act_n^{(2j-1)}(E)
\Big|
\nonumber
\\
&&
\leq
\left(
\frac{ 2^{78}\morse^{7}}{s_0^{12}\b^{7} \mathtt r_2^{3/2}}
+\frac{s_*^{12}\morse}{\tilde\checco}
\right)
\frac{2^{22}\morse^{1/2}}{\b^2 s_0}
\checco
\,.
\end{eqnarray*}
Combining this estimate with \eqref{damietta2}
we get \eqref{rosetta}.
\eproof

\begin{lemma}
Let $0\leq j\leq N.$
 For $\hat\act\in \hat D_{r_0}$ 
 and\footnote{Recalling \eqref{autunno2+} and \eqref{gruffalo}
$\mathcal E^{2j}(\mathtt r,\tilde\checco)
:=
 (\bar E_{2j} ,R_0^2-2\morse)_{\mathtt r}
\cap \{\Re  E>\bar E_{2j} + \tilde \checco\}.$}
 \begin{equation}\label{bassora+}
   E\in \mathcal E^{2j}(r_0^2\b/2^{10}\morse,\tilde\checco)\,,\qquad
 \forall\, \tilde\checco\geq \frac{2^5\morse}{\b s_*^2}\checco
\end{equation}
 we have
 \begin{equation}
\Big|
\partial_E \act_n^{(2j),\pm}(E,\hat\act)
-
\partial_E \bar \act_n^{(2j),\pm}(E)
\Big|
\leq
\left[
\left(
\frac{18}{r_0^2}
+
\frac{2^{193} \morse^{27}}
{ \b^{28}  s_0^{45}s_*^2} 
\right)
\frac{ 2^{57}\morse^7}{\b^{15/2}  s_0^{13}}
\ln\left(
4+\frac{\b}{2\tilde\checco}
\right)
+
\frac{2^{21}\morse^2}{\b^{5/2}s_0^3 \tilde\checco}
\right]
\checco
\,.
\label{damietta+}
\end{equation}
\end{lemma}
\proof
For
$\hat\act\in \hat D_{r_0}$ and
$\breve E\in \O^{(2j),-}_{r_0^2/2^8\morse}$ (defined in \eqref{vltava}) with
$\Re \breve E>1$, namely 
\begin{equation}\label{vltava+}
\breve E\in 
\O_\sharp:=
\left(1,\frac{2 R_0^2}{\bar E_{2j}-\bar E_{2j-1}}-1\right)_{r_0^2/2^8\morse}
\cap\{ \Re \breve E>1 \}\,.
\end{equation}
we have,
deriving \eqref{valeriabar} (recall \eqref{acrobat})
and \eqref{collemaggio+},
that
\begin{eqnarray}
&&\Big|
\frac{\sqrt 2\sqrt{E_{2j}- E_{2j-1}}}{\sa_{2j}-\sa_{2j-1}}
\partial_E \act_n^{(2j),-}\big(\l_{2j}^{-1}(\breve E)\big)
-
\frac{\sqrt 2\sqrt{\bar E_{2j}- \bar E_{2j-1}}}{\bar \sa_{2j}-\bar \sa_{2j-1}}
\partial_E \bar \act_n^{(2j),-}\big(\bar \l_{2j}^{-1}(\breve E)\big)
\Big|
\nonumber
\\
&&=
|\partial_{\breve E}\breve\act_n^{(2j),-}(\breve E,\hat\act)
-
\partial_{\breve E}\breve{\bar\act}_n^{(2j),-}(\breve E)
|
\nonumber
\\
&&
\leq
\checco
\left(
\frac{18}{r_0^2}
+
\frac{2^{193} \morse^{27}}
{ \b^{28}  s_0^{45}s_*^2} 
\right)
\frac{ 2^{55}\morse^7}{\b^7  s_0^{13}}
\ln\left(
4+\frac{1}{\Re\breve E-1}
\right)
\,.
\label{afrodite+}
\end{eqnarray}
Making the substitution 
$\breve E=\bar\l_{2j}(E)$, \eqref{afrodite+} becomes (recalling \eqref{acrobat})
\begin{eqnarray}
&&\Big|
\frac{\sqrt{E_{2j}- E_{2j-1}}}{\sa_{2j}-\sa_{2j-1}}
\partial_E \act_n^{(2j-1),-}\big(\l_{2j}^{-1}(\bar\l_{2j}(E))\big)
-
\frac{\sqrt{\bar E_{2j}- \bar E_{2j-1}}}{\bar \sa_{2j}-\bar \sa_{2j-1}}
\partial_E \bar \act_n^{(2j-1),-}(E)
\Big|
\nonumber
\\
&&\leq
\checco
\left(
\frac{18}{r_0^2}
+
\frac{2^{193} \morse^{27}}
{ \b^{28}  s_0^{45}s_*^2} 
\right)
\frac{ 2^{55}\morse^7}{\b^7  s_0^{13}}
\ln\left(
4+\frac{1}{\Re\bar\l_{2j}(E)-1}
\right)
\,,
\label{afrodite3+}
\end{eqnarray}
for
 $\hat\act\in \hat D_{r_0}$
 and 
\begin{equation}\label{venere+}
 E\in \bar\l_{2j}^{-1}(\O_\sharp)
 \supseteq
 \left(
 \bar E_{2j}, R_0^2+\frac{\bar E_{2j}}{2}+\frac{3\bar E_{2j-1}}{2}
 \right)_{\frac{r_0^2\b}{2^9 \morse}}
 \cap \{ \Re E>\bar E_{2j}\}
 \supset
 \mathcal E^{2j}(r_0^2\b/2^{10}\morse,0)
 \,,
 \end{equation}
recalling \eqref{acrobat2} and \eqref{ladispoli}.
\\
For $E$ as in \eqref{bassora+}, by \eqref{acrobat} and \eqref{ladispoli}
we get
$$
\Re\bar\l_{2j}(E)-1=
\frac{2\tilde\checco}{\bar E_{2j}-\bar E_{2j-1}}
\geq 
\frac{2\tilde\checco}{\b}\,,
$$
then,
by \eqref{afrodite+}, we get
\begin{eqnarray}
&&\Big|
\frac{\sqrt{E_{2j}- E_{2j-1}}}{\sa_{2j}-\sa_{2j-1}}
\partial_E \act_n^{(2j-1),-}\big(\l_{2j}^{-1}(\bar\l_{2j}(E))\big)
-
\frac{\sqrt{\bar E_{2j}- \bar E_{2j-1}}}{\bar \sa_{2j}-\bar \sa_{2j-1}}
\partial_E \bar \act_n^{(2j-1),-}(E)
\Big|
\nonumber
\\
&&\leq
\checco
\left(
\frac{18}{r_0^2}
+
\frac{2^{193} \morse^{27}}
{ \b^{28}  s_0^{45}s_*^2} 
\right)
\frac{ 2^{55}\morse^7}{\b^7  s_0^{13}}
\ln\left(
4+\frac{\b}{2\tilde\checco}
\right)
\,,
\label{afrodite3++}
\end{eqnarray}
Note that
by \eqref{baghdad}, \eqref{bassora+}, \eqref{genesis}, \eqref{harlock}
\begin{equation}\label{hormuz+}
E\in  \mathcal E^{2j}(r_0^2\b/2^{10}\morse,\tilde\checco)
\quad
\Longrightarrow
\quad
\l_{2j}^{-1}(\bar\l_{2j}(E))
\in
\mathcal E^{2j}(r_0^2/2^8,\tilde\checco/2)\,.
\end{equation}
By  \eqref{hormuz+},
 \eqref{baghdad},  \eqref{suzanne2} and Cauchy 
 estimates\footnote{
 Indeed we have by \eqref{suzanne2} used with
 $\tilde\checco\rightsquigarrow \tilde\checco/4$ and Cauchy 
 estimates
 $$
 \sup_{\mathcal E^{2j}(r_0^2/2^8,\tilde\checco/2)
\times \hat D_{r_0}}
 |\partial_{EE}\act_n^{(2j),\pm}(E,\hat\act)|\leq 
\frac{2^{17}}{\sqrt\b \max\{r_0^2,\tilde\checco\}}\left(
\frac{\morse}{\b s_0^3}+
\ln \frac{\b^3 s_0^6}{\tilde\checco \morse^2}\right)
 $$} 
 we have 
that 
\begin{eqnarray}
&&\Big|
\partial_E \act_n^{(2j),-}\big(\l_{2j}^{-1}(\bar\l_{2j}(E))\big)
-
\partial_E \act_n^{(2j),-}(E)
\Big|
\nonumber
\\
&&\leq
\sup_{\mathcal E^{2j}(r_0^2/2^8,\tilde\checco/2)
\times \hat D_{r_0}}
|\partial_{EE} \act_n^{(2j),-}| 
\frac{12\checco\morse}{\b}
\nonumber
\\
&&\leq
\frac{2^{21}\morse}{\b^{3/2} \min\{r_0^2,\tilde\checco\}}\left(
\frac{\morse}{\b s_0^3}+
\ln \frac{\b^3 s_0^6}{\tilde\checco \morse^2}\right)
\checco\,.
\label{medina+}
\end{eqnarray}

Arguing as in \eqref{beirut} we get
\begin{eqnarray}\label{beirut+}
\left|\sqrt{E_{2j}- E_{2j-1}}-\sqrt{\bar E_{2j}- \bar E_{2j-1}}\right|
\leq
8\frac{\sqrt\morse}{\b}\checco
\,. 
\end{eqnarray}
Arguing as in \eqref{palmira}
we get
\begin{equation}\label{palmira+}
\left|
\frac{ \sqrt{E_{2j}- E_{2j-1}}}{\sa_{2j}-\sa_{2j-1}}
-
\frac{\sqrt{\bar E_{2j}- \bar E_{2j-1}}}{\bar \sa_{2j}-\bar \sa_{2j-1}}
\right|
\leq
\frac{2^6\morse^{3/2}}{s_0^4\b^2}\checco
\,.
\end{equation}
By \eqref{ladispoli}
\begin{equation}\label{pasargade+}
\left|
\frac{\bar \sa_{2j}-\bar \sa_{2j-1}}{\sqrt{\bar E_{2j}- \bar E_{2j-1}}}
\right|
\leq \frac{2\pi}{\sqrt\b}\,.
\end{equation}
For $E$ as in \eqref{bassora+}
we have  by \eqref{hormuz+},
 \eqref{palmira+} and
 \eqref{suzanne2} (used with  $\tilde\checco\rightsquigarrow \tilde\checco/2$) that
\begin{eqnarray}
&&\left|
\left(\frac{ \sqrt{E_{2j}- E_{2j-1}}}{\sa_{2j}-\sa_{2j-1}}
-
\frac{\sqrt{\bar E_{2j}- \bar E_{2j-1}}}{\bar \sa_{2j}-\bar \sa_{2j-1}}
\right)
\partial_E \act_n^{(2j),-}\big(\l_{2j}^{-1}(\bar\l_{2j}(E))\big)
\right|
\nonumber
\\
&\leq&
\frac{2^{17}\morse^{3/2}}{s_0^4\b^{5/2}}\left(
\frac{\morse}{\b s_0^3}+
\ln \frac{\b^3 s_0^6}{\tilde\checco \morse^2}\right)\checco
\,.
\label{krak+}
\end{eqnarray} 
Then
we have for $E$ as in \eqref{bassora+}
\begin{eqnarray}
& &
\Big|
\partial_E \act_n^{(2j),-}\big(\l_{2j}^{-1}(\bar\l_{2j}(E))\big)
-
\partial_E \bar \act_n^{(2j),-}(E)
\Big|
\nonumber
\\
&\stackrel{\eqref{pasargade+}}\leq&
\frac{2\pi}{\sqrt\b}
\left|\frac{\sqrt{\bar E_{2j}- \bar E_{2j-1}}}{\bar \sa_{2j}-\bar \sa_{2j-1}}
\Big(
\partial_E \act_n^{(2j),-}\big(\l_{2j}^{-1}(\bar\l_{2j}(E))\big)
-
\partial_E \bar \act_n^{(2j),-}(E)
\Big)
\right|
\nonumber
\\
&\stackrel{\eqref{krak+},\eqref{afrodite3++}}\leq&
\frac{2^{17}\morse^{3/2}}{s_0^4\b^{3}}\left(
\frac{\morse}{\b s_0^3}+
\ln \frac{\b^3 s_0^6}{\tilde\checco \morse^2}\right)\checco
+
\checco
\left(
\frac{18}{r_0^2}
+
\frac{2^{193} \morse^{27}}
{ \b^{28}  s_0^{45}s_*^2} 
\right)
\frac{ 2^{55}\morse^7}{\b^{15/2}  s_0^{13}}
\ln\left(
4+\frac{\b}{2\tilde\checco}
\right)
\nonumber
\\
&\stackrel{\eqref{harlock}}\leq&
\checco
\left(
\frac{18}{r_0^2}
+
\frac{2^{193} \morse^{27}}
{ \b^{28}  s_0^{45}s_*^2} 
\right)
\frac{ 2^{56}\morse^7}{\b^{15/2}  s_0^{13}}
\ln\left(
4+\frac{\b}{2\tilde\checco}
\right)\,.
\label{antiochia+}
\end{eqnarray}
Recollecting for $E$ as in \eqref{bassora+} we have
\begin{eqnarray*}
& &\Big|
\partial_E \act_n^{(2j),-}(E)
-
\partial_E \bar \act_n^{(2j),-}(E)
\Big|
\\
&\stackrel{\eqref{medina+}}\leq&
\Big|
\partial_E \act_n^{(2j),-}\big(\l_{2j}^{-1}(\bar\l_{2j}(E))\big)
-
\partial_E \bar \act_n^{(2j),-}(E)
\Big|
\\
& &
+
\frac{2^{21}\morse}{\b^{3/2} \min\{r_0^2,\tilde\checco\}}\left(
\frac{\morse}{\b s_0^3}+
\ln \frac{\b^3 s_0^6}{\tilde\checco \morse^2}\right)
\checco
\\
&\stackrel{\eqref{antiochia+}}\leq&
\checco
\left(
\frac{18}{r_0^2}
+
\frac{2^{193} \morse^{27}}
{ \b^{28}  s_0^{45}s_*^2} 
\right)
\frac{ 2^{56}\morse^7}{\b^{15/2}  s_0^{13}}
\ln\left(
4+\frac{\b}{2\tilde\checco}
\right)
+
\frac{2^{21}\morse}{\b^{3/2} \min\{r_0^2,\tilde\checco\}}\left(
\frac{\morse}{\b s_0^3}+
\ln \frac{\b^3 s_0^6}{\tilde\checco \morse^2}\right)
\checco
\\
&\stackrel{\eqref{harlock}}\leq&
\left[
\left(
\frac{18}{r_0^2}
+
\frac{2^{193} \morse^{27}}
{ \b^{28}  s_0^{45}s_*^2} 
\right)
\frac{ 2^{57}\morse^7}{\b^{15/2}  s_0^{13}}
\ln\left(
4+\frac{\b}{2\tilde\checco}
\right)
+
\frac{2^{21}\morse^2}{\b^{5/2}s_0^3 \tilde\checco}
\right]
\checco
\end{eqnarray*}
proving  \eqref{damietta+}.
\eproof

Let us define the rectangle
\begin{equation}\label{retto}
\mathcal R^{(2j)}(r,\tilde\checco)
:=\{  \bar E^{(2j)}_- +\tilde\checco
< \Re E <  \bar E^{(2j)}_+ -\tilde\checco\,,\ \ 
|\Im E|<r\}\,.
\end{equation}

\begin{corollary}\label{nebraska}
Let $1\leq j< N.$
Let 
\begin{equation}\label{bassora++}
 \frac{2^5\morse}{\b s_*^2}\checco
 \leq \tilde\checco \leq
 \frac{\mathtt r_2\morse}{2^5}\,.
\end{equation}
Then 
$\partial_E \act_n^{(2j)}(E,\hat\act)$
and
$\partial_E \bar \act_n^{(2j)}(E)$
are holomorphic 
for  $E\in\mathcal R^{(2j)}(\mathtt r_4/4,\tilde\checco)$
and $\hat\act\in\hat D_{r_0}$
with estimates
 \begin{equation}
\sup_{\mathcal R^{(2j)}(\mathtt r_4/4,\tilde\checco)
\times\hat D_{r_0}}
\Big|
\partial_E \act_n^{(2j)}(E,\hat\act)
-
\partial_E \bar \act_n^{(2j)}(E)
\Big|
\leq
\left(
\frac{ 2^{75}\morse^{6}}{s_0^{21/2}\b^{6} \mathtt r_2^{3/2}}
+
\frac{\b}{\tilde\checco}
\right)
\frac{ 2^{25}\morse^{5/2}}{\b^{4}  s_0^{9/2}\mathtt r_2}
\checco
\,.
\label{rosetta+}
\end{equation}
Then\footnote{By Cauchy estimates.}
 \begin{equation}
\sup_{\mathcal R^{(2j)}(\mathtt r_4/4,\tilde\checco)
\times\hat D_{r_0}}
\Big|
\partial_{E\hat\act} \act_n^{(2j)}(E,\hat\act)
\Big|
\leq
\left(
\frac{ 2^{75}\morse^{6}}{s_0^{21/2}\b^{6} \mathtt r_2^{3/2}}
+
\frac{\b}{\tilde\checco}
\right)
\frac{ 2^{26}\morse^{5/2}}{\b^{4}  s_0^{9/2}\mathtt r_2 r_0}
\checco
\,.
\label{rosetta+*}
\end{equation}
\end{corollary}
\proof
First note that
\eqref{bassora++} implies
both \eqref{bassora} and   \eqref{bassora+}. 
Moreover,
recalling  \eqref{autunno}, \eqref{retto}, \eqref{ladispoli},  \eqref{tropicana}
and \eqref{harlock}, we have
$$
\mathcal R^{(2j)}(\mathtt r_4/4,\tilde\checco)\ 
\subset\ 
\mathcal E^{2j-1}(\mathtt r_4,\tilde\checco)
\cap \{\Re  E\geq\bar E^{(2j-1)}_- + \mathtt r_2\morse/2^7\}
$$
and, by \eqref{acri} and
\eqref{autunno2},
$$
\mathcal R^{(2j)}(\mathtt r_4/4,\tilde\checco)\ 
\subset\ 
\mathcal E^{2j}(r_0^2\b/2^{10}\morse,\tilde\checco)\,.
$$
Therefore we can apply estimates \eqref{damietta} and 
\eqref{damietta+},
recalling \eqref{citola}, \eqref{noce2}
and the fact that by \eqref{alcafone}
$N\leq 4\sqrt{\morse/\b s_0^3},$
 obtaining
 \begin{eqnarray*}
&&
\sup_{\mathcal R^{(2j)}(\mathtt r_4/4,\tilde\checco)
\times\hat D_{r_0}}
\Big|
\partial_E \act_n^{(2j)}(E,\hat\act)
-
\partial_E \bar \act_n^{(2j)}(E)
\Big|
\nonumber
\\
&&
\qquad
\leq
\left(
\frac{ 2^{78}\morse^{7}}{s_0^{12}\b^{7} \mathtt r_2^{3/2}}
+\frac{s_*^{12}\morse}{\tilde\checco}
\right)
\frac{2^{22}\morse^{1/2}}{\b^2 s_0}
\checco
\\
&&
\qquad\quad
+
\frac{8\morse^{1/2}}{\b^{1/2}s_0^{3/2}}
\left[
\left(
\frac{18}{r_0^2}
+
\frac{2^{193} \morse^{27}}
{ \b^{28}  s_0^{45}s_*^2} 
\right)
\frac{ 2^{57}\morse^7}{\b^{15/2}  s_0^{13}}
\ln\left(
4+\frac{\b}{2\tilde\checco}
\right)
+
\frac{2^{21}\morse^2}{\b^{5/2}s_0^3 \tilde\checco}
\right]
\checco
\\
&&\qquad
\stackrel{\eqref{tropicana},\eqref{badedas}}\leq
\left[
\frac{ 2^{75}\morse^{5}}{s_0^{17/2}\b^{5} \mathtt r_2^{3/2}}
+
\frac{ 2^{16}\morse^{4}}{\b^4  s_0^{10}\mathtt r_2}
\ln \frac{\b}{\tilde\checco}
+
\frac{\b}{\tilde\checco}
\right]
\frac{ 2^{25}\morse^{5/2}}{\b^{4}  s_0^{9/2}\mathtt r_2}
\checco
\,,  
\end{eqnarray*}
where in the last inequality we have used \eqref{harlock} and that
$$
\frac{ 2^{16}\morse^{4}}{\b^4  s_0^{10}\mathtt r_2}
\ln \frac{\b}{\tilde\checco}
\leq
\frac{ 2^{75}\morse^{6}}{s_0^{21/2}\b^{6} \mathtt r_2^{3/2}}
+
\frac{\b}{\tilde\checco}\,.
$$
The last estimates can be proved as 
\eqref{compleannoTommy}
substituting  $a$ and $b$ in \eqref{compleanno Betta} with
$$
a:=\frac{ 2^{16}\morse^{4}}{\b^4  s_0^{10}\mathtt r_2}
\,,\qquad
b:=\frac{ 2^{75}\morse^{6}}{s_0^{21/2}\b^{6} \mathtt r_2^{3/2}}\,.
$$
\eproof

\begin{corollary}\label{alaska}
  For  $i=0,2N,$ and  $\hat\act\in \hat D_{r_0}$ 
 and\footnote{Recalling \eqref{autunno2+}, \eqref{gruffalo},
 \eqref{tracollo} and \eqref{tracollo1} we have
$\mathcal E^{2N}(\mathtt r,\tilde\checco)
=
\mathcal E^{0}(\mathtt r,\tilde\checco)
=
 (\bar E_{2N} ,R_0^2-2\morse)_{\mathtt r}
\cap \{\Re  E>\bar E_{2N} + \tilde \checco\}.$}
 \begin{equation}\label{bassora2N}
   E\in \mathcal E^{i}(r_0^2\b/2^{10}\morse,\tilde\checco)\,,\qquad
 \forall\, \tilde\checco\geq \frac{2^5\morse}{\b s_*^2}\checco\,,
\end{equation}
we have
\begin{equation}\label{rosetta2N}
\Big|
\partial_E \act_n^{(i)}(E,\hat\act)
-
\partial_E \bar \act_n^{(i)}(E)
\Big|
\leq
\left(
\frac{ 2^{56}\morse^{13/2}}{\mathtt r_2 \b^{8}  s_0^{29/2}}
\ln\left(
4+\frac{\b}{2\tilde\checco}
\right)
+
\frac{2^{24}\morse^{5/2}}{\b^{3}s_0^{9/2} \tilde\checco}
\right)
\checco\,.
\end{equation}
Then\footnote{By Cauchy estimates.}
\begin{equation}\label{rosetta2N*}
\Big|
\partial_{E\hat\act} \act_n^{(i)}(E,\hat\act)
\Big|
\leq
\left(
\frac{ 2^{56}\morse^{13/2}}{\mathtt r_2 \b^{8}  s_0^{29/2}}
\ln\left(
4+\frac{\b}{2\tilde\checco}
\right)
+
\frac{2^{24}\morse^{5/2}}{\b^{3}s_0^{9/2} \tilde\checco}
\right)
\frac{2\checco}{r_0}\,.
\end{equation}
\end{corollary}
\proof
By estimate 
\eqref{damietta+},
recalling \eqref{citola}, \eqref{ancomarzio4}
and the fact that by \eqref{alcafone}
$N\leq 4\sqrt{\morse/\b s_0^3},$ we get
\begin{eqnarray*}
&&
\Big|
\partial_E \act_n^{(i)}(E,\hat\act)
-
\partial_E \bar \act_n^{(i)}(E)
\Big|
\\
&\leq&
\frac{8\morse^{1/2}}{\b^{1/2}s_0^{3/2}}
\left[
\left(
\frac{18}{r_0^2}
+
\frac{2^{193} \morse^{27}}
{ \b^{28}  s_0^{45}s_*^2} 
\right)
\frac{ 2^{57}\morse^7}{\b^{15/2}  s_0^{13}}
\ln\left(
4+\frac{\b}{2\tilde\checco}
\right)
+
\frac{2^{21}\morse^2}{\b^{5/2}s_0^3 \tilde\checco}
\right]
\checco
\\
&\stackrel{\eqref{tropicana},\eqref{badedas}}\leq&
\frac{8\morse^{1/2}}{\b^{1/2}s_0^{3/2}}
\left[
\frac{1}{\mathtt r_2}
\frac{ 2^{53}\morse^6}{\b^{15/2}  s_0^{13}}
\ln\left(
4+\frac{\b}{2\tilde\checco}
\right)
+
\frac{2^{21}\morse^2}{\b^{5/2}s_0^3 \tilde\checco}
\right]
\checco\,.
\end{eqnarray*}

\eproof


\subsection{The energy as a function of the actions}

Lemma \ref{lorien} also holds (with slightly modified constants) for the perturbed
action:
\begin{lemma}\label{lorien+}
 We have that 
 \begin{equation}\label{november+}
\partial_E  \act_n^{(2j-1),\pm}(E,\hat\act)
	\geq \frac{s_0}{8\sqrt\morse}\,,
	\qquad 
	\forall\,
 E^{(2j-1)}_-(\hat\act)<E< E^{(2j-1)}_+(\hat\act)\,,
 \ \ \hat\act\in \hat D
\end{equation}
and
\begin{equation}\label{rain+}
\partial_E  \act_n^{(2j),\pm}(E,\hat\act)
\geq \frac{1}{16\sqrt{E- E_{2j-1}}}
\sqrt{\frac{\b s_0^3}{ \morse}}\,,
\qquad
\forall\,
E^{(2j-1)}_-(\hat\act)<E< E^{(2j-1)}_+(\hat\act)\,,
 \ \ \hat\act\in \hat D\,.
\end{equation}
As a consequence
\begin{equation}\label{galadriel+}
\min_{1\leq i\leq 2N-1 } 
\inf_{\hat\act \in \hat D}
\inf_{E\in ( E^{(i)}_-(\hat\act), E^{(i)}_+(\hat\act))}
\partial_E  \act_n^{(i)}(E,\hat\act)
\geq \frac{\sqrt\b s_0^{3/2}}{64\morse}\,.
\end{equation}
Moreover
\begin{equation}\label{galadriel2+}
\partial_E  \act_n^{(2N)}(E,\hat\act)\,,
\
-\partial_E  \act_n^{(0)}(E,\hat\act)
\geq \frac{1}{4\sqrt{E+3\morse/2}}
\,,
\ \ 
\forall\,
 E_{2N}(\hat\act)<E<R_0^2-2\morse\,,\  \hat\act\in\hat D\,.
\end{equation}
\end{lemma}
\proof
First we note that, since
$$
\Gm(\sa_{2j-1}+\sa)-\Gm(\sa_{2j-1})
=
\Gm(\sa_{2j-1}+\sa)- E_{2j-1}
\leq \frac{2\morse}{s_0^2}\sa^2
$$
for every $\sa$, then
\begin{equation}\label{lothlorien+}
\Sa_{2j}(E)-\sa_{2j-1}\,,\ 
\sa_{2j-1}-\Sa_{2j-1}(E)\,\geq
\frac{s_0}{\sqrt{2\morse}}\sqrt{E- E_{2j-1}}\,.
\end{equation}
Therefore\footnote{Recall
 \eqref{lamento2}.}
\begin{eqnarray}
\partial_E \act_n^{(2j-1),+}(E)
&=&
\frac{1}{2\pi}
\int_{\sa_{2j-1}}^{\Sa_{2j}(E )}
\frac1{\sqrt{E-\Gm(\sa)}}\Big(
1
+ \tilde b\big(E-\Gm(\sa),\sa \big)
\Big)\, d\sa\,,
\nonumber
\\
&\stackrel{\eqref{3holesquater}}\geq&
\frac{1}{8}
\int_{\sa_{2j-1}}^{\Sa_{2j}(E )}
\frac{1}{\sqrt{E- E_{2j-1}}}\, d\sa
\geq \frac{s_0}{16\sqrt\morse}
\label{turnit}
\end{eqnarray}
for $ E_{2j-1}(\hat\act)<E< E_{2j}(\hat\act).$
The estimates for 
$\partial_E  \act_n^{(2j-1),-}$
is analogous. 
\\
For $\sa_{2j-1}\leq \sa\leq \sa_{2j}$ 
we have
$$
E-\Gm(\sa)\leq E- E_{2j-1}\,,
$$
then,
recalling \eqref{ancomarzio2} and \eqref{lamento2}, we get
\begin{eqnarray}
\partial_E \act_n^{(2j),-}(E)
&=&
\frac{1}{2\pi}
\int_{\sa_{2j-1}}^{\sa_{2j}}
\frac1{\sqrt{E-\Gm(\sa)}}\Big(
1
+ \tilde b\big(E-\Gm(\sa),\sa \big)
\Big)\, d\sa\,,
\nonumber
\\
&\stackrel{\eqref{3holesquater}}\geq&
\frac{1}{8}
\int_{\sa_{2j-1}}^{\sa_{2j}}
\frac{1}{\sqrt{E- E_{2j-1}}}\, d\sa
\end{eqnarray}
and \eqref{rain+}
follows\footnote{The estimates for 
$\partial_E\act_n^{(2j),+}$
is analogous.}   by \eqref{draghetto} and \eqref{octoberx}.
\\
 \eqref{galadriel+} 
 follows by \eqref{november+},
\eqref{rain+}, \eqref{harlock} and \eqref{noce2}.
\\
Finally \eqref{galadriel2+} directly follows by
 \eqref{lamento3}, \eqref{piove}, \eqref{fragile}, \eqref{ciccio2}.\eproof

By Lemma \ref{lorien+}, \eqref{noce}, \eqref{noce2}
and \eqref{ancomarzio4} we have that,
for $\hat \act\in\hat D$ and $1\leq i\leq 2N,$ the functions
$\act_n^{(i)}(E,\hat\act),-\act_n^{(0)}(E,\hat\act)$
are strictly increasing increasing w.r.t. $E$ and, therefore,
invertible with inverse functions
\beqa{mortazza+}
&&\mathtt E^{(i)}(\hat\act,\cdot): 
\act_n\in  [ \acci^{(i)}_-(\hat\act), \acci^{(i)}_+(\hat\act)]
\mapsto  \mathtt E^{(i)}(\hat\act,\act_n)\in 
[  E^{(i)}_-(\hat\act),  E^{(i)}_+(\hat\act)]
\,,
\qquad
1\leq i<2N\,,
\nonumber
\\
&&\mathtt E^{(2N)}(\hat\act,\cdot): 
\act_n\in   [ \acci^{(2N)}_-(\hat\act),\acci^{(2N)}_+(\hat\act)]\mapsto  
[  E^{(2N)}_-(\hat\act),R_0^2-2\morse]\,,
\nonumber
\\
&&
\mathtt E^{(0)}(\hat\act,\cdot): 
\act_n\in   [\acci^{(0)}_-(\hat\act),\acci^{(0)}_+(\hat\act)]\mapsto  
[  E^{(0)}_-(\hat\act),R_0^2-2\morse]\,,
\eeqa
where
\begin{eqnarray}\label{torpignattara+}
&& \acci^{(i)}_\pm(\hat\act):=  \act_n^{(i)}(  E^{(i)}_\pm(\hat\act),\hat\act)\,,
\qquad {\rm except \  for} 
\\
&&
 \acci^{(2N)}_+(\hat\act):=  \act_n^{(2N)}(R_0^2-2\morse,\hat\act)\,,\qquad
 \acci^{(0)}_-(\hat\act):=  \act_n^{(0)}(R_0^2-2\morse,\hat\act) 
\,.
\end{eqnarray}
Note that actually
\begin{equation}\label{torpignattara2+}
 \acci^{(2j-1)}_-:=0\,,\qquad
\forall\, 1\leq j\leq N\,.
\end{equation}
Set
\begin{equation}\label{bigjim}
\Bu^i :=\Big\{ \act=(\hat \act,\act_n)\ |\ \hat \act\in\hat D,\ \  \acci^{(i)}_-(\hat \act)<\act_n
<\acci^{(i)}_+(\hat \act)
\Big\}\subseteq  \hat D\times \R\subseteq\R^n\,.
\end{equation}
By construction
\begin{equation}\label{mezzemaniche}
 \act_n^{(i)}\big(\mathtt E^{(i)}(\act),\hat\act)\big)=\act_n
 \qquad
 {\rm on}\ \ \ \Bu^i
 \,.
\end{equation}

For $\loge>0$ we define
\beqano
&&
\arr{\acci^{(2j-1)}_- (\hat\act,\loge):=0\\
\acci^{(2j-1)}_+(\hat\act,\loge):=
\act_n^{(2j-1)}\big(E^{(2j-1)}_+(\hat \act)-\loge,\hat \act\big)
}
\phantom{AAA}{\ (1\leq j\leq N)} \,,
\nonumber
\\
&&
\arr{\acci^{(2j)}_-(\hat\act,\lambda):=
\act_n^{(2j)}\big(E^{(2j)}_-(\hat \act)+\loge,\hat \act\big)
\\
\acci^{(2j)}_+(\hat\act,\lambda):=
\act_n^{(2j)}\big(E^{(2j)}_+(\hat \act)-\loge,\hat \act\big)}
\phantom{AAAAAAAa}(1\leq j< N) 
\nonumber
\\
&&
\arr{\acci^{(0)}_-(\hat\act,\lambda):=
\act_n^{(0)}\big(R_{0}^2-2\morse,\hat \act\big)
\\
\acci^{(0)}_+(\hat\act,\lambda):=\act_n^{(0)}\big(E^{(0)}_-(\hat \act)+\loge,\hat \act\big)}
\\
&&
\arr{\acci^{(2N)}_-(\hat\act,\lambda):=\act_n^{(2N)}\big(E^{(2N)}_-(\hat \act)+\loge,\hat \act\big)
\\
\acci^{(2N)}_+(\hat\act,\lambda):=
\act_n^{(2N)}\big(R_{0}^2-2\morse,\hat \act\big)}
\eeqano
and
\begin{equation}\label{playmobilk}
\Bu^i (\loge):=\Big\{ \act=(\hat \act,\act_n)\ |\ \hat \act\in\hat D,\ \  \acci^{(i)}_-(\hat \act,\loge)<\act_n
<\acci^{(i)}_+(\hat \act,\loge)
\Big\}\subseteq  \hat D\times \R\subseteq\R^n \,.
\end{equation}
Note that
$$
\acci^{(i)}_\pm(\hat\act,0)=\acci^{(i)}_\pm(\hat\act)\,,
\qquad
\qquad
\Bu^i (0)=\Bu^i \,.
$$
Set 
\begin{eqnarray}
\mathbb E^{2j-1}(\loge,\hat\act)
&:=&
(E^{(2j-1)}_-(\hat\act),E^{(2j-1)}_+(\hat\act)-\loge)\,,
\qquad
1\leq  j\leq N\,,
\nonumber
\\
\mathbb E^{2j}(\loge,\hat\act)
&:=&
(E^{(2j)}_-(\hat\act)+\loge,E^{(2j)}_+(\hat\act)-\loge)\,,
\qquad
1\leq  j< N\,,
\nonumber
\\
\mathbb E^i(\loge,\hat\act)
&:=&
(E^{(i)}_-(\hat\act)+\loge,R_0^2-2\morse)\,,
\qquad\quad
  i=0,2N\,.
  \label{croccanti}
\end{eqnarray}
Note  that by construction
\begin{equation}\label{entangled}
\forall\, I^*=(\hat\act^*,\act_n^*)\in \Bu^i (\loge) \qquad \exists\,!\ E^*\in 
\mathbb E^i(\loge,\hat\act^*) \ \ 
{\rm s.t.}\ \ 
\act_n^{(i)}(E^*,\hat\act^*)=\act_n^*\,.
\end{equation}
Define the following domains of $\C^n$:
\begin{equation}\label{texas}
\mathcal D^i(\loge,r_0):=
\{
(E,\hat\act)\ \ E\in \mathbb E^i_\loge(\loge,\hat\act^*)\,,
 |\hat\act-\hat\act^*|<r_0\ \ 
{\rm with}
 \ \hat\act^*\in \hat D
\}\,.
\end{equation}

\begin{lemma}
We have, 
for\footnote{$\mathtt r_4$ defined in \eqref{acri}.} 
$0<\loge\leq  \mathtt r_4/2$
\begin{eqnarray}
\sup_{\mathcal D^i(\loge,r_0)}
|\partial_{E}\act_n^{(i)}(E,\hat\act)|
&\leq&
c_{1,0}\ln\frac{\morse}{\loge}\,,
\nonumber
\\
\sup_{\mathcal D^i(\loge,r_0)}
|\partial_{EE}\act_n^{(i)}(E,\hat\act)|
&\leq&
c_{2,0}\frac{1}{\loge}\,,
\nonumber
\\
\sup_{\mathcal D^i(\loge,r_0/2)}
|\partial_{\hat\act}\act_n^{(i)}(E,\hat\act)|
&\leq&
c_{0,1}\checco\,,
\nonumber
\\
\sup_{\mathcal D^i(\loge,r_0/2)}
|\partial_{E\hat\act}\act_n^{(i)}(E,\hat\act)|
&\leq&
c_{1,1}\checco\ln\frac{\morse}{\loge}\,,
\nonumber
\\
\sup_{\mathcal D^i(\loge,r_0/2)}
|\partial_{\hat\act\hat\act}\act_n^{(i)}(E,\hat\act)|
&\leq&
c_{0,2}\checco\,,
\label{spoleto}
\end{eqnarray}
for suitable constants $c_{i,j}$.
\end{lemma}
\proof
We omit the details.
\eproof

\begin{lemma}\label{florida}
 Let $0<\loge\leq  \mathtt r_4/2$ and set
 \begin{equation}\label{sibari}
 \rho=\rho(\loge):=
 \left\{\begin{array}{l}
\displaystyle \min\left\{\frac{\b s_0^{3} \loge }{2^{16}c_{2,0}(c_{0,1}+1)\morse^2},\frac{r_0}{2}\right\}\ \ \ \ 
\qquad\ \quad{\rm if}	 \ \ \ 
\loge_0
\leq \loge\leq \frac{\morse}{2}
\,,
\\ 
\displaystyle
\min\left\{\frac{\sqrt\b s_0^{5/2} \loge }{2^{16}c_{2,0}(c_{0,1}+1)\morse^{3/2}}
\ln\frac{\morse}{\loge},\frac{r_0}{2}\right\}\ \ \ \
\ {\rm if}\ \ \  
0<\loge\leq \loge_0
\,,
\end{array}
\right.
\end{equation}
where
$$
\loge_0:=
\morse\exp\left(-\frac{2^{97}\morse^{19/2}}{s_0^{15}\b^{19/2}}\right)\,.
$$
For every $0\leq i\leq 2N,$ the function 
 $\mathtt E^{(i)}(\act)$ has analytic extension on the complex 
 $\rho$-neighborhood  
 of $\Bu^i (\loge),$ namely
 $\Bu^i_{\rho} (\loge).$ 
 Finally
 \begin{equation}\label{alabama}
\Xi \big(\Bu^i_{\rho} (\loge)\big)\subset
 \mathcal D^i(\loge,\rho)\,,
 \qquad
 \text{where}\qquad
 \Xi(I):=\big(\mathtt E^{(i)}(\act),\hat\act\big)\,.
\end{equation}
\end{lemma}
\proof
We consider first the case $1\leq i<2N.$
Set
$$
\mathcal F(E,I):=
\act_n^{(i)}(E,\hat\act)-\act_n\,.
$$
Fix $I^*\in \Bu^i (\loge).$ By \eqref{entangled} there exist unique
$E^*=\mathtt E^{(i)}(\act^*)\in 
\mathbb E^i(\loge,\hat\act^*)$ 
such that 
$\mathcal F(E^*,\act^*)=0$. We want to apply the Implicit Function 
Theorem to $\mathcal F$  in order to find a holomorphic function
\begin{equation}\label{alabama2}
\mathtt E^{(i)}: \{|I-I^*|\leq \rho\}\ \to\ \{|E-E^*|\leq c\loge\}\,,
\qquad c	\leq 1/2\,, 
\end{equation} 
($c$ to be chosen below see \eqref{argo})
such that
$$
\mathcal F (\mathtt E^{(i)}(\act),\act)=0\,.
$$
This is possible since
$$
\partial_E \mathcal F(E^*,\act^*)=
\partial_{E}\act_n^{(i)}(E^*,\hat\act^*)\neq 0
$$
by Lemma \ref{lorien+}.
By a quantitative version of the Implicit Function 
Theorem we have to check that
\begin{eqnarray}\label{clitunno1}
\sup_{|I-I^*|\leq \rho}|\act_n^{(i)}(E^*,\hat\act)-I_n|
&\leq& 	
\frac{c\loge}{2}\d_*\,,
\\
\sup_{|E-E^*|\leq c\loge\,,\    |I-I^*|\leq \rho}
|
\partial_{E}\act_n^{(i)}(E,\hat\act)-\partial_{E}\act_n^{(i)}(E^*,\hat\act^*)
|
&\leq&
\frac{\d_*}2\,,
\label{clitunno2}
\end{eqnarray}
where
$$
\d_*:=|\partial_{E}\act_n^{(i)}(E^*,\hat\act^*)|\,.
$$
Recalling \eqref{spoleto} we have that
\begin{eqnarray*}
|\act_n^{(i)}(E^*,\hat\act)-I_n|
&\leq&
|\act_n^{(i)}(E^*,\hat\act)-\act_n^{(i)}(E^*,\hat\act^*)|
+
|I_n^*-I_n|
\\
&\leq&
(c_{0,1}+1)\rho
\leq \frac{c\loge}{2}\d_*
\end{eqnarray*}
by \eqref{sibari}.
Then, taking
\begin{equation}\label{micene}
\rho\leq \frac{c\loge \d_*}{2(c_{0,1}+1)}
\end{equation}
we satisfy
\eqref{clitunno1}.
Moreover by \eqref{spoleto} and since $c\leq 1/2$ we have
\begin{eqnarray*}
|
\partial_{E}\act_n^{(i)}(E,\hat\act)-\partial_{E}\act_n^{(i)}(E^*,\hat\act^*)
|
&\leq&
\frac{2c_{2,0}}{\loge}c\loge+
c_{1,1} \checco(\ln\frac{2\morse}{\loge} )\rho
\\
&\leq&
2c c_{2,0}
+
c_{1,1} \checco(\ln\frac{2\morse}{\loge} )
\frac{c\loge \d_*}{2(c_{0,1}+1)}
\leq
\frac{\d_*}{2}
\end{eqnarray*}
if
$$
c\leq \frac{\d_*}{8c_{2,0}}
\qquad
{\rm and}
\qquad
c\leq \frac{c_{0,1}+1}{4c_{1,1} \checco \loge \ln(2\morse/\loge)}\,.
$$
The second restriction on $c$ is always satisfied 
for $\checco$ small enough by \eqref{genesis}.
The first  restriction is satisfied if we take\footnote{Even if we did not evaluate
$c_{2,0}$, it is very large!}
\begin{equation}\label{argo}
c\leq\min\left\{
\frac12\,,\, 
\frac{\d_*}{8c_{2,0}}
\right\}\,,
\quad \text{which is implied choosing}\quad
c:= 
\frac{\sqrt\b s_0^{3/2}}{2^9c_{2,0}\morse}
\,,
\end{equation}
since
\begin{equation}\label{tirinto}
\d_*\geq \frac{\sqrt\b s_0^{3/2}}{64\morse}
\end{equation}
by
\eqref{galadriel+}.

On the other hand for $\loge$ small enough we could always take 
$c:=1/2$ (but for simplicity we take $c$ as in \eqref{argo} in any case).
Indeed,
recalling lemmata \ref{giga2} and \ref{giga2+}, 
if
$$
\loge\leq \morse 
\exp\left(-\frac{2^{97}\morse^{19/2}}{s_0^{15}\b^{19/2}}\right)\,,
$$
we have 
the better (w.r.t. \eqref{tirinto}) estimates from below\footnote{Recall in particular
\eqref{maracaibo}, \eqref{trocadero}, \eqref{otello}, \eqref{badedas}.} 
$$
\d_*\geq \frac{s_0}{2^6\sqrt\morse}\ln\frac{\morse}{\loge}
\geq
\frac{2^{91}\morse^{9}}{s_0^{14}\b^{19/2}}
\,.
$$
Then, taking $\rho$ as in \eqref{sibari} and recalling \eqref{tirinto}
we have that \eqref{micene} is satisfied.
\\
Finally \eqref{alabama} follows by \eqref{alabama2}.

We omit the details of the case $i=0,2N.$ 
\eproof

\subsection{Derivatives of the energy}

We have 
\begin{eqnarray}
&&\partial_{\act_n} \mathtt E^{(i)} 
=
\frac{1}{\partial_E \act_n^{(i)}}\,,\qquad
\partial_{\hat \act} \mathtt E^{(i)} 
=
-\frac{\partial_{\hat \act} \act_n^{(i)}}{\partial_E \act_n^{(i)}}\,,
\qquad
\partial_{\act_n \act_n} \mathtt E^{(i)} 
=
-\frac{\partial_{EE} \act_n^{(i)}}{(\partial_E \act_n^{(i)})^3}\,,
\nonumber
\\
&&\partial_{\act_n \hat \act} \mathtt E^{(i)} 
=
\frac{\partial_{EE} \act_n^{(i)} \partial_{\hat \act}\act_n^{(i)}}{(\partial_E \act_n^{(i)})^3}-
\frac{\partial_{E\hat \act} \act_n^{(i)} }{(\partial_E \act_n^{(i)})^2}
\,,
\nonumber
\\
&&\partial_{\hat \act \hat \act} \mathtt E^{(i)} 
=
-\frac{\partial_{\hat \act\hat \act} \act_n^{(i)}}{\partial_E \act_n^{(i)}}
+\frac{2\partial_{E\hat \act} \act_n^{(i)}\partial_{\hat \act}\act_n^{(i)}}{(\partial_E \act_n^{(i)})^2}
-\frac{\partial_{EE} \act_n^{(i)} (\partial_{\hat \act}\act_n^{(i)})^2}{(\partial_E \act_n^{(i)})^3}\,,
\label{daitarn3}
\end{eqnarray}
where $\mathtt E^{(i)}$ and  $\act_n^{(i)}$ are evaluated in  $\big(\act_n^{(i)}(E,\hat \act),\hat \act\big)$
and $(E,\hat \act),$ respectively\footnote{
Or, which is equivalent, in $\act$ and
$\big(\mathtt E^{(i)}(\act),\hat \act\big),$
respectively.}.

\subsection{The cosine case}

Let us introduce the parameter $\ch>0$.
Let us consider the case in which the unperturbed potential is $(-\ch)$-cosine, namely 
$$
\FO(\sa)= -\ch\cos\sa\,.
$$ 
Given $s_0>0,$
we obviously have that this function is 
$(\morse,\b,s_0)$-Morse-non-degenerate
according to Definition \ref{morso}
with, recalling \eqref{lontra},
\begin{equation}\label{nibelungo}
\morse=\ch\cosh s_0\,,\qquad 
\b=\ch\,.
\end{equation}
We have only two critical points, namely $N=1$: the minimum
$\bar  \sa_1=0$ and the maximum  $\bar  \sa_2=\pi,$ with corresponding critical energies $\bar  E_1=-\ch$ and 
$\bar  E_2=\ch,$ respectively.
The functions $\bar \Sa_i,$ $i=1,2$ (defined in \eqref{andria})
are $\bar\Sa_1(E)=-\arccos(-E/\ch)$ and  $\bar\Sa_2(E)=\arccos(-E/\ch).$

{\sl We will put the apex $\empty^\star$ (instead of the bar $\ \bar{\empty}$)  to mean that we are consider exactly the $(-\ch)$-cosine case.}

Set
$$
\bm E :=E/\ch\,.
$$
Then recalling \eqref{sunday} we have that the action variable
corresponding to the $-\ch\cos$ potential is
\begin{equation}\label{eralogico2}
\act_n^{(1),\star}(E)=\frac{2\sqrt\ch}{\pi}\int_0^{\arccos(-\bm E)}\sqrt{\bm E+\cos \sa}d\sa\,,
\end{equation}
with inverse function  (recall \eqref{mortazza})
$$
\mathtt E^{(1),\star}:(0, 4\sqrt {2\ch}/\pi)\ \to\ (-\ch,\ch)\,.
$$ 
The derivative of $\act_n^{(1),\star}$ is (recall also \eqref{cippone})
\begin{eqnarray*}
\partial_E \act_n^{(1),\star}(E) &=&
\frac{1}{\pi\sqrt\ch}\int_0^{\arccos(-\bm E)}\frac{d\sa}{\sqrt{\bm E+\cos \sa}}
=\frac{2}{\pi\sqrt\ch}\int_0^1 \frac{dy}{\sqrt{1-\bm E+2\bm E y^2 - ( \bm E+1)y^4}}
\\
&=&
\frac{2}{\pi\sqrt\ch}\int_0^1 \frac{dy}{\sqrt{1-y^4 -\bm E(1-y^2)^2}}
\ > 0\,,
\end{eqnarray*}
making the change of variables  $\sa=\arccos\big(( \bm E+1)y^2 -\bm E\big).$
Moreover
$$
\partial_{EE} \act_n^{(1),\star}(E)
=\frac{1}{\pi\ch^{3/2}}\int_0^1 \frac{(1-y^2)^2}{\big(1-y^4 -\bm E(1-y^2)^2\big)^{3/2}}
dy\ >0\,.
$$
In particular $\partial_{E} \act_n^{(1),\star}(E)$ is an increasing function.
Analogously we have that
\begin{equation}\label{thatsall}
\partial_{EEE} \act_n^{(1),\star}(E)
=\frac{3}{2\pi\ch^{5/2}}\int_0^1 \frac{(1-y^2)^4}{\big(1-y^4 -\bm E(1-y^2)^2\big)^{5/2}}
dy\ >0\,.
\end{equation}
Similarly all the derivatives of $\act_n^{(1),\star}$ are positive functions.
\\
Note that (by the Lebesgue's theorem)
\begin{equation}\label{settenani}
\lim_{\bm E\to -1^+}\partial_{E} \act_n^{(1),\star}(E)=\frac{1}{\sqrt{2\ch}}\,,\qquad
\lim_{\bm E\to -1^+}\partial_{EE} \act_n^{(1),\star}(E)=\frac{1}{8\sqrt 2\ch^{3/2}}\,.
\end{equation}
The fact that 
 $c_{-\cos}$ (defined in  \eqref{vana0})
 is equal to the first limit in \eqref{settenani}, namely
\begin{equation}\label{senzafonte}
\inf_{-1<\bm E<1}
\partial_{E} \act_n^{(1),\star}(E)=\frac{1}{\sqrt{2\ch}}\,,
\end{equation}
  follows since 
 $\partial_{E} \act_n^{(1),\star}(E)$ is increasing.
By direct calculation (or general arguments, recall\eqref{maracaibo}),
$$
\lim_{\bm E\to 1^-}\partial_{EE} \act_n^{(1),\star}(E)=+\infty\,.
$$
Since $\partial_{EE} \act_n^{(1),\star}>0$, by
\eqref{thatsall} and \eqref{settenani} we have
\begin{equation}\label{storione}
\inf_{-1<\bm E<1} \partial_{EE} \act_n^{(1),\star}(E)\geq 
\frac{1}{8\sqrt 2\ch^{3/2}}>0\,.
\end{equation}
We also have, recalling \eqref{daitarn3},
\begin{equation}\label{pisolo}
\partial_{\act_n} \mathtt E^{(1),\star} (\act_n)=
\frac{1}{\partial_{E} \act_n^{(1),\star}\big(\mathtt E^{(1),\star} (\act_n)\big)}>0\,,\qquad
\partial_{\act_n \act_n} \mathtt E^{(1),\star} (\act_n)=
-\frac{\partial_{EE} \act_n^{(1),\star}\big(\mathtt E^{(1),\star} (\act_n)\big)}
{\Big(\partial_{E} \act_n^{(1),\star}\big(\mathtt E^{(1),\star} (\act_n)\big)\Big)^3}<0
\end{equation}
and, by a direct calculation\footnote{The function is increasing.}
and \eqref{settenani},
\begin{equation}\label{biancaneve}
\inf_{0<\act_n<4\sqrt{2\ch}/\pi}
-\partial_{\act_n \act_n} \mathtt E^{(1),\star} (\act_n)
=
\lim_{\bm E\to -1^+} \frac{\partial_{EE} \act_n^{(1),\star}(E)}
{\big( \partial_{E} \act_n^{(1),\star}(E)\big)^3}=
\frac14\,.
\end{equation}

\medskip
Recalling \eqref{sunday} we get
\begin{equation}\label{eralogico}
\act_n^{(2),\star}(E)=\frac{\sqrt\ch}{\pi}\int_0^{\pi}\sqrt{\bm E+\cos\sa}d\sa=
-\act_n^{(0),\star}(E)\,,
\end{equation}
with inverse functions 
$$
\mathtt E^{(2),\star}:(2\sqrt{2\ch}/\pi,+\infty)\ \to\ (\ch,+\infty)\,,\qquad
\mathtt E^{(0),\star}:(-\infty,-2\sqrt{2\ch}/\pi)\ \to\ (\ch,+\infty)\,.
$$ 
The derivatives of $\act_n^{(2),\star}$ are
\begin{equation}\label{ducabianco}
\partial_E \act_n^{(2),\star}(E)=
\frac{1}{2\pi\sqrt\ch}\int_0^{\pi}\frac{d\sa}{\sqrt{\bm E+\cos \sa}}
> 0\,,\qquad
\partial_{EE} \act_n^{(2),\star}(E)=-
\frac{1}{4\pi\ch^{3/2}}\int_0^{\pi}\frac{d\sa}{(\bm E+\cos \sa)^{3/2}}
< 0\,.
\end{equation}
Note also that $\partial_{EEE} \act_n^{(2),\star}>0.$
We have
\begin{equation}\label{pauperes}
 E\geq 2\ch\qquad \Longrightarrow\qquad
\partial_E \act_n^{(2),\star}(E)\leq \frac{1}{\sqrt{2 E}}\,,\qquad
-\partial_{EE} \act_n^{(2),\star}(E)\geq 
\frac{1}{8\sqrt 2  E^{3/2}}\,.
\end{equation}
We get
\begin{equation}\label{pisolo2}
\partial_{\act_n} \mathtt E^{(2),\star} (\act_n)=
\frac{1}{\partial_{E} \act_n^{(2),\star}\big(\mathtt E^{(2),\star} (\act_n)\big)}>0\,,\qquad
\partial_{\act_n \act_n} \mathtt E^{(2),\star} (\act_n)=
-\frac{\partial_{EE} \act_n^{(2),\star}\big(\mathtt E^{(2),\star} (\act_n)\big)}
{\Big(\partial_{E} \act_n^{(2),\star}\big(\mathtt E^{(2),\star} (\act_n)\big)\Big)^3}>0
\end{equation}
and
$$
\inf_{\act_n>2\sqrt{2\ch}/\pi}\partial_{\act_n \act_n} \mathtt E^{(2),\star} (\act_n)>0\,,
$$
since, by direct calculation, 
$$
\lim_{\bm E\to 1^+} - \frac{\partial_{EE} \act_n^{(2),\star}(E)}
{\big( \partial_{E} \act_n^{(2),\star}(E)\big)^3}=+\infty
$$
and by \eqref{ducabianco}
$$
\lim_{\bm E\to +\infty} - \frac{\partial_{EE} \act_n^{(2),\star}(E)}
{\big( \partial_{E} \act_n^{(2),\star}(E)\big)^3}=2\,.
$$
By the previous limits, \eqref{pisolo2} and \eqref{eralogico} we get 
\begin{equation}\label{biancaneve2}
\inf_{\act_n>4\sqrt 2/\pi}\left|\partial_{\act_n \act_n} \mathtt E^{(2),\star} (\act_n)\right|
=\inf_{\act_n<-4\sqrt{2\ch}/\pi}\left|\partial_{\act_n \act_n} \mathtt E^{(0),\star} (\act_n)\right|
	\geq c_{\star\star}>0\,,
\end{equation}
for a suitable $c_{\star\star}>0$
(that can be explicitly evaluated!).




\subsection{The cosine-like case}

\nl
An important class of Morse non--degenerate functions, as we will shortly show,  is the following.

\begin{lemma}\label{sibillini}
Let $s_0,\ch>0.$
 If $G$ satisfies
  \begin{equation*}
 {\modulo}G(\sa)+\ch\cos \sa {\modulo}_{s_0}\leq \cgot\ch\,,
 \qquad \text{for some}\qquad 
 0<\cgot\leq \frac{1}{4}\min\{1,s_0^2\}\,,
\end{equation*}
 then it is 
 $(\morse,\b,s_0)$--Morse--non--degenerate\footnote{According to Definition
 \ref{morso}.}
 with
 $$
 \b=\ch/4\,,  \qquad
 \morse=\ch(\cgot+\cosh s_0)\leq \ch(\frac14+\cosh s_0)\,.
 $$
 Moreover $G$
 has only two non--degenerate critical points
(a maximum and a minimum).
\end{lemma}
\proof
We have, by Cauchy estimates,
$$
\frac{1}{\ch}( |G'(\sa)|+|G''(\sa)|)
 \geq |\sin \sa| +|\cos \sa| 
 -\frac{\cgot}{s_0}-2\frac{\cgot}{s_0^2}
 \geq 1-\frac{\cgot}{s_0}-2\frac{\cgot}{s_0^2}\geq \frac14\,.
$$
We can choose $\morse$ as above since
${\modulo}\cos \sa{\modulo}_{s_0}=\cosh s_0.$
Regarding the last sentence
we note that for $\sa\in(-\pi,\pi]$
we have only two critical points,
a minimum in $(-\pi/6,\pi/6)$
and a maximum in $(-\pi,-5\pi/6)\cup
(5\pi/6,\pi].$ Indeed we have that, 
setting $g(\sa):=\ch^{-1}G(\sa)+\cos \sa$,
$\ch^{-1}G'(\sa)=\sin \sa+g'(\sa),$ so that
\begin{equation}\label{coratella}
\ch^{-1}G'(\sa)=\sin \sa+g'(\sa)
\geq \sin \sa-\cgot/s_0\geq \sin \sa -1/4\,.
\end{equation}
This implies that 
$G'(\pi/6)\geq \ch/4,$ 
$G'(-\pi/6)\leq -\ch/4.$
Then, by continuity, there exists a critical
point of $G$ in $(-\pi/6,\pi/6).$
Moreover such point is a minimum and there are no other critical points in $(-\pi/6,\pi/6)$ since there $G$
is strictly convex:
$$
\ch^{-1}G''(\sa)=\cos \sa+g''(\sa)
\geq \sqrt 3/2-2\cgot/s_0^2\geq \sqrt 3/2 -1/2
>0\,.
$$
Similarly in $(-\pi,-5\pi/6)\cup
(5\pi/6,\pi]$ there is only one critical point,
which is a maximum.
Finally, by \eqref{coratella},
$G'(\sa)\geq \ch/4$ for $\sa\in[\pi/6,5\pi/6]$
and, analogously, 
$G'(\sa)\leq -\ch/4$ for 
$\sa\in[-5\pi/6,-\pi/6];$
so that there are no other critical points.
\qed

We set\footnote{
This is exactly the value of  $\mathtt r_2$ in \eqref{tropicana}
in the case of a
$(2\ch  e^{s_0}, \ch/4,s_0)$-Morse-non-degenerate
function.}

\begin{equation}\label{varsavia}
\mathtt r_2^\star=\mathtt r_2^\star(\ch,s_0,r_0):=
\min\left\{
\frac{s_0^{49}}{2^{304} e^{30s_0}},\
\frac{r_0^2}{2^{11}\ch  e^{s_0}}\right\}
\end{equation}
and
\begin{equation}\label{stoccolma}
\checco_\diamond^\star=\checco_\diamond^\star(\ch,s_0,r_0):=
\frac{s_0^{27/2}(\mathtt r_2^\star)^4\ch}{2^{135}s_*^{12}e^{6s_0}}\,.
\end{equation}
Note that
\begin{equation}\label{oslo}
\checco_\diamond^\star(\ch,s_0,r_0)
\leq \frac{1}{16}
\checco_\diamond(\ch e^{s_0},\ch,s_0,r_0)
=
\frac{  s_0^{15}}{ 2^{124}e^{9s_0}}
\min
\left\{
r_0^2\,, \ \frac{r_0^3}{\sqrt\ch e^{s_0/2}}\,,
\ 
\frac{\ch s_0^{75}}{2^{321}e^{44 s_0}}
\right\}
\,.
\end{equation}

\begin{definition}\label{sorellastre}
 Given $r_0,s_0,\ch>0,$ we say that a holomorphic function $\Gm:\hat D_{r_0}\times \T_{s_0}\to\C$  is
$(-\ch)$-cosine--like
if\footnote{Note that $\cosh s_0\leq e^{s_0}.$}
 \begin{equation}
 {\modulo}\Gm(\sa,\hat\act)+\ch\cos \sa {\modulo}_{\hat D,r_0,s_0} \leq 
 \checco_\diamond^\star(\ch,s_0,r_0)
\leq
 \frac1{16}\checco_\diamond(2\ch e^{s_0},\ch/4,s_0,r_0)
\,.
 \label{psico}
\end{equation}
\end{definition}

\begin{proposition}\label{squalo}
Assume that
$\Gm$ is
 $(-\ch)$-cosine-like\footnote{ According to Definition
\ref{sorellastre}.}.
Then for every $\hat\act\in\hat D_{r_0}$
the function $\Gm(\sa,\hat\act)$ is 
$(2\ch  e^{s_0}, \ch/4,s_0)$-Morse-non-degenerate\footnote{
According to Definition \ref{morso}.}
and, therefore, by \eqref{galadriel+}
\begin{equation}\label{vanacoseno}
\inf_{\hat\act\in \hat D}\inf_{E_1(\hat\act)<E<E_2(\hat\act)} \partial_E \act_n^{(1)}(E,\hat\act)  \geq
\frac{ s_0^{3/2}}{2^8e^{s_0}\sqrt\ch}
\,.
\end{equation}
Moreover
\begin{equation}\label{cardellino}
\inf_{\hat\act\in \hat D}\inf_{E_1(\hat\act)<E<E_2(\hat\act)} \frac{\partial_{EE} \act_n^{(1)}(E,\hat\act)}{\Big(\partial_E \act_n^{(1)}(E,\hat\act)\Big)^3}
\, \geq \,
\frac{1}{16}\,,
\qquad
\inf_{\hat\act\in \hat D}\inf_{E_2(\hat\act)<E<R_0^2-2\morse}
\frac{-\partial_{EE} \act_n^{(2)}(E,\hat\act)}{\Big(\partial_E \act_n^{(2)}(E,\hat\act)\Big)^3}
\, \geq \,
2\,,
\end{equation}
with $c_{\star\star}$ defined in \eqref{biancaneve2}.
\end{proposition}

\begin{remark}
 Imposing a stronger condition in \eqref{psico} 
 and using Lemma \ref{giga2}
 (in particular \eqref{maracaibo} and \eqref{otello})
 we can prove that \eqref{vanacoseno} holds with
 $1/4\sqrt\ch$
 on the right hand side.
\end{remark}

\proof
Set
\begin{equation}\label{mosca}
\checco^\star:={\modulo} \Gm+\ch\cos \sa^\star{\modulo}_{\hat D,r_0,s_0}
\stackrel{\eqref{psico}}\leq 
\checco_\diamond^\star(\ch,s_0,r_0)
\leq
 \frac1{16}\checco_\diamond(2\ch e^{s_0},\ch/4,s_0,r_0)\,.
\end{equation}
Then
the Morse-non-degeneracy of $\Gm$ follows by Lemma \ref{sibillini}
(with $G\rightsquigarrow\Gm$).

Let us now prove the first estimate in
\eqref{vanacoseno}.

We note that by Lemma \ref{bratislava}
for
$$
E_2(\hat\act)
  -\frac{s_0^4\ch(\mathtt r_2^\star)^2}{2^{37}e^{s_0}}\leq E<
  E_2(\hat\act)
$$
we have
\begin{equation}\label{tallin}
\frac{\partial_{EE} \act_n^{(1)}(E,\hat\act)}{\Big(\partial_E \act_n^{(1)}(E,\hat\act)\Big)^3}
\geq 
\frac{s_0^2}{2^{20}  e^{2 s_0}\mathtt r_2^\star}
\stackrel{\eqref{varsavia}}\geq 1\,.
\end{equation}

Let us now consider the case
\begin{equation}\label{cracovia}
E_1(\hat\act)<
E<
E_2(\hat\act)
 -\frac{s_0^4\ch(\mathtt r_2^\star)^2}{2^{37}e^{s_0}}
 \stackrel{\eqref{stoccolma}}\leq
 E_2(\hat\act)
 -2^{15}e^{s_0} \checco_\diamond^\star
  \end{equation}
By \eqref{mosca} and \eqref{october}
we have that\footnote{Note that in this case
$\bar E_1=-\ch,$ $\bar E_2=\ch.$}
$$
|E_1(\hat\act)+\ch|\,,\ \ 
|E_2(\hat\act)-\ch|\leq 2 \checco^\star\leq 2 \checco_\diamond^\star\,.
$$
Then we have that, if $E$ satisfies \eqref{cracovia} then it also satisfies
\begin{equation}\label{cracoviabis}
-\ch-4 \checco_\diamond^\star<
E
 <\ch - 2^{14}e^{s_0} \checco_\diamond^\star\,.
  \end{equation}
   Recalling \eqref{acri} we set
$$
\mathtt r_4^\star
:=
\frac{\mathtt r_2^\star \ch e^{s_0}}{2^4}
\,,\qquad
\tilde\checco^\star:=\frac{s_0^4\ch(\mathtt r_2^\star)^2}{2^{38}e^{s_0}}\geq
2^{14} e^{s_0} \checco_\diamond^\star\,.
$$
Since\footnote{This is exactly condition \eqref{bassora} for the present case.}
 \begin{equation*}
2^8 e^{s_0}\checco^\star
  \leq
  \tilde\checco^\star
  \leq \mathtt r_4^\star/4
  \,,
\end{equation*}
we can apply Corollary \ref{sarabanda3}
obtaining\footnote{$\act_n^{(2j-1),\star}$ defined in \eqref{eralogico2}.} 
 \begin{eqnarray}
&&
\sup_{\mathcal E^{2j-1}(\mathtt r_4^\star/4,\tilde\checco^\star)
\times\hat D_{r_0}}
\Big|
\partial_E \act_n^{(2j-1)}(E,\hat\act)
-
\partial_E  \act_n^{(2j-1),\star}(E)
\Big|
\nonumber
\\
&&
\qquad
\leq
\left(
\frac{ 2^{99}e^{7 s_0}}{s_0^{12} (\mathtt r_2^\star)^{3/2}}
+\frac{2s_*^{12}e^{s_0}\ch}{\tilde\checco^\star}
\right)
\frac{2^{25}e^{s_0}}{ s_0\ch^{3/2}}
\checco^\star
\nonumber
\\
&&
\qquad
=
\left(
\frac{ 2^{99}e^{7 s_0}}{s_0^{12} (\mathtt r_2^\star)^{3/2}}
+\frac{s_*^{12}e^{2s_0}
2^{39}}{s_0^4(\mathtt r_2^\star)^2}
\right)
\frac{2^{25}e^{s_0}}{ s_0\ch^{3/2}}
\checco^\star
\nonumber
\\
&&
\qquad
\leq
\frac{2^{65}s_*^{12}e^{2s_0}}{s_0^5(\mathtt r_2^\star)^2\ch^{3/2}}
\checco^\star
\,.
\label{atene}
\end{eqnarray}
By Cauchy estimate 
we get (noting that $\mathtt r_4^\star/8\geq\tilde\checco^\star/2$)
 \begin{equation}
\sup_{\mathcal E^{2j-1}(\mathtt r_4^\star/8,\tilde\checco^\star/2)
\times\hat D_{r_0}}
\Big|
\partial_{EE} \act_n^{(2j-1)}(E,\hat\act)
-
\partial_{EE}  \act_n^{(2j-1),\star}(E)
\Big|
\leq
\frac{2^{66}s_*^{12}e^{2s_0}}{s_0^5(\mathtt r_2^\star)^2\ch^{3/2}\tilde\checco^\star}
\checco^\star
\,.
\label{atene2}
\end{equation}
Now we note that if $E$ satisfies \eqref{cracovia}, then it also satisfies 
\eqref{cracoviabis} and, therefore,
$$
E\in \mathcal E^{2j-1}(\mathtt r_4^\star/8,\tilde\checco^\star/2)\,,
$$
then, for $\hat\act\in\hat D,$ by \eqref{senzafonte}, \eqref{storione}, \eqref{vanacoseno}, \eqref{atene}, \eqref{atene2}
we have\footnote{For $a,a^\star,b,b^\star>0$
we have
$$
\frac{a}{b^3}\geq \frac{a^\star}{(b^\star)^3}
\left(1
-48 \frac{|b-b^\star|}{b^\star}
\right)
-\frac{|a-a^\star|}{b^3}
\,.
$$
Use it with $a=\partial_{EE} \act_n^{(1)}$,
$a^\star=\partial_{EE} \act_n^{(1),\star},$
$b=\partial_E \act_n^{(1)},$
$b^\star=\partial_E \act_n^{(1),\star}.$
}
\begin{eqnarray*}
\frac{\partial_{EE} \act_n^{(1)}}{\Big(\partial_E \act_n^{(1)}\Big)^3}
\geq
\frac14
\left(
1-
\frac{2^{72}s_*^{12}e^{2s_0}}{s_0^5(\mathtt r_2^\star)^2\ch}
\checco^\star
\right)
-
\frac{2^{90}s_*^{12}e^{5s_0}}{s_0^{19/2}(\mathtt r_2^\star)^2\tilde\checco^\star}
\checco^\star
\geq\frac{1}{16}\,,
\end{eqnarray*}
where the last inequality follows by \eqref{mosca} and \eqref{stoccolma}.
Recalling \eqref{tallin} this conclude the proof of the first estimate in
\eqref{cardellino}.

We omit the proof of the second inequality in \eqref{cardellino}.
\eproof



\section{The theorem on the universal  analytic behavior  of actions}

Let $r_0, s_0,\morse,\b,\ch>0$ and $\hat D$ a domain of $\R^{n-1}.$
Consider a real analytic function on $\T_{s_0}\times \hat D_{r_0}.$
We will make the following assumptions on $\Gm.$

There exists a real--analytic function $\sa\to\FO(\sa)\in\hol^1_{s_0}$ such that the following assumptions hold\footnote{$\checco_\diamond$ defined in \eqref{genesis}, while
 $(\morse,\b)$-Morse-non-degeneracy in
Definition \ref{morso}.}
\begin{flalign}\label{A1}
&
{\modulo}\Gm  - \FO {\modulo}_{\hat D,{r_0},s_0}
\leq 
\checco_\diamond=\checco_\diamond(\morse,\b,s_0,r_0)
\tag{A1}
\\ \ \nonumber\\
\label{A2}
&\FO\in\hol^1_{s_0}\  \text{is} \ 
(\morse,\b)\text{--}Morse\text{--}non\text{--}degenerate
\tag{A2}
\end{flalign} 
In alternative to \equ{A2}     
we will later consider  the following
{\sl stronger}\footnote{Recall \eqref{nibelungo} and \eqref{psico} .}  
assumption:
$\Gm$ is $(-\ch)$-cosine--like
according to Definition \ref{sorellastre}, namely
\begin{flalign}\label{A2'}
&
 {\modulo}\Gm(\sa,\hat\act)+\ch\cos \sa {\modulo}_{\hat D,r_0,s_0} \leq 
 \checco_\diamond^\star(\ch,s_0,r_0)
\,. \tag{A2$'$}
\end{flalign}

\begin{theorem}\label{glicemiak} {\bf (Universal  analytic behavior  of actions)}\\
{\bf Part I.}
Assume that $\Gm$ satisfies \eqref{A1} and \eqref{A2}.

\noi
{\rm (i)} {\bf (Universal behavior at critical energies)}\\
There exist
real-analytic functions $\phi^{(i)}_\pm(z,\hat \act),$ $\psi^{(i)}_\pm(z,\hat \act)$, 
with holomorphic extension 
on\footnote{$\mathtt r_2$ 
was defined in \eqref{tropicana}.} 
$|z|< \mathtt r_2$, $\hat \act\in \hat D_{r_0}$, 
with estimates
\begin{eqnarray}\label{pappagallok}
&&\sup_{|z|<\mathtt r_2, \,\hat \act\in \hat D_{r_0}}
|\phi^{(i)}_\pm|
\leq
\frac{ 2^{84}\morse^{8}}{s_* s_0^{13}\b^{17/2} }
\,,
\qquad
\sup_{|z|<\mathtt r_2, \,\hat \act\in \hat D_{r_0}}
|\psi^{(i)}_\pm|
 < \frac{2^9\sqrt\morse}{\b s_*^{3/2}}\,,
\nonumber
\\
&&
\sup_{|z|<\mathtt r_2, \,\hat \act\in \hat D_{r_0/2}}
|\partial_{\hat \act}\phi^{(i)}_\pm|
\leq
M_\f
\,,
\qquad
\sup_{|z|<\mathtt r_2, \,\hat \act\in \hat D_{r_0/2}}
|\partial_{\hat \act}\psi^{(i)}_\pm|
< M_\psi \checco\,,
\end{eqnarray}
where $M_\f,M_\psi$ where defined in \eqref{trocadero*}.
Moreover\footnote{For $0\leq i\leq 2N $ except
$i=0,2N $ and the $+$ sign.
}
\begin{equation}\label{LEGOk}
\act_n^{(i)}\Big(E_\pm^{(i)}(\hat \act)\mp \morse z, \,\hat \act\Big)=
\phi^{(i)}_\pm(z,\hat \act) +\psi^{(i)}_\pm(z,\hat \act)\ z \log z\, 
 \,,\qquad \text{for}\ \ \ 0<z<\mathtt r_2\,,\ \ 
 \hat \act\in \hat D\,.
\end{equation}

\noi
{\rm (ii)} {\bf (Analyticity at minimal energies)}\\
In the case of relative minimal critical energies (i.e., $i=2j-1$) 
\begin{equation}\label{lamponek}
\psi^{(2j-1)}_-=0\,,
\end{equation}
while in the other cases\footnote{
Namely $\psi^{(2j)}_\pm$  and $\psi^{(2j-1)}_+.$}
\beq{ciofecak}
|\psi^{(i)}_\pm(0,\hat \act)| 
\geq 
\frac{s_0}{32 \sqrt{\morse}}
\,.
\eeq

\noi
{\rm (iii)} {\bf (Perturbative behavior away from critical energies)}\\
Let 
 \begin{equation}\label{bassoraTH}
  \frac{2^5\morse}{\b s_*^2}\checco
  \leq
  \tilde\checco
  \leq \mathtt r_4=
  \frac{\mathtt r_2\morse}{2^5}
  \,.
\end{equation}
For $1\leq j\leq N,$
we have that the functions 
$
\partial_E \act_n^{(2j-1)}(E,\hat\act),$
$
\partial_E \bar \act_n^{(2j-1)}(E)
$ are holomorphic 
on\footnote{Recall \eqref{autunno2}.}
 $\mathcal E^{2j-1}(\mathtt r_4/4,\tilde\checco)
\times\hat D_{r_0}$ with 
 \begin{equation}
\sup_{\mathcal E^{2j-1}(\mathtt r_4/4,\tilde\checco)
\times\hat D_{r_0}}
\Big|
\partial_E \act_n^{(2j-1)}(E,\hat\act)
-
\partial_E \bar \act_n^{(2j-1)}(E)
\Big|
\leq
\left(
\frac{ 2^{78}\morse^{7}}{s_0^{12}\b^{7} \mathtt r_2^{3/2}}
+\frac{s_*^{12}\morse}{\tilde\checco}
\right)
\frac{2^{22}\morse^{1/2}}{\b^2 s_0}
\checco
\,.
\label{rosettaTH}
\end{equation}
For $1\leq j< N,$
the functions
$\partial_E \act_n^{(2j)}(E,\hat\act)$
and
$\partial_E \bar \act_n^{(2j)}(E)$
are holomorphic 
on\footnote{Recall \eqref{retto}.}  $\mathcal R^{(2j)}(\mathtt r_4/4,\tilde\checco)\times\hat D_{r_0}$
with estimates
 \begin{equation}
\sup_{\mathcal R^{(2j)}(\mathtt r_4/4,\tilde\checco)
\times\hat D_{r_0}}
\Big|
\partial_E \act_n^{(2j)}(E,\hat\act)
-
\partial_E \bar \act_n^{(2j)}(E)
\Big|
\leq
\left(
\frac{ 2^{75}\morse^{6}}{s_0^{21/2}\b^{6} \mathtt r_2^{3/2}}
+
\frac{\b}{\tilde\checco}
\right)
\frac{ 2^{25}\morse^{5/2}}{\b^{4}  s_0^{9/2}\mathtt r_2}
\checco
\,.
\label{rosetta+TH}
\end{equation}
 Finally for  $i=0,2N,$ 
 the functions
$\partial_E \act_n^{(i)}(E,\hat\act)$
and
$\partial_E \bar \act_n^{(i)}(E)$
are holomorphic 
on\footnote{Recall \eqref{autunno2+}.}  
$\mathcal E^{i}(r_0^2\b/2^{10}\morse,\tilde\checco)$
$\times\hat D_{r_0}$
 with estimates
\begin{equation}\label{rosetta2NTH}
\Big|
\partial_E \act_n^{(i)}(E,\hat\act)
-
\partial_E \bar \act_n^{(i)}(E)
\Big|
\leq
\left(
\frac{ 2^{56}\morse^{13/2}}{\mathtt r_2 \b^{8}  s_0^{29/2}}
\ln\left(
4+\frac{\b}{2\tilde\checco}
\right)
+
\frac{2^{24}\morse^{5/2}}{\b^{3}s_0^{9/2} \tilde\checco}
\right)
\checco\,.
\end{equation}

\noi
{\rm (iv)}
{\bf (Estimates on the derivatives of the actions)}\\
We have
\begin{equation}\label{vana}
\min_{1\leq i\leq 2N-1 } 
\inf_{\hat\act \in \hat D}
\inf_{E\in ( E^{(i)}_-(\hat\act), E^{(i)}_+(\hat\act))}
\partial_E  \act_n^{(i)}(E,\hat\act)
\geq \frac{\sqrt\b s_0^{3/2}}{64\morse}
=:\frac{1}{\bar c}\,.
\end{equation}

\begin{equation}\label{moldavater}
\partial_E  \act_n^{(2N)}(E,\hat\act)\,,
\
-\partial_E  \act_n^{(0)}(E,\hat\act)
\geq \frac{1}{4\sqrt{E+3\morse/2}}
\,,
\ \ 
\forall\,
 E_{2N}(\hat\act)<E<R_0^2-2\morse\,,\  \hat\act\in\hat D\,.
\end{equation}

\noi
{\rm (v)} {\bf (Estimates on the derivatives of the 
energies)}\\
Take $0<\loge\leq \mathtt r_4/2.$
Recalling the definition of
$\Bu^i (\loge)$ in \eqref{playmobilk}
and of $\rho(\loge)$ in \eqref{sibari}
we have that the following estimates hold uniformly
in $I\in \Bu^i_{\rho(\loge)} (\loge)$
and for every $1\leq i<2N,$
\begin{eqnarray}
&&
\left|\partial_{\act_n \act_n} \mathtt E^{(i)} \right|
\leq
\bar c^3 c_{2,0}\frac{1}{\loge}\,,
\nonumber
\\
&&\left|\partial_{\act_n \hat \act} \mathtt E^{(i)} \right|
\leq
\left(\bar c^3 c_{2,0}c_{0,1}\frac{1}{\loge}
+\bar c^2 c_{1,1}\ln\frac{\morse}{\loge}\right)
\checco
\,,
\nonumber
\\
&&\left|\partial_{\hat \act \hat \act} \mathtt E^{(i)} \right|
\leq
\left(
\bar c\, c_{0,2}
+
2\bar c^2 c_{0,1}
c_{1,1}\checco\ln\frac{\morse}{\loge}
+
\bar c^3  c_{2,0} c_{0,1}^2 \frac{\checco}{\loge}
\right)
\checco
\,,
\label{ofena}
\end{eqnarray}
where the above constants are defined in
\eqref{spoleto} and \eqref{vana}.
The same estimates hold for the case
$i=0,2N,$ with $8 R_0$ instead of $\bar c.$

\noindent
{\bf Part II (cosine-like case)}

\noindent
If $\Gm$  satisfies 
\equ{A2'} 
then, 
\begin{equation}\label{malandrino} 
\inf_{\hat\act\in \hat D}\inf_{E_1(\hat\act)<E<E_2(\hat\act)} \partial_E \act_n^{(1)}(E,\hat\act)  \geq
\frac{ s_0^{3/2}}{2^8e^{s_0}\sqrt\ch}
\end{equation}
and
\begin{equation}\label{thor}
\inf_{E_1(\hat\act)<E<E_2(\hat\act),\,\hat \act\in \hat D}
-\partial^2_{\act_n} {\mathtt E}^{(1)}
\, \geq \,
\frac{1}{16}\,,
\qquad
\inf_{E_2(\hat\act)<E<R_0^2-2\morse,\,\hat \act\in \hat D} 
\partial^2_{\act_n} {\mathtt E}^{(2)}
\, \geq \,
2\,.
\end{equation}

\end{theorem}

\subsection{Proof of Theorem \ref{glicemiak}}

\noi
{\rm (i)} follows by Lemmata 
\ref{giga2} and \eqref{giga2+},
noting that $\mathtt r_2<\mathtt r_3$
(defined in \eqref{tropicana} and 
\eqref{tropicana+} respectively)
and using \eqref{noce2} (and \eqref{alcafone}, \eqref{harlock}).

\noi
{\rm (ii)}
\eqref{lamponek}
follows by Lemma \ref{giga}.
\eqref{ciofecak} follows by 
\eqref{otello} and \eqref{otello+}.

\noi
{\rm (iii)}
follows by Corollaries \ref{sarabanda3}, \ref{nebraska} and  \ref{alaska}.

\noi
{\rm (iv)} \eqref{vana} is
\eqref{galadriel+}.
\eqref{moldavater} is \eqref{galadriel2+}.

\noi
{\rm (iv)} 
\eqref{ofena} follows by
\eqref{daitarn3}, \eqref{spoleto}, \eqref{alabama}
and \eqref{vana} (	\eqref{moldavater} in the case
$i=0,2N$).

\noi
{\bf Part II}\\
\eqref{malandrino} is \eqref{vanacoseno}.
\eqref{thor} is \eqref{cardellino} (recall \eqref{daitarn3}).

\section{Appendix}
\subsection{Technical lemmata on holomorphic functions}
\label{tucumano}

\begin{lemma}
Let $f(\z)$ be a holomorphic function on the domain
$C_*\cap \{ |z|<r\}$ for some $r>0.$
Assume that $f'(\z)=a(\z)+b(\z)\ln \z,$ for some  functions
$a(\z),b(\z)$ holomorphic in $\{ |z|<r\}$.
 Then  $f(\z)=\phi(\z)+\chi(\z)\z \ln \z$ where
 $\phi(\z),\chi(\z)$ are holomorphic functions in $\{ |z|<r\}$
 with 
 \begin{equation}\label{ilota}
 \chi(\z)=\frac{1}{\z}\int_0^\z b(z)dz\,,\qquad\qquad
 \phi(\z)=c+\int_0^\z (a(\z)-\chi(\z))dz\,,
\end{equation}
 for some constant $c.$
\end{lemma}
\proof
We first note\footnote{Writing
$b(\z)=\sum_{n\geq 0} b_n\z^n$, we get
 $\chi(\z)=\sum_{n\geq 0} \frac{b_n}{n+1}\z^n$.}
 that $\chi(\z)$ is 
holomorphic in $\{ |z|<r\}$.
It is immediate to see that
$\chi,\phi$ must satisfy the equations
$
\z \chi'+\chi=b
$
and
$\phi'+\chi=a.$
Since every solution of the homogeneous equation
$\z g'+g=0$ has the form
$g(\z)=\const/\z,$ we have that the only   solution
 of the
inhomogeneous  equation which is continuous at $\z=0$ is the
one defined in \eqref{ilota}.
The formula for $\phi$ is obvious. 
\eproof

\begin{lemma}
Let $f(\z)$ be  holomorphic on
$C_*\cap \{ |z|<r\}$ and continuous on 
$ \{ |z|<r\}$, 
 for some $r>0.$
 Then 
 $f(\z)$ is holomorphic on 
$ \{ |z|<r\}$.
\end{lemma}
\proof
It directly follows by the following well known result in complex analysis\footnote{See, e.g.,  Theorem 16.8 of Rudin's book, {\sl Real and Complex Analysis}.}:
 {\it Suppose $\O$ is a region, $L$ is a straight line or a circular are, $\O\setminus L$ is the union of 
two regions $\O_1$ and $\O_2$, $f$ 
is continuous in $\O$, and $f$ 
is holomorphic in $\O_1$ and in $\O_2$. Then $f$ is holomorphic in $\O$.}
\eproof

\begin{lemma}\label{fettuccine1}
 Assume that $\chi(\sa)$
 is a holomorphic function on the complex ball
 $B_r=B_r(0)$ for some radius $r>0,$
 with $\chi(0)=0$ and
 $\sup_{B_{r}}|\chi|\leq \checco.$
 Then the function $\phi(\sa):=
 \chi(\sa)/\sa$ for $\sa\neq 0$
 and $\phi(0):=\chi'(0),$ is holomorphic
 in $B_r$ with 
 $\sup_{B_{r/2}}|\phi|\leq 2\checco/r.$
\end{lemma}
 \proof
 By Cauchy estimates 
 $\sup_{B_{r/2}}|\chi'|\leq 2\checco/r,$
 then, since $\chi(0)=0,$
 $|\chi(\sa)|\leq \frac{2\checco}{r}|\sa|$
 for every $\sa\in B_{r/2}.$
 \eproof

\begin{lemma}\label{fettuccine2}
 Let $g,\phi$ be holomorphic 
 on $B_r$ with 
 $\sup_{B_r}|g|\leq M_1$,
 $\sup_{B_r}|\phi|\leq\epsilon\leq 1/4$
 and
 $\inf_{B_{r/2}}|g|=:g_0>0$.
 Then 
 $$
 |g(\sa+\sa \phi(\sa))-g(\sa)|
 \leq
 \frac{2M_1}{g_0}\epsilon |g(\sa)|\,,
 \qquad
 \forall\, |\sa|< \frac{r}{2}\,.
 $$
\end{lemma}
 \proof
If $|\sa|< r/2$ then
$|\sa+\sa \phi(\sa)|<3r/4$ and, by Cauchy
estimates,
\begin{eqnarray*}
&&
 |g(\sa+\sa \phi(\sa))-g(\sa)|
 \leq \sup_{B_{3r/4}}|g'|\epsilon \frac{r}{2}
 \leq
 2M_1\epsilon
 \leq
 \frac{2M_1\epsilon}{g_0}|g(\sa)|\,. 
\end{eqnarray*}
 
 \eproof

 \begin{lemma}\label{fettuccine3}
 Let $f,\chi$ be holomorphic on 
 $B_r$, satisfying $f(0)=\chi(0)=0,$
 $f'(0)\neq 0,$ $\sup_{B_r}|f|\leq M_2,$
 $\sup_{B_r}|\chi|\leq \checco\leq r/8.$
Then for every
 $
 \rho\leq\frac{r^2|f'(0)|}{16 M_2}
 $
 \begin{equation}\label{pollo}
 |f(\sa+\chi(\sa))-f(\sa)|\leq
 \frac{2}{\rho}\checco|f(\sa)|\,,
 \qquad\qquad
 \forall\, \sa\in B_\rho
 \,. 
 \end{equation}
\end{lemma}
 \proof
 First we note that
 $
 \rho\leq r/4
 $
 since by Cauchy estimates
 $|f'(0)|\leq M_2/r.$
 Set
 $g(\sa):=
 f(\sa)/\sa$ for $\sa\neq 0$
 and $g(0):=f'(0)$
 and set also
 $\phi(\sa):=
 \chi(\sa)/\sa$ for $\sa\neq 0$
 and $\phi(0):=\chi'(0).$
 $g$ and $\phi$ are holomorphic on 
 $B_r.$ By Lemma \ref{fettuccine1}
 $$
 \sup_{B_{r/2}}|\phi|\leq 
 \frac{2\checco}{r}\leq \frac14\,,
 \qquad
  \sup_{B_{r/2}}|g|\leq 
 \frac{2 M_2}{r}
 \,.
 $$
 We have (recalling that $\rho\leq r/4$)
 $$
 \inf_{B_\rho}|g|\geq |g(0)|-
 \sup_{B_\rho}|g'|\rho
 \geq |f'(0)|-  \frac{8 M_2}{r^2} \rho
 \geq \frac12|f'(0)|\,.
 $$
 By Lemma \ref{fettuccine2}
 (with $r\rightsquigarrow 2\rho$) we  have
  $$
 |g(\sa+\sa \phi(\sa))-g(\sa)|
 \leq
 \frac{8M_2}{r |f'(0)|}
 \frac{2\checco}{r}
  |g(\sa)|
  =
  \frac{1}{\rho}\checco  |g(\sa)|
  \,,
 \qquad
 \forall\, \sa\in B_\rho\,.
 $$
 Moreover, by the last estimate, we get
  $$
 |g(\sa+\sa \phi(\sa))
 \phi(\sa)
 |
 \leq
 |g(\sa+\sa \phi(\sa))
 |
 \frac{2\checco}{r}
 \leq
  \frac{2\checco}{r}
  (1+\frac{\checco}{\rho})|g(\sa)|
 \leq
  \frac{1}{\rho}\checco  |g(\sa)|
  \,,
 \qquad
 \forall\, \sa\in B_\rho\,,
 $$ 
This shows that
  $$
 |g(\sa+\sa \phi(\sa))
 (1+\phi(\sa))
 -g(\sa)|
 \leq
  \frac{2}{\rho}\checco  |g(\sa)|
  \,,
 \qquad
 \forall\, \sa\in B_\rho\,,
 $$ 
 which is equivalent to 
 \eqref{pollo} (dividing by $|\sa|$).
 \eproof



\begin{thebibliography}{99}

\footnotesize

\bibitem{AKN}
V.~I. Arnold, V.~V. Kozlov, and A.~I. Neishtadt.
Mathematical aspects of classical and celestial mechanics,\emph{
volume~3 of Encyclopaedia of Mathematical Sciences}.
Springer-Verlag, Berlin, third edition, 2006.
[Dynamical systems. III], Translated from the Russian original by E. Khukhro.

%

\bibitem{BFS2} D. Bambusi, A. Fus\`e, M. Sansottera. Exponential stability in the perturbed central force problem,
Regular and Chaotic Dynamics volume 23, pages 821-- 841(2018), doi
https://doi.org/10.1134/S156035471807002X

\bibitem{BC} L. Biasco, and L. Chierchia. KAM Theory for secondary tori, 
 arXiv:1702.06480v1 [math.DS] (21 Feb 2017)


\bibitem{BClin} L. Biasco, and L. Chierchia. 
On the measure of Lagrangian invariant  tori in nearly--integrable mechanical systems. 
Rend. Lincei Mat. Appl. {\bf 26} (2015), 1--10 

\bibitem{BCnonlin} L. Biasco, and L. Chierchia. 
On the topology of nearly--integrable Hamiltonians at simple resonances. To appear in Nonlinearity (2020)









\bibitem{MNT} A.G. Medvedev, A.I.  Neishtadt, D.V. Treschev, 
Lagrangian tori near resonances of near--integrable Hamiltonian systems, 
Nonlinearity, {\bf 28}:7 (2015), 2105--2130

%
%



%




\end{thebibliography}
\end{document}